\documentclass{memo-l}
\usepackage{graphicx,graphics,amssymb,amscd,amsmath,amsthm,amssymb,latexsym,mathrsfs,verbatim,ifthen}
\usepackage[dvips]{color}
\usepackage{xy}
\xyoption{all}
\usepackage{multirow}
\usepackage{tabls}
\input epsf

\numberwithin{section}{chapter}
\numberwithin{equation}{chapter}
\numberwithin{figure}{chapter}

\newcommand{\reals}{{\mathbb R}}
\newcommand{\integers}{{\mathbb Z}}
\newcommand{\complex}{{\mathbb C}}
\newcommand{\posint}{{\mathbb N}}
\newcommand{\rational}{{\mathbb Q}}
\newcommand{\A}{{\mathcal A}}
\newcommand{\Z}{{\mathcal Z}}
\renewcommand{\P}{{\mathcal P}}
\newcommand{\Q}{{\mathcal Q}}
\newcommand{\M}{{\mathcal M}}
\newcommand{\N}{{\mathcal N}}
\newcommand{\R}{{\mathcal R}}
\renewcommand{\S}{{\mathcal S}}
\renewcommand{\L}{{\mathscr L}}
\newcommand{\C}{{\mathcal C}}
\newcommand{\WL}{{\sf WL}}
\newcommand{\FL}{{\sf FL}}
\newcommand{\CL}{{\sf CL}}

\newcommand{\filter}{{\rm V}}

\newcommand{\turnup}[1]{\rotatebox[origin=c]{90}{\ensuremath#1}}
\newcommand{\turndown}[1]{\rotatebox[origin=c]{270}{\ensuremath#1}}
\newcommand{\turndl}[1]{\rotatebox[origin=c]{225}{\ensuremath#1}}
\newcommand{\turnd}[1]{\rotatebox[origin=c]{-90}{\ensuremath#1}}

\newcommand{\turnnegsixty}[1]{\rotatebox[origin=c]{-60}{\ensuremath#1}}
\newcommand{\turnonetwenty}[1]{\rotatebox[origin=c]{120}{\ensuremath#1}}
\newcommand{\turnaround}[1]{\rotatebox[origin=c]{180}{\ensuremath#1}}

\newcommand{\Cat}{{\rm Cat}}
\newcommand{\Nar}{{\rm Nar}}

\newcommand{\Weak}{{\sf Weak}}
\newcommand{\Abs}{{\sf Abs}}
\newcommand{\Int}{{\sf Int}}
\newcommand{\rk}{{\sf rk}}
\newcommand{\Fix}{{\sf Fix}}
\newcommand{\Mov}{{\sf Mov}}
\newcommand{\im}{{\rm im}}
\newcommand{\Iso}{{\sf Iso}}
\newcommand{\sint}{\textstyle{\int}}
\newcommand{\mins}{{\sf mins}}
\newcommand{\nz}{{\sf nz}}
\newcommand{\Shi}{{\sf Shi}}

\newcommand{\Pge}{\Phi_{\geq -1}}
\newcommand{\abs}[1]{ \left|#1 \right| }

\newcommand{\set}[1]{\left\{ #1\right\} }

\newcommand{\qbin}[2]{\genfrac{[}{]}{0pt}{}{#1}{#2}}

\newtheorem{theorem}{Theorem}[section]
\newtheorem{proposition}[theorem]{Proposition}
\newtheorem{conjecture}[theorem]{Conjecture}
\newtheorem{corollary}[theorem]{Corollary}
\newtheorem{lemma}[theorem]{Lemma}
\newtheorem{carterslemma}[theorem]{Carter's Lemma}
\newtheorem{shiftinglemma}[theorem]{The Shifting Lemma}
\newtheorem{genshiftinglemma}[theorem]{Generalized Shifting Lemma}
\newtheorem*{exchangeproperty}{Exchange Property}
\newtheorem*{deletionproperty}{Deletion Property}
\newtheorem{subwordproperty}[theorem]{The Subword Property}
\newtheorem{dualtitslemma}[theorem]{Dual Tits' Lemma}

\theoremstyle{definition}

\newtheorem*{update}{Update}

\newtheorem{disclaimer}[theorem]{Disclaimer}

\newtheorem{example}[theorem]{Example}
\newtheorem{notation}[theorem]{Notation}
\newtheorem{definition}[theorem]{Definition}
\newtheorem{altdefinition}[theorem]{Alternate Definition}
\newtheorem{problem}[theorem]{Open Problem}

\begin{document}

\frontmatter

\title{Generalized Noncrossing Partitions and Combinatorics of Coxeter Groups}

\author{Drew Armstrong}
\address{Department of Mathematics, Cornell University, Ithaca, New York 14853}
\curraddr{School of Mathematics, University of Minnesota, Minneapolis, Minnesota 55455}
\email{armstron@math.umn.edu}
\thanks{This work was supported in part by NSF grant DMS-0603567.}
	
\date{\today}

\maketitle

\cleardoublepage
\thispagestyle{empty}
\vspace*{13.5pc}
\begin{center}
  For Moira.
\end{center}
\cleardoublepage

\setcounter{page}{7}
\tableofcontents


\chapter*{Abstract}
This memoir is a refinement of the author's PhD thesis --- written at Cornell University (2006). It is primarily a desription of new research but we have also included a substantial amount of background material. At the heart of the memoir we introduce and study a poset $NC^{(k)}(W)$ for each finite Coxeter group $W$ and each positive integer $k$. When $k=1$, our definition coincides with the generalized noncrossing partitions introduced by Brady and Watt \cite{brady-watt:kpione} and Bessis \cite{bessis:dual}. When $W$ is the symmetric group, we obtain the poset of classical $k$-divisible noncrossing partitions, first studied by Edelman \cite{edelman:kdiv}.

In general, we show that $NC^{(k)}(W)$ is a graded join-semilattice whose elements are counted by a generalized ``Fuss-Catalan number'' $\Cat^{(k)}(W)$ which has a nice closed formula in terms of the degrees of basic invariants of $W$. We show that this poset is locally self-dual and we also compute the number of multichains in $NC^{(k)}(W)$, encoded by the zeta polynomial. We show that the order complex of the poset is shellable (hence Cohen-Macaulay) and we compute its homotopy type. Finally, we show that the rank numbers of $NC^{(k)}(W)$ are polynomials in $k$ with nonzero rational coefficients alternating in sign. This defines a new family of polynomials (called ``Fuss-Narayana'') associated to the pair $(W,k)$. We observe some interesting properties of these polynomials.

In the case that $W$ is a classical Coxeter group of type $A$ or $B$, we show that $NC^{(k)}(W)$ is isomorphic to a poset of ``noncrossing'' set partitions in which each block has size divisible by $k$. This motivates our general use of the term ``$k$-divisible noncrossing partitions'' for the poset $NC^{(k)}(W)$. In types $A$ and $B$ we prove ``rank-selection'' and ``type-selection'' formulas refining the enumeration of multichains in $NC^{(k)}(W)$. We also describe bijections relating multichains of classical noncrossing partitions to ``$k$-divisible'' and ``$k$-equal'' noncrossing partitions. Our main tool is the family of Kreweras complement maps.

Along the way we include a comprehensive introduction to related background material. Before defining our generalization $NC^{(k)}(W)$, we develop from scratch the theory of the generalized noncrossing partitions $NC^{(1)}(W)$ as defined by Brady and Watt \cite{brady-watt:kpione} and Bessis \cite{bessis:dual}. This involves studying a finite  Coxeter group $W$ with respect to its generating set $T$ of {\em all} reflections, instead of the usual Coxeter generating set $S$. This is the first time that this material has appeared together.

Finally, it turns out that our poset $NC^{(k)}(W)$ shares many enumerative features in common with the generalized nonnesting partitions of Athanasiadis \cite{athanasiadis:cat,athanasiadis:nar} and the generalized cluster complexes of Fomin and Reading \cite{fomin-reading}. We give a basic introduction to these topics and we make several conjectures relating these three families of ``Fuss-Catalan objects''.


\chapter*{Acknowledgements}
This work has taken several years and gone through several phases to arrive in its present form. The following people have contributed along the way; they deserve my sincerest thanks.

First, I thank Lou Billera, my advisor at Cornell University, under whose guidance the original draft \cite{armstrong} was written. I thank my wife Heather for her thorough reading of the manuscript and for helpful mathematical discussions. Thanks to the anonymous referee for many valuable suggestions. And thanks to Hugh Thomas who made several contributions to my understanding; in particular, he suggested the proof of Theorem \ref{th:ELNCk}. 

I have benefitted from valuable conversations with many mathematicians; their suggestions have improved this memoir in countless ways. Thanks to: (in alphabetical order) Christos Athanasiadis, David Bessis, Tom Brady, Ken Brown, Fr\'ed\'eric Chapoton, Sergey Fomin, Christian Krattenthaler, Cathy Kriloff, Jon McCammond, Nathan Reading, Vic Reiner and Eleni Tzanaki. I would especially like to thank Christian Krattenthaler for his close attention to the manuscript and for his subsequent mathematical contributions \cite{krattenthaler1,krattenthaler2,kratt:typeD,kratt-muller}.

Thank you to the American Institute of Mathematics, Palo Alto, and to the organizers --- Jon McCammond, Alexandru Nica and Vic Reiner --- of the January 2005 conference \cite{armstrong:braids} at which I was able to meet many of the above-named researchers for the first time. Thanks also to John Stembridge for the use of his software packages \verb|posets| and \verb|coxeter| for \verb|Maple| \cite{stembridge}.

Finally, I thank the people who introduced me to noncrossing partitions. My first research experience was under the supervision of Roland Speicher at Queen's University. His papers on noncrossing partitions, as well as those of Paul Edelman and Rodica Simion, inspired me to look further. In particular, I thank Paul Edelman for his encouragement and advice.
\mainmatter


\chapter{Introduction}
\label{sec:introduction}

The main purpose of this memoir is to communicate new results in the theory of finite Coxeter groups. However, we have also included a great deal of background material, which we hope will serve as a reference and a foothold for further progress in the subject of ``Coxeter-Catalan combinatorics''. In this introductory chapter we will sketch our motivations and give an outline of the rest of the memoir.

\section{Coxeter-Catalan Combinatorics}

An element of $GL(\reals^n)$ is a {\sf reflection} if it sends some nonzero vector $\alpha\in\reals^n$ to its negative $-\alpha$ and fixes the orthogonal hyperplane $\alpha^\perp$ pointwise. If $W$ is a finite group with a faithful representation $\rho: W\hookrightarrow GL(\reals^n)$ generated by reflections, we call the pair $(W,\rho)$ a {\sf finite real reflection group}. The main motivation for the study of these groups comes from Lie theory. However, the general notion of symmetries generated by reflections is fundamental in geometry and it predates the notion of ``group''. The subject of finite real reflection groups developed through the nineteenth century --- notably in the works of M\"obius, Jordan, Schl\"afli, Killing, Cartan, and Weyl --- and reached its definitive form with the complete classification \cite{coxeter} by Coxeter in 1935, using the {\sf Coxeter diagrams}. Coxeter referred to these groups as ``reflection groups'' throughout his life, but they have commonly been known as {\sf finite Coxeter groups} since Tits used the term in 1961~\cite{tits}. For notes on the history of reflection groups, we refer to Bourbaki~\cite[Chapter 26]{bourbaki} and Coxeter \cite[Chapter 11]{coxeter:polytopes}.

We say that a finite reflection group $(W,\rho)$ is {\sf reducible} if we have $W=W'\times W''$ where $W'$ and $W''$ are proper subgroups generated by reflections of $(W,\rho)$. In this case, $\rho$ decomposes as a direct sum $\rho'\oplus\rho''$ of $W$-representations such that $(W',\rho'|_{W'})$ and $(W'',\rho''|_{W''})$ are themselves finite reflection groups. If $(W,\rho)$ cannot be written in this way, it is {\sf irreducible}.

It turns out that finite reflection groups are completely reducible; thus the problem of classification is to enumerate all of the possible irreducible cases. This has been done \cite{coxeter}. Following the standard Cartan-Killing notation, there are eight families of finite irreducible Coxeter groups denoted by the letters $A$ through $I$. ($B$ and $C$ denote isomorphic groups, although we will see in Section \ref{sec:rootsystems} that they possess nonisomorphic ``crystallographic root systems''.) A subscript in the notation indicates the {\sf rank} of the reflection group, which is the dimension of the representation.

\begin{figure}[ht]
\vspace{.1in}
\begin{center}
\scalebox{0.9}{\input{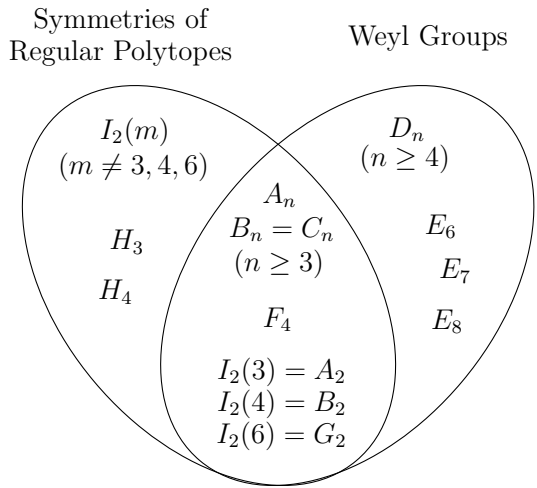}}
\end{center}
\caption{The finite irreducible Coxeter groups}
\label{fig:cox}
\end{figure}

The finite irreducible Coxeter groups fall into two overlapping classes: the {\sf Weyl groups} (or the ``crystallographic'' groups), which stabilize a lattice in $\reals^n$; and the {\sf groups of symmetries of regular polytopes} (see Figure \ref{fig:cox}).
This classification is an important landmark in contemporary mathematics; not only does it contain the classification of regular polytopes in Euclidean space, but it is also strongly related to the classification of semisimple Lie algebras. Due to the importance of Lie theory, it has become standard to refer to a Coxeter group by the letter $W$ for ``Weyl''. For our purposes, however, we will usually not require the crystallographic hypothesis (see Section \ref{sec:crystallography}).

In the past ten years, a new perspective in this field has emerged and this perspective has recently begun to synthesize several combinatorial, algebraic and geometric topics. At the center of this synthesis is a generalized ``Catalan number'' $\Cat(W)$ defined for each finite Coxeter group $W$, with the property that the classical Catalan number $\Cat(A_{n-1})=\frac{1}{n+1}\binom{2n}{n}$ corresponds to the symmetric group of type $A_{n-1}$.\footnote{{\bf Warning:} Throughout this memoir, we will use the Cartan-Killing symbol $A_{n-1}$ to denote the symmetric group on $n$ letters. This is shorthand for $W(A_{n-1})$, ``the Weyl group of type $A_{n-1}$''. We hope that this notation will cause no confusion with the alternating groups.} We refer to $\Cat(W)$ as the {\sf Coxeter-Catalan number} of the group $W$. Its explicit formula is
\begin{equation}
\label{eq:catalan}
\Cat(W):=\frac{1}{\abs{W}}\prod_{i=1}^n (h+d_i),
\end{equation} 
where $h$ is the Coxeter number and $d_1,d_2,\ldots,d_n$ are the degrees of  $W$, arising from its ring of polynomial invariants (see Section \ref{sec:invariant}). The number $\Cat(W)$ has been discovered independently in several different areas and wherever it appears it is accompanied by a wealth of new combinatorics.

We will now briefly describe the three main streams of thought that have converged into this {\sf Coxeter-Catalan combinatorics}. Typically, the ``type $A$'' combinatorics (corresponding to the symmetric group) was observed previously and it is the general ``type $W$'' perspective that is new. In the following sketches we use some undefined terminology. The reader unfamiliar with Coxeter theory may wish to skip these now and return after reading the background material in Chapter \ref{sec:background}. A more thorough introduction to these ideas is given in Chapter \ref{sec:fusscatcomb}.

\bigskip
\noindent {\bf Stream 1: Ideals and Antichains.}  It seems that the formula \eqref{eq:catalan} was first written down by Djokovi\'c in 1980~\cite{djokovic} in the case that $W$ is a Weyl group, although he did not recognize a connection with the Catalan numbers (personal communication). For him, the formula $\Cat(W)$ counted the number of conjugacy classes of elements of order $h+1$ in the semisimple Lie group corresponding to $W$. (More generally, he showed that
\begin{equation*}
\prod_{i=1}^n \frac{d_i+\ell-1}{d_i}
\end{equation*}
is the number of conjugacy classes of elements with order dividing $\ell$.) Haiman later rediscovered Djokovi\'c's result in an equivalent form --- in particular he showed that $\Cat(W)$ counts the number of $W$-orbits in the quotient $Q/(h+1)Q$ of the root lattice $Q$ corresponding to $W$ \cite[Theorem 7.4.4]{haiman} --- and he conjectured \cite[Conjecture 7.3.1]{haiman} that this should be the multiplicity of the sign representation in some ``type $W$ ring of diagonal coinvariants''. This conjecture was subsequently verified by Gordon~\cite{gordon}.

Soon after, Shi observed the numbers $\Cat(W)$ in his study of affine Weyl groups and hyperplane arrangements~\cite{shi}. Let $W$ be a Weyl group of rank $n$. If $(\cdot,\cdot)$ is the inner product on $\reals^n$ and $\alpha\in\reals^n$ is a root of $W$, define the hyperplanes
\begin{equation*}
H^i_\alpha:=\{x\in\reals^n:(x,\alpha)=i\}.
\end{equation*}
The collection of $H^0_\alpha$ where $\alpha$ ranges over the roots of $W$ is called the {\sf Coxeter arrangement} of $W$. To define the {\sf Catalan arrangement} we include the additional hyperplanes $H^{+1}_\alpha$, $H^{-1}_\alpha$ for all roots $\alpha$. If $\A$ is a collection (``arrangement'') of hyperplanes in $\reals^n$ then the connected components of the complement $\reals^n\setminus \cap_{H\in\A}H$ are called {\sf chambers} of $\A$. In this language, Shi's result is equivalent to the fact that the Catalan arrangement has $\abs{W}\cdot\Cat(W)$ chambers. Since the group $W$ acts simply transitively on the chambers of the Coxeter arrangement, there is a natural bijection between group elements and chambers; consequently, each of the $\abs{W}$ chambers of the Coxeter arrangement is divided into $\Cat(W)$ chambers by the Catalan arrangement. However, Shi's proof was case-by-case and he did not guess the formula \eqref{eq:catalan}.

At around the same time, Postnikov had the idea to study antichains in the root poset of $W$. If $\Phi$ is a crystallographic root system with positive roots $\Phi^+$ and simple roots $\Pi$, there is a useful partial order on $\Phi^+$ defined by setting $\alpha\leq\beta$ for $\alpha,\beta\in\Phi^+$ whenever $\beta-\alpha$ is in the positive span of the simple roots $\Pi$. If the positive roots $\Phi^+$ are taken to define a ``positive'' half-space for each hyperplane $H_\alpha^0$, then the intersection of these positive half-spaces is called the {\sf positive cone}. Postnikov independently noticed that antichains in the root poset are in bijection with ``positive'' regions of the Catalan arrangement, and he conjectured the formula \eqref{eq:catalan} (see \cite[Remark 2]{reiner}).

Later on, two independent theories gave uniform proofs of formula \eqref{eq:catalan} for counting the positive Catalan chambers, and explained both the connection to Haiman's and Djokovi\'c's work and the connection to Postnikov's ideas.
On the one hand, Athanasiadis calculated the characteristic polynomial $\phi(t)$ of the Catalan hyperplane arrangement by showing that $\phi(t)=\chi(t-h)$, where $h$ is the Coxeter number and $\chi(t)$ is the characteristic polynomial of the Coxeter arrangement. He did this for the classical types in his thesis (the type $A$ case was also considered by Postnikov and Stanley~\cite{postnikov-stanley}), and later proved a more general uniform theorem \cite[Theorem 1.2]{athanasiadis:cat}. Combining this with the well-known formula
\begin{equation}
\chi(t)=\prod_{i=1}^n (t-d_i+1)
\end{equation}
of Orlik and Solomon \cite{orlik-solomon} and Zaslavsky's theorem \cite{zaslavsky} (which implies that the number of chambers of the Catalan arrangement is equal to $(-1)^n\phi(-1)$) establishes the formula.

On the other hand, Cellini and Papi~\cite{cellini-papi} gave a bijective proof that the positive chambers in the Catalan arrangement are counted by \eqref{eq:catalan}. They defined two bijections: one from the positive Catalan chambers to the set of antichains in the root poset, and another from these antichains to $W$-orbits in the ``finite torus'' $Q/(h+1)Q$. Haiman's result \cite[Theorem 7.4.4]{haiman} then gives the enumeration.

A more detailed introduction to these ideas is given later in Section \ref{sec:nonnesting}. Note the important fact that all of the ideas in Stream 1 depend on a {\em crystallographic} root system. We will see that the other two streams do not have this restriction.

\bigskip
\noindent{\bf Stream 2: Cluster Combinatorics.}
Through their study of total positivity in Lie groups, Fomin and Zelevinsky have introduced the subject of {\sf cluster algebas}. Axiomatically, a cluster algebra can be defined as a commutative ring whose generators are grouped into {\sf clusters} of equal cardinality. These algebras occur ``in nature'' as the homogeneous coordinate rings of certain flag varieties and Grassmannians (see the survey \cite{fomin-zelevinsky:survey} by Fomin and Zelevinsky). A fundamental result in the subject is the finite type classification~\cite{fomin-zelevinsky:finitetype}, in which the finite type cluster algebras are shown to be  classified by the Dynkin diagrams; that is, they correspond to crystallographic root systems. However, we will see that the essential combinatorial structure does not depend on crystallography.

Let $\Phi$ be a crystallographic root system (for definition see Section \ref{sec:rootsystems}) with positive roots $\Phi^+$, simple roots $\Pi$ and Weyl group $W$. The generators of the type $W$ cluster algebra are in bijection with the set of {\sf almost-positive roots}
\begin{equation}
\label{eq:almostpos}
\Pge:=\Phi^+\cup (-\Pi).
\end{equation}
The relationships among the generators are recorded in the {\sf cluster complex} $\Delta(W)$, which is a pure, flag simplicial complex on $\Pge$ whose maximal simplices are the clusters of the corresponding algebra. This complex, in principle, depends on the algebra, but Fomin and Zelevinsky have given a purely combinatorial construction; they describe a binary relation on $\Pge$, called {\sf compatibility}, that determines when two almost-positive roots occur together in a cluster~\cite{fomin-zelevinsky:ysystems}.

When $W$ is the symmetric group $A_{n-1}$ the notion of compatibility is familiar: the almost-positive roots of $A_{n-1}$ correspond to the diagonals of a convex $(n+2)$-gon, and two diagonals are compatible precisely when they don't cross. In this case, the maximal compatible sets of diagonals are triangulations. That is, the complex $\Delta(W)$ is a generalization of the classical (simplicial) {\sf associahedron}.

Fomin and Zelevinsky showed case-by-case that the complex $\Delta(W)$ has $\Cat(W)$ maximal faces \cite[Proposition 3.8]{fomin-zelevinsky:ysystems}. When $W$ is a noncrystallographic finite Coxeter group, there is no associated cluster algebra, but the combinatorial definition of $\Delta(W)$ as a flag complex on $\Pge$ does generalize (see \cite[Section 5.3]{fomin-reading:survey}). However, the only known geometric realization of $\Delta(W)$ as a convex polytope (due to Chapoton, Fomin and Zelevinsky~\cite{chapoton-fomin-zelevinsky}) does not obviously extend to the noncrystallographic case, and it is an open problem to give a geometric realization of $\Delta(W)$ when $W$ is noncrystallographic\footnote{Recent work by Hohlweg, Lange and Thomas \cite{hohlweg-lange,HLT} gives a polytopal construction of $\Delta(W)$ that extends to the noncrystallographic types.}.

In a development parallel to the theory of cluster algebras, Reading has defined a class of lattice quotients of the weak order on $W$, which he calls {\sf Cambrian lattices}~\cite{reading}. Corresponding to each graph orientation of the Coxeter diagram he constructs a lattice with $\Cat(W)$ elements, and he conjectured \cite[Conjecture 1.1]{reading} that the Hasse diagram of this lattice is isomorphic to the $1$-skeleton of (the dual sphere of) the Fomin-Zelevinsky associahedron $\Delta(W)$.\footnote{Reading and Speyer have now proved this conjecture \cite{reading-speyer}.} This conjecture generalizes the known fact in type $A$ that the {\sf Tamari lattice} is an orientation of the $1$-skeleton of the classical associahedron. These structures have also been observed by Thomas. In his study of {\sf trim lattices}~\cite{thomas}, he constructed a family of lattices, called {\sf pre-Cambrian}, and he conjectured \cite[Conjecture 3]{thomas} that these are isomorphic to the Cambrian lattices of Reading.

In Section \ref{sec:cluster}, we will give a more thorough introduction to associahedra, polygon dissections and cluster complexes.

\bigskip
\noindent{\bf Stream 3: Noncrossing Partitions.} 
\label{sec:stream3}
We say a partition of the set $\{1,2,\ldots,n\}$ is {\sf noncrossing} if there do not exist $1\leq a<b<c<d\leq n$ with $a$ and $c$ together in a block and $b$ and $d$ together in a different block. The set $NC(n)$ of noncrossing partitions of $\{1,2,\ldots,n\}$ can be thought of as a poset (partially-ordered set) under refinement, and moreover it is a {\sf lattice} (each pair of elements has a least upper bound and a greatest lower bound).

The systematic study of the noncrossing partitions began with Kreweras~\cite{kreweras} in 1972, but their early history is not so clear. Stanley gives a short historical bibliography \cite[pages 261--262]{stanley:ec2} in which he attributes the first appearance of noncrossing partitions to Becker \cite{becker2} in 1948. It is true that Becker had been studying ``the general theory of rhyme'' (what we would call set partitions) \cite{becker1,becker2,becker3}, but we find no explicit mention of ``planar rhyme schemes'' (his name for noncrossing partitions) until a presentation given to the Washington DC meeting of the AMS, October 27, 1951 \cite{becker3}. A more overlooked reference is Motzkin's note \cite{motzkin} of 1948 (in which he also introduced the famous {\sf Motzkin numbers}) whose final sentence seems to suggest that he was considering noncrossing partitions. (Thanks to Dave Callan and Len Smiley \cite{callan-smiley} for the Motzkin reference.) Whatever their origin, the literature on noncrossing partitions is extensive, and they have long been a favorite object in algebraic combinatorics. For an account of the history of noncrossing partitions, see the surveys by Simion~\cite{simion} and McCammond~\cite{mccammond:noncrossing}. For a modern perspective, see our historical sketch at the end of Section \ref{sec:classicalNC}.

While initially studied for its own sake, the lattice of noncrossing partitions has recently found two surprising applications. {\sf Free probability} is a branch of operator algebras with close ties to physics. Essentially, it is a noncommutative analogue of probability theory in which the property of ``independence'' is replaced by the notion of ``freeness''. The subject was invented around 1985 by Voiculescu in an effort to understand certain von Neumann algebras called free group factors, but its most natural formulation involves ensembles of random matrices. In the early 1990's, Speicher showed that Voiculescu's theory could be encoded in the language of noncrossing partitions. It is an observation of Rota from the 1960's that the classical convolution of random variables is --- in some sense --- the same as M\"obius inversion on the lattice of set partitions. Speicher proved an analogous result in a different setting: he showed that by restricting attention to the lattice of {\em noncrossing} set partitions, one obtains Voiculescu's {\sf free convolution} of random variables. For further details, see the surveys by Speicher~\cite{speicher:survey} and Biane~\cite{biane:probability}. 

The other recent application involves a convergence of the noncrossing partitions with the theory of Coxeter groups. In the early 2000's, Brady and Bessis independently defined a beautiful algebraic generalization of the noncrossing partitions. There is a poset $NC(W)$ for each finite Coxeter group $W$, such that $NC(A_{n-1})$ is isomorphic to the classical noncrossing partitions $NC(n)$. Much of the theory of these posets was developed independently by Brady and Watt~\cite{brady,brady-watt:kpione} and Bessis~\cite{bessis:dual}, the latter of whom gave a case-by-case proof that the number of elements of $NC(W)$ is equal to \eqref{eq:catalan} \cite[Proposition 5.2.1]{bessis:dual}.

In this context, the poset $NC(W)$ is defined as an example of a {\sf Garside structure}. 
In his 1969 thesis at Oxford (see~\cite{garside}), Garside developed a new approach to the study of the word and conjugacy problems in Artin's braid groups. While Garside himself never wrote another paper in mathematics, his approach has been quite influential. Based on the ideas in Garside's thesis, a Garside structure for a group $G$ is a labelled partially-ordered set (with the lattice property) that is used to encode algorithmic solutions to the word and conjugacy problems in the group. The poset $NC(W)$ was constructed as a Garside structure for the Artin group (generalized braid group) corresponding to $W$. For an introduction to Garside structures, see McCammond~\cite{mccammond:garside}. Prior to the work of Brady-Watt and Bessis, some important special cases had been considered: Biane~\cite{biane:crossings} (working in free probability) showed that $NC(n)$ is related to the Cayley graph of the symmetric group $A_{n-1}$ and Reiner~\cite{reiner} generalized many of the classical $NC(n)$ results to the type $B$ case, using a geometrically-motivated definition of $NC(B_n)$ that later turned out to agree with the general $NC(W)$.

In the first half of Chapter \ref{sec:classical} we develop the theory of classical noncrossing partitions in detail, and we give a historical sketch at the end of Section \ref{sec:classicalNC}. The topic of algebraic noncrossing partitions $NC(W)$ is a central theme in this memoir, and it is the focus of Chapters \ref{sec:background} and \ref{sec:kdiv}.

\bigskip

\section{Noncrossing Motivation}
\label{sec:motivation}

In Chapter \ref{sec:kdiv} we will define and study a generalization of the lattice of noncrossing partitions. For each positive integer $k$ and finite Coxeter group $W$, we will define a poset $NC^{(k)}(W)$ --- called the poset of {\sf $k$-divisible noncrossing partitions} --- with the property that $NC^{(1)}(W)$ is isomorphic to $NC(W)$. We will see that this poset has beautiful enumerative and structural properties. Furthermore, it turns out that $NC^{(k)}(W)$ is closely related to other structures recently studied by Athanasiadis \cite{athanasiadis:cat,athanasiadis:nar} and Fomin and Reading \cite{fomin-reading}. Together these three topics comprise a generalization of the Coxeter-Catalan combinatorics for each positive integer $k$. The central enumerative feature of these extended theories is a ``generalized Coxeter-Catalan number'' $\Cat^{(k)}(W)$, which we call the {\sf Fuss-Catalan} number. It has an explicit formula generalizing \eqref{eq:catalan}:\begin{equation*}
\Cat^{(k)}(W)=\frac{1}{\abs{W}}\prod_{i=1}^n (kh+d_i).
\end{equation*}
We will refer to the combinatorics surrounding this number as the {\sf Fuss-Catalan combinatorics of $W$}. After developing the theory of $NC^{(k)}(W)$, we will discuss the Fuss-Catalan combinatorics more broadly in Chapter \ref{sec:fusscatcomb}, offering several conjectures relating the ``Fuss-Catalan objects''.

However, our original inspiration comes from pure combinatorics. The germ of this memoir began as undergraduate research at Queen's University under the supervision of Roland Speicher. Though this memoir bears little resemblance to the original,  it will be valuable in all that follows to keep in mind the motivating example of ``classical noncrossing partitions''.

\begin{figure}
\vspace{.1in}
\begin{center}
\scalebox{0.7}{
\input{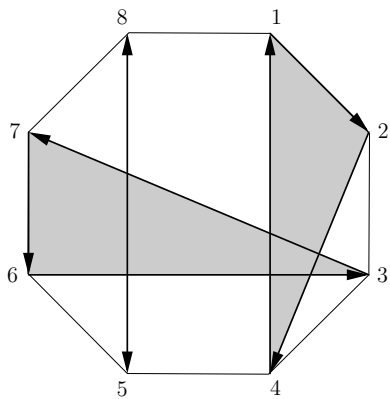}
}
\end{center}
\caption{The cycle diagram of the permutation $(124)(376)(58)$}
\label{fig:cycles}
\end{figure}

Given a permutation $\pi$ of the set $[n]=\{1,2,\ldots,n\}$, we define its {\sf cycle diagram} by labelling the vertices of a convex $n$-gon clockwise by $1,2,\ldots,n$ and drawing a directed edge from vertex $i$ to vertex $j\neq i$ whenever $\pi(i)=j$. For example, Figure \ref{fig:cycles} displays the cycle diagram of the permutation $(124)(376)(58)$ of the set $[8]$. Notice that we can identify each cycle of $\pi$ with the convex hull of its vertices. {\em If each of the cycles of $\pi$ is oriented clockwise and the convex hulls of its cycles are mutually disjoint, we say that $\pi$ is a {\sf noncrossing permutation}. Furthermore, we say that the corresponding partition of $[n]$ by $\pi$-orbits is a {\sf noncrossing partition}}. (Thus, the permutation $(124)(376)(58)$ and the partition $\{\{1,2,4\},\{3,7,6\},\{5,8\}\}$ are {\sf crossing}.)

In the 1990's, several researchers independently came up with the following classification of noncrossing permutations (see Section \ref{sec:classicalNC}). Consider the symmetric group $A_{n-1}$ of permutations of $[n]$, generated by the set $T$ of all transpositions, and consider the Cayley graph $(A_{n-1},T)$. {\em It turns out that a permutation $\pi\in A_{n-1}$ is noncrossing if and only if it lies on a geodesic between the identity $1$ and the $n$-cycle $(12\cdots n)$ in the (left or right) Cayley graph $(A_{n-1},T)$.} Of course, a different choice of $n$-cycle would lead to a different labelling of the vertices of the $n$-gon and a different notion of ``noncrossing''. (Later, we will replace ``$n$-cycles'' in the group $A_{n-1}$ with ``Coxeter elements'' in the finite Coxeter group $W$.)

\begin{figure}
\vspace{.1in}
\begin{center}
\input{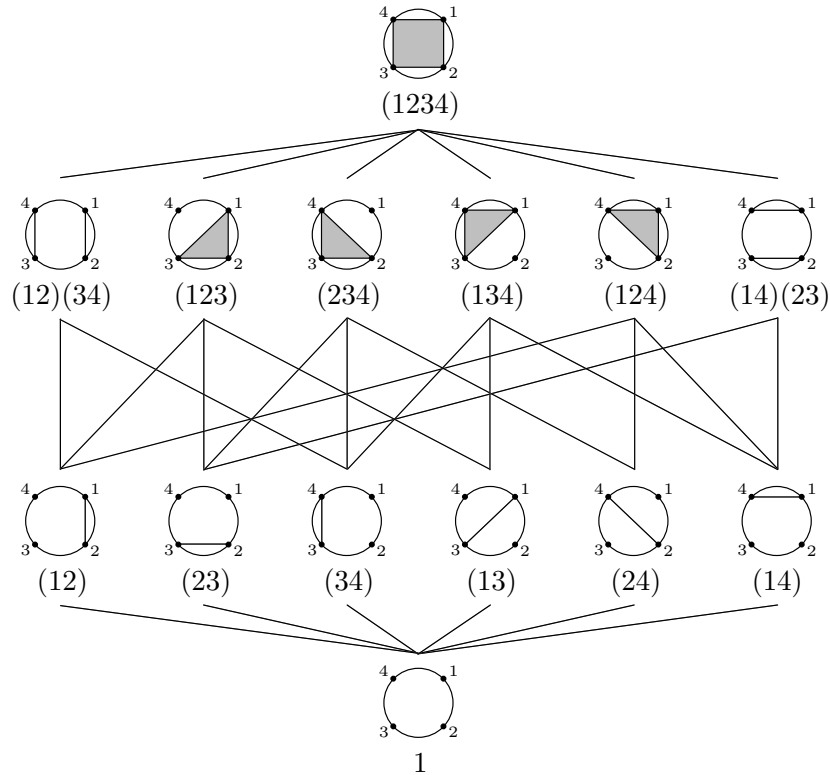}
\end{center}
\caption{The interval between $1$ and $(1234)$ in the Cayley graph of $A_3$ with respect to transpositions is isomorphic to the lattice of noncrossing partitions of the set $[4]$}
\label{fig:noncrossing}
\end{figure}

Moreover, the structure of $(A_{n-1},T)$ is closely related to the well-known ``refinement order'' on noncrossing partitions. Given noncrossing partitions $\P$ and $\Q$ of $[n]$, we say that $\P$ {\sf refines} $\Q$ --- and write $\P\leq\Q$ --- if each block of $\P$ is fully contained in some block of $\Q$. It turns out that this partial order coincides with the Cayley graph structure on $(A_{n-1},T)$. {\em The interval in the Cayley graph $(A_{n-1},T)$ between $1$ and $(12\cdots n)$ is isomorphic as a poset to the refinement order on noncrossing partitions of $[n]$}. In Figure \ref{fig:noncrossing} we display this isomorphism for $n=4$. 

\begin{figure}
\vspace{.1in}
\begin{center}
\input{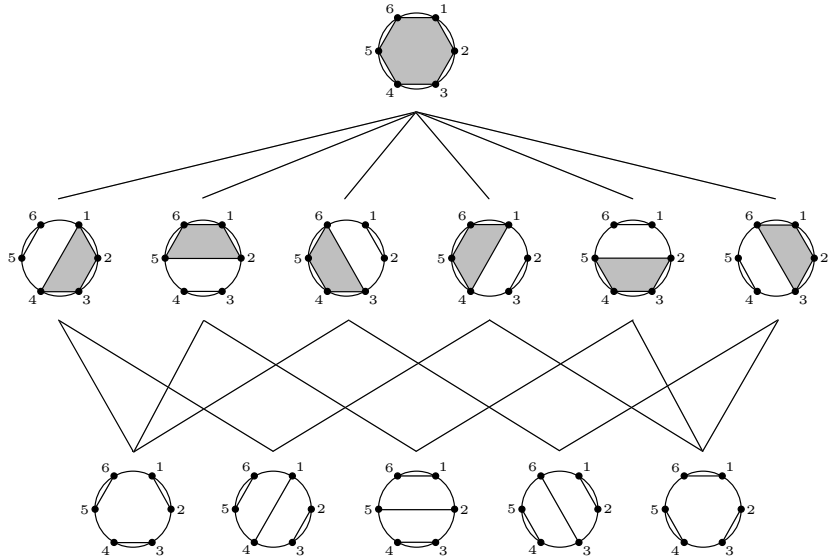}
\end{center}
\caption{$2$-divisible noncrossing partitions of the set $[6]$}
\label{fig:2divNC}
\end{figure}

The classical noncrossing partitions have a natural generalization, first studied by Edelman \cite[Section 4]{edelman:kdiv}: We say that a noncrossing set partition is {\sf $k$-divisible} if each of its blocks has cardinality divisible by $k$. Note that this can occur only when the underlying set has cardinality divisible by $k$. We let $NC^{(k)}(n)$ denote the collection of $k$-divisible noncrossing partitions of $[kn]$, partially ordered by refinement. Figure \ref{fig:2divNC} displays the Hasse diagram of the poset $NC^{(2)}(3)$. Edelman showed that this poset contains
\begin{equation*}
\frac{1}{n}\binom{(k+1)n}{n-1}
\end{equation*}
elements, which is the Fuss-Catalan number $\Cat^{(k)}(A_{n-1})$ of type $A$.

The poset $NC^{(k)}(W)$ was later considered by Stanley in connection with parking functions and quasisymmetric functions~\cite{stanley:parking}. However, the literature on $k$-divisible noncrossing partitions is not extensive. In this memoir we hope to focus more attention on Fuss-Catalan generalizations by placing $NC^{(k)}(n)$ in a more general and more natural setting.

In particular, we will show that our poset $NC^{(k)}(A_{n-1})$ is isomorphic to Edelman's $NC^{(k)}(n)$ and hence we will extend the notion of ``$k$-divisibility'' to all finite Coxeter groups. Moreover, we hope to demonstrate that reflection groups provide the ``correct'' setting for these ideas since many combinatorial arguments become simpler and more motivated in this algebraic language.

\section{Outline of the Memoir}
\label{sec:outline}

\noindent{\bf Chapter \ref{sec:background}:}\,
We begin with a quick introduction to finite Coxeter systems and root systems from a basic level, following the approach in Humphreys \cite{humphreys}. We have included this material for the reader who may be interested in noncrossing partitions but who is not familiar with Coxeter theory. The reader with a background in Coxeter theory may wish to skip these sections, or refer back to them as necessary. 

In Section \ref{sec:reducedwords} we briefly review the combinatorial approach to classical Coxeter systems via the theory of reduced $S$-words. This material is classical and it may be found in Bj\"orner and Brenti \cite{bjorner-brenti}. Sections \ref{sec:absolute}, \ref{sec:shifting} and \ref{sec:NC} describe a more recent approach to the subject, using the theory of reduced $T$-words. Much of this theory about the pair $(W,T)$ has not appeared together before. Here we will define and describe basic properties of the lattice $NC(W)$ of noncrossing partitions corresponding to the finite Coxeter group $W$. We have based this introduction on results of Carter \cite{carter}, Brady and Watt \cite{brady-watt:partialorder,brady-watt:kpione} and Bessis \cite{bessis:dual}.

In the final Section \ref{sec:invariant}, we will motivate the notions of {\sf degrees} and {\sf exponents} for the group $W$. These integers are important to the combinatorics of the finite Coxeter groups and they will appear in most of our enumerative formulas related to the poset $NC^{(k)}(W)$. Here again we follow Humphreys \cite{humphreys}.

\bigskip

\noindent{\bf Chapter \ref{sec:kdiv}:}\, 
This chapter contains the definition and development of the poset $NC^{(k)}(W)$ from a Coxeter theory point of view. All of this material is original except where stated otherwise.

First we develop the basic notions of {\sf multichains} and {\sf minimal sequences} that lead to the main definition of a {\sf delta sequence} (Definition \ref{def:mchaindsequence}). The poset $NC^{(k)}(W)$ is defined as the (dual) componentwise order on delta sequences. In Section \ref{sec:basics} we explore the basic structural properties of $NC^{(k)}(W)$, including a characterization of the principal order ideals of $NC^{(k)}(W)$ in terms of parabolic subgroups of $W$.

Then we consider enumerative questions in Section \ref{sec:fusscat}. Here we define the {\sf Fuss-Catalan number} $\Cat^{(k)}(W)$ and show that it counts the elements of $NC^{(k)}(W)$. We also define the {\sf Fuss-Narayana} numbers as the ``rank numbers'' of $NC^{(k)}(W)$ and show that these are polynomials in $k$. We compute the Fuss-Narayana numbers for all finite types (Figure \ref{fig:fussnar}) and observe several intriguing properties of these polynomials. This suggests several topics for further research.

In the final two sections \ref{sec:iterated} and \ref{sec:shelling}, we study more subtle structural and enumerative properties of $NC^{(k)}(W)$. We define a generalization $(NC^{(k)}(W))^{(\ell)}$ of the $k$-divisible noncrossing partitions, and show that this is isomorphic to the ``$kl$-divisible noncrossing partitions'' $NC^{(k\ell)}(W)$ (Theorem \ref{th:nckl=nclk}). The techniques used in this proof have a ``homological'' feel and may be interesting in themselves. In turn, we obtain an enumeration formula for multichains in $NC^{(k)}(W)$: the number of $\ell$-multichains in $NC^{(k)}(W)$ is given by the zeta polynomial
\begin{equation*}
\Z(NC^{(k)}(W),\ell)=\Cat^{(k\ell)}(W).
\end{equation*}
As a corollary, we find that the number of maximal chains in $NC^{(k)}(W)$ is equal to $n!(kh)^n/\abs{W}$, where $n$ is the rank and $h$ is the Coxeter number of $W$.

Finally, in joint work with Hugh Thomas, we construct a shelling of the order complex of $NC^{(k)}(W)$ (Theorem \ref{th:ELNCk}), based on the shelling of $NC(W)$ by Athanasiadis, Brady and Watt \cite{athanasiadis-brady-watt}. This implies that the order complex is homotopy equivalent to a wedge of $(n-1)$-dimensional spheres. We show that the number of these spheres is equal to $(-1)^n\Cat^{(-k)}(W)$.

Throughout this chapter there are suggestions for further research. One such problem is to study the action of a certain dihedral group on the homology groups of $NC^{(k)}(W)$ (Problem \ref{prob:homologyrep}).

\bigskip

\noindent{\bf Chapter \ref{sec:classical}:}\, 
Here we illustrate and motivate the material of Chapter \ref{sec:kdiv} by examining the classical types in detail. We give combinatorial realizations of the posets $NC^{(k)}(A_{n-1})$ and $NC^{(k)}(B_n)$ as posets of ``$k$-divisible'' set partitions. However, we are unable at this time to give a similar characterization of $NC^{(k)}(D_n)$.\footnote{This problem has now been solved by Krattenthaler \cite{kratt:typeD}.} This chapter contains a mix of new material and a survey of previous work.

The first two sections serve as an introduction to the classical theory of noncrossing partitions. Section \ref{sec:classicalNC} introduces the subject of noncrossing set partitions $NC(n)$ as initiated by Germain Kreweras in 1972 \cite{kreweras}. We describe the relationship between the lattice $NC(n)$ of noncrossing partitions under refinement and the Cayley graph $(A_{n-1},T)$ of the symmetric group $A_{n-1}$ with respect to transpositions $T$. At the end of the section we outline a ``modern history of noncrossing partitions'' from the mid 1990's to the present (mid 2000's). In Section \ref{sec:classicalkreweras}, we develop the theory of the ``Kreweras complement maps'', which are an important tool in the remainder of the chapter. These were introduced by Kreweras \cite{kreweras} and generalized by Nica and Speicher \cite{nica-speicher:ntuples}. We take this work to its logical conclusion by defining the fully general Kreweras maps.

In the second half of the chapter we describe new research. In Section \ref{sec:classicalNCk}, we apply the Kreweras complement maps to prove that our poset $NC^{(k)}(A_{n-1})$ is isomorphic to the poset $NC^{(k)}(n)$ of $k$-divisible noncrossing  partitions of $[kn]$, first studied by Edelman \cite{edelman:kdiv}. The main construction in this proof (that of ``shuffle partitions''), leads to several new bijective results on multichains in the lattice of noncrossing partitions. Then, in Section \ref{sec:typeA}, we translate a result of Edelman to count multichains in $NC^{(k)}(A_{n-1})$ by ``rank-jump vector'' (Theorem \ref{th:edelman}). We also generalize a result of Kreweras to count multichains in $NC^{(k)}(A_{n-1})$ whose bottom element has a prescribed ``parabolic type'' (Theorem \ref{th:kltype}). Either of these results can be used to obtain the Fuss-Narayana numbers for type $A$.

In Section \ref{sec:typeB} we extend all of our above type $A$ (symmetric group) results to type $B$ (the signed symmetric group), based on the lattice $\widetilde{NC}(2n)$ of type $B$ noncrossing partitions defined by Reiner \cite{reiner}. In this case, the poset $\widetilde{NC}\hspace{0in}^{(k)}(2n)\cong NC^{(k)}(B_n)$ of $k$-divisible type $B$ noncrossing partitions has not been considered previously. In the case that $k$ is even and $n$ is odd, the poset $\widetilde{NC}\hspace{0in}^{(k)}(2n)$ is not amenable to our current techniques. This is a mystery (Problem \ref{prob:typeBmystery}).

Finally, in Section \ref{sec:typeD} we compute the Fuss-Narayana numbers for type $D$, based on a result of Athanasiadis and Reiner \cite{athanasiadis-reiner}. As mentioned, we do not currently know a combinatorial realization of the poset $NC^{(k)}(D_n)$ in terms of $k$-divisible set partitions. We suggest this as an interesting problem\footnotemark[\value{footnote}] (Problem \ref{prob:typeD}).

\bigskip

\noindent{\bf Chapter \ref{sec:fusscatcomb}:}\, 
In the final chapter of the memoir we survey the other two subjects that --- together with our poset $NC^{(k)}(W)$ --- currently comprise the {\sf Fuss-Catalan combinatorics of finite Coxeter groups}. The first is the theory of ``nonnesting partitions'', whose Fuss-Catalan version is due to Athanasiadis \cite{athanasiadis:cat,athanasiadis:nar}; the second is the theory of cluster complexes, whose Fuss-Catalan version is due to Fomin and Reading \cite{fomin-reading}. While considering the three ``Fuss-Catalan families'', mysterious coincidences emerge, and we make several conjectures.\footnote{Since an early version of this memoir was circulated in 2005, some of our conjectures have now been proven by Krattenthaler \cite{krattenthaler1,krattenthaler2,kratt:typeD,kratt-muller} and Tzanaki \cite{tzanaki:faces}.}

In Section \ref{sec:nonnesting} we discuss the theory of nonnesting partitions as introduced by Postnikov (see \cite[Remark 2]{reiner}). We motivate the algebraic definition of nonnesting partitions as {\sf antichains in the root poset}, and show how this corresponds in type $A$ to the notion of ``nonnesting set partitions''. We show how the refinement order $NN(W)$ on nonnesting partitions and the lattice $NC(W)$ of noncrossing partitions may be directly compared as induced subposets of the lattice $\L(W)$ of parabolic subgroups of $W$.

Finally, we discuss the {\sf Shi hyperplane arrangement} and its relation to nonnesting partitions due to Cellini and Papi \cite{cellini-papi}. Athanasiadis \cite{athanasiadis:cat,athanasiadis:nar} has defined a class $NN^{(k)}(W)$ of {\sf geometric multichains of nonnesting partitions}, which he proves are in bijection with positive chambers in the {\sf extended Shi arrangement}. We define maps from $NN^{(k)}(W)$ and $NC^{(k)}(W)$ to the lattice $\L(W)$, and conjecture that the images are equidistributed in a nice way (Conjecture \ref{conj:NCk=NNk}).

We begin Section \ref{sec:cluster} by introducing the ideas of polygon triangulations and the {\sf  associahedron}. We discuss how the associahedron has been extended to all finite Coxeter groups by Fomin and Zelevinsky \cite{fomin-zelevinsky:ysystems} as part of their theory of {\sf cluster algebras}. Then we discuss general polygon dissections and the {\sf generalized cluster complex} $\Delta^{(k)}(W)$ of Fomin and Reading \cite{fomin-reading}. We observe that the $h$-vector of this complex coincides with our Fuss-Narayana numbers, which follows from explicit formulas computed by Fomin and Reading.

In Section \ref{sec:triangles}, we discuss the {\sf Chapoton triangles} introduced by Fr\'ed\'eric Chapoton \cite{chapoton:one,chapoton:two}. These are two-variable generating functions related to the three families of Coxeter-Catalan objects. We extend Chapoton's definitions to the three Fuss-Catalan families (Definition \ref{def:triangles}) and we propose that a conjecture of Chapoton from the $k=1$ case holds in general (Conjecture \ref{conj:triangles}). Then we define the {\sf dual triangles} and conjecture a formula for the dual $F$-triangle\footnote{This is now proved by Krattenthaler \cite{krattenthaler1}.} (Conjecture \ref{conj:dualF}).

Finally, in Section \ref{sec:future} we suggest some problems for future research and speculate on the future of this subject. The Fuss-Catalan combinatorics of finite Coxeter groups is a new topic and it offers many exciting open problems.

\bigskip

\begin{disclaimer}
In this memoir, we will deal exclusively with the case of finite Coxeter groups --- that is, finite real reflection groups. It may happen that some of the topics discussed can be extended meaningfully to other classes of groups. There are two promising directions of generalization. On one hand, the notion of noncrossing partitions has been extended to certain ``well-generated'' {\sf complex reflection groups} by Bessis \cite{bessis:complex}, and on the other hand, Brady, Crisp, Kaul and McCammond \cite{mccammond-etal}, as well as Thomas (personal communication), have considered noncrossing partitions related to a certain class of {\em infinite} Coxeter groups, the {\sf affine Weyl groups}. Many of our techniques and results may extend to the complex and infinite cases but we leave this for future consideration.
\end{disclaimer}


\chapter{Coxeter Groups and Noncrossing Partitions}
\label{sec:background}

Here we provide an introduction to the theory of finite reflection groups, followed by a development of the algebraic theory of noncrossing partitions. In our treatment of the basic theory of Coxeter groups, we follow Humphreys~\cite{humphreys}. Another excellent reference for reflection groups is Kane \cite{kane}. For the classical theory of reduced words, we refer to Bj\"orner and Brenti~\cite{bjorner-brenti}. For general poset theory, we refer to Stanley \cite{stanley:ec1}. Our introduction to the theory of dual Coxeter systems is based on results from Carter~\cite{carter}, Brady and Watt~\cite{brady-watt:kpione} and Bessis~\cite{bessis:dual}.

\section{Coxeter Systems}
\label{sec:coxsystems}
Suppose that a group $W$ is generated by a finite set $S$ with one relation of the form $(ss')^{m(s,s')}=1$ ($m(s,s')\in\{1,2,\ldots,\infty\}$) for each pair of generators $(s,s')\in S\times S$. In particular, when $m(s,s')=2$ the generators $s$ and $s'$ commute, and $m(s,s')=\infty$ means there is no relation between $s$ and $s'$. If the numbers $m(s,s')$ satisfy the properties
\begin{eqnarray*}
& m(s,s')=m(s',s)& \quad\text{and}\\
& m(s,s')=1 \Leftrightarrow s=s', &
\end{eqnarray*}
then we say that the symmetric array $m: S\times S\to \{1,2,\ldots,\infty\}$ is a {\sf Coxeter matrix}, and the group $W$ is a {\sf Coxeter group}. It turns out in this case that $m(s,s')$ is precisely the order of the element $ss'$ in the group. Notice that a Coxeter group is generated by involutions, since $m(s,s)=1$ for all $s\in S$. This is intended to model the property of being ``generated by reflections''.  The cardinality of $S$ is called the {\sf rank} of the group.

If there is a partition of the generators $S=S'\sqcup S''$ such that elements of $S'$ commute with elements of $S''$, then we have $W=\langle S'\rangle\times \langle S''\rangle$, where $\langle S'\rangle$ and $\langle S''\rangle$ are themselves Coxeter groups. In this case we say that $W$ is {\sf reducible}. Otherwise, it is {\sf irreducible}.

\begin{figure}[ht]
\vspace{.1in}
\scalebox{.9}{
\input{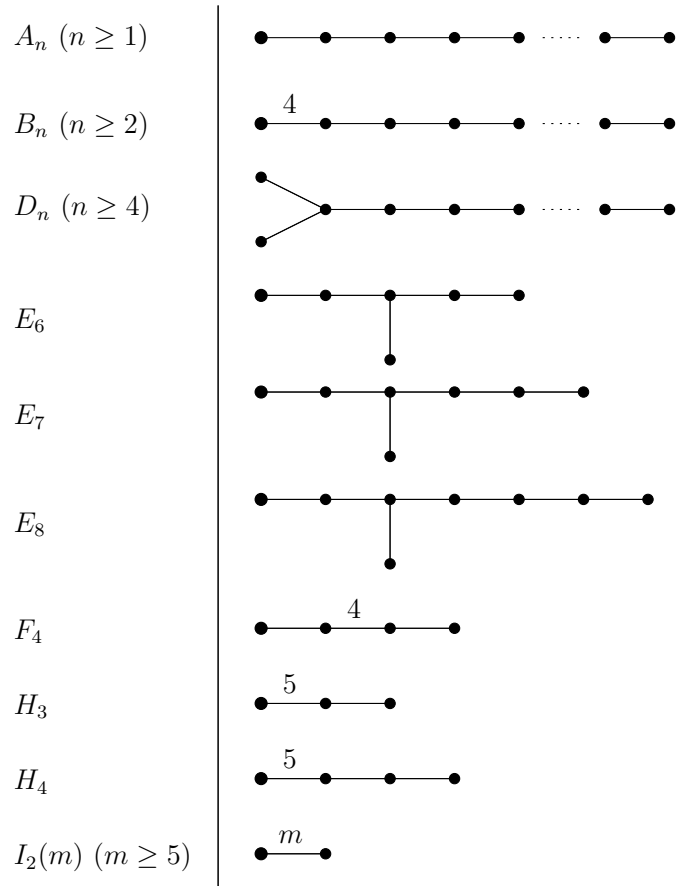}
}
\caption{Coxeter diagrams of the finite irreducible Coxeter groups}
\label{fig:coxdiagrams} 
\end{figure}

A convenient way to encode a Coxeter group is via its {\sf Coxeter diagram}, which is a graph with one vertex $v_s$ for each generator $s\in S$. If $m(s,s')\geq 3$ then we connect $v_s$ and $v_{s'}$ by an edge, and if $m(s,s')\geq 4$ then we label this edge by the number $m(s,s')$. Because two generators commute precisely when they are not connected by an edge, it is easy to read off the irreducible factors of $W$: they correspond to the connected components of the graph. Figure \ref{fig:coxdiagrams} displays the Coxeter diagrams of the finite irreducible Coxeter groups. It is interesting to note that the nonbranching graphs correspond to groups of symmetries of regular polytopes, and the Coxeter diagrams with edge labels in the set $\{4,6\}$ correspond to Weyl groups.

The term ``Coxeter group'' is standard and we will frequently use it, but we should mention that this notation is ambiguous. It is possible for a finite group $W$ to be a Coxeter group with respect to two different generating sets $S$ and $S'$ such that $S$ and $S'$ are not related by a group automorphism of $W$. For example, the dihedral group of order $12$ can be realized as a Coxeter group in two ways:
\begin{eqnarray*}
& \big\langle s,t : s^2=t^2=(st)^6=1\big\rangle, \\
& \big\langle a,b,c : a^2=b^2=c^2=(ab)^3=(ac)^2=(bc)^2=1\big\rangle .
\end{eqnarray*}
The first of these is irreducible of rank $2$ and it is denoted as $I_2(6)$ or $G_2$ (depending on whether we wish to emphasize the fact that it is a Weyl group --- see Section \ref{sec:rootsystems}). The second is the {\em reducible} Coxeter group $A_1\times A_2$, which has rank $3$. For this reason, it is more precise when discussing $W$ to refer to the {\sf Coxeter system} $(W,S)$, including both the group $W$ and the generating set $S$.

Of course, defining a Coxeter group in this way obscures the fact that it is an essentially geometric object. The abstract definition in terms of a presentation was inspired by the finite case, in which each Coxeter group has a faithful representation as a group generated by Euclidean reflections. In this case, we will see, the numbers $m(s,s')$ determine the dihedral angles between reflecting hyperplanes.

In general, we can construct a {\sf geometric representation} for each Coxeter system $(W,S)$, whether $W$ is finite or not, although the generators will not in general be Euclidean reflections. We say that a linear transformation is a {\sf linear reflection} if it fixes a hyperplane $H$ and sends some nonzero vector $\alpha$ to its negative (if we are in a Euclidean space and $H=\alpha^\perp$, this reduces to the usual notion of reflection). Let $V$ be the vector space over $\reals$ spanned by the abstract symbols $\{\alpha_s: s\in S\}$ and define a symmetric bilinear form $B$ on $V$ by setting
\begin{equation*}
B(\alpha_s,\alpha_{s'}):=-\cos\frac{\pi}{m(s,s')}
\end{equation*}
for all pairs of generators $(s,s')\in S\times S$, with the convention that $B(\alpha_s,\alpha_{s'})=-1$ if $m(s,s')=\infty$. (We note that the form $B$ has an important analogue in the theory of Lie algebras, called the {\sf Killing form}.) For each $s\in S$, we then define a linear reflection on $V$ by setting
\begin{equation*}
\sigma_s\lambda:=\lambda-2B(\alpha_s,\lambda)\alpha_s.
\end{equation*}
It turns out that the linear map defined on the generators by $s\mapsto\sigma_s$ is a faithful representation $\sigma:W\hookrightarrow GL(V)$ that preserves the bilinear form $B$. It is also {\sf essential}; that is, there is no nontrivial subspace of $V$ that is fixed pointwise by $W$ (see \cite[Sections 5.3 and 5.4]{humphreys}). Hence $W$ is isomorphic to a group of linear transformations of $V$ generated by the linear reflections $\sigma_s$. We call the elements $\{ \sigma_s :s\in S\}$ the {\sf simple reflections} of the geometric representation.

In this memoir, we are interested in the case when $W$ is a {\em finite} group. In this case the form $B$ will be positive definite, hence an inner product. A proof of the following theorem is given, for example, in Humphreys \cite[Section 6.4]{humphreys}.
\begin{theorem}
\label{theorem:finiteW}
The form $B$ is positive definite if and only if $W$ is finite.
\end{theorem}

The definiteness of the form, in turn, places strict constraints on the integers $m(s,s')$. This is the essential insight that allowed Cartan and Killing to classify the semisimple Lie algebras, and later allowed Coxeter to classify the finite real reflection groups \cite{coxeter} using precisely the diagrams in Figure \ref{fig:coxdiagrams}.

Thus for a finite Coxeter group $W$ we may identify the inner product space $(V,B)$ with Euclidean space $\reals^n$, where $n$ is the rank of $(W,S)$. Since the inner product $B$ is also $W$-invariant, the geometric representation $\sigma:W\hookrightarrow GL(V)$ is orthogonal and the generating linear reflections $\sigma_s$ become Euclidean reflections. It follows that the dihedral angle between the fixed hyperplanes of the simple reflections $\sigma_s$ and $\sigma_{s'}$ is equal to $\pi/m(s,s')$ for all $(s,s')\in S\times S$.

Moreover, the fact that $B$ is nondegenrate (which follows from positive definiteness)  implies that the geometric representation of an irreducible finite Coxeter system is irreducible as a representation (it has no invariant subspaces) \cite[Section 6.3]{humphreys}. In summary, each finite Coxeter system $(W,S)$ has an associated faithful, essential, orthogonal representation generated by reflections, and this representation is irreducible whenever the system is. Conversely, every finite reflection group arises in this way.

It was known since the mid-nineteenth century that a finite reflection group acting on $\reals^3$ must have angles between reflecting hyperplanes given by $\pi/m$, where $m$ is in $\{2,3,4,5\}$. However, the reasoning was purely geometric, since the abstract concept of a group did not crystallize until near the end of the century. M\"obius understood the situation in $\reals^3$ around 1852, and Schl\"afli had classified the regular polytopes and determined their symmetry groups by about 1850, although his work was ignored for a long time. M\"obius and Schl\"afli understood that these symmetry groups are generated by reflections, but the problem of determining all finite groups generated by reflections was not solved until much later. Goursat gave the answer for $\reals^4$ in 1889 and Cartan solved the general crystallographic case in his work on Lie algebras. Finally, Coxeter completed the classification in 1935, using some ideas of Witt (see Bourbaki \cite[Chapter 26]{bourbaki} and references therein).

\section{Root Systems}
\label{sec:rootsystems}
Historically, the main impetus for the classification of reflection groups came from Lie theory. The classification of finite Coxeter groups is closely related to the classification of semisimple Lie algebras, which was obtained by Cartan and Killing before 1894. Their classification depended on a combinatorial structure called a root system.

Let $\Phi$ be a finite spanning set in the Euclidean space $V$ with inner product $(\cdot,\cdot)$, and for each $\alpha\in\Phi$ let
\begin{equation*}
t_\alpha(\lambda)=\lambda-2\frac{(\lambda,\alpha)}{(\alpha,\alpha)}\alpha
\end{equation*}
denote the reflection in the hyperplane $\alpha^\perp$. We say that $\Phi$ is a {\sf (finite) root system} if the following two properties hold:
\begin{eqnarray}
\label{eq:RS1}
& \Phi\cap \reals\alpha=\{\alpha,-\alpha\}\quad\text{for all}\,\, \alpha\in\Phi &\quad\text{and}\\
\label{eq:RS2}
& t_\alpha\Phi=\Phi\quad \text{for all}\,\,\alpha\in\Phi.&
\end{eqnarray}
(Note that many authors use the term ``root system'' to refer to a more restrictive concept, which we will call a ``crystallographic root system''. See \ref{sec:crystallography} below.) It is easy to see that the group generated by the $t_\alpha$ is finite since it injects into the group of permutations of $\Phi$. We will denote this group by $W(\Phi)=\langle t_{\alpha}:\alpha\in\Phi\rangle$. Hence each root system gives rise to a finite reflection group. Moreover, we will see that every finite reflection group arises in this way and that the correspondence is unique up to the lengths of the $\alpha\in\Phi$.

Let $\sigma:W\hookrightarrow GL(V)$ be the geometric representation of a finite Coxeter system $(W,S)$ as in Section \ref{sec:coxsystems}, where $V$ is the real vector space spanned by the symbols $\{ \alpha_s:s\in S\}$. To simplify notation, we will write $w(\alpha_s)$ instead of $\sigma(w)(\alpha_s)$. Now let $\Phi(W)$ denote the set of images of the vectors $\alpha_s$ under the action of $W$:
\begin{equation*}
\Phi(W)  =  \left\{ w(\alpha_s): w\in W, s\in S\right\}.
\end{equation*}
Notice that $\Phi(W)$ consists of unit vectors, since $B(\alpha_s,\alpha_s)=-\cos(\pi/1)=1$ for all $s\in S$, and the action of $W$ preserves the inner product $B$. Thus, the first property of root systems \eqref{eq:RS1} is satisfied. To verify the second property \eqref{eq:RS2}, we use the following elementary fact.
\begin{lemma}
\label{lemma:reflections}
The reflection $t_{w(\alpha_s)}$ is given by $\sigma(wsw^{-1})$.
\end{lemma}
\begin{proof}
First note that $wsw^{-1}$ sends $w(\alpha_s)$ to its negative since
\begin{equation*}
wsw^{-1}(w(\alpha_s))=ws(\alpha_s)=w(-\alpha_s)=-w(\alpha_s).
\end{equation*}
Then we must show that $wsw^{-1}$ fixes $w(\alpha_s)^\perp$ pointwise. But if $\alpha$ is in $w(\alpha_s)^\perp$ then $w^{-1}(\alpha)$ is in $\alpha_s^\perp$ since $w^{-1}$ is an orthogonal transformation. Hence
\begin{equation*}
wsw^{-1}(\alpha)=ws(w^{-1}(\alpha))=w(w^{-1}(\alpha))=\alpha.
\end{equation*}
\end{proof}
This implies that each $t_{w(\alpha_s)}$ acts as a permutation on the set $\Phi(W)$, and we conclude that $\Phi(W)$ is indeed a root system. Furthermore, this lemma shows that  the reflection group generated by the root system is
\begin{equation*}
W(\Phi(W))=\big\langle \sigma(wsw^{-1}):w\in W, s\in S\big\rangle=\sigma(W),
\end{equation*}
so that $W(\Phi(W))\cong W$. If we define isomorphism of root systems up to orthogonal transformations and the lengths of the roots, it is also true that $\Phi(W(\Phi))\cong \Phi$. This sets up a bijection between reflection groups and isomorphism classes of root systems.

\begin{figure}
\vspace{.1in}
\begin{center}
\input{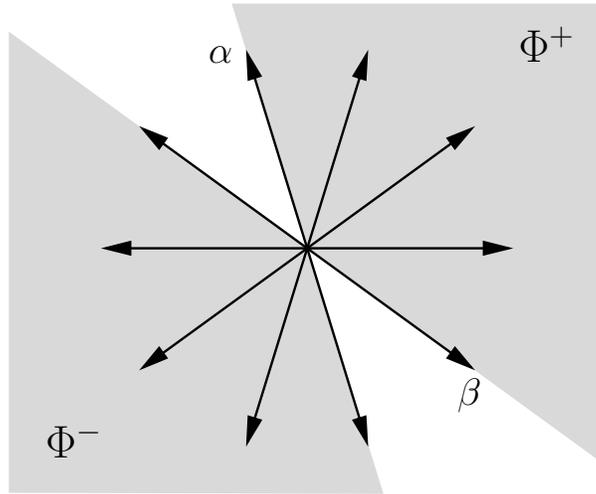}
\end{center}
\caption{The root system $\Phi(I_2(5))$}
\label{fig:rootsystem}
\end{figure}

The structure of a root system is described in terms of ``positive systems'' and ``simple systems''. Let $\Phi$ be a root system in the Euclidean space $V$. Any hyperplane in $V$ that does not intersect $\Phi$ partitions $\Phi$ into two sets,
\begin{equation*}
\Phi=\Phi^+\sqcup\Phi^-,
\end{equation*}
which we call a {\sf positive system} $\Phi^+$ and a {\sf negative system} $\Phi^-=-\Phi^+$ for $\Phi$. The entire root system is contained in the cone generated by $\Phi^+$, which consists of the {\sf positive cone} (the positive span of $\Phi^+$) and the {\sf negative cone} (the negative span of $\Phi^+$). Let $\Pi$ denote the set of roots generating the extremal rays of the positive cone. It turns out that $\Pi$ is a {\sf simple system} for $\Phi$, in the sense that every root $\alpha\in\Phi$ can be expressed as a linear combination from $\Pi$ in which the coefficients are all nonnegative or all nonpositive (obvious), and that $\Pi$ is a vector space basis for $V$ (not obvious). In general, simple systems and positive systems uniquely determine each other \cite[Section 1.3]{humphreys}. Figure \ref{fig:rootsystem} displays the root system $\Phi(I_2(5))$ for the dihedral group of order $10$, with positive and negative cones shaded. The corresponding simple system is $\Pi=\{\alpha,\beta\}$.

If we think of $\Phi$ together with its Coxeter group $W=W(\Phi)$, the simple systems of $\Phi$ have a very natural interpretation: they correspond to the Coxeter generating sets of $W$. Consider the geometric representation of the finite Coxeter system $(W,S)$ and let $\{\alpha_s : s\in S\}\subseteq \Phi(W)$ be called the {\sf simple roots}. Then we have the following result \cite[Section 5.4]{humphreys}.
\begin{theorem}
The simple roots $\{\alpha_s:s\in S\}$ are a simple system for $\Phi(W)$.
\end{theorem}

\subsection{Addendum --- Crystallography}
\label{sec:crystallography}

We should note that our use of the term ``root system'' is slightly different than that usually encountered in Lie theory. The main definitions are the following.

We say that a discrete additive subgroup $L$ of a real vector space $V$ is a {\sf lattice}. If $L$ is isomorphic to the free $\integers$-module $\integers^n$, where $n$ is the dimension of $V$, we say that $L$ has {\sf full rank}.

\begin{definition}
A finite real reflection group $W$ acting on $V$ is called a {\sf Weyl group} if it stabilizes a full rank lattice $L\subseteq V$. In general, a group that stabilizes a lattice is called {\sf crystallographic}.
\end{definition}

\begin{definition}
We say that a root system $\Phi\subseteq V$ is {\sf crystallographic} if it satisfies conditions \eqref{eq:RS1} and \eqref{eq:RS2} above, as well as the following:
\begin{equation}
\label{eq:RS3}
2\frac{(\alpha,\beta)}{(\alpha,\alpha)}\in\integers\quad\text{for all}\,\,\alpha,\beta\in\Phi.
\end{equation}
\end{definition}
In contrast to the noncrystallographic case, let $\Pi\subseteq\Phi$ be a simple system in a crystallographic root system $\Phi$. Then every root can be expressed uniquely as a $\integers$-linear (not just $\reals$-linear) combination of simple roots. (As before the coefficients are all nonnegative or all nonpositive.)

We will see that Weyl groups are precisely those finite reflection groups that arise from crystallographic root systems. First, suppose that $W(\Phi)$ is the reflection group generated by a crystallographic root system $\Phi\subseteq V$. If we define the set
\begin{equation*}
Q=\integers\Phi=\left\{ \sum n_i \alpha_i : n_i\in\integers \text{ and } \alpha_i\in\Phi \text{ for all } i\right\},
\end{equation*}
then condition \eqref{eq:RS3} and the fact that $\Phi$ spans $V$ guarantee that $Q\subseteq V$ is a full rank lattice, called the {\sf root lattice} of $\Phi$.\footnote{The same construction for noncrystallographic $\Phi$ yields a {\em dense} --- in particular, not discrete --- additive subgroup of $V$. This is the topic of {\sf quasicrystals} \cite{quasicrystal}.} Condition \eqref{eq:RS2} tells us that $Q$ is $W$-equivariant, hence $W$ is a Weyl group.

Conversely, let $W$ be a finite reflection group with its geometric representation $\sigma:W\hookrightarrow V$ and suppose that $W$ stabilizes a full rank lattice $L\subseteq V$. If we compute the trace of an element $\sigma(w)$ with respect to a $\integers$-linear basis for $L$, we find that it must be an integer for all $w\in W$, and this translates to the necessary condition that  $m(s,s')$ is in the set $\{2,3,4,6\}$ for all $(s,s')\in S\times S'$, $s\neq s'$.

Now we have a slight problem. If $\Phi=\Phi(W)$ is the root system constructed above from $\sigma$, then $\Phi$ is rarely crystallographic, even when $W$ is known to be a Weyl group. However, since $W$ satisfies the property that $m(s,s')\in\{2,3,4,6\}$ for all $(s,s')\in S\times S$, $s\neq s'$, it is possible to modify the lengths of the roots $\alpha_s$ to create a crystallographic root system $\Phi'$ that is isomorphic to $\Phi(W)$. First, set $\lambda_s=c_s\alpha_s$ for some scalars $c_s\in\reals $, $s\in S$. Suppose we are able to choose the scalars $c_s$ so that the following properties hold:
\begin{eqnarray*}
m(s,s')=3 & \Rightarrow & c_s=c_{s'},\\
m(s,s')=4 & \Rightarrow & c_s=\sqrt{2}c_{s'}\,\,\text{or}\,\, \sqrt{2}c_s=c_{s'},\\
m(s,s')=6 & \Rightarrow & c_s=\sqrt{3}c_{s'}\,\,\text{or}\,\, \sqrt{3}c_s=c_{s'}.
\end{eqnarray*}
Then we will have $\sigma_s(\lambda_{s'})=\lambda_{s'}+d(s,s')\lambda_s$ where the $d(s,s')$ are integers, so that $\{\lambda_s:s\in S\}$ generates a crystallographic root system $\Phi'$ satisfying $W=W(\Phi')$. The only difficulty here is to choose the scalars $c_s$ in a consistent way, but since the Coxeter diagram of a finite reflection group is a forest, this is always possible. Note, however, that the choice of scalars might not be unique. This is the case for the Weyl group of type $B_n$ which has two nonisomorphic crystallographic root systems, called $B_n$ and $C_n$. (The stronger notion of isomorphism for noncrystallographic root systems does not allow scaling of roots.)

Finally, note that the root system $\Phi(I_2(5))$ in Figure \ref{fig:rootsystem} can {\em not} be made crystallographic by changing root lengths, since the Coxeter diagram of $I_2(5)$ contains the label $5\not\in\{2,3,4,6\}$. The fact that $I_2(5)$ is noncrystallographic is equivalent to the fact that the Euclidean plane can not be tiled by regular pentagons.

We will usually not work with crystallographic root systems or Weyl groups, but the notions will be important in Chapter \ref{sec:fusscatcomb}, especially Section \ref{sec:nonnesting}.

\section{Reduced Words and Weak Order}
\label{sec:reducedwords}
The idea of a Coxeter system is very important in modern mathematics. Part of the beauty of this subject comes from the fact that it shows up frequently in unexpected places, and lies at the intersection of algebra, geometry, and combinatorics. We described above some of the algebraic and geometric aspects of the theory. From this point on, we will be concerned primarily with the combinatorial side of Coxeter systems, which is based on the study of reduced words.

Let $(W,S)$ be a (possibly infinite) Coxeter system with finite generating set $S$, and consider the word length $\ell_S:W\to\integers$ on $W$ with respect to $S$. That is, $\ell_S(w)$ is the minimum integer $r$ such that there exists an expression $w=s_1s_2\cdots s_r$ with $s_1,s_2,\ldots,s_r\in S$. We call such a minimal expression $s_1s_2\cdots s_r$ a {\sf reduced S-word} for $w$, and we refer to $\ell_S$ as the {\sf standard length} on $(W,S)$.

We have already seen two different ways to define a Coxeter system: in terms of its Coxeter presentation, or as a group generated by reflections with an associated root system. Much of the theory of Coxeter systems can also be expressed in the language of reduced $S$-words. For example, define the Exchange Property and the Deletion Property as follows. Here $\hat{s}$ denotes the fact that an occurrence of the symbol $s$ has been deleted from a word.
\begin{exchangeproperty}
Let $w=s_1s_2\cdots s_r$ be a reduced $S$-word and consider $s\in S$. If $\ell_S(sw)<\ell_S(w)$ then $sw=s_1\cdots \hat{s_i}\cdots s_r$ for some $1\leq i\leq r$.
\end{exchangeproperty}
\begin{deletionproperty}
If $w=s_1s_2\cdots s_r$ and $\ell_S(w)< r$, then we have \linebreak
$w=s_1\cdots\hat{s_i}\cdots\hat{s_j}\cdots s_r$ for some $1\leq i<j\leq r$.
\end{deletionproperty}
It turns out that either of these properties is sufficient to characterize a Coxeter system. This result appears as Theorem 1.5.1 in~\cite{bjorner-brenti}.
\begin{theorem}
If $W$ is a group with a generating set $S$ of involutions, the following are equivalent:
\begin{enumerate}
\item $(W,S)$ is a Coxeter system.
\item $(W,S)$ satisfies the Exchange Property.
\item $(W,S)$ satisfies the Deletion Property.
\end{enumerate}
\end{theorem}
A powerful way to encode information about reduced words is to consider $W$ as a partially ordered set, or a {\sf poset}. The word length $\ell_S$ naturally induces a partial order on $W$ in the following way. For all $\pi,\mu$ in $W$, we have the triangle inequality $\ell_S(\mu)\leq \ell_S(\pi)+\ell_S(\pi^{-1}\mu)$. Whenever this inequality is an {\em equality}, we set $\pi\leq_S\mu$.
\begin{definition}
Define the {\sf (right) weak order} on $W$ by setting
\begin{equation*}
\pi\leq_S\mu\quad\Longleftrightarrow\quad \ell_S(\mu)=\ell_S(\pi)+\ell_S(\pi^{-1}\mu)
\end{equation*}
for all $\pi,\mu$ in $W$, and denote this poset by $\Weak(W)$.
\end{definition}
It is easy to verify that this relation satisfies the reflexive, transitive and antisymmetric properties of a partial order, and that the identity $1\in W$ satisfies $1\leq_S w$ for all $w\in W$. A poset $P$ is called {\sf graded} if there exists a rank function $\rk:P\to\integers$ with the property that every unrefinable chain $x=z_0\leq z_2\leq\cdots\leq z_r=y$ between $x,y\in P$ has the same length $r=\rk(y)-\rk(x)$. It is easy to see that the weak order is graded with rank function $\ell_S$.

In general, $\Weak(W)$ is also a {\sf meet-semilattice} in the sense that every pair of elements $\pi,\mu\in W$ has a greatest lower bound, or a {\sf meet} $\pi\wedge\mu\in W$. The existence of meets can be proved using the Exchange Property, and has many consequences for the structure of $W$ (see \cite[Chapter 3]{bjorner-brenti}). If $W$ is finite, it is also true that $\Weak(W)$ has a maximum element called $w_\circ$, and by a basic property of posets \cite[Proposition 3.3.1]{stanley:ec1} this implies that every pair $\pi,\mu\in W$ also has a least upper bound, or {\sf join} $\pi\vee\mu\in W$. In this case, we say that $\Weak(W)$ is a {\sf lattice}\footnote{We hope no confusion will arise with the other use of the term ``lattice'' as a discrete additive subgroup of $\reals^n$.}. 

Note that the weak order also has a natural interpretation in terms of reduced $S$-words: by definition, we have $\pi\leq_S\mu$ if and only if there exists a reduced $S$-word for $\mu$ such that $\pi$ occurs as a {\em prefix}. That is, there exists a reduced $S$-word $\mu=s_1 s_2\cdots s_r$ for $\mu$, such that $\pi=s_1s_2\cdots s_k$ for some $k\leq r$. Of course, we could also define a {\em left} weak order on $W$ in terms of {\em suffixes} of reduced $S$-words. The left and right weak orders do not coincide but they are isomorphic via the map $w\mapsto w^{-1}$ (since this map reverses reduced words, it switches prefixes with suffixes). Hence, we will usually not make a distinction between the left and right weak orders, referring to $\Weak(W)$, simply, as {\em the} weak order on $W$.

\begin{figure}
\vspace{.1in}
\begin{center}
\input{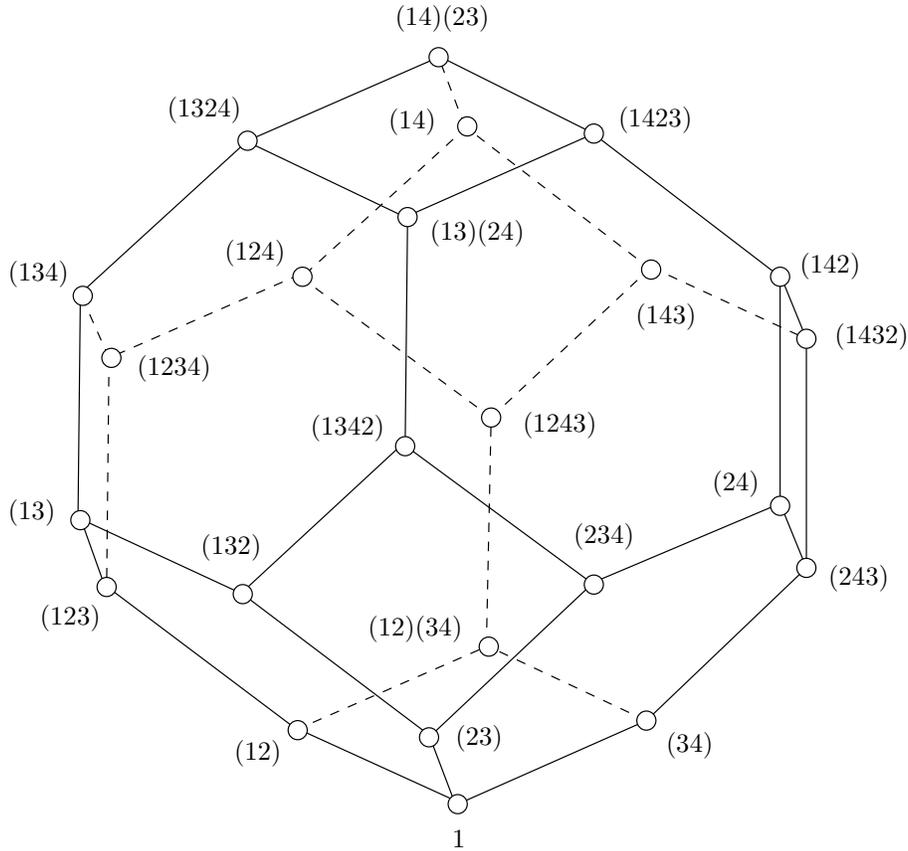}
\end{center}
\caption{The Hasse diagram of $\Weak(A_3)$; $1$-skeleton of the type $A_3$ permutohedron; and Cayley graph of $A_3$ with respect to the adjacent transpositions $S=\{(12),(23),(34)\}$}
\label{fig:weak(a3)}
\end{figure}

When $(W,S)$ is a finite Coxeter system, there is also a nice geometric interpretation of the poset $\Weak(W)$. Consider the {\sf Coxeter arrangement} of $(W,S)$, which is the set $\A$ of reflecting hyperplanes for the geometric representation $\sigma:W\hookrightarrow GL(V)$. The connected components of the complement $V\setminus(\cup_{H\in\A}  H)$ are called {\sf chambers}, and each of these chambers is a fundamental domain for the action of $W$. Each of the hyperplanes in $\A$ has a positive side, corresponding to the direction of its positive root, and the intersection of these positive half-spaces is called the {\sf fundamental chamber}, or the {\sf positive chamber}.

Now if we select a generic point $v$ in the interior of the fundamental chamber, its $W$-orbit consists of $\abs{W}$ points, one in each chamber, and they are naturally in bijection with the elements of $W$. The convex hull of these points is called the {\sf permutohedron} of type $W$, and the $1$-skeleton of the permutohedron is isomorphic to the Hasse diagram of the weak order $\Weak(W)$. The {\sf Hasse diagram} is an oriented graph on the vertex set of poset elements. We draw an edge from $\pi$ to $\mu$ in $W$ if $\pi\leq_S\mu$ and there does not exist $\varphi$ such that $\pi\leq_S\varphi\leq_S\mu$ (in this case, we say that $\mu$ {\sf covers} $\pi$), and we illustrate the orientation on edges by directing them upward in the diagram.

We explain this situation with an example.

\begin{example}
Consider the Coxeter system of type $A_3$. The corresponding group is isomorphic to the symmetric group $\frak{S}_4$, and it acts as the symmetries of a regular tetrahedron in $\reals^3$ by permuting vertices. If we label the vertices of the tetrahedron by the symbols $\{1,2,3,4\}$, then one choice of Coxeter generating set for $A_3$ is the adjacent transpositions $S=\{(12),(23),(34)\}$. In this case, the longest element is $w_\circ=(14)(23)$ with length $\ell_S(w_\circ)=6$. Figure \ref{fig:weak(a3)} displays the Hasse diagram of $\Weak(A_3)$, distorted to emphasize the fact that it is isomorphic as a graph to the $1$-skeleton of the permutohedron\footnote{An easy way to construct the permutohedron of type $A_{n-1}$ is to take as vertices the point $(1,2,\ldots,n)\in\reals^n$ and all permutations of its coordinates. Dividing out by the invariant line on which the coordinates sum to $0$ yields a convex polytope in $\reals^{n-1}$.}.
\end{example}

It is also worth mentioning that the Hasse diagram of $\Weak(W)$ coincides with the {\sf Cayley graph} of $W$ with respect to $S$ (that is, where we connect two vertices $\pi,\mu\in W$ by an edge if $\mu=\pi s$ for some $s\in S$). Indeed the only way this could fail is if two elements with the same standard length are connected by an edge in the Cayley graph. But consider the {\sf sign representation} $w\mapsto \det(\sigma(w))$. Since all of the generators $s\in S$ are reflections, we have $\det(\sigma(s_1s_2\cdots s_k))=(-1)^k$ for all $s_1,s_2,\ldots,s_k\in S$. Hence if $\mu=\pi s$ for some $s\in S$, then $\pi$ and $\mu$ must differ in length.

We remarked above that the left and right weak orders are isomorphic, but that they do not coincide. A common way to remove the ``sidedness'' of the weak order is to define $\pi\leq_B\mu$ whenever $\pi$ is an {\em arbitrary} subword (not necessarily a prefix) of a reduced $S$-word for $\mu$. The resulting order is called the {\sf Bruhat order} (or {\sf strong order}) on $W$ and it {\sf extends} both the left and right weak orders. That is, $\pi\leq_S \mu$ in either the left or right weak order implies $\pi\leq_B \mu$. This order arises as the inclusion order on closures of Bruhat cells for the corresponding semisimple Lie group, and we will not discuss it further here.

The combinatorics of reduced $S$-words and associated partial orders on $W$ has been extensively studied. For more on this beautiful topic, see the text of Bj\"orner and Brenti~\cite{bjorner-brenti}. 

\section{Absolute Order}
\label{sec:absolute}
Now we move away from the standard theory of Coxeter groups to describe more recent work on algebraic noncrossing partitions. The essential new idea here is that we substitute for the Coxeter generators $S$ a larger generating set.

Consider a finite Coxeter system $(W,S)$ with root system $\Phi(W)$ and geometric representation $\sigma:W\hookrightarrow GL(V)$. In Lemma \ref{lemma:reflections}, we showed that the reflections orthogonal to the roots all have the form $\sigma(wsw^{-1})$ for some $s\in S$ and $w\in W$. That is, not only are the group elements $wsw^{-1}\in W$ involutions, but they also act as reflections. This is certainly not true of every involution in $W$: for example, the involution $\sigma((16)(25)(34))$ in the geometric representation of $A_5$ has determinant $-1$, but its eigenvalues are $1,1,-1,-1,-1$, so it is {\em not} a reflection. (Recall that a reflection is an orthogonal transformation with all but one eigenvalue equal to $1$, and the last equal to $-1$.) It turns out that {\em every} reflection in the geometric representation has the form $\sigma(wsw^{-1})$ \cite[Section 1.14]{humphreys}, and this inspires the following definition.

\begin{definition}
The conjugate closure of the Coxeter generating set,
\begin{equation*}
T:=\left\{ wsw^{-1}: s\in S, w\in W\right\},
\end{equation*}
is called the generating set of {\sf reflections}.
\end{definition}

In general, we will call the elements of $T$ the ``reflections'' of $W$, even when we are not explicitly considering the geometric representation. Consider the motivating example of the symmetric group.

\begin{example}
For the symmetric group $\frak{S}_n$ on $n$ letters (the Coxeter group of type $A_{n-1}$) the standard choice of Coxeter generating set $S$ is the set of adjacent transpositions; hence $T$ is the generating set of {\em all} transpositions. We take this to motivate the general use of the letter $T$.
\end{example}

The systematic study of a Coxeter group $W$ with respect to the generators $T$ is fairly recent. The earliest reference we can find is Carter \cite{carter}, from 1972. In 2001, Bessis defined the following terminology \cite[Version 1]{bessis:dual}.

\begin{definition}
\label{def:dualCoxeter}
If $W$ is a finite Coxeter group with set of reflections $T$, we call the pair $(W,T)$ a {\sf dual Coxeter system}.
\end{definition}

Note that a dual Coxeter system is {\em not} a Coxeter system. For instance, if we naively try to form a ``Coxeter diagram'' for the symmetric group $A_2$ with respect to the transpositions $T=\{(12),(23),(13)\}$ it will be a triangle; but this diagram does not appear in Figure \ref{fig:coxdiagrams}, contradicting the fact that $A_2$ is finite. That is, there are more relations among the elements of $T$ than just the pairwise ones, and the presentation of $W$ with respect to $T$ is more complicated than the standard Coxeter presentation. However, the structure of the pair $(W,T)$ is very rich, and the notion of a dual Coxeter system has some advantages over the usual notion of a ``classical'' Coxeter system. One immediate advantage is the fact that a finite Coxeter system has a unique set of reflections $T$, but many equivalent choices of simple reflections $S$. The fact that $T$ is closed under conjugation determines the essential character of this theory.

Focusing now on $T$, we can mimic the classical theory of reduced $S$-words. Let $\ell_T:W\to\integers$ denote the word length on $W$ with respect to the generating set $T$. We will call this the {\sf absolute length} on $W$. Clearly we have $\ell_T(w)\leq\ell_S(w)$ for all $w\in W$ since $S\subseteq T$ (indeed, $\ell_T$ has occasionally been known as the ``short length'' on $W$).  If $w=t_1t_2\cdots t_r$ with $t_1,t_2,\ldots,t_r\in T$ and $\ell_T(w)=r$, then we call $t_1t_2\cdots t_r$ a {\sf reduced $T$-word} for $w$. Note that the absolute length naturally induces a partial order on $W$, just as the standard length induces the weak order on $W$.

\begin{figure}
\vspace{.1in}
\begin{center}
\input{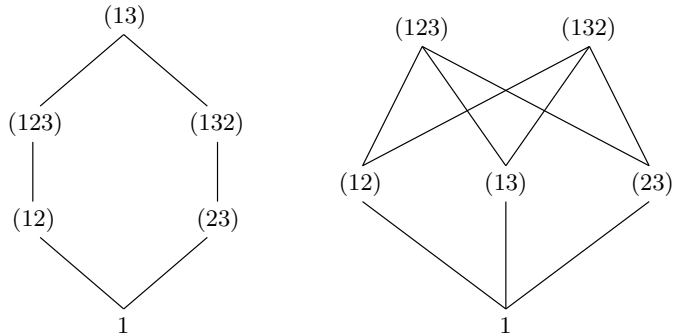}
\end{center}
\caption{Weak order versus absolute order on $A_2$}
\label{fig:weakvsabsolute}
\end{figure}

\begin{definition}
\label{def:absorder}
Define the {\sf absolute order} on $W$ by setting
\begin{equation*}
\pi\leq_T\mu\quad\Longleftrightarrow\quad \ell_T(\mu)=\ell_T(\pi)+\ell_T(\pi^{-1}\mu)
\end{equation*}
for all $\pi$, $\mu$ in $W$, and denote this poset by $\Abs(W)$.
\end{definition}
Again, it is straightforward to see that $\Abs(W)$ is a graded poset with rank function $\ell_T$, and the identity $1\in W$ is the unique minimum element. However, in constrast with the weak order, the absolute order does {\em not} in general have a maximum element, even when $W$ is finite. For example,  Figure \ref{fig:weakvsabsolute} compares the Hasse diagrams of $\Weak(A_2)$ and $\Abs(A_2)$. Notice that both of the $3$-cycles in $A_2$ are maximal elements of $\Abs(A_2)$, and that the Hasse diagram is isomorphic to the Cayley graph with respect to the generating set $T=\{(12),(13),(23)\}$. In general, the Hasse diagram of $\Abs(W)$ is isomorphic to the Cayley graph of $W$ with respect to $T$, again because each of the generators has determinant $-1$.

Like the weak order, the absolute order also has an important geometric interpretation. For each root $\alpha\in\Phi(W)$, let $t_{\alpha}$ again denote the reflection in the hyperplane $\alpha^\perp$ orthogonal to $\alpha$. The essential lemma about the geometry of reduced $T$-words is the following standard but nontrivial result, proved by Carter in his 1972 paper~\cite{carter}. Actually, he stated the result only in the case of Weyl groups, but his argument works in general. Since the result is central to our topic and since it has such intrinsic interest, we feel that it deserves a name.

\begin{carterslemma}
\label{lemma:carter}
The $T$-word $w=t_{\alpha_1}t_{\alpha_2}\cdots t_{\alpha_r}$ for $w\in W$ is reduced if and only if the set $\{\alpha_1,\alpha_2,\ldots,\alpha_r\}$ is linearly independent.
\end{carterslemma}

This has two immediate consequences: first, a reflection may not be repeated in a reduced $T$-word; and second, we have $\ell_T(w)\leq n$ for all $w\in W$, where $n$ is the rank of the reflection group. Note in Figure \ref{fig:weakvsabsolute} that every element of $A_2$ has absolute length less than or equal to $2$.

Now, given $w\in W$, consider the transformation $\sigma(w)-\sigma(1)$ acting on $V$, which we will denote simply by $w-1$. The kernel of $w-1$ is the {\sf fixed space} of $w$ (the space fixed pointwise by $\sigma(w)$) and the image of $w-1$ is called the {\sf moved space}.
\begin{definition}
For all $w\in W$, set
\begin{equation*}
\Fix(w)=\ker(w-1)\quad\text{and}\quad \Mov(w)=\im(w-1).
\end{equation*}
\end{definition}
An elementary fact in linear algebra tells us that $\im(w-1)$ and the kernel of the adjoint $\ker(w-1)^*$ are orthogonal complements in $V$. But $w$ is an orthogonal transformation, so we have $(w-1)^*=w^*-1=w^{-1}-1$. Since the kernels of $w-1$ and $w^{-1}-1$ are clearly the same, it follows that $\Mov(w)=\Fix(w)^\perp$ for all $w\in W$.

It turns out that the geometry of the absolute order is best expressed in terms of moved spaces. The next theorem was known in an equivalent form to Carter, and was first expressed in the language of moved spaces by Brady and Watt~\cite{brady-watt:kpione}.

\begin{theorem}
\label{th:basicMov}
For all $\pi,\mu$ in $W$, we have
\begin{enumerate}
\item $\ell_T(\pi)=\dim \Mov(\pi)$.
\item $\pi\leq_T\mu\Rightarrow \Mov(\pi)\subseteq \Mov(\mu)$.
\item If $t\in T$, then $\Mov(t)\subseteq \Mov(\pi)\Rightarrow t\leq_T\pi$.
\label{basicMov3}
\end{enumerate}
\end{theorem}

\begin{proof} Suppose that $\pi=t_{\alpha_1}t_{\alpha_2}\cdots t_{\alpha_r}$ is a reduced $T$-word for $\pi\in W$. To prove $(1)$, we will show that $\{\alpha_1,\alpha_2,\ldots,\alpha_r\}$ is a basis for $\Mov(\pi)$. By Carter's Lemma, the set $\{\alpha_1,\alpha_2,\ldots,\alpha_r\}$ is linearly independent, so that the intersection of hyperplanes $\cap_{i=1}^r \alpha_i^\perp$ has codimension $r$. But certainly this intersection is contained in the fixed space of $w$, so that $\dim\Fix(\pi)\geq n-r$, or $\dim\Mov(\pi)\leq r$. To finish the proof, we will show that $\{\alpha_1,\alpha_2,\ldots,\alpha_r\}$ is contained in the moved space $\Mov(\pi)=\im(\pi-1)$.

To show that $\alpha_1$ is in $\im(\pi-1)$, consider a vector $x$ in $(\cap_{i=2}^r \alpha_i^\perp) \setminus \alpha_1^\perp$, which exists by linear independence. Then $\pi(x)$ is just equal to $t_{\alpha_1}(x)=x-2\frac{(x,\alpha_1)}{(\alpha_1,\alpha_1)}\alpha_1$, and so $c\alpha_1=\pi(x)-x$ for some constant $c\neq 0$. Hence $\alpha_1$ is in $\im(\pi-1)$. Now assume by induduction that $\{\alpha_1,\alpha_2,\ldots,\alpha_{k-1}\}$ is contained in $\im(\pi-1)$ for some $2\leq k\leq r$. By linear independence, there exists a vector $y\in (\cap_{i=k+1}^r \alpha_i^\perp)\setminus \alpha_k^\perp$ (where we understand $\cap_{i=r+1}^r \alpha_i^\perp$ to be the whole space $V$), which then satisfies
\begin{equation*}
\pi(y)=t_{\alpha_1}t_{\alpha_2}\cdots t_{\alpha_k}(y)=y+c_1\alpha_1+c_2\alpha_2+\cdots +c_k\alpha_k
\end{equation*}
for some scalars $c_1,c_2,\ldots, c_k$ with $c_k\neq 0$. It follows that $c_k\alpha_k=\pi(y)-y -\sum_{i=1}^{k-1} c_i\alpha_i$, so that $\alpha_k$ is in $\im(\pi-1)$ and this completes the proof of $(1)$.

To prove $(2)$, suppose that $\pi\leq_T\mu$ for some $\pi$, $\mu$ in $W$. That is, there exists a reduced $T$-word $\mu=t_{\alpha_1}t_{\alpha_2}\cdots t_{\alpha_r}$ for $\mu$ and $1\leq k\leq r$ so that the prefix $t_{\alpha_1}t_{\alpha_2}\cdots t_{\alpha_k}$ (necessarily reduced) is a $T$-word for $\pi$. In this case, since $\{\alpha_1,\ldots,\alpha_k\}$ is a basis for $\Mov(\pi)$ and $\{\alpha_1,\ldots,\alpha_r\}$ is a basis for $\Mov(\mu)$, it follows that $\Mov(\pi)\subseteq\Mov(\mu)$.

Finally, to prove $(3)$, consider $t_\alpha\in T$ and $\pi\in W$, with $\Mov(t_\alpha)\subseteq\Mov(\pi)$. By the triangle inequality, we have $\ell_T(\pi)\leq\ell_T(t_\alpha)+\ell_T(t_\alpha\pi)=1+\ell_T(t_\alpha\pi)$, or $\ell_T(t_\alpha\pi)\geq\ell_T(\pi)-1$. To show $t_\alpha\leq_T\pi$, then, it is sufficient to prove that $\ell_T(t_\alpha)\leq\ell_T(\pi)-1$.

So consider a reduced $T$-word $\pi=t_{\alpha_1}t_{\alpha_2}\cdots t_{\alpha_r}$ for $\pi$. Note that $\Mov(t_\alpha)$ is just the line spanned by $\alpha$, and since $\Mov(t_\alpha)$ is contained in $\Mov(\pi)$, which has basis $\{\alpha_1,\ldots,\alpha_r\}$, it follows that the set $\{\alpha\}\cup\{\alpha_1,\ldots,\alpha_r\}$ is not linearly independent. By Carter's Lemma this implies that the $T$-word $t_\alpha \pi=t_\alpha t_{\alpha_1}\cdots t_{\alpha_r}$ is not reduced. But since the generators $T$ all have determinant $-1$, the length of any two $T$-words for $t_\alpha\pi$ must differ by a multiple of $2$. That is, $\ell_T(t_\alpha\pi)\leq (r+1)-2=\ell_T(\pi)-1$.

\end{proof}

We might hope that the map $w\mapsto \Mov(w)$ is a poset inclusion into the lattice of subspaces of $V$ under containment. Unfortunately, this is not the case. While this map does preserve order, it is not in general injective, and it is not true that $\Mov(\pi)\subseteq\Mov(\pi)$ implies $\pi\leq_T\mu$. For example, consider the two $3$-cycles in $A_2$. It is easy to show that each of these has moved space $V$, but they are not comparable in $\Abs(A_2)$ (see Figure \ref{fig:weakvsabsolute}).

The essential obstacle here is the fact that, unlike the weak order, the absolute order  does not have a unique maximum element. If we remove this difficulty, however, by restricting to an interval of $\Abs(W)$, then the map $w\mapsto\Mov(w)$ becomes a {\sf poset isomorphism}. That is, it is injective, and both directions preserve order. This property was first explained by Brady and Watt in~\cite{brady-watt:kpione}, and it is based on the following general result which they proved in the earlier~\cite{brady-watt:partialorder}. 

\begin{lemma}[\cite{brady-watt:partialorder}]
\label{lemma:brady-watt}
Let $w$ be any element of the orthogonal group $O(V)$, and suppose that $M$ is any subspace of the moved space $\Mov(w)=\im(w-1)$. Then there exists a unique orthogonal transformation $\pi\in O(V)$ with the properties $M=\Mov(\pi)$ and
\begin{equation*}
\dim\Mov(w)=\dim\Mov(\pi)+\dim\Mov(\pi^{-1}w).
\end{equation*}
\end{lemma}

\begin{proof}[Proof Sketch]
The proof relies on the fact that $w-1$ acts invertibly on the moved space $\Mov(w)=\im(w-1)$. Certainly $\im(w-1)$ is closed under the action of $w-1$, and we have already shown that the kernel $\Fix(w)=\ker(w-1)$ is orthogonal to $\Mov(w)$, hence their intersection is trivial. In this case, let $U$ denote the unique subspace of $\Mov(w)$ with the property that $(w-1)U=M$.

Now define a map $\pi$ by setting $\pi=w$ on $U$ and $\pi=1$ on $\Fix(w)$ and extend linearly. In a series of lemmas, Brady and Watt show that this map is well-defined (that is, $V=\Fix(w)\oplus U$), and it is the unique map with the desired properties.
\end{proof}

Actually, Brady and Watt used a slightly cleaner language by defining a partial order on the orthogonal group $O(V)$, setting $A\leq_O B$ whenever
\begin{equation*}
\dim\Mov(B)=\dim\Mov(A)+\dim\Mov(A^{-1}B).
\end{equation*}
By Theorem \ref{th:basicMov} $(1)$, this partial order agrees with the absolute order $\leq_T$ when restricted to the discrete subgroup $W$ of $O(V)$.

\begin{theorem}[\cite{brady-watt:kpione}]
\label{theorem:hardMov}
Consider $\pi$ and $\mu$ in $W$, and suppose that there exists $w\in W$ with $\pi\leq_T w$ and $\mu\leq_T w$. In this case, we have
\begin{equation*}
\pi\leq_T\mu\quad\Longleftrightarrow\quad\Mov(\pi)\subseteq\Mov(\mu).
\end{equation*}
\end{theorem}

\begin{proof}
We already know that $\pi\leq_T\mu$ implies $\Mov(\pi)\subseteq\Mov(\mu)$. So suppose $\Mov(\pi)\subseteq\Mov(\mu)$ and let $A\in O(V)$ denote the unique orthogonal transformation with the properties $\Mov(A)=\Mov(\pi)$ and $A\leq_O\mu$, guaranteed by Lemma \ref{lemma:brady-watt}. If we can show that $A=\pi$ then we are done since $\pi\leq_O\mu$ is the same as $\pi\leq_T\mu$. Now, since $\mu\leq_T w$, we have $A\leq_O\mu\leq_O w$ which implies $A\leq_O w$ by transitivity; and since $\pi\leq_T w$, we have $\Mov(\pi)\subseteq\Mov(w)$. That is, $A$ satisties the properties of Lemma \ref{lemma:brady-watt} with respect to $w$. But $\pi$ also satisfies these properties, hence $A=\pi$ by uniqueness.
\end{proof}

Thus the intervals of $\Abs(W)$ have a nice geometric interpretation: each is isomorphic via the map $w\mapsto\Mov(w)$ to a poset of subspaces of $V$ under inclusion. (The fact that this map is injective follows directly from the antisymmetric property of partial order.) Below we will study the intervals of $\Abs(W)$, showing that they are always self-dual. Beyond that, we will isolate a natural subclass of intervals that have further nice properties, including the lattice property (see Section \ref{sec:NC}). These are the so called ``lattices of noncrossing partitions'', and they will be the main focus of this memoir.

Although it seems very natural to consider a Coxeter group $W$ with respect to its generating set of reflections $T$, the idea is more recent than one might expect. It seems that Carter was the first to systematically study the group $W$ from this point of view in the early 1970's, for example in~\cite{carter}, and the idea of a dual Coxeter system is very recent, first considered explicitly in 2001 by Bessis~\cite[Version 1]{bessis:dual}, where he developed many properties of dual Coxeter systems in parallel with the classical theory. Certainly there will be more to say in this direction. For example, a classification of dual Coxeter systems in terms of some sort of ``dual'' Exchange and Deletion properties (see Section \ref{sec:reducedwords}) is not yet known.

\section{Shifting and Local Self-Duality}
\label{sec:shifting}
We mentioned above that the weak order on $W$ has an essentially ``sided'' nature: one may speak of the left weak order or the right weak order. The fact that the generating set $T$ is closed under conjugation, however, means that the absolute order is immune from these sorts of difficulties. In this section, we will see that the fact that $T$ is closed under conjugation has many nice consequences for the structure of the absolute order.

To begin with, consider $w$ in $W$ with $\ell_T(w)=r$ and let $\sigma$ be an arbitrary element of $W$. If $w=t_1 t_2\cdots t_r$ is a reduced $T$-word for $w$, then since $T$ is closed under conjugation it follows that
\begin{equation*}
\sigma w\sigma^{-1}=(\sigma t_1\sigma^{-1})(\sigma t_2\sigma^{-1})\cdots (\sigma t_r \sigma^{-1})
\end{equation*}
is a $T$-word for $\sigma w\sigma^{-1}$. Now, if there were a shorter $T$-word for $\sigma w\sigma^{-1}$, then, conjugating by $\sigma^{-1}$, we would obtain a $T$-word for $w$ with length less than $r$, contradicting the fact that $\ell_T(w)=r$. It follows that $\ell_T(\sigma w\sigma^{-1})=r$, and we conclude that the absolute length is invariant under conjugation. 

It is also easy to see that the poset structure of the absolute order is invariant under conjugation. Indeed, the map $w\mapsto\sigma w\sigma^{-1}$ is bijective, with inverse $w\mapsto \sigma^{-1}w\sigma$. Then suppose that $\pi\leq_T\mu$, or $\ell_T(\mu)=\ell_T(\pi)+\ell_T(\pi^{-1}\mu)$. Since $\ell_T$ is invariant under conjugation, this is the same as 
\begin{equation*}
\ell_T(\sigma\mu\sigma^{-1})=\ell_T(\sigma\pi\sigma^{-1})+\ell_T(\sigma\pi^{-1}\mu\sigma^{-1}),
\end{equation*}
which is equivalent to $\sigma\pi\sigma^{-1}\leq_T\sigma\mu\sigma^{-1}$. We say that $w\mapsto\sigma w\sigma^{-1}$ is a {\sf poset automorphism} of $\Abs(W)$.

A nice way to encapsulate the properties of reduced $T$-words is the following lemma, which we will use frequently in our study of $\Abs(W)$.

\begin{shiftinglemma}
\label{shifting}
If $w=t_1 t_2\cdots t_r$ is a reduced $T$-word for $w\in W$ and $1<i<r$, then the two expressions
\begin{eqnarray}
\label{eq:leftshift}
w= t_1 t_2\cdots t_{i-2} t_i (t_i t_{i-1} t_i) t_{i+1}\cdots t_r & \text{and}\\
\label{eq:rightshift}
w= t_1 t_2\cdots t_{i-1}(t_i t_{i+1} t_i) t_i t_{i+2} \cdots t_r &
\end{eqnarray}
are also reduced $T$-words for $w$.
\end{shiftinglemma}

\begin{proof} Both of these are $T$-words for $w$ and they have length $r=\ell_T(w)$, so they are reduced.
\end{proof}

In other words, if a reflection $t$ occurs in the $i$-th place of some reduced $T$-word for $w$ with $1<i<r=\ell_T(w)$, then there exist two other reduced $T$-words for $w$ in which $t$ appears in the $(i-1)$-th place and the $(i+1)$-th place, respectively. That is, we may effectively ``shift'' any symbol $t$ in a reduced $T$-word to the left or to the right, at the expense of conjugating the shifted-over symbol by $t$. 

We observed above that reduced $T$-words contain no repetition, as a consequence of Carter's Lemma. We can now obtain the same result in a purely algebraic way using the Shifting Lemma. Indeed, if a reduced $T$-word contains two occurrences of the symbol $t$, then the two occurrences may be shifted until they are adjacent, at which point they will cancel, contradicting the minimality of length. Contrast this with the case of reduced $S$-words which {\em may} contain repetition.

The Shifting Lemma also implies that there is no difference between the ``left'' absolute order and the ``right'' absolute order, since a $T$-word for $\pi\in W$ occurs as a prefix of a reduced $T$-word for $\mu\in W$ if and only if it occurs as a suffix of some other reduced $T$-word for $\mu$ (obtained by shifting). More generally, we get a characterization of the absolute order in terms of arbitrary subwords.

\begin{subwordproperty}
\label{subword}
We have $\pi\leq_T\mu$ in $\Abs(W)$ if and only if $\pi$ occurs as an arbitrary subword of some reduced $T$-word for $\mu$.
\end{subwordproperty}

\begin{proof}
If $\pi\leq_T\mu$, then $\ell_T(\mu)=\ell_T(\pi)+\ell_T(\pi^{-1}\mu)$ implies the existence of a reduced $T$-word for $\mu$ in which $\pi$ occurs as a prefix. Conversely, suppose that $\mu=t_1 t_2\cdots t_r$ is a reduced $T$-word and that $\pi=t_{i_1} t_{i_2}\cdots t_{i_k}$ for some $1\leq i_1<i_2<\cdots <i_k\leq r$. By repeatedly shifting the symbols $t_{i_1},t_{i_2},\ldots,t_{i_k}$ to the left, we obtain a reduced $T$-word for $\mu$ of the form
\begin{equation*}
\mu=t_{i_1}t_{i_2}\cdots t_{i_k} t'_{k+1}\cdots t'_r,
\end{equation*}
and this implies that $\pi\leq_T\mu$.
\end{proof}

Warning: Given $\pi\leq_T\mu$ and a reduced $T$-word for $\mu$, note that it is {\em not} generally true that $\pi$ occurs as a subword; the lemma merely states that there {\em exists} some reduced $T$-word for $\mu$ containing $\pi$ as a subword. We note that the absolute order shares some features in common with both the weak order and the Bruhat order on $W$. The definition of $\Abs(W)$ in terms of absolute length $\ell_T$ exactly mimics the definition of $\Weak(W)$ in terms of standard length $\ell_S$. However, the subword characterization of the absolute order is more similar to the subword characterization of the Bruhat order. With the generating set $T$, there is no need to distinguish between the two perspectives.

Now we will study the structure of the intervals in $\Abs(W)$. Recall that a poset $P$ is called {\sf locally self-dual} if, for each interval $[x,y]=\{z\in P: x\leq z\leq y\}$, there exists a bijection $[x,y]\to [x,y]$ that is order-reversing. Such a map is called an {\sf anti-automorphism} of the poset $[x,y]$. The technique of shifting furnishes us with a family of anti-automorphisms that demonstrate the locally self-duality of the absolute order.

\begin{definition}
\label{def:algkreweras}
For all pairs of elements $(\mu,\nu)\in W\times W$, define the map $K_\mu^\nu:W\to W$ by setting $K_\mu^\nu(w):=\mu w^{-1}\nu$.
\end{definition}

Clearly the map $K_\mu^\nu$ acts as a permutation on $W$ for all $\mu$ and $\nu$. In the special case that $\mu\leq_T\nu$, it also behaves well with respect to absolute order.

\begin{lemma}
\label{lemma:kreweras}
Consider $\mu\leq_T\nu$ in $\Abs(W)$. Then the following hold.
\begin{enumerate}
\item $K_\mu^\nu$ is an anti-automorphism of the interval $[\mu,\nu]\subseteq\Abs(W)$.
\label{krewitem:one}
\item For all $\pi\in[\mu,\nu]$ we have
\begin{equation*}
\ell_T(K_\mu^\nu(\pi))=\ell_T(\mu)+\ell_T(\nu)-\ell_T(\pi).
\end{equation*}
\label{krewitem:two}
\end{enumerate}
\end{lemma}

\begin{figure}
\vspace{.1in}
\begin{center}
\input{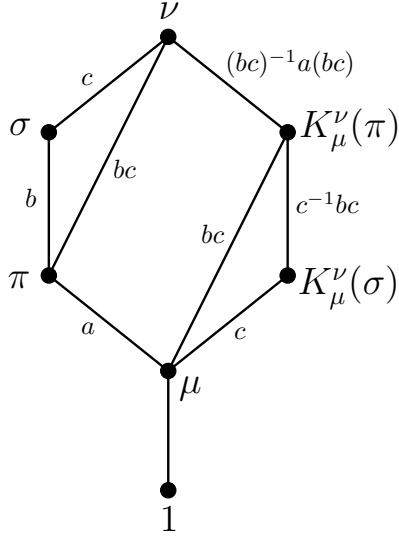}
\end{center}
\caption{$K_\mu^\nu$ is an anti-automorphism of the interval $[\mu,\nu]\subseteq\Abs(W)$}
\label{fig:orderreversing}
\end{figure}

\begin{proof}
To prove \eqref{krewitem:one}, note that the map $K_\mu^\nu$ is invertible, with $(K_\mu^\nu)^{-1}=K_\nu^\mu$. If we can show that $K_\mu^\nu$ takes $[\mu,\nu]$ to itself and that it reverses order on this interval, then the fact that the inverse reverses order will follow. Indeed, given $\pi\leq_T\sigma$ in $[\mu,\nu]$, we know that $(K_\mu^\nu)^{-1}(\pi)$ and $(K_\mu^\nu)^{-1}(\sigma)$ are comparable because $K_\mu^\nu$ has finite order and hence $(K_\mu^\nu)^{-1}=(K_\mu^\nu)^r$ for some $r\geq 1$. If we had $(K_\mu^\nu)^{-1}(\pi)\leq_T(K_\mu^\nu)^{-1}(\sigma)$, then applying the order-reversing $K_\mu^\nu$ on both sides leads to the contradiction $\sigma\leq_T\pi$. Hence we must have $(K_\mu^\nu)^{-1}(\sigma)\leq_T(K_\mu^\nu)^{-1}(\pi)$.

Now consider $\mu\leq_T\pi\leq_T\sigma\leq_T\nu$, and set $a=\mu^{-1}\pi$, $b=\pi^{-1}\sigma$ and $c=\sigma^{-1}\nu$. By assumption, we have $\nu=\mu abc$ with $\ell_T(\nu)=\ell_T(\mu)+\ell_T(a)+\ell_T(b)+\ell_T(c)$. But it is also true that $\nu=\mu c(c^{-1}b c)((bc)^{-1}a(bc))$, and
\begin{equation*}
\ell_T(\nu)=\ell_T(\mu)+\ell_T(c)+\ell_T(c^{-1}bc)+\ell_T((bc)^{-1}a(bc)),
\end{equation*}
since $\ell_T$ is invariant under conjugation. Following definitions, and applying the Subword Property, this last fact is equivalent to the relations
\begin{equation*}
\mu\leq_T K_\mu^\nu(\sigma)\leq_T K_\mu^\nu(\pi)\leq_T\nu,
\end{equation*}
which proves the result. Figure \ref{fig:orderreversing} illustrates this situation, with lines representing reduced $T$-words.

Now recall that $\Abs(W)$ is a graded poset with rank function $\ell_T$. Because $K_\mu^\nu$ reverses order on $[\mu,\nu]$, the two elements $\pi$ and $K_\mu^\nu(\pi)$ must have complementary rank within $[\mu,\nu]$, proving \eqref{krewitem:two}.
\end{proof}

This lemma shows that the absolute order $\Abs(W)$ is a locally self-dual poset since each interval $[\mu,\nu]$ has a self-duality given by $K_\mu^\nu$. In addition, the following easy corollary allows us to compare the maps $K_\mu^\nu$ on {\em different} intervals. We will need this later.

\begin{corollary}
\label{cor:kreweras}
For all $\mu\leq_T\pi\leq_T\sigma\leq_T\nu$, we have
$K_\mu^\sigma(\pi)\leq_T K_\mu^\nu(\pi)$.
\end{corollary}

\begin{proof}
First note that  $K_\mu^\nu(\pi)=K_\mu^\sigma(\pi)\,K_1^\nu(\sigma)$ in the group $W$. Then Lemma \ref{lemma:kreweras} $(2)$ implies that
\begin{equation*}
\ell_T(K_\mu^\nu(\pi))=\ell_T(K_\mu^\sigma(\pi))+\ell_T(K_1^\nu(\sigma)),
\end{equation*}
and we conclude by Definition \ref{def:absorder} that $K_\mu^\sigma(\pi)\leq_T K_\mu^\nu(\pi)$.
\end{proof}

We call $K_\mu^\nu$ the {\sf Kreweras complement} on $[\mu,\nu]$ since, in the case of the symmetric group, it is equivalent to a combinatorial construction of Kreweras (see Section \ref{sec:classicalkreweras}). Directly below, we will examine a class of intervals $[\mu,\nu]$ in $\Abs(W)$ that have the lattice property. In this case, the notation ``complement'' is motivated since the map $K_\mu^\nu$ turns out to be a {\sf lattice complement}. That is, for all $\pi$ in $[\mu,\nu]$ we have $\pi\vee K_\mu^\nu(\pi)=\nu$ and $\pi\wedge K_\mu^\nu(\pi)=\mu$. The Kreweras complement is our most important tool for understanding the structure of the absolute order. It will play an essential role throughout Chapters \ref{sec:kdiv} and \ref{sec:classical}.

\section{Coxeter Elements and Noncrossing Partitions}
\label{sec:NC}
We noted above that the absolute order $\Abs(W)$ does not in general have a maximum element. It is natural then to consider the set of maximal elements in more detail. It turns out that the maximal elements do not play equivalent roles; there is a very special class among them.

Suppose that $(W,S)$ is a finite Coxeter system and denumerate the simple generators by $S=\{s_1,s_2,\ldots,s_n\}$.

\begin{definition}
A {\sf standard Coxeter element} is any element of the form
\begin{equation*}
c=s_{\sigma(1)}s_{\sigma(2)}\cdots s_{\sigma(n)},
\end{equation*}
where $\sigma$ is some permutation of the set $\{1,2,\ldots,n\}$. A {\sf Coxeter element} is any conjugate of a standard Coxeter element in $W$.
\end{definition}

(Note that many authors use the term ``Coxeter element'' to refer to what we have called a ``standard Coxeter element''.) The corresponding ``Coxeter transformations'' in the geometric representation were first used by Cartan and Killing around 1890 in their study of Lie algebras, but they did not consider them as elements of a group; that insight had to wait for Weyl (see \cite[Chapter 26]{bourbaki}). The elements are named for Coxeter who showed in 1951 that they have a remarkable role~\cite{coxeter2} to play in the invariant theory of $W$, as described below in Section \ref{sec:invariant}. The Coxeter elements are also the motivating examples of Springer's regular elements~\cite{springer}. We say an element of a finite reflection group is {\sf regular} if it possesses an eigenvector that lies in none of the reflecting hyperplanes. This is an important geometric perspective, but we will not use it here.

Again consider the example of the symmetric group $A_{n-1}$ on $n$ letters, generated by the set of adjacent transpositions $S$. In this case, note that any standard Coxeter element is an $n$-cycle, so that the Coxeter elements form the conjugacy class of $n$-cycles. This is an important motivating example to keep in mind. However, it can be misleading: {\em In the case of the symmetric group, the Coxeter elements are precisely the maximal elements of the absolute order; in all other cases, the Coxeter elements are a proper subclass of maximal elements}.

The following basic properties of Coxeter elements follow from Humphreys \cite[Section 3.16]{humphreys}.

\begin{lemma} 
\label{lemma:basicCox}
Let $(W,S)$ be a finite Coxeter system with reflections $T$ and geometric representation $\sigma:W\hookrightarrow GL(V)$.
\begin{enumerate}
\item Any two standard Coxeter elements are conjugate. Hence the Coxeter elements form a single conjugacy class in $W$.
\item The moved space of any Coxeter element is $V$.
\item If $c\in W$ is a Coxeter element then we have $t\leq_T c$ for all $t\in T$.
\label{basicCox3}
\end{enumerate}
\end{lemma}

\begin{proof}
Property $(1)$ is Proposition 3.16 in Humphreys. If $c\in W$ is a standard Coxeter element, Lemma 3.16 in Humphreys implies that $\Fix(c)=\{0\}$. Then for any $w\in W$ we have $\Fix(wcw^{-1})=\sigma(w)\Fix(c)=\{0\}$ and  $\Mov(wcw^{-1})=\Fix(wcw^{-1})^\perp$, proving $(2)$. Property $(3)$ follows from our Theorem \ref{th:basicMov} $(3)$ and the fact that $\Mov(t)\subseteq V$ for all $t\in T$.
\end{proof}

Since any two simple generating sets for $W$ are conjugate, Lemma \ref{lemma:basicCox} $(1)$ shows that the set of Coxeter elements is independent of the choice of simple generating set used to define them. Furthermore, notice that Lemma \ref{lemma:basicCox} $(2)$ and Theorem \ref{th:basicMov} imply that the Coxeter elements are among the maximal elements of the absolute order, as we claimed above.

In his study of dual Coxeter systems~\cite{bessis:dual} Bessis showed that the Coxeter elements play a fundamental role in the structure of the absolute order. To describe his results, we need the notion of a ``parabolic subgroup''. Part of the importance of these subgroups is that they allow results about finite Coxeter groups to be proved by induction. Let $(W,S)$ be a finite Coxeter system.

\begin{definition}
For any subset $I\subseteq S$ of the simple generators, define $W_I:=\langle I\rangle$ to be the {\sf standard parabolic subgroup} generated by $I$. More generally, for any $w\in W$ and $I\subseteq S$, we call $wW_I w^{-1}$ a {\sf parabolic subgroup} of $W$. If $c$ is a Coxeter element in $wW_Iw^{-1}$, we call it a {\sf parabolic Coxeter element} of $W$.
\end{definition}

Given $I\subseteq S$, it is an elementary fact that the pair $(W_I,I)$ is itself a Coxeter system and its geometric representation agrees with the geometric representation of $(W,S)$, restricted to a subspace $V_I\subseteq V$ \cite[Sections 1.10 and 5.5]{humphreys}. Thus, for $w\in W$, the pair $(wW_Iw^{-1},wIw^{-1})$ is also a Coxeter system with geometric representation in $wV_I\subseteq V$. We observe that the absolute length on $wW_Iw^{-1}$ agrees with the absolute length on $W$, restricted to $wW_Iw^{-1}$.

\begin{proposition}
\label{prop:restriction}
Let $T'$ denote the set of reflections for the Coxeter system $(wW_Iw^{-1},wIw^{-1})$. Then we have $\ell_T(\pi)=\ell_{T'}(\pi)$ for all $\pi\in wW_Iw^{-1}$.
\end{proposition}

\begin{proof}
By Theorem \ref{th:basicMov}, we have $\ell_T(\pi)=\dim\Mov(\pi)$. But the moved space of $\pi$ in the geometric representation of $(W,S)$ coincides with the moved space in the geometric representation of $(wW_Iw^{-1},wIw^{-1})$, hence $\ell_T(\pi)=\ell_{T'}(\pi)$.
\end{proof}

Now we will present two important theorems from Bessis~\cite{bessis:dual}. The first is a combinatorial lemma \cite[Proposition 1.6.1]{bessis:dual} that is ``dual'' to a result of Tits on classical Coxeter systems. For $w$ in $W$, let $R_T(w)$ denote the set of reduced $T$-words for $w$.

\begin{dualtitslemma}
Given $w\in W$, the set $R_T(w)$ is connected under shifts of the form \eqref{eq:leftshift} and \eqref{eq:rightshift}  if and only if $w$ is a parabolic Coxeter element.
\end{dualtitslemma}

The classical Tits' Lemma is stated in terms of reduced $S$-words, and it says that the set $R_S(w)$ of reduced $S$-words for a given $w\in W$ is connected under {\sf braid moves}. A braid move is just an application of one of the defining Coxeter relations. For this reason, Bessis referred to moves of the form \eqref{eq:leftshift} and \eqref{eq:rightshift} as {\sf dual braid moves}. Tits' use of the lemma was to solve the word problem in $W$ with respect to $S$; Bessis' dual lemma is part of an efficient solution for the word problem with respect to $T$.\footnote{Bessis originally referred to this lemma as the Dual Matsumoto Property, but he now suggests the name Dual Tits' Lemma (personal communication).} The fact that the dual lemma holds only for parabolic Coxeter elements highlights the significance of these elements.

The second result \cite[Lemma 1.4.3]{bessis:dual} characterizes the parabolic Coxeter elements in terms of the absolute order.

\begin{theorem}
\label{th:coxelement}
An element $w$ in $W$ is a parabolic Coxeter element if and only if there exists some Coxeter element $c\in W$ with $w\leq_T c$.
\end{theorem}

This says that the set of parabolic Coxeter elements in $W$ forms an {\sf order ideal} in $\Abs(W)$. That is, if $\mu$ is a parabolic Coxeter element and $\pi\leq_T\mu$, then $\pi$ is also a parabolic Coxeter element. We now define the object that will occupy our attention for the main part of this memoir.
\begin{definition}
\label{def:NC(W)}
Relative to a Coxeter element $c$ in $W$, define the poset of {\sf noncrossing partitions} 
\begin{equation*}
NC(W,c):=[1,c] =\{w\in W: 1\leq_T w\leq_T c\}
\end{equation*}
as the interval in $\Abs(W)$ between $1$ and $c$.
\end{definition}

\begin{figure}
\vspace{.1in}
\begin{center}
\input{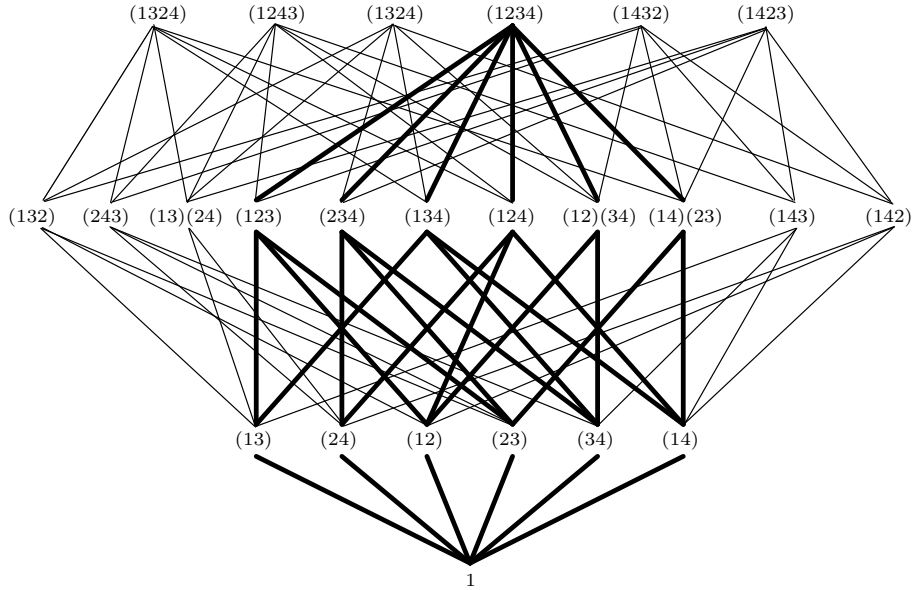}
\end{center}
\caption{$NC(A_3)$ as an interval in $\Abs(A_3)$}
\label{fig:nc(a3)}
\end{figure}

Some comments: The definition of the noncrossing partitions $NC(W,c)$ appears to depend on the Coxeter element $c$. However, the fact that Coxeter elements form a conjugacy class (Lemma \ref{lemma:basicCox} $(1)$), together with the fact that conjugation by a fixed group element $w \in W$ is an automorphism of $\Abs(W)$, implies that $NC(W,c)\cong NC(W,c')$ for all Coxeter elements $c,c'\in W$. Thus the isomorphism type of $NC(W,c)$ is independent of $c$, and we will often write $NC(W)$ when the Coxeter element is understood.

The notation ``noncrossing partitions'' for an interval of the form $[1,c]$ may cause some initial confusion, since it is not a priori apparent that this structure has anything whatsoever to do with partitions, let alone the notion of ``crossing''. The notation comes from the field of algebraic combinatorics, in which a lattice of ``noncrossing'' set partitions has been well-studied for some time. This structure was introduced by Kreweras in 1972 \cite{kreweras}, and it turns out that his poset is isomorphic to $NC(A_{n-1})$. See the Introduction or Chapter \ref{sec:classical} for a thorough discussion of this topic.

\begin{example}
Figure \ref{fig:nc(a3)} displays the Hasse diagram of $\Abs(A_3)$, which is isomorphic to the Cayley graph of $(A_3,T)$. We have highlighted the interval $NC(A_3)$ with respect to the Coxeter element $c=(1234)$. Observe that this coincides with the lattice of noncrossing set partitions from Figure \ref{fig:noncrossing}. Since the maximal elements of $\Abs(A_3)$ are all Coxeter elements, all maximal intervals in $\Abs(A_3)$ are isomorphic to $NC(A_3)$. Compare this with the Hasse diagram of $\Weak(A_3)$ in Figure \ref{fig:weak(a3)}, which has height $6$ and $3$ {\sf atoms} (elements covering the minimum element), whereas $NC(A_3)$ has height $3$ and $6$ atoms.
\end{example}

Let us consider some basic properties of the poset $NC(W)$. As well as being an interval in the Cayley graph of $(W,T)$, there is an important geometric ``representation'' of the poset. This follows immediately from Theorem \ref{theorem:hardMov}.

\begin{proposition}
\label{prop:Movembedding}
The map $\pi\mapsto \Mov(\pi)$ is a poset embedding of $NC(W)$ into the lattice of subspaces of $V$ under inclusion.
\end{proposition}

Recall from Lemma \ref{lemma:kreweras} that the map $K_1^c$ is an anti-automorphism of the interval $NC(W,c)$. Hence the poset is self-dual.

\begin{notation}
We call $K_1^c$ {\em the} Kreweras complement on $NC(W,c)$. When it will cause no confusion, we will tend to drop the subscript $1$ or the superscript $c$, and write the Kreweras complement simply as $K:NC(W)\to NC(W)$.
\end{notation}

In Section \ref{sec:classicalkreweras}, we will discuss a pictorial interpretation of this map.

Now, subintervals of the poset $NC(W)$ also have an important interpretation, following from Bessis' result, Theorem \ref{th:coxelement}.

\begin{proposition}
\label{prop:subintervals}
Every interval in $NC(W)$ is isomorphic to $NC(W')$, where $W'$ is a parabolic subgroup of $W$.
\end{proposition}

\begin{proof}
Given $\pi\leq_T\mu\leq_T c$, so that $\pi$ and $\mu$ are in $NC(W,c)$, notice that the composition of Kreweras complements $K_1^{\pi^{-1}\mu}\circ K_1^\mu$ restricts to an isomorphism of intervals $[\pi,\mu]\cong [1,\pi^{-1}\mu]$. Then if $W'$ is a parabolic subgroup in which $\pi^{-1}\mu$ is a Coxeter element (whose existence is guaranteed by Theorem \ref{th:coxelement}) and since absolute length restricts to parabolic subgroups (Proposition \ref{prop:restriction}), we have
\begin{equation*}
[\pi,\mu]\cong [1,\pi^{-1}\mu] \cong NC(W').
\end{equation*}
\end{proof}

This is an important observation that allows the use of induction to prove results about $NC(W)$. The next theorem was surprisingly difficult to establish in a uniform way. The first uniform proof was recently given in 2005 by Brady and Watt \cite[Theorem 7.8]{brady-watt:lattice}.
\begin{theorem}
The poset $NC(W)$ is a lattice.
\end{theorem}

This property was suspected from the beginning, and it is central to the definition of $NC(W)$ as a Garside structure \cite{bessis:dual,brady-watt:kpione}. In the classical types $A$, $B$ and $D$, the lattice property follows easily from known combinatorial realizations (see Chapter \ref{sec:classical}). Bessis first verified the exceptional types by computer.

Finally, the fact that $NC(W)$ is a lattice allows us to prove two results relating the lattice structure to the group structure of $NC(W)$. Let $\wedge$ and $\vee$ denote the meet (greatest lower bound) and join (least upper bound) in $NC(W)$, respectively. The first result is a lemma that we will need later.

\begin{lemma}\hspace{.1in}
\label{lemma:joins}
\begin{enumerate}
\item Suppose that $w=t_1 t_2 \cdots t_r$ is a reduced $T$-word for $w\in NC(W)$. Then we have $w=t_1\vee t_2\vee\cdots\vee t_r$.
\item If we have $\pi$, $\mu$ and $\pi\mu$ in $NC(W)$, with either $\pi\leq_T\pi\mu$ or $\mu\leq_T \pi\mu$, it follows that $\pi\mu=\pi\vee\mu$.
\end{enumerate}
\end{lemma}

\begin{proof}
To prove $(1)$, suppose that $w'$ and $w't$ are in $NC(W)$ with $\ell_T(w't)=\ell_T(w')+1$, for some $w'\in W$ and $t\in T$. We will show in this case that $w't=w'\vee t$ and the result follows by induction.

By assumption we have $w'\leq_T w't$, and by the Subword Property we have $t\leq_T w't$, hence $w'\vee t\leq_T w't$ since join $\vee$ is the least upper bound. We will be done if we can show that $\ell_T(w'\vee t)\geq\ell_T(w't)=\ell_T(w')+1$.

By Theorem \ref{th:basicMov} $(2)$, notice that $\Mov(w')\cup\Mov(t)\subseteq\Mov(w'\vee t)$, since $w'$ and $t$ are both below $w'\vee t$. If $t$ were less than $w'$, we would have a reduced $T$-word for $w'$ containing the symbol $t$. Then $w't$ would contain repetition and could not be reduced. This contradition implies that $w'$ and $t$ are incomparable, and by Proposition \ref{prop:Movembedding} their moved spaces are also incomparable. That is, $\Mov(w')\cup\Mov(t)$ properly contains $\Mov(w')$ so that
\begin{equation*}
\Mov(w')\subsetneq \Mov(w')\cup\Mov(t)\subseteq \Mov(w'\vee t),
\end{equation*}
and we conclude that $\dim\Mov(w'\vee t)\geq\dim\Mov(w')+1$, or $\ell_T(w'\vee t)\geq\ell_T(w')+1$, which completes the proof.

To prove $(2)$, consider $\pi$, $\mu$ and $\pi\mu$ in $NC(W)$ with $\pi\leq_T\pi\mu$. Then there exists a reduced $T$-word $\pi\mu=t_1t_2\cdots t_r$ with $\pi=t_1t_2\cdots t_k$ for some $1\leq k\leq r$. In this case $\mu=t_{k+1}\cdots t_r$ is also a reduced $T$-word. Apply $(1)$ and the associativity of join. The proof in the case $\mu\leq_T\pi\mu$ is similar.
\end{proof}

The second result shows how the lattice property interacts with the Kreweras complement maps.

\begin{theorem}
\label{theorem:kreweras}
Given $\mu\leq_T\pi\leq_T\nu$ in $NC(W)$, we have
\begin{enumerate}
\item $K_\mu^\nu(\pi)=\mu\vee K_1^\nu(\pi)=\nu\wedge K_\mu^c(\pi)$.
\label{krewitem:one}
\item The map $K_\mu^\nu$ is a {\sf lattice complement} on $[\mu,\nu]$. That is,
\begin{equation*}
\pi\wedge K_\mu^\nu(\pi)=\mu \quad\text{ and }\quad \pi\vee K_\mu^\nu(\pi)=\nu.
\end{equation*}
\label{krewitem:two}
\end{enumerate}
\end{theorem}

\begin{proof}
Applying Lemma \ref{lemma:joins}, we immediately have
\begin{equation*}
K_\mu^\nu(\pi)=\mu(\pi^{-1}\nu)=\mu\vee (\pi^{-1}\nu)=\mu\vee K_1^\nu(\pi).
\end{equation*}
Now set $\mu'=K_1^c(\mu)$, $\pi'=K_1^c(\pi)$ and $v'=K_1^c(\nu)$. Since $\nu'\leq_T\pi'\leq_T\mu'$, it makes sense to consider $K_{\nu'}^{\mu'}(\pi')$ and it is easy to see that $K_{\nu'}^{\mu'}(\pi')=K_1^c(K_\mu^\nu(\pi))$. Finally, applying the anti-automorphism $(K_1^c)^{-1}$ to equation $K_{\nu'}^{\mu'}(\pi')=\nu'\vee K_1^{\mu'}(\pi')$ gives
\begin{eqnarray*}
K_\mu^\nu(\pi) &=& (K_1^c)^{-1}(K_{\nu'}^{\mu'}(\pi'))\\
&=& (K_1^c)^{-1}(\nu'\vee K_1^{\mu'}(\pi'))\\
&=& (K_1^c)^{-1}(\nu')\wedge (K_1^c)^{-1}(K_1^{\mu'}(\pi'))\\
&=& \nu\wedge K_\mu^c(\pi),
\end{eqnarray*}
proving \eqref{krewitem:one}.

To show \eqref{krewitem:two}, we apply \eqref{krewitem:one} and properties of meet to get
\begin{eqnarray*}
\pi\wedge K_\mu^\nu(\pi) &=& \pi\wedge(\nu\wedge K_\mu^c(\pi))\\
&=& \nu\wedge(\pi\wedge K_\mu^c(\pi))\\
&=& \nu\wedge K_\mu^\pi(\pi)\\
&=& \nu\wedge\mu=\mu.
\end{eqnarray*}
The proof of $\pi\vee K_\mu^\nu(\pi)=\nu$ is similar.
\end{proof}

In the 1990's it became popular to generalize classical combinatorial objects in the setting of reflection groups, and this was successively done for Kreweras' noncrossing set partitions by Reiner~\cite{reiner} and Biane~\cite{biane:crossings}. Reiner asked in~\cite[Remark 2]{reiner} whether the construction could be generalized uniformly to all finite reflection groups. Without apparent knowledge of Biane's work~\cite{biane:crossings} or Reiner's question, the same definition (Definition \ref{def:NC(W)}) was given nearly simultaneously by Brady and Watt~\cite{brady-watt:kpione} and Bessis~\cite{bessis:dual}. This happened around 2000 and the subject has quickly taken off. (See the historical comments at the end of Section \ref{sec:classicalNC}).

 In~\cite{bessis:dual} Bessis gave numerological evidence for a sort of ``duality'' between $NC(W)$ and the weak order $\Weak(W)$, but he was not able to formalize this relationship. In the next section we will encounter some of the {\sf numerology} inherent in this subject. We use the term here in a non-superstitious way to describe suggestive and unexplained (or only dimly-understood) enumerative coincidences. A general paradigm in this subject is the use of numerological evidence to suggest new uniform theories.

\section{Invariant Theory and Catalan Numbers}
\label{sec:invariant}
We begin with an ``empirical'' observation. Let $(W,S)$ be a finite Coxeter system of rank $\abs{S}=n$, and consider the length generating polynomials of the standard and absolute lengths on $W$:
\begin{equation*}
P(W,S,q)=\sum_{w\in W} q^{\ell_S(w)}\quad\text{and}\quad P(W,T,q)=\sum_{w\in W} q^{\ell_T(w)}.
\end{equation*}
It is a surprising fact that each of these polynomials has a nice factorization. Considering some examples, one might observe that
\begin{equation*}
P(W,S,q)=\prod_{i=1}^n \frac{q^{d_i}-1}{q-1}\quad\text{and}\quad P(W,T,q)=\prod_{i=1}^n (1+(d_i-1)q),
\end{equation*}
where the sequence of integers $\{d_1,d_2,\ldots,d_n\}$ is the same in both cases. We will soon see that these integers arise as degrees of generators in the ring of polynomial invariants of $W$, hence they are called the {\sf degrees} of $W$. We will always number the degrees in weakly increasing order $d_1\leq d_2\leq\cdots\leq d_n$. If $W=W'\times W''$ is a reducible Coxeter group, and $W'$ and $W''$ have reflections $T'$ and $T''$, respectively, it is evident from the definition of the polynomial $P(W,T,q)$ that
\begin{equation*}
P(W,T,q)=P(W,T',q)\, P(W'',T'',q),
\end{equation*}
and so the degree sequence of $W$ is the concatenation of the degree sequences of its factors $W'$ and $W''$. The complete list of degrees for the finite irreducible Coxeter groups is displayed in Figure \ref{fig:degrees}.

\begin{figure}
\vspace{.1in}
\begin{center}
\begin{tabular}{|c|l|}
\hline
$W$ & $d_1,\ldots,d_n$\\
\hline
$A_n$ & $2,3,\ldots, n+1$\\
$B_n$ & $2,4,6,\ldots,2n$\\
$D_n$ & $2,4,6,\ldots,2(n-1),n$\\
$E_6$ & $2,5,6,8,9,12$\\
$E_7$ & $2,6,8,10,12,14,18$\\
$E_8$ & $2,8,12,14,18,20,24,30$\\
$F_4$ & $2,6,8,12$\\
$H_3$ & $2,6,10$\\
$H_4$ & $2,12,20,30$\\
$I_2(m)$ & $2,m$\\
\hline
\end{tabular}
\end{center}
\caption{Degrees of the finite irreducible Coxeter groups}
\label{fig:degrees}
\end{figure}

Simple observations lead to some interesting properties of these integers. As is standard, we will denote the number of reflections in $W$ by $N=\abs{T}$.  Putting $q=1$ in $P(W,T,q)$ we get $\prod_{i=1}^n d_i=\abs{W}$. Differentiating $P(W,T,q)$ by $q$ and setting $q=0$ yields $\sum_{i=1}^n d_i =N+n$. Also, observing the coefficients of $P(W,T,q)$ we see that the rank numbers of $\Abs(W)$ are given by elementary symmetric combinations of the integers $d_1-1,d_2-1,\ldots,d_n-1$. In particular, the number of maximal elements of the absolute order is equal to $\prod_{i=1}^n (d_i-1)$. 

One can check that the length generating polynomial for a general finite group $G$ with generating set $H$ does {\em not} factor in this way. So a finite Coxeter group for some reason comes equipped with a sequence of special integers. What is the nature these numbers? Where do they come from? In what other groups does this phenomenon occur? Much of the work on Coxeter groups in the second half of the twentieth century was aimed at answering these questions, and still questions remain.

In discussing this topic, we will necessarily present some results out of chronological order, since the web of ideas is highly connected. The seminal work in this field \cite{coxeter2} was published in 1951 by Coxeter. In 1950, Chevalley had given an address at the International Congress of Mathematicians at Cambridge, where he proposed a new method for computing the Poincar\'e polynomials of the exceptional compact simple Lie groups. In general, if $X$ is a topological space with finite dimensional homology groups $H_i(X,\reals)$, the {\sf Poincar\'e series} of $X$ is the generating function of the Betti numbers,
\begin{equation*}
P(X,q)=\sum_{i\geq 0} \dim H_i(X,\reals)\, q^i.
\end{equation*}
When $X$ is a simple Lie group, Chevalley observed that this series has a factorization in terms of special integers. He had essentially shown that
\begin{equation*}
P(X,q)=\prod_{i=1}^n (1+q^{2d_i-1}),
\end{equation*}
where the $d_1,d_2,\ldots,d_n$ are what we have called the degrees of the corresponding Weyl group. Coxeter was in attendance at this talk and he immediately recognized the numbers $\{d_1,d_2,\ldots,d_n\}$, having encountered them previously in his own work. This coincidence was the inspiration for the paper~\cite{coxeter2} and the beginning of a long line of investigation. Using modern notation, Coxeter's insight was the following.

\begin{theorem}[\cite{coxeter2}]
\label{theorem:Coxeter}
Let $c\in W$ be a Coxeter element, where $W$ is a finite Coxeter group of rank $n$. The eigenvalues of $c$ are given by
\begin{equation*}
\omega^{d_1-1}, \omega^{d_2-1},\ldots, \omega^{d_n-1},
\end{equation*}
where $\{d_1,d_2,\ldots,d_n\}$ are the degrees of $W$, $\omega=e^{2\pi i/h}$ and $h$ is the order of $c$ as an element of $W$.
\end{theorem}

Accordingly, the integers $\{m_1,m_2,\ldots,m_n\}$, where $m_i=d_i-1$, are called the {\sf exponents} of the group (Coxeter's notation). The order $h$ of a Coxeter element is today called the {\sf Coxeter number} of $W$.

Thus, the  degrees of $W$ are related to the word length on $W$ (with respect to $S$ or $T$); the homology of the corresponding Lie group (of course, only when $W$ is a Weyl group); and the eigenvalues of a Coxeter element. Why call them ``degrees'' then? This name comes from yet another definition of these integers in the subject of {\sf invariant theory}.

Consider a finite dimensional $K$-linear representation $\rho:G\to GL(V)$ of a finite group $G$, with $K\in\{\reals,\complex\}$, and let $S=S(V^*)$ denote the symmetric algebra of the dual space $V^*$, which is the algebra of {\sf polynomial functions} on $V$. That is, relative to any fixed basis $\{x_1,x_2,\ldots,x_n\}$ of $V$, we may identify $S$ with the ring of polynomials $K[dx_1,dx_2,\ldots,dx_n]$ in the coordinate functions. There is a natural action of $G$ on $S$ induced by the {\sf contragredient action} of $G$ on $V^*$: $(g\cdot f)(v)=f(g^{-1}\cdot v)$ for all $g\in G$, $v\in V$ and $f\in V^*$. Define the {\sf ring of polynomial invariants} of $\rho:G\to GL(V)$ as
\begin{equation*}
S^G=\left\{ f\in S: g\cdot f=f \text{ for all } g\in G\right\}.
\end{equation*}

The motivating example of an invariant ring is the ring $\Lambda^n$ of {\sf symmetric polynomials} in $x_1,x_n,\ldots,x_n$, whose study goes back at least to Newton. We say that  a polynomial $f\in K[x_1,x_2,\ldots,x_n]$ is {\sf symmetric} if it satisfies
\begin{equation*}
f(x_1,x_2,\ldots,x_n)=f(x_{\sigma(1)},x_{\sigma(2)},\ldots,x_{\sigma(n)})
\end{equation*}
for all permutations $\sigma$ of the set $\{1,2,\ldots,n\}$. If $\rho: \frak{S}_n\to GL(V)$ is the {\sf permutation representation} of the symmetric group $\frak{S}_n$, in which $\frak{S}_n$ permutes a basis of $V$, then $\Lambda^n$ is isomorphic to the invariant ring of $(\frak{S}_n,\rho)$. Note, however, that this is {\em not} the geometric representation of $\frak{S}_n$ as a Coxeter group (Section \ref{sec:coxsystems}), since, for instance, it is not irreducible. To obtain the geometric representation, we divide by the $1$-dimensional invariant subspace spanned by $x_1+x_2+\cdots +x_n$. This accounts for the sometimes confusing index discrepancy $\frak{S}_n=A_{n-1}$.

It was well-known since Burnside's 1911 text \cite{burnside} that the field of fractions of $S^G$ has transcendence degree $n$ over $K$. Thus, $S^G$ possesses a set of $n$ algebraically independent, homogeneous elements, called {\sf basic invariants}. In the case of a finite real reflection group $\rho:W\hookrightarrow GL(V)$, Chevalley and Coxeter (see \cite[Section 6]{coxeter2}) had both studied the ring of invariants and shown the stronger condition that {\em the basic invariants are a generating set for $S^W$}. That is, if $\{f_1,f_2,\ldots,f_n\}$ is a set of basic invariants, then $S^W=\reals[f_1,f_2,\ldots,f_n]$. One may ask for which other groups and representations this happens. In 1954, Shephard and Todd showed that this condition is sufficient to characterize finite reflection groups, where they allowed both real and complex reflection groups (discussed below). Their result \cite[Theorem 5.1]{shephard-todd} builds on and generalizes the work of Chevalley, Coxeter and Racah.

\begin{theorem}[\cite{shephard-todd}]
\label{theorem:shephard-todd1}
Let $G$ be a finite group with $n$-dimensional unitary representation $\rho:G\to GL(V)$. The following are equivalent.
\begin{enumerate}
\item $G$ is generated by (pseudo-) reflections.
\item $G$ possesses a set of basic invariants of degrees $d_1,d_2,\ldots,d_n$ such that \begin{equation*}
\prod_{i=1}^n d_i=\abs{G}.
\end{equation*}
\item $G$ possesses a set of basic invariants $f_1,f_2,\ldots,f_n$ such that the ring of invariants is polynomial, $S^G=\complex[f_1,f_2,\ldots,f_n]$.
\end{enumerate}
\end{theorem}

Because every finite dimensional complex representation of a finite group is unitary, this theorem says that, from the perspective of invariant theory, reflection groups are the {\em best} groups. Under these conditions Shephard and Todd also showed that the sequence of {\sf degrees of basic invariants} $d_1,d_2,\ldots,d_n$ is uniquely determined, even though the basic invariants themselves are not. We will discuss the idea of a pseudo-reflection at the end of this section, but for now assume that the representation is real (orthogonal). In this case, the degrees of basic invariants agree with the ``degrees'' in Figure \ref{fig:degrees}. In fact, this is the usual way to define these numbers, and it explains the name.

Now suppose that $W$ is a real reflection group with geometric representation $\rho:W\hookrightarrow GL(\reals^n)$. In general, $W$ preserves the standard inner product on $V$ (we have $(w\cdot u,w\cdot v)=(u,v)$ for all $w\in W$ and $u,v\in\reals^n$); and, since $W$ is real, it preserves the homogeneous degree $2$ polynomial $dx_1^2+\cdots +dx_n^2$. Hence $2$ is a degree of $W$. Moreover, if $W$ had an invariant $1$-form this would imply a fixed subspace, contradicting the fact that the geometric representation is essential. So the smallest degree is always $d_1=2$.

Thinking again of a Coxeter element $c$ and exponents $m_1,m_2,\ldots,m_n$, we know that the eigenvalues come in conjugate pairs, so that $m_i+m_{n-i+1}=h$ for all $i\leq i\leq n$. Taking $i=1$ or $i=n$, we find that the highest degree is always equal to the Coxeter number $d_n=m_n+1=m_n+m_1=h$. In addition, Springer \cite{springer} proved that the centralizer of a Coxeter element $c$ is just  the cyclic group of order $h$ generated by $c$. This proves that the conjugacy class of Coxeter elements in $W$ has size $\abs{W}/h=\prod_{i=1}^{n-1} d_i$. (In general, this is less than the number of maximal elements $\prod_{i=1}^n (d_i-1)$ in the absolute order.) Finally, since $c$ has no eigenvalue equal to $w^0=1$ ($1$ is not a degree), the numbers $h-m_i$ are a permutation of the exponents. Hence
\begin{equation*}
N=\sum_{i=1}^n m_i = \sum_{i=1}^n h-m_i = nh-N,
\end{equation*}
or $N=nh/2$. Here we needed the fact that $W$ is real, since in the general complex case, Coxeter's Theorem \ref{theorem:Coxeter} may fail.

The classification theorem \ref{theorem:shephard-todd1}  explains the definition of the degrees, but we have not yet explained the connection with the absolute length generating function, which was our original motivation. This result was first proved by Shephard and Todd \cite[Theorem 5.3]{shephard-todd} in a case-by-case way, and the first uniform proof was given by Solomon \cite{solomon}.

\begin{theorem}[\cite{shephard-todd}]
\label{theorem:shephard-todd2}
If $G$ is a finite group with $n$-dimensional unitary representation $\rho:G\to GL(V)$ generated by (pseudo-) reflections, let $\Fix(g)$ denote the fixed space $\ker(\rho(g)-\rho(1))$ for all $g\in G$. Then we have
\begin{equation}
\label{eq:shephard-todd}
\sum_{g\in G} q^{n-\dim\Fix(g)}=\prod_{i=1}^n (1+(d_i-1)q).
\end{equation}
\end{theorem}

From here, our observed factorization of the polynomial $P(W,T,q)$ is an easy corollary of Theorem \ref{th:basicMov} $(1)$. We do not know who was first to observe this fact.

We note that Shephard and Todd's Theorems \ref{theorem:shephard-todd1} and \ref{theorem:shephard-todd2} both apply in the slightly broader context of {\sf finite unitary groups generated by pseudo-reflections} (often called {\sf complex reflection groups}). A {\sf pseudo-reflection} is a unitary transformation of finite order on an $n$-dimensional vector space $V$, with exactly $n-1$ eigenvalues equal to $1$. That is, it fixes a hyperplane, and acts as multiplication by a root of unity on the orthogonal copy of $\complex$. In~\cite{shephard-todd} Shephard and Todd gave a full classification of these groups, generalizing Coxeter's classification in the orthogonal case. Remarkably, there are not many additional examples: the full classification contains 4 infinite classes and 34 exceptional groups. A.M.~Cohen has suggested a classification scheme generalizing Cartan-Killing, in which the letters $A$ through $N$ are used to denote the groups~\cite{cohen}.

All of the proofs provided in Shephard and Todd~\cite{shephard-todd} are case-by-case, building on previous work. This is typical for a subject in which there is a complete classification available. The program to provide uniform proofs of these results occupied several subsequent generations of mathematicians and it continues today. An interesting question is whether the degrees can be defined in a purely combinatorial way, or whether they are essentially geometric, or topological. Barcelo and Goupil discuss these issues in a 1994 survey~\cite{barcelo-goupil} with many historical references.

Now a new chapter has opened in the numerology of Coxeter groups.
\begin{definition}
If $W$ is a finite Coxeter group, define the {\sf Coxeter-Catalan number} associated to $W$ by
\begin{equation}
\Cat(W):=\prod_{i=1}^n \frac{h+d_i}{d_i}=\frac{1}{\abs{W}}\prod_{i=1}^n (h+d_i).
\end{equation}
\end{definition}

\begin{figure}
\vspace{.1in}
\begin{center}
\begin{tabular}{|c|c|c|c|c|c|c|c|c|c|c|c|}
\hline
$A_{n-1}$ &
$B_n$ & $D_n$ & $I_2(m)$ & $H_3$ & $H_4$ & $F_4$ &$E_6$ & $E_7$ & $E_8$ \\
\hline
&&&&&&&&&\\[-.1in]
$\textstyle\frac{1}{n} \binom{2n}{n-1}$
&$\textstyle\binom{2n}{n}$
&$\textstyle\frac{3n-2}{n}\binom{2n-2}{n-1}$
&$\textstyle m+2$
&$\textstyle 32$
&$\textstyle 280$
&$\textstyle 105$
&$\textstyle 833$
&$\textstyle 4160$
&$\textstyle 25080$\\[.05in]
\hline
\end{tabular}
\end{center}
\caption{The numbers $\Cat(W)$ for the finite irreducible Coxeter groups} 
\label{fig:cox-cat}
\end{figure}

It is also possible to define a Coxeter-Catalan number for certain ``well-generated'' complex reflection groups (see \cite{bessis:complex}), but we will consider only the real case here.

Figure \ref{fig:cox-cat} displays the complete list of Coxeter-Catalan numbers for finite irreducible Coxeter groups. These are readily computed from the degrees in Figure \ref{fig:degrees}. In the introduction we explained some of the provenance of the number $\Cat(W)$. It was first written down by Djocovi\'c in 1980~\cite{djokovic} and it is known to count several seemingly unrelated objects. Recently, interest in ``Coxeter-Catalan combinatorics'' has exploded, and there has been a workshop at the American Institute of Mathematics devoted exclusively to this topic (see~\cite{armstrong:braids}). Of immediate interest is the following theorem.

\begin{theorem}
The cardinality of the set $NC(W)$ is $\Cat(W)$.
\end{theorem}

\begin{proof}
This formula was first conjectured by Postnikov, based on theorems of Reiner in the classical types (see \cite[Remark 2]{reiner}), before the uniform definition of $NC(W)$ was known. The exceptional types were later verified by computer~\cite{bessis:dual}. As of this writing {\em a uniform proof is not known}.
\end{proof}

Work on the Coxeter-Catalan combinatorics continues on several fronts and one of the aims of this memoir is to give an introduction to this program. Another of our goals is to generalize the theory of Coxeter-Catalan combinatorics to the theory of ``Fuss-Catalan combinatorics'' (see Section \ref{sec:fusscat}). In Chapter \ref{sec:kdiv} we will define and study a ``Fuss-Catalan'' generalization of the lattice of noncrossing partitions $NC(W)$ for each positive integer $k$. It turns out that these ``$k$-divisible noncrossing partitions'' are intimately related to other structures recently considered by Athanasiadis, Fomin and Reading. In Chapter \ref{sec:fusscatcomb} we will explore some of these relationships.


\chapter{$k$-Divisible Noncrossing Partitions}
\label{sec:kdiv}

Here we define and study a poset $NC^{(k)}(W)$ --- called the poset of ``$k$-divisible noncrossing partitions'' --- for each finite Coxeter group $W$ and each extended positive integer $k\in\posint\cup\{\infty\}$. When $k=1$, our definition coincides with $NC(W)$. 

Whereas Chapter \ref{sec:background} was largely a survey, the material in this chapter is new. Throughout this chapter let $(W,S)$ denote a fixed finite Coxeter system with rank $n=\abs{S}$ and set of reflections $T$. We also fix a Coxeter element $c\in W$.

\section{Minimal Factorizations}
\label{sec:minimal}
Essentially, this chapter is a study of phenomena that occur in the $k$-fold direct product of the absolute order. In general, if P is a poset, the {\sf $k$-fold direct product poset} is the set
\begin{equation*}
P^k=\{(x_1,x_2,\ldots,x_k): x_1,x_2,\ldots,x_k\in P\}
\end{equation*}
together with the componentwise partial order. We will denote the $k$-fold direct product of the absolute order by $\Abs^k(W)$ and for notational convenience we will often denote an element $(w_1,w_2,\ldots,w_k)$ simply by $(w)_k$.

Note that every element $(w)_k\in\Abs^k(W)$ determines an element of $W$ by multiplication; we define the {\sf multiplication map} $m:\Abs^k(W)\to W$ by
\begin{equation*}
m(w_1,w_2,\ldots,w_k):=w_1w_2\cdots w_k,
\end{equation*}
and we say that $(w)_k\in\Abs^k(W)$ is a {\sf $k$-factorization} of $m(w)_k\in W$. Less formally, we may sometimes refer to the word $w_1w_2\cdots w_k$ as a $k$-factorization of $w$. In general, the triangle inequality for $\ell_T$ implies that
\begin{equation*}
\ell_T(w)\leq \ell_T(w_1)+\cdots+\ell_T(w_k).
\end{equation*}
When equality holds, we give the factorization a special name.\begin{definition}
The sequence $(w_1,w_2,\ldots,w_k)\in\Abs^k(W)$ is called a {\sf minimal $k$-factorization} of $w\in W$ if
\begin{enumerate}
\item $w=w_1w_2\cdots w_k$
\item $\ell_T(w)=\ell_T(w_1)+\ell_T(w_2)+\cdots +\ell_T(w_k)$
\end{enumerate}
\end{definition}

Note that we have deliberately allowed a minimal factorization to contain copies of the identity $1\in W$. This introduces some redundancy; that is,  if we permute the entries of a minimal factorization $(w)_k$ without changing the relative order of the non-identity entries then the element $m(w)_k\in W$ is unchanged. For example, the sequences $(w,1,1)$, $(1,w,1)$ and $(1,1,w)$ are all minimal $3$-factorizations for a given $w\in W$. Another consequence is the fact that there exist minimal $k$-factorizations of each $w\in W$ for arbitrarily large values of $k$: simply augment the sequence with copies of the identity.

The concept of a minimal factorization is intended to generalize the definition of a reduced $T$-word. If $w=t_1 t_2\cdots t_k$ is a reduced $T$-word for $w$ then $(t_1,t_2,\ldots,t_k)$ is clearly a minimal $k$-factorization of $w$. Further, any minimal factorization for $w$ that does {\em not} contain copies of the identity can be thought of as a coarsening, or a bracketing, of a reduced $T$-word for $w$. Indeed, if $(w_1,w_2,\ldots,w_k)$ is a minimal factorization of $w\in W$ with $w_i\neq 1$ for all $1\leq i\leq k$, then we have $k\leq_T\ell_T(w)$, with equality if and only if $w_i\in T$ for all $1\leq i\leq k$.

Like reduced $T$-words, minimal factorizations have a corresponding notion of shifting. This is a useful rephrasing of the Shifting Lemma \ref{shifting}.

\begin{genshiftinglemma}
\label{shifting2}
Given a minimal $k$-factorization $(w)_k$ of $w\in W$ and an integer $1<i<k$, the two sequences
\begin{eqnarray*}
&&(w_1,w_2,\ldots,w_{i-2},w_i,w_i^{-1}w_{i-1}w_i,w_{i+1},\ldots,w_k)\quad\text{and}\\
&&(w_1,w_2,\ldots,w_{i-1},w_i w_{i+1}w_i^{-1},w_i,w_{i+2},\ldots,w_k)
\end{eqnarray*}
are also minimal $k$-factorizations for $w$.
\end{genshiftinglemma}

\begin{proof}
Each is clearly a factorization of $w$. The fact that they are minimal follows from the fact that the absolute length $\ell_T:W\to\integers$ is invariant under conjugation (see Section \ref{sec:shifting}).
\end{proof}

That is, if the element $w'$ occurs in the $i$-th place of some minimal $k$-factorization for $w\in W$ with $1<i<k$, then there exist other minimal $k$-factorizations of $w$ in which $w'$ occurs in the $(i+1)$-th and the $(i-1)$-th place, respectively.

Finally, given a minimal $k$-factorization $(w_1,w_2,\ldots,w_k)$ of $w\in W$, note that we have $w_i\leq_T w$ for all $1\leq i\leq k$. This follows easily from the subword property \ref{subword}.

\section{Multichains and Delta Sequences}
\label{sec:mchains}
Now we will apply the idea of minimal factorizations to the lattice of noncrossing partitions $NC(W)$ with respect to the fixed Coxeter element $c\in W$.

First consider a finite positive integer $k\in\posint$ and let $NC^k(W)$ denote the $k$-fold direct product poset of $NC(W)$. If $W$ has Coxeter generators $S=\{s_1,s_2,\ldots,s_n\}$, observe that the $k$-fold direct product $W^k$ is itself a Coxeter group with respect to the generators
\begin{equation*}
S^k:=\{s_{i,j}=(1,1,\ldots,1,s_j,1,\ldots,1): 1\leq i,j\leq n\},
\end{equation*}
where $s_j$ occurs in the $i$-th entry of $s_{i,j}$. (The Coxeter diagram of $W^k$ is just the disjoint union of $k$ copies of the Coxeter diagram for $W$.) We will denote the reflections of $W^k$ by $T^k$.  Hence, we may consider the lattice of noncrossing partitions $NC(W^k)$ of the Coxeter system $(W^k,S^k)$ with respect to the Coxeter element $(c,c,\ldots,c)\in W^k$. The following is immediate.

\begin{lemma}
$NC^k(W)=NC(W^k)$
\end{lemma}

We will confuse these two structures whenever we please. In our study of the poset $NC^k(W)$ we are interested in the following two related families of elements, which play ``dual'' roles.

\begin{definition}\hspace{.1in}
\label{def:mchaindsequence}
\begin{enumerate}
\item We say $(\pi)_k:=(\pi_1,\pi_2,\ldots,\pi_k)\in NC^k(W)$ is a {\sf $k$-multichain} if
\begin{equation*}
\pi_1\leq_T\pi_2\leq_T\cdots\leq_T\pi_k.
\end{equation*}
\item We say $[\delta]_k:=[\delta_0;\delta_1,\delta_2,\ldots,\delta_k]\in NC^{k+1}(W)$ is a {\sf $k$-delta sequence} if
\begin{equation*}
c=\delta_0\delta_1\cdots\delta_k\quad\text{and}\quad\ell_T(c)=\sum_{i=0}^k\ell_T(\delta_i).
\end{equation*}
\end{enumerate}
\end{definition}

Some explanation is in order: The notion of a multichain is standard in the combinatorics literature. In general, a {\sf chain} in a poset $P$ is a set of strictly increasing elements $x_1<x_2<\cdots<x_k$, and the ``multi-'' prefix indicates that we allow repetition of elements. Note that what we have called a ``$k$-multichain'' is sometimes called a ``multichain of length $k-1$'' (see \cite[Chapter 3]{stanley:ec1}). By definition, a $k$-delta sequence $[\delta]_k$ is nothing but a minimal $(k+1)$-factorization of the Coxeter element $c$. Indeed, this condition is enough to guarantee that $\delta_i$ is in $NC(W)$ for all $0\leq i\leq k$.

We introduced the notions of multichain and delta sequence together because they are equivalent, in the sense that a delta sequence is the sequence of ``successive differences'' of a multichain. We formalize this relationship with the following two maps.

\begin{definition}
For each $k$-multichain $(\pi)_k$, define a corresponding $k$-delta sequence $\partial(\pi)_k$ by
\begin{equation}
\label{partial}
\partial(\pi)_k:=[\pi_1;\pi_1^{-1}\pi_2,\pi_2^{-1}\pi_3,\ldots,\pi_{k-1}^{-1}\pi_k,\pi_k^{-1}c],
\end{equation}
and given any $k$-delta sequence, define a corresponding $k$-multichain $\sint[\delta]_k$ by
\begin{equation}
\label{sint}
\sint[\delta]_k:=(\delta_0,\delta_0\delta_1,\delta_0\delta_1\delta_2,\ldots,\delta_0\delta_1\cdots\delta_{k-1}).
\end{equation}
\end{definition}

Note that we might instead have defined $\partial$ using Kreweras complements (Definition \ref{def:algkreweras}); if we suppress the subscript $1$ in $K_1^\pi$, we get 
\begin{equation*}
\partial(\pi)_k=[K^{\pi_1}(1),K^{\pi_2}(\pi_1),K^{\pi_3}(\pi_2),\ldots,K^{\pi_k}(\pi_{k-1}),K^c(\pi_k)].
\end{equation*}
Then since $\ell_T(K^\nu(\mu))=\ell_T(\nu)-\ell_T(\mu)$ for all $\mu\leq_T\nu$ (Lemma \ref{lemma:kreweras} $(2)$), we see that $\partial(\pi)_k$ is indeed a delta sequence.
The fact that $\sint[\delta]_k$ is a multichain follows from the subword property.

\begin{figure}
\vspace{.1in}
\begin{center}
\input{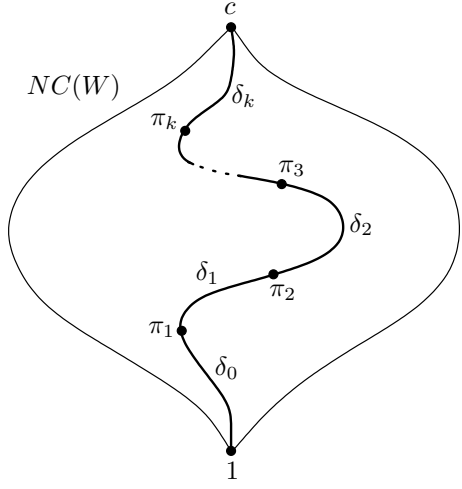}
\caption{A multichain / delta sequence pair in $NC(W)$}
\label{fig:multichain}
\end{center}
\end{figure}

To make sense of these definitions, consider the diagram in Figure \ref{fig:multichain}. The curved line represents a geodesic in the Cayley graph from $1$ to the Coxeter element $c$, which is the same as a reduced $T$-word for $c$. The dots are the elements of the multichain and the lines between dots represent reduced $T$-words for the elements of the delta sequence. It is easy to verify that $\partial$ and $\sint$ are reciprocal bijections. Hence multichains and delta sequences determine each other. We will use them interchangeably.

The notation for ``delta sequences'' is meant to suggest an analogy with sequences of homology groups or sequences of differential maps, although the analogy is vague at this point. Theorem \ref{th:nckl=nclk} in Section \ref{sec:iterated} will clarify the analogy.

Now we will examine some elementary properties of delta sequences. Notice the use of the semi-colon in the notation $[\delta]_k=[\delta_0;\delta_1,\ldots,\delta_k]$. We have chosen this notation\footnote{Thanks to Krattenthaler \cite{krattenthaler1} for the suggestion.} to indicate that the first element $\delta_0$ plays a special role, separate from the other elements $(\delta)_k=(\delta_1,\ldots,\delta_k)$. Indeed, we need not even specify the first element $\delta_0$ since it is implicity determined by the other $k$ elements. This leads to an alternate definition of delta sequences.

\begin{lemma}
There is a bijection between $k$-delta sequences $[\delta]_k$ and pairs $(\pi,(\delta)_k)$ --- where $(\delta)_k$ is a minimal $k$-factorization of $\pi\in NC(W)$ --- given by
\begin{eqnarray*}
[\delta_0;\delta_1,\ldots,\delta_k] &\mapsto &(\delta_0^{-1}c,(\delta)_k), \\
(\pi,(\delta)_k) &\mapsto &[c\pi^{-1};\delta_1,\ldots,\delta_k].  
\end{eqnarray*}
\end{lemma}

\begin{proof}
The two maps are mutual inverses. We need to show that they are well-defined.

First, suppose that $[\delta]_k$ is a $k$-delta sequence and set $\pi=\delta_0^{-1}c$, so that $\pi=\delta_1\cdots\delta_k$. Since $\delta_0\leq_T c$, it follows that $\pi=K_1^c(\delta_0)\leq_T c$ because the Kreweras complement $K_1^c$ maps $NC(W)$ to itself. Hence $\pi\in NC(W)$. And because $\ell_T(\pi)=\ell_T(c)-\ell_T(\delta_0)$, we have $\ell_T(\pi)=\sum_{i=1}^k \ell_T(\delta_k)$, which implies that $(\delta)_k$ is a minimal $k$-factorization of $\pi$.

Conversely, consider an arbitrary element $\pi\in NC(W)$ with minimal factorization $(\delta)_k$. If we set $\delta_0=c\pi^{-1}=(K_1^c)^{-1}(\pi)$, then we have $c=\delta_0\delta_1\cdots\delta_k$. And $\ell_T(\pi)=\ell_T(c)-\ell_T(\delta_0)$ because $\delta_0=c\pi^{-1}\leq_T c$. It follows that $\ell_T(c)=\sum_{i=0}^k\ell_T(\delta_i)$, hence $[\delta]_k$ is a $k$-delta sequence.
\end{proof}

\begin{altdefinition}
When no confusion can result, we will refer to a minimal $k$-factorization $(\delta)_k$ of an arbitrary element $\pi\in NC(W)$ as a {\sf $k$-delta sequence}, in which case the element $\delta_0$ is implicitly understood.
\end{altdefinition}

Thus a $k$-delta sequence can be thought of as a minimal $(k+1)$-factorization $[\delta]_k$ of the Coxeter element or as a minimal $k$-factorization $(\delta)_k$ of an arbitrary element $\pi\in NC(W)$. For consistency, we will always use the notation $[\delta]_k$ to refer to a sequence with $k+1$ entries and the notation $(\delta)_k$ to refer to a sequence with $k$ entries.

Notice that the notion of a multichain is common to all posets, but the definition of a delta sequence depends on the group structure of $W$. The poset $NC(W)$ has both a group structure and a lattice structure. In this case, we can say more.

\begin{lemma}
\label{lemma:delta}
Consider an arbitrary sequence $(\delta)_k\in NC^k(W)$ and let $K$ denote the Kreweras complement $\pi\mapsto \pi^{-1}c$. The following are equivalent.
\begin{enumerate}
\item $(\delta)_k$ is a delta sequence.
\label{item:one}
\item $\delta_j\leq_T K(\delta_i)$ for all $1\leq i<j\leq k$.
\label{item:two}
\item $\delta_{i_1}\delta_{i_2}\cdots\delta_{i_m}=\delta_{i_1}\vee\delta_{i_2}\vee\cdots\vee\delta_{i_m}$ for all $1\leq i_1<i_2<\cdots <i_m\leq k$.
\label{item:three}
\end{enumerate}
\end{lemma}

\begin{proof}
We will show that \eqref{item:one}$\Rightarrow$\eqref{item:two}$\Rightarrow$\eqref{item:three}$\Rightarrow$\eqref{item:one}.

First, consider indices $1\leq i<j\leq k$. To prove \eqref{item:one}$\Rightarrow$\eqref{item:two} we will use the shifting lemma \ref{shifting2}. Indeed, since $c=\delta_0\delta_1\cdots\delta_k$ is a minimal factorization, we may shift the symbol $\delta_i$ to the left to obtain a new minimal factorization 
\begin{equation*}
c=\delta_i\delta_0'\delta_1'\cdots\delta_{i-1}'\delta_{i+1}\cdots\delta_k,
\end{equation*}
where $\delta_p'=(\delta_i^{-1}\delta_p\delta_i)$. But now we see that $\delta_i^{-1}c=\delta_0'\cdots\delta_{i-1}'\delta_{i+1}\cdots\delta_k$ is a minimal factorization for $\delta_i^{-1}c=K_1^c(\delta_i)$, which must contain the symbol $\delta_j$ since $i<j$. 
By the subword property, we get $\delta_j\leq_T \delta_i^{-1}c$, as desired.

Now suppose that \eqref{item:two} holds, and fix $1\leq i_1<i_2<\cdots<i_m\leq k$. Since $\delta_{i_2}\leq_T K(\delta_{i_1})$, there exists a minimal factorization of the  form $c=\delta_{i_1}K(\delta_{i_1})=\delta_{i_1}\delta_{i_2}w'$, hence $\delta_{i_1}\delta_{i_2}$ is in $NC(W)$ and Lemma \ref{lemma:joins} implies that $\delta_{i_1}\delta_{i_2}=\delta_{i_1}\vee\delta_{i_2}$. Then, since $\delta_{i_3}\leq_T K(\delta_{i_1})$ and $\delta_{i_3}\leq_T K(\delta_{i_2})$, we have $\delta_{i_3}\leq_T K(\delta_{i_1})\wedge K(\delta_{i_2})=K(\delta_{i_1}\vee\delta_{i_2})$ by the definition of meet and the fact that $K$ is an anti-automorphism. Repeating the above argument, we get $\delta_{i_1}\delta_{i_2}\delta_{i_3}=\delta_{i_1}\vee\delta_{i_2}\vee\delta_{i_3}$, and \eqref{item:three} follows by induction.

Finally, suppose that \eqref{item:three} holds and define a sequence $(\pi)_k$ in $NC^k(W)$ by setting
\begin{equation*}
\pi_i:=K^{-1}(\delta_i\delta_{i+1}\cdots\delta_k)=c(\delta_i\delta_{i+1}\cdots\delta_k)^{-1}
\end{equation*}
for all $1\leq i\leq k$. By \eqref{item:three} and properties of join, $(\pi)_k$ is a $k$-multichain. Hence $\partial(\pi)_k=[\delta]_k=[\delta_0;\delta_1,\ldots,\delta_k]$ is a $k$-delta sequence, proving \eqref{item:one}.
\end{proof}

Now we consider the case where $k=\infty$. Let $NC^\infty(W)$ denote the set of semi-infinite sequences
\begin{equation*}
NC^\infty(W):=\{(w_1,w_2,w_3,\ldots):w_i\in NC(W)\text{ for all }i\in\posint\}
\end{equation*}
together with the componentwise partial order. In this case, the definitions of $\infty$-multichain and $\infty$-delta sequence extend without difficulty, because a multichain may contain repetition and a delta sequence may contain copies of the identity. In particular, an $\infty$-delta sequence may contain only finitely many entries not equal to $1$. The reciprocal bijections $\partial$ and $\sint$ also extend without difficulty.  The reason we consider the case $k=\infty$ separately is because $W^\infty$ is an {\em infinitely generated Coxeter group}, and in this case the noncrossing partitions $NC(W^\infty)$ are not defined. To get around this difficulty, we set
\begin{equation*}
NC(W^\infty):=NC^\infty(W).
\end{equation*}

Now we may regard $NC(W^k)$ as an induced subposet of $NC(W^\ell)$ for all extended positive integers $k,\ell\in\posint\cup\{\infty\}$ with $k\leq\ell$. In order to facilitate this point of view, we fix notation.

\begin{notation}
\label{notation:extend}
Given any $k$-multichain $(\pi)_k$ and $k$-delta sequence $(\delta)_k$, we set $\pi_r=c$ and $\delta_r=1$ for all $r>k$.
\end{notation}

Bessis and Corran have independently considered a particular combinatorial realization of delta sequences, which they call {\sf derived sequences} \cite[Definition 8.4]{bessis-corran}. We have borrowed the differential notations $\partial$ and $\sint$ from them.

\section{Definition of $k$-Divisible Noncrossing Partitions}
\label{sec:kdivdef}
This section is the heart of the memoir. Here we will define two mutually dual posets $NC_{(k)}(W)$ and $NC^{(k)}(W)$ for each finite Coxeter group $W$ and each extended positive integer $k\in\posint\cup\{\infty\}$. In the case $k=1$, we will see that both of these coincide with the lattice of noncrossing partitions $NC(W)$.

Later, in Chapter \ref{sec:classical}, we will see that the poset $NC^{(k)}(W)$ has a remarkable combinatorial interpretation for the classical groups $A_{n-1}$ and $B_n$: In these cases, $NC^{(k)}(W)$ is isomorphic to a poset of ``noncrossing'' set partitions, each of whose blocks has cardinality divisible by $k$. The general definition of $NC^{(k)}(W)$ was inspired by these two cases; for this reason, we will refer to $NC^{(k)}(W)$ as the poset of {\sf $k$-divisible noncrossing partitions} of $W$, even when $W$ is not of classical type. Hopefully this nomenclature will cause no more confusion than does the notation ``noncrossing partitions'' for $NC(W)$.

\begin{definition}
\label{def:NCk}
Suppose that $k\in\posint\cup\{\infty\}$ is an extended positive integer.
\begin{enumerate}
\item Let $NC_{(k)}(W)$ denote the set of $k$-delta sequences $(\delta)_k\in NC(W^k)$ together with the componentwise partial order
\begin{equation*}
(\delta)_k\leq(\varepsilon)_k\quad\Longleftrightarrow\quad (\delta)_k\leq_{T^k} (\varepsilon)_k.
\end{equation*}
\item Let $NC^{(k)}(W)$ denote the set of $k$-multichains $(\pi)_k\in NC(W^k)$ together with the partial order
\begin{equation*}
(\pi)_k\leq (\mu)_k \quad\Longleftrightarrow\quad \partial(\mu)_k\leq_{T^k}\partial(\pi)_k.
\end{equation*}
\end{enumerate}
\end{definition}

By definition, $NC_{(k)}(W)$ is just the induced subposet of $NC(W^k)$ consisting of $k$-delta sequences. That is, we have $(\delta)_k\leq (\varepsilon)_k$ in $NC_{(k)}(W)$ if and only if $\delta_i\leq_T\varepsilon_i$ for all $1\leq i\leq k$. Notice that we place no {\em a priori} restriction on the elements $\delta_0$ and $\varepsilon_0$, but we will soon find that $(\delta)_k\leq (\varepsilon)_k$ implies $\varepsilon_0\leq_T\delta_0$.

The poset of $k$-divisible noncrossing partitions $NC^{(k)}(W)$ consists of multichains; thus we might think of a $k$-multichain $(\pi)_k\in NC(W^k)$ as an ``algebraic $k$-divisible noncrossing partition''. However, the partial order on $NC^{(k)}(W)$ is defined {\em in terms of the dual order on corresponding delta sequences}: we have $(\pi)_k\leq (\mu)_k$ in $NC^{(k)}(W)$ if and only if
\begin{equation*}
\mu_i^{-1}\mu_{i+1}\leq_T \pi_i^{-1}\pi_{i+1}\quad\text{ for all } 1\leq i\leq k.
\end{equation*}

By definition, $NC_{(k)}(W)$ and $NC^{(k)}(W)$ are just two ways of looking at the same structure, and the two posets are dual to each other via the reciprocal anti-isomorphisms
\begin{equation*}
\partial:NC^{(k)}(W)\to NC_{(k)}(W)\quad\text{ and }\quad\sint:NC_{(k)}(W)\to NC^{(k)}(W).
\end{equation*}

\begin{figure}
\vspace{.1in}
\begin{center}
\input{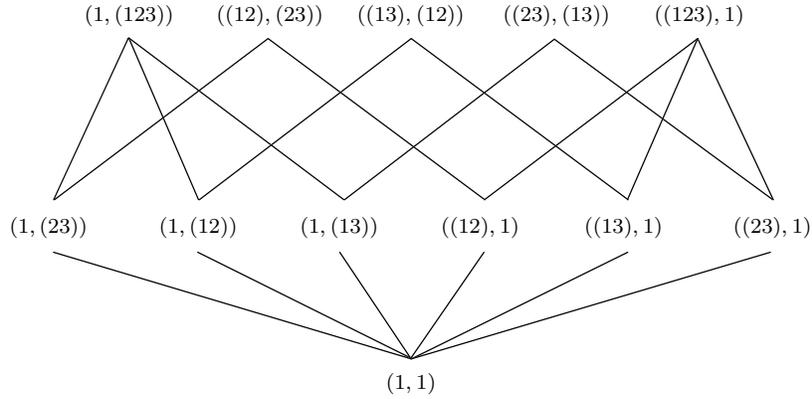}
\end{center}
\caption{The Hasse diagram of $NC_{(2)}(A_2)$}
\label{fig:nc_2(a2)}
\end{figure}

\begin{figure}
\vspace{.1in}
\begin{center}
\input{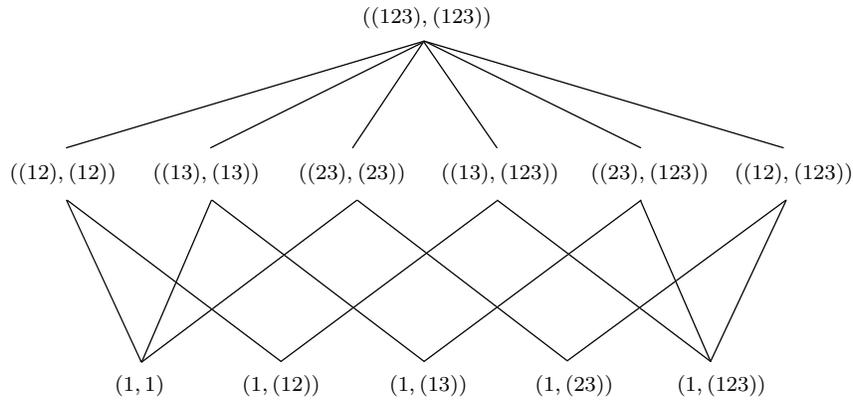}
\end{center}
\caption{The Hasse diagram of $NC^{(2)}(A_2)$}
\label{fig:nc^2(a2)}
\end{figure}

\begin{example}
If we fix the Coxeter element $(123)$ in the symmetric group $A_2$, Figures \ref{fig:nc_2(a2)} and \ref{fig:nc^2(a2)} display the Hasse diagrams of $NC_{(2)}(A_2)$ and $NC^{(2)}(A_2)$, respectively. Observe the duality between these  posets via the anti-isomorphisms $\partial$ and $\sint$. Unlike the noncrossing partitions $NC(W)$, the poset $NC_{(2)}(A_2)$ has no maximum element and the poset $NC^{(2)}(A_2)$ has no minimum element.
\end{example}

In general, the posets $NC_{(k)}(W)$ and $NC^{(k)}(W)$ {\em are not lattices} because (for instance) any two distinct minimal elements of $NC^{(k)}(W)$ have no lower bound.

The term ``$k$-divisible noncrossing partitions'' can be used equally well to refer to either of these posets since they are anti-isomorphic, but in formal situations we will reserve the term for the poset of multichains $NC^{(k)}(W)$. For consistency, we will use the notation $\leq$ with no subscript for both partial orders $NC_{(k)}(W)$ and $NC^{(k)}(W)$.

Finally, note that the construction of $P_{(k)}$ and $P^{(k)}$ can be carried out whenever the poset $P$ arises from the word length on a group. For example, consider a group $G$ with generating set $H$ and word length $\ell_H$. Define a partial order $P(G)$ on $G$ by setting
\begin{equation*}
g\leq_H g'\quad\Longleftrightarrow\quad \ell_H(g')=\ell_H(g)+\ell_H(g^{-1}g').
\end{equation*}
In this case we may mimic \ref{def:NCk} to define the posets $P^{(k)}(G)$ and $P_{(k)}(G)$. For instance, if $W$ is a finite Coxeter group with generators $S$, then $P(W)$ is the right weak order. It may be interesting to consider $P^{(k)}(W)$ in this case.

\section{Basic Properties of $k$-Divisible Noncrossing Partitions}
\label{sec:basics}
Now we will develop some basic theory. First, we examine the cases $k=1$ and $k=2$, for which the definition of $NC^{(k)}(W)$ reduces to known constructions.

\subsection{Small values of $k$}
It is clear from the definition that $NC_{(1)}(W)$ is isomorphic to $NC(W)$ via the trivial identification $\iota(\delta):=\delta$, and so the composition $\iota\,\circ\,\partial$ is an anti-isomorphism from $NC^{(1)}(W)$ to $NC(W)$. Recalling that $NC(W)$ is self-dual via the Kreweras complement anti-automorphism $K$, this implies the following.

\begin{lemma}
\label{lemma:nc1}
The map $\iota(\pi)=\pi$ is a poset isomorphism
\begin{equation*}
NC^{(1)}(W)\cong NC(W).
\end{equation*}
\end{lemma}

The $2$-divisible noncrossing partitions also have a familiar interpretation. Given any poset $P$, let $\Int(P)$ denote the {\sf interval poset} of $P$: that is, the set of intervals
\begin{equation*}
\Int(P):=\left\{ [x,y]\subseteq P: x\leq y\right\}
\end{equation*}
under the inclusion partial order. This is a standard object in combinatorics, and it is related to the concept of the {\sf incidence algebra} \cite[Chapter 3.6]{stanley:ec1}.

\begin{lemma}
\label{lemma:intervals}
The map $(\pi_1,\pi_2)\mapsto [K_1^{\pi_2}(\pi_1),\pi_2]$ is a poset isomorphism
\begin{equation*}
NC^{(2)}(W)\cong \Int(NC(W)).
\end{equation*}
\end{lemma}

\begin{proof}
The map is bijective since $K_1^{\pi_2}$ is a bijection on the interval $[1,\pi_2]$, and it is clearly invertible with inverse $[\mu_1,\mu_2]\mapsto ((K_1^{\mu_2})^{-1}(\mu_1),\mu_2)$. We need to show that both the map and its inverse preserve order.

So consider $(\pi_1,\pi_2)$ and $(\mu_1,\mu_2)$ in $NC^{(2)}(W)$. By definition, we will have $(\pi_1,\pi_2)\leq(\mu_1,\mu_2)$ if and only if $K_1^{\mu_2}(\mu_1)\leq_T K_1^{\pi_2}(\pi_1)$ and $K_1^c(\mu_2)\leq_T K_1^c(\pi_2)$. But since $K_1^c(\mu_2)\leq_T K_1^c(\pi_2)$ is equivalent to $\pi_2\leq_T\mu_2$, we conclude that
\begin{equation*}
(\pi_1,\pi_2)\leq (\mu_1,\mu_2)\quad\Longleftrightarrow\quad[K_1^{\pi_2}(\pi_1),\pi_2]\subseteq [K_1^{\mu_2}(\mu_1),\mu_2],
\end{equation*}
as desired.
\end{proof}

Thus the poset $NC^{(k)}(W)$ of $k$-divisible noncrossing partitions is simultaneously a generalization of both $NC(W)$ and $\Int(NC(W))$. Combining with the results of Chapter \ref{sec:classical}, this represents a two-fold generalization of Exercise 3.68(c) in Stanley~\cite{stanley:ec1}.

Note that Lemma \ref{lemma:intervals} allows us to classify the minimal elements of $NC^{(2)}(W)$. Clearly the minimal elements of $\Int(NC(W))$ are the ``singleton intervals'' $[\pi,\pi]$ for $\pi\in NC(W)$. Applying the isomorphism with $NC^{(2)}(W)$, we see that the minimal elements of $NC^{(2)}(W)$ are precisely the $2$-multichains of the form $(1,\pi)$ with $\pi\in NC(W)$; hence they are in natural bijection with the elements of $NC(W)$. One can observe this phenomenon in Figure \ref{fig:nc^2(a2)}. In general, we will soon show that the minimal elements of $NC^{(k+1)}(W)$ are in bijection with the elements of $NC^{(k)}(W)$ for all $k\in\posint$.

\subsection{Graded Semilattice}
\label{sec:semilattice}
We have observed that the poset $NC^{(k)}(W)$ is not, in general, a lattice since it does not have a minimum element. It turns out, however, that this is the only obstruction to this poset being a lattice, because any two elements in $NC^{(k)}(W)$ {\em do have a least upper bound}. To prove this, we will use the following lemma.

\begin{lemma}
\label{lemma:orderideal}
$NC_{(k)}(W)$ is an order ideal in the lattice $NC(W^k)$.
\end{lemma}

\begin{proof}
We will suppose that $k<\infty$. The proof for $k=\infty$ is analogous.

Let $(\delta)_k$ be a delta sequence, and consider an arbitrary sequence $(\varepsilon)_k\in NC(W^k)$ with the property that $\varepsilon_i\leq_T\delta_i$ for all  $1\leq i\leq k$. We wish to show in this case that $(\varepsilon)_k$ is also a delta sequence.

Since $\varepsilon_i\leq_T\delta_i$ for all $1\leq i\leq k$, there exist elements $\varepsilon_1',\varepsilon_2',\ldots,\varepsilon_k'$ with the property that $(\varepsilon_i,\varepsilon_i')$ is a minimal factorization of $\delta_i$ for all $1\leq i\leq k$. Then since $(\delta_1,\delta_2,\ldots,\delta_k)$ is a $k$-delta sequence by assumption, we conclude that the sequence
\begin{equation*}
(\varepsilon_1,\varepsilon_1',\varepsilon_2,\varepsilon_2',\ldots,\varepsilon_k,\varepsilon_k')
\end{equation*}
is a $2k$-delta sequence. Finally, using the shifting lemma \ref{shifting2}, we shift the elements $\varepsilon_i$ to the left to obtain a new $2k$-delta sequence
\begin{equation*}
(\varepsilon_1,\varepsilon_2,\ldots,\varepsilon_k,\varepsilon_1'',\varepsilon_2'',\ldots,\varepsilon_k''),
\end{equation*}
which implies that $(\varepsilon_1,\ldots,\varepsilon_k)$ is a $k$-delta sequence.
\end{proof}

Now recall that $NC(W^k)$ is a graded lattice with rank function
\begin{equation*}
\rk(w_1,w_2,\ldots,w_k)=\ell_T(w_1)+\ell_T(w_2)+\cdots+\ell_T(w_k).
\end{equation*}
This is just a property of the direct product of posets. If a poset $P$ contains meets but not joins, we say that $P$ is a {\sf meet-semilattice}. On the other hand, if $P$ possesses joins but not meets, we call $P$ a {\sf join-semilattice}. It turns out that $NC_{(k)}(W)$ and $NC^{(k)}(W)$ are both graded semilattices.

\begin{theorem}\hspace{.1in}
\label{th:semilattice}
\begin{enumerate}
\item $NC_{(k)}(W)$ is a graded meet-semilattice with rank function
\begin{equation*}
\rk(\delta)_k=n-\ell_T(\delta_0).
\end{equation*}
\item $NC^{(k)}(W)$ is a graded join-semilattice with rank function
\begin{equation*}
\rk(\pi)_k=\ell_T(\pi_1).
\end{equation*}
\end{enumerate}
\end{theorem}

\begin{proof}
Note that an order ideal $I\subseteq P$ in a graded lattice $P$ is a meet-semilattice. The existence of meets follows since, for $x$ and $y$ in $I$, the meet $x\wedge y$ in $P$ is below both $x$ and $y$, hence it must be in $I$. If $x$ and $y$ had a greater lower bound in $I$, this would be a greater lower bound in $P$. The fact that $I$ is graded follows because every interval in $I$ is also an interval in $P$. By Lemma \ref{lemma:orderideal}, $NC_{(k)}(W)$ is a graded meet-semilattice, and the rank function is equal to the rank function in $NC(W^k)$:
\begin{equation*}
\rk(\delta)_k=\sum_{i=1}^k \ell_T(\delta_i)=\ell_T(c)-\ell_T(\delta_0)=n-\ell_T(\delta_0).
\end{equation*}

The fact that $NC^{(k)}(W)$ is a join-semilattice then follows since $NC^{(k)}(W)$ is anti-isomorphic to $NC_{(k)}(W)$. Given $(\pi)_k\in NC^{(k)}(W)$, suppose that $\partial(\pi)_k=(\delta)_k\in NC_{(k)}(W)$, so that $\pi_1=\delta_0$. Then the rank function of $NC^{(k)}(W)$, initialized so that minimal elements have rank $0$, must be
\begin{equation*}
\rk(\pi)_k=n-\rk(\partial(\pi)_k)=n-(n-\ell_T(\delta_0))=\ell_T(\delta_0)=\ell_T(\pi_1).
\end{equation*}
\end{proof}

In particular, notice that the height of $NC^{(k)}(W)$ is equal to $n=\rk(W)$, which is the same as the height of $NC(W)$. That is, as $k$ increases, the poset $NC^{(k)}(W)$ {\em becomes ``wider'' but not ``taller''}.

We say that a poset $P$ is {\sf bounded} if it possesses a maximum element $\hat{1}$ and a minimum element $\hat{0}$. An elementary result in poset theory says that finite bounded semilattices are lattices \cite[Proposition 3.3.1]{stanley:ec1}. Thus if we formally adjoin a maximum element $\hat{1}$ to $NC_{(k)}(W)$ and a minimum element $\hat{0}$ to $NC^{(k)}(W)$, for $k<\infty$, we conclude that the resulting posets
\begin{equation*}
NC_{(k)}(W)\cup\{\hat{1}\}\quad\text{ and }\quad NC^{(k)}(W)\cup\{\hat{0}\}
\end{equation*}
are lattices. Observe that this is the case when $k=2$ and $W=A_2$.

Theorem \ref{th:semilattice} also allows us to characterize the minimal elements of $NC^{(k)}(W)$. Since $NC^{(k)}(W)$ is graded, these are just the elements $(\pi)_k$ with $\rk(\pi)_k=\ell_T(\pi_1)=0$, or $\pi_1=1\in W$. That is, the minimal elements of $NC^{(k)}(W)$ are precisely the $k$-multichains of the form $(1,\pi_2,\ldots,\pi_k)$. This implies the following.

\begin{corollary}
\label{cor:minimal}
For $k\geq 2$, the forgetful map $(\pi_1,\pi_2,\ldots,\pi_k)\mapsto (\pi_2,\ldots,\pi_k)$ is a bijection from the minimal elements of $NC^{(k)}(W)$ to the set $NC^{(k-1)}(W)$.
\end{corollary}

\begin{proof}
We have $\pi_1=1$, and the $(k-1)$-multichain $(\pi_2,\ldots,\pi_k)$ is unrestricted since $1\leq_T\pi$ for all $\pi\in NC(W)$.
\end{proof}

This bijection could be used to induce a partial order on the minimal elements of $NC^{(k)}(W)$. We do not know if this perspective is useful.

\subsection{Intervals and Order Ideals}
Because $NC^{(k)}(W)$ is anti-isomorphic to an order ideal $NC_{(k)}(W)\subseteq NC(W^k)$, every interval in $NC^{(k)}(W)$ is isomorphic to an interval in $NC(W^k)$, and we can describe these easily.

\begin{lemma}
Every interval in $NC^{(k)}(W)$ is isomorphic to
\begin{equation*}
NC(W_1)\times NC(W_2)\times\cdots\times NC(W_k)
\end{equation*}
for some parabolic subgroups $W_1,W_2,\ldots, W_k$ of $W$.
\end{lemma}

\begin{proof}
This follows since every interval in $NC(W^k)$ is isomorphic to $NC(W')$ for some parabolic subgroup $W'\subseteq W^k$ (Proposition \ref{prop:subintervals}), and parabolic subgroups of $W^k$ have the form $W_1\times W_2\times\cdots \times W_k$, where each $W_i$ is a parabolic subgroup of $W$. Finally, it is easy to show that
\begin{equation*}
NC(W_1\times\cdots\times W_k)=NC(W_1)\times\cdots\times NC(W_k).
\end{equation*}
\end{proof}

But notice that not all $k$-tuples of parabolics are possible here since the height of any interval in $NC^{(k)}(W)$ is bounded by $n=\rk(W)$. In fixed cases, it is possible to classify all $k$-tuples $W_1,W_2,\ldots,W_k$ of parabolics of $W$ with the property that
\begin{equation*}
\rk(W_1)+\rk(W_2)+\cdots+\rk(W_k)\leq \rk(W),
\end{equation*}
and Krattenthaler has used this idea in his computation of the generalized $M$-triangle for the exceptional types~\cite{krattenthaler1,krattenthaler2} (see Chapter \ref{sec:fusscatcomb}).

This takes care of the intervals in $NC^{(k)}(W)$. After this, we might be interested in the order ideals and order filters of the poset. If $P$ is a poset,  the {\sf principal order ideal} $\Lambda(x)$ and the {\sf principal order filter} $\filter(x)$ generated by $x\in P$ are defined as
\begin{eqnarray*}
\Lambda(x)&:=&\{ y\in P:y\leq x\}\quad\text{and}\\
\filter(x)&:=&\{ y\in P:x\leq y\}.
\end{eqnarray*}
That is, $\Lambda(x)$ is the set of elements below $x$, and $\filter(x)$ is the set of elements above $x$. The notation is meant to suggest the ``shape'' of these objects. 
Because $NC^{(k)}(W)$ contains a maximum element (the Coxeter element $(c,c,\ldots,c)\in W^k$), its principal order filters are just intervals and we already understand these. The principal order ideals, however, are more complicated. We will show that every principal order ideal in $NC^{(k)}(W)$ is isomorphic to $NC^{(k)}(W')$ for some parabolic subgroup $W'\subseteq W$.

The proof will require the following alternate characterization of the partial order $NC^{(k)}(W)$.

\begin{lemma}
\label{lemma:horizontal}
Consider $k$-multichains $(\pi)_k$ and $(\mu)_k$. We have $(\pi)_k\leq (\mu)_k$ in $NC^{(k)}(W)$ if and only if
\begin{equation*}
(\pi)_k\leq_{T^k}(\mu)_k\quad\text{and}\quad\mu_1\pi_1^{-1}\geq_T \mu_2\pi_2^{-1}\geq_T \cdots \geq_T \mu_k\pi_k^{-1}.
\end{equation*}
\end{lemma}

\begin{proof}
First we will show that $(\pi)_k\leq (\mu)_k$ in $NC^{(k)}(W)$ implies that $(\pi)_k\leq_{T^k}(\mu)_k$ (that is, componentwise) in $NC(W^k)$. So suppose that $(\pi)_k\leq (\mu)_k$ and consider delta sequences $\partial(\pi)_k=(\delta)_k$ and $\partial(\mu)_k=(\varepsilon)_k$, so that $\varepsilon_i\leq_T\delta_i$ for all $1\leq i\leq k$. We will be done if we can show that
\begin{equation*}
\delta_0\delta_1\cdots\delta_i\leq_T\varepsilon_0\varepsilon_1\cdots\varepsilon_i
\end{equation*}
for all $0\leq i\leq k-1$, since $\pi_i=\delta_0\cdots\delta_{i-1}$ and $\mu_i=\varepsilon_0\cdots\varepsilon_{i-1}$.

To do this, fix $0\leq i\leq k-1$ and consider the minimal $i$-factorization $(\delta_1,\delta_2,\ldots,\delta_i)$. Because $\varepsilon_j\leq_T \delta_i$ for all $1\leq j\leq i$, there exist elements $\varepsilon_1',\ldots,\varepsilon_i'$ such that
\begin{equation*}
(\varepsilon_1,\varepsilon_1',\varepsilon_2,\varepsilon_2',\ldots,\varepsilon_i,\varepsilon_i')
\end{equation*}
is a minimal $2i$-factorization of $\delta_1\cdots\delta_i$. Shifting the symbols $\varepsilon_j$ to the left (Lemma \ref{shifting2}), we obtain a new minimal $2i$-factorization of $\delta_1\cdots\delta_i$ of the form
\begin{equation*}
(\varepsilon_1,\varepsilon_2,\ldots,\varepsilon_i,\varepsilon_1'',\varepsilon_2'',\ldots,\varepsilon_i'')
\end{equation*}
and we conclude that $\varepsilon_1\cdots\varepsilon_i\leq_T\delta_1\cdots\delta_i$. Finally, if we apply the anti-automorphism $(K_1^c)^{-1}$ to both sides of this inequality, we get $\delta_0\delta_1\cdots\delta_i\leq_T\varepsilon_0\varepsilon_1\cdots\varepsilon_i$, as desired.

Now suppose that $(\pi)_k\leq(\mu)_k$, so that $(\pi)_k\leq_{T^k}(\mu)_k$ as we have just shown. In particular, fixing $1\leq i\leq k-1$, we get a square of the form
\begin{equation}
\label{eq:pimusquare}
\begin{array}{ccc}
\mu_i & \leq & \mu_{i+1}\\
\turnup{\leq} & & \turnup{\leq}\\
\pi_i & \leq & \pi_{i+1}
\end{array}
\end{equation}
and the following sequence of statements make sense:
\begin{equation}
\label{eq:statements}
\begin{array}{rcll}
\mu_i^{-1}\mu_{i+1} & \leq_T & \pi_i^{-1}\pi_{i+1}, & \qquad\text{apply $(K_1^{\mu_{i+1}})^{-1}$}\\
\mu_i & \geq_T  & \mu_{i+1}\pi_{i+1}^{-1}\pi_i, & \qquad\text{apply $K_{\pi_i}^{\mu_{i+1}}$}\\
\pi_i\mu_i^{-1}\mu_{i+1} & \leq_T  & \pi_{i+1}, & \qquad\text{apply $(K_1^{\mu_{i+1}})^{-1}$}\\
\mu_i\pi_i^{-1} & \geq_T & \mu_{i+1}\pi_{i+1}^{-1}. & 
\end{array}
\end{equation}
Hence $\mu_i\pi_i^{-1}\geq_T\mu_{i+1}\pi_{i+1}^{-1}$ as desired. Conversely, suppose that $(\pi)_k\leq_{T^k}(\mu)_k$ and that $\mu_i\pi_i^{-1}\geq_T\mu_{i+1}\pi_{i+1}^{-1}$ holds for all $1\leq i\leq k-1$. Fixing $1\leq i\leq k-1$, we have a square of the form \eqref{eq:pimusquare}, and hence it makes sense to run the statements \eqref{eq:statements} in reverse, proving $\mu_i^{-1}\mu_{i+1}\leq_T \pi_i^{-1}\pi_{i+1}$, as desired.
\end{proof}

Now we can classify the principal order ideals in $NC^{(k)}(W)$. Recall from Theorem \ref{th:coxelement} that every element $\pi\in NC(W)$ is a Coxeter element in some parabolic subgroup of $W$. Let us denote this parabolic subgroup by $W_\pi$ (we will see in Section \ref{sec:parabolic} that $W_\pi$ is unique).

\begin{theorem}
\label{th:orderideal}
Let $(\pi)_k$ be in $NC^{(k)}(W)$ with $\partial(\pi)_k=(\delta)_k$, and consider the parabolic subgroup $W'=W_{\pi_1}=W_{\delta_0}$ of $W$. Then we have:
\begin{enumerate}
\item The principal order ideal $\Lambda((\pi)_k)$ in $NC^{(k)}(W)$ is isomorphic to $NC^{(k)}(W')$ via the map
\begin{equation*}
(\mu_1,\mu_2,\ldots,\mu_k)\mapsto (\mu_1,\mu_2\pi_2^{-1}\pi_1,\mu_3\pi_3^{-1}\pi_1,\ldots,\mu_k\pi_k^{-1}\pi_1).
\end{equation*}
\label{idealitem:one}
\item The principal order filter $\filter((\delta)_k)$ in $NC_{(k)}(W)$ is isomorphic to $NC_{(k)}(W')$ via the map
\begin{equation*}
(\varepsilon_1,\varepsilon_2,\ldots,\varepsilon_k)\mapsto(\varepsilon_1\pi_2^{-1}\pi_1,\pi_1^{-1}\pi_2\varepsilon_2\pi_3^{-1}\pi_1,\ldots,\pi_1^{-1}\pi_{k-1}\varepsilon_{k-1}\pi_k^{-1}\pi_1,\pi_1^{-1}\pi_k\varepsilon_k).
\end{equation*}
\label{idealitem:two}
\end{enumerate}
\end{theorem}

\begin{proof}
We will prove only \eqref{idealitem:one}. Then \eqref{idealitem:two} follows from the commutative square:
\begin{equation*}
\xymatrix{                                                                    
\filter((\delta)_k) \ar[r] \ar[d]_{\int} & NC_{(k)}(W') \\                           
\Lambda((\pi)_k) \ar[r] & NC^{(k)}(W') \ar[u]_{\partial}\\                                                         
} 
\end{equation*}

Let $I$ denote the map in the statement of \eqref{idealitem:one}. First, we will show that $I$ sends $\Lambda((\pi)_k)$ to $NC^{(k)}(W')$, and that it preserves order. Suppose that $(\mu)_k$ is in $\Lambda((\pi)_k)$ so that $(\mu)_k\leq(\pi)_k$ in $NC^{(k)}(W)$. By Lemma \ref{lemma:horizontal}, we know that
\begin{equation}
\label{eq:bmultichain}
\pi_1\mu_1^{-1}\geq_T \pi_2\mu_2^{-1}\geq_T \cdots\geq_T \pi_k\mu_k^{-1}
\end{equation}
and $(\mu)_k\leq_{T^k}(\pi)_k$. In particular, since $\mu_1\leq_T\pi_1$, we have $\pi_1\geq_T (K_1^{\pi_1})^{-1}(\mu_1)=\pi_1\mu_1^{-1}$. Hence $\pi_1$ is above every element of the multichain \eqref{eq:bmultichain} and we may apply $K_1^{\pi_1}$ to get
\begin{equation*}
\mu_1\leq_T \mu_2\pi_2^{-1}\pi_1\leq_T\cdots\leq_T \mu_k\pi_k^{-1}\pi_1.
\end{equation*}
That is, $I((\mu)_k)$ is a multichain, and every entry of $I((\mu)_k)$ is below $\pi_1$ so it is in $NC^{(k)}(W')$ as desired.

To see that $I$ preserves order, consider $(\mu)_k\leq (\sigma)_k$ in $\Lambda((\pi)_k)$. For fixed $1\leq i\leq k$, Lemma \ref{lemma:horizontal} says that $\mu_i\leq_T\sigma_i\leq_T\pi_i$, hence $(K_1^{\pi_i})^{-1}(\sigma_i)\leq_T (K_1^{\pi_i})^{-1}(\mu_i)$, or $\pi_i\sigma_i^{-1}\leq_T \pi_i\mu_i^{-1}$. Applying $K_1^{\pi_1}$ to this inequality gives $\mu_i\pi_i^{-1}\pi_1\leq_T \sigma_i\pi_i^{-1}\pi_1$. Now consider $I((\mu)_k)$ with $i$-th entry $\mu_i\pi_i^{-1}\pi_1$ and $I((\sigma)_k)$ with $i$-th entry $\sigma_i\pi_i^{-1}\pi_1$ and note that
\begin{equation*}
(\sigma_i\pi_i^{-1}\pi_1)(\mu_i\pi_i^{-1}\pi_1)^{-1}=\sigma_i\mu_i^{-1}.
\end{equation*}
Lemma \ref{lemma:horizontal}, together with the relations
\begin{equation*}
\sigma_1\mu_1^{-1}\geq_T \sigma_2\mu_2^{-1}\geq_T \cdots \geq_T \sigma_k\mu_k^{-1}
\end{equation*}
(because $(\mu)_k\leq (\sigma)_k$) and the fact that absolute length restricts well to parabolic subgroups (Proposition \ref{prop:restriction}), implies that $I((\mu)_k)\leq I((\sigma)_k)$ in $NC^{(k)}(W')$, as desired. 

Now $I$ is clearly invertible as a map on $W^k$, with $I^{-1}$ given by
\begin{equation*}
(\mu_1,\mu_2,\ldots,\mu_k)\mapsto (\mu_1,\mu_2\pi_1^{-1}\pi_2,\mu_3\pi_1^{-1}\pi_3,\ldots,\mu_k\pi_1^{-1}\pi_k).
\end{equation*}
To complete the proof, we must show that $I^{-1}$ sends $NC^{(k)}(W')$ to $\Lambda((\pi)_k)$ and preserves order. If $(\mu)_k$ is in $NC^{(k)}(W')$, then we have $\mu_i\leq_T \pi_1$ for all $1\leq i\leq k$, and we may apply the map $(K_1^{\pi_1})^{-1}$ to get
\begin{equation*}
\pi_1\geq_T \pi_1\mu_1^{-1} \geq_T \pi_1\mu_2^{-1}\geq_T \cdots\geq_T \pi_1\mu_k^{-1}.
\end{equation*}
Then using Corollary \ref{cor:kreweras} and the fact that $\pi_1\leq_T\pi_i\leq_T \pi_{i+1}$, we conclude that
\begin{equation*}
\mu_i\pi_1^{-1}\pi_i=K_1^{\pi_i}(\pi_1\mu_i^{-1})\leq_T K_1^{\pi_i}(\pi_1\mu_{i+1}^{-1})\leq_T K_1^{\pi_{i+1}}(\pi_1\mu_{i+1}^{-1})=\mu_{i+1}\pi_1^{-1}\pi_{i+1}
\end{equation*}
for all $1\leq i\leq k-1$, which proves that $I^{-1}((\mu)_k)$ is a $k$-multichain. To verify that $I^{-1}((\mu)_k)$ is in $\Lambda((\pi)_k)$, we must show that $(\pi)_k\leq I^{-1}((\mu)_k)$, or $\pi_i^{-1}\pi_{i+1}\leq_T (\mu_i\pi_1^{-1}\pi_i)^{-1}(\mu_{i+1}\pi_1^{-1}\pi_{i+1})$ for fixed $1\leq i\leq k-1$. But clearly 
\begin{equation*}
\pi_1\mu_{i+1}^{-1}=(K_1^{\pi_1})^{-1}(\mu_{i+1})\leq_T \pi_1\leq_T \pi_i,
\end{equation*}
and all of this is below $\pi_{i+1}$ so we may apply $K_1^{\pi_{i+1}}$ to the inequality $\pi_1\mu_{i+1}^{-1}\leq_T\pi_i$ to get $\pi_i^{-1}\pi_{i+1}\leq_T \mu_{i+1}\pi_1^{-1}\pi_{i+1}$. Finally, starting with the fact that $\pi_1\mu_{i+1}^{-1}\leq_T \pi_1\mu_i^{-1}\leq_T\pi_i$, the following sequence of statements makes sense:
\begin{equation*}
\begin{array}{rcll}
\pi_1\mu_{i+1}^{-1} & \leq_T & \pi_1\mu_i^{-1},  \\
\mu_{i+1}\pi_1^{-1}\pi_i & \geq_T & \mu_i\pi_1^{-1}\pi_i,  \\
(\mu_{i+1}\pi_1^{-1}\pi_{i+1})(\pi_i^{-1}\pi_{i+1})^{-1} & \geq_T & \mu_i\pi_1^{-1}\pi_i,   \\
\pi_i^{-1}\pi_{i+1} & \leq_T & (\mu_i\pi_1^{-1}\pi_i)^{-1}(\mu_{i+1}\pi_1^{-1}\pi_{i+1}).
\end{array}
\end{equation*}
(Between the first and second lines, we apply $K_1^{\pi_i}$, and we apply $K_1^{\mu_{i+1}\pi_1^{-1}\pi_{i+1}}$ between the third and fourth.) The proof that $I^{-1}$ preserves order is similar to the above verification for $I$.
\end{proof}

The style of this tricky proof emphasizes the fact that delta sequences are homological in spirit, and they have their own sort of ``diagram chasing''. We made the effort because this result is extremely useful. It allows one to use the method of induction to prove results about the $k$-divisible noncrossing partitions $NC^{(k)}(W)$. We will use this later in Theorem \ref{th:ELNCk}.

\subsection{Meta-Structure}
\label{sec:meta}
Recall that $NC(W^\infty)$ is the lattice of semi-infinite sequences
\begin{equation*}
\{(w_1,w_2,w_3,\ldots): \,\,\,w_i\in NC(W)\text{ for all } i\in\posint\},
\end{equation*}
under the componentwise absolute order. If we fix some place-holder element $w'\in W$, then for every choice of index set $I=\{i_1,i_2,\ldots,i_k\}$ with $i_1<\cdots<i_k$, we can identify $NC(W^k)$ with the induced subposet of $NC(W^\infty)$ consisting of sequences $(w_1,w_2,w_3,\ldots)$ where $w_j=w'$ for $j\not\in I$. Moreover, if $k\leq\ell<\infty$, this perspective allows us to think of $NC(W^k)$ as an induced subposet of $NC(W^\ell)$ in $\binom{\ell}{k}$ distinct ways (one for each choice of index subset).

It turns out that this broader perspective also applies to the $k$-divisible noncrossing partitions. For example, consider Figure \ref{fig:nc_2(a2)} and Figure \ref{fig:nc^2(a2)}. Notice that $NC_{(2)}(A_2)$ contains two isomorphic copies of $NC_{(1)}(A_2)=NC(A_2)$\,: the first corresponds to delta sequences $(\delta_1,\delta_2)$ where $\delta_1=1$, and the second corresponds to delta sequences $(\delta_1,\delta_2)$ with $\delta_2=1$. Moreover, by duality, $NC^{(2)}(A_2)$ contains two isomorphic copies of $NC^{(1)}(A_2)=NC(A_2)$. In general,

\begin{lemma}
\label{lemma:contained}
For integers $1\leq k\leq \ell<\infty$, the poset $NC^{(\ell)}(W)$ contains $\binom{\ell}{k}$ distinct isomorphic copies of $NC^{(k)}(W)$, embedded as order filters.
\end{lemma}

\begin{proof}
We will prove the analogous result for the poset $NC_{(k)}(W)$ of $k$-delta sequences. For each subset $I$ of the indices $\{1,2,\ldots,\ell\}$ with $\abs{I}=k$, it is clear that the induced subposet of $NC_{(\ell)}(W)$ consisting of the delta sequences
\begin{equation*}
\left\{ (\delta)_\ell \in NC_{(\ell)}(W): \delta_j=1 \text{ for all } j\not\in I\right\}
\end{equation*}
is isomorphic to $NC_{(k)}(W)$. Furthermore, since $1\leq_T w$ for all $w\in W$, it is clear that this subposet is an order ideal. Now apply the anti-isomorphism $\sint$.
\end{proof}

In Chapter \ref{sec:classical}, we will see that this result yields new information about the classical $k$-divisible noncrossing partitions.

In summary, for every finite subset $I$ of the positive integers $\posint$, there is a corresponding poset of $\abs{I}$-divisible noncrossing partitions $NC^{(\abs{I})}(W,I)$, and we have
\begin{equation*}
NC^{(\abs{I})}(W,I)\subseteq NC^{(\abs{J})}(W,J)
\end{equation*}
as induced subposets if and only if $I\subseteq J$. The poset of finite subsets of $\posint$ under inclusion is called the {\sf infinite boolean lattice}. So the collection of posets 
\begin{equation*}
\left\{NC^{(\abs{I})}(W,I):I\subseteq\posint, \abs{I}<\infty \right\}
\end{equation*} forms an infinite boolean lattice under poset inclusion as order filters. It is possible to define many non-trivial induction and restriction maps on this lattice, and there may be interesting structure to discover here.

\subsection{Cover Relations}
Now we give a characterization of the cover relations in $NC^{(k)}(W)$ in terms of the lattice structure on $NC(W)$. Let $\wedge$ and $\vee$ denote the meet and join in $NC(W)$, respectively.

The key result is the following.

\begin{lemma}
\label{lemma:mcjoin}
If $(\pi)_k\leq (\mu)_k$ in $NC^{(k)}(W)$, then we have
\begin{equation*}
(\mu)_k= (\pi_1\vee w,\pi_2\vee w,\ldots,\pi_k\vee w),
\end{equation*}
where $w=\pi_1^{-1}\mu_1$
\end{lemma}

\begin{proof}
Consider $(\pi)_k\leq (\mu)_k$ in $NC^{(k)}(W)$ and recall from Lemma \ref{lemma:horizontal} that $(\pi)_k\leq_{T^k} (\mu)_k$. To show that $\mu_i=\pi_i\vee w$ for all $1\leq i\leq k$, we proceed by induction on $i$. First, Lemma \ref{lemma:joins} implies that
\begin{equation*}
\mu_1=\pi_1(\pi_1^{-1}\mu_1)=\pi_1\vee\pi_1^{-1}\mu_1=\pi_1\vee w.
\end{equation*}
Now suppose that $\mu_j=\pi_j\vee w$ for all $1\leq j\leq i< k$. We wish to show that $\mu_{i+1}=\pi_{i+1}\vee w$. Since $\pi_{i+1}\leq_T \mu_{i+1}$ and $w=K_1^{\mu_1}(\pi_1)\leq_T \mu_1\leq_T \mu_{i+1}$, we see that $\mu_{i+1}$ is an upper bound for $\pi_{i+1}$ and $w$, hence $\mu_{i+1}\geq_T \pi_{i+1}\vee w$. We are done if we can show that $\mu_{i+1}\leq_T \pi_{i+1}\vee w$.

So let $\partial(\pi)_k=(\delta)_k$ and $\partial(\mu)_k=(\varepsilon)_k$ denote the corresponding delta sequences. Using the fact that $\varepsilon_i\leq_T\delta_i$ and repeatedly applying Lemma \ref{lemma:joins}, we have
\begin{align*}
\mu_{i+1} &=\varepsilon_0\varepsilon_1\cdots\varepsilon_i\leq_T \varepsilon_0\varepsilon_1\cdots\varepsilon_i\vee \varepsilon_i^{-1}\delta_i =(\varepsilon_0\varepsilon_1\cdots\varepsilon_{i-1})\vee(\varepsilon_i\vee\varepsilon_i^{-1}\delta_i)\\
 &= \mu_i\vee\delta_i = (\pi_i\vee w)\vee \delta_i=(\pi_i\vee\delta_i)\vee w=\pi_{i+1}\vee w.
\end{align*}
\end{proof}

This shows that to go ``up'' in $NC^{(k)}(W)$ from a given multichain $(\pi)_k$, we should join a fixed element $w\in NC(W)$ to each entry of $(\pi)_k$. However, the lemma says nothing about {\em which} $w$ may be joined to a given $(\pi)_k$. Certainly some choices are not allowed. For example, consider the $2$-multichain $(1,(12))$ in $NC^{(2)}(A_2)$ (Figure \ref{fig:nc^2(a2)}). Joining with the element $(23)$ we get $(1\vee (23),(12)\vee (23))=((23),(123))$, which of course is a multichain, but it is {\em not} above $(1,(12))$.

In view of Lemmas \ref{lemma:horizontal} and \ref{lemma:mcjoin}, we will have $(\pi_1,\ldots,\pi_k)\leq (\pi_1\vee w,\ldots,\pi_k\vee w)$ in $NC^{(k)}(W)$ if and only if
\begin{equation}
\label{eq:wrelations}
\left(K_1^{\pi_1\vee w}\right)^{-1}(\pi_1)\geq_T \left(K_1^{\pi_2\vee w}\right)^{-1}(\pi_2)\geq_T \cdots\geq_T \left(K_1^{\pi_k\vee w}\right)^{-1}(\pi_k).
\end{equation}
However, it is quite difficult in general to classify elements $w$ satisfying \eqref{eq:wrelations}. (As an interesting side note, if $NC(W)$ were a {\sf distributive lattice} \cite[Chapter 3.4]{stanley:ec1}, the relations \eqref{eq:wrelations} would hold {\em for all} $w\in NC(W)$.) Instead, we will describe such $w$ implicitly by classifying the reflections they are built from.

If $P$ is a poset containing $x$ and $y$, we say that {\sf $y$ covers $x$}, and we write $x\prec y$, whenever $x<y$ and there does not exist any $z$ such that $x<z<y$. These {\sf cover relations} are familiar since they correspond to the edges in the Hasse diagram of $P$. To describe a partial order, it is sufficient to specify its set of cover relations (the partial order is then the {\sf transitive closure} of these relations).

For example, the cover relations in $NC(W)$ are easy to describe. Because the Hasse diagram of $\Abs(W)$ is isomorphic to the Cayley graph of $W$ with respect to $T$, there is an edge joining $\pi$ and $\mu$ in $NC(W)$ precisely when $\mu=\pi t$ for some $t\in T$. If, in addition, $\ell_T(\mu)=\ell_T(\pi)+1$, then we have $\pi\prec\mu$. That is, {\em the edges in the Hasse diagram of $NC(W)$ are naturally labelled by reflections}.

More generally, we can also describe the cover relations in $NC_{(k)}(W)$ since this is an order ideal in the lattice of noncrossing partitions $NC(W^k)$. There is an edge in the Hasse diagram of $NC_{(k)}(W)$ joining delta sequences $(\delta)_k$ and $(\varepsilon)_k$ precisely when $(\delta)_k^{-1}(\varepsilon)_k$ is in $T^k$; that is, if
\begin{equation*}
(\varepsilon)_k=(\delta_1,\delta_2,\ldots,\delta_i t,\ldots,\delta_k)
\end{equation*}
for some $1\leq i\leq k$ and $t\in T$. (Recall that the reflections $T^k$ in $W^k$ have the form $(1,1,\ldots,t,\ldots,1)$ for some $t\in T$.) If $\ell_T(\varepsilon_i)=\ell_T(\delta_i)+1$ we have $(\delta)_k\prec (\varepsilon)_k$. If, on the other hand, $\ell_T(\varepsilon_i)=\ell_T(\delta_i)-1$, we have $(\varepsilon)_k\prec (\delta)_k$. In either case, we call this a {\sf cover relation of index $i$}.

These observations, together with Lemma \ref{lemma:mcjoin}, allow us to characterize the cover relations in the $k$-divisible noncrossing partitions.

\begin{theorem}
\label{theorem:covers}
Let $(\pi)_k$ be in $NC^{(k)}(W)$. Then for each reflection $t\in T$ with $t\leq_T \pi_i^{-1}\pi_{i+1}$ there exists a cover relation
\begin{equation}
\label{eq:tcover}
(\pi)_k\prec (\pi_1\vee t',\pi_2\vee t',\ldots,\pi_k\vee t'),
\end{equation}
where $t'=(\pi_1^{-1}\pi_{i+1})\,t\,(\pi_1^{-1}\pi_{i+1})^{-1}$. Moreover, every cover relation in $NC^{(k)}(W)$ has this form. 
\end{theorem}

\begin{proof}
Since $NC_{(k)}(W)$ is an order ideal in $NC(W^k)$ we understand its cover relations. They have precisely the form
\begin{equation}
\label{eq:cover}
(\varepsilon)_k=(\delta_1,\delta_2,\ldots,\delta_i t,\ldots,\delta_k)\prec (\delta_1,\delta_2,\ldots,\delta_k),
\end{equation}
where $\ell_T(\delta_i t)=\ell_T(\delta_i)-1$. Note that the condition $\ell_T(\delta_i t)=\ell_T(\delta_i)-1$ is equivalent to $t\leq_t\delta_i$ by definition of the absolute order.
That is, we get one cover relation of index $i$ for each reflection below $\delta_i$.

Now suppose that $\sint(\delta)_k=(\pi)_k$ and $\sint(\varepsilon)_k=(\mu)_k$, so that, in particular, $\delta_i=\pi_i^{-1}\pi_{i+1}$. Then since $\pi_1=c(\delta_1\cdots\delta_k)^{-1}$ and $\mu_1=c(\delta_1\cdots\delta_i t\cdots\delta_k)^{-1}$, we have
\begin{align*}
\pi_1^{-1}\mu_1 & = (\delta_1\cdots\delta_k) c^{-1} c (\delta_1\cdots\delta_i t\cdots\delta_k)^{-1}\\
 & = (\delta_1\cdots\delta_i)\,t\,(\delta_1\cdots\delta_i)^{-1}\\
 & = (\pi_1^{-1}\pi_{i+1})\,t\,(\pi_1^{-1}\pi_{i+1})^{-1},
\end{align*}
and this is a reflection since $T$ is closed under conjugation. The result now follows from Lemma \ref{lemma:mcjoin}.
\end{proof}

Since the set of cover relations determines the partial order, Theorem \ref{theorem:covers} represents a complete characterization of the poset $NC^{(k)}(W)$ in terms of joins. That is, given $(\pi)_k$, we have $(\pi_1,\ldots,\pi_k)\leq (\pi_1\vee w,\ldots,\pi_k\vee w)$ if and only if $w$ has a reduced $T$-word $w=t_1't_2'\cdots t_r'$ where the reflections $t_1',t_2',\ldots,t_r'$ come from a chain of cover relations $(\pi)_k\prec\cdots\prec (\pi_1\vee w,\ldots,\pi_k\vee w)$.

There is also a nice description of the ``index'' of the cover relation \eqref{eq:tcover}. Notice that if $t\leq_T \pi_i^{-1}\pi_{i+1}$ for some $1\leq i\leq k$, {\em then this $i$ is  unique}. This is because $(\pi_1,\pi_1^{-1}\pi_2,\ldots,\pi_{k-1}^{-1}\pi_k,\pi_k^{-1}c)$ is a minimal factorization of the Coxeter element. If $\pi_i^{-1}\pi_{i+1}$ has a reduced $T$-word containing the symbol $t$, then none of the other entries may contain $t$, because reduced $T$-words may not contain repetition. By duality with $NC_{(k)}(W)$, we say that this cover relation has {\sf index $i$}. It is also easy to see that the index of the cover relation \eqref{eq:tcover} is equal to the smallest $i$ such that $t'\leq_T \pi_{i+1}$, or equivalently, $\pi_{i+1}=\pi_{i+1}\vee t'$. For example, the cover relation $(1,(12))\prec ((12),(12))$ in $NC^{(2)}(A_2)$ has index $1$ because it ``stays the same'' in the $2$nd position, and the cover relation $(1,(12))\prec ((13),(123))$ has index $2$ because all entries change (it stays the same in the ``third position''; recall Notation \ref{notation:extend}).

We will return to the idea of cover relations in Section \ref{sec:shelling} when we talk about edge-labellings of posets.

\subsection{Automorphisms}
\label{sec:automorphisms}

To end this section, we examine the group of poset automorphisms of $NC_{(k)}(W)$, which is isomorphic to the group of automorphisms of $NC^{(k)}(W)$. We say that $F$ is an {\sf automorphism} of the finite poset $P$ if $F$ is a self-bijection that preserves order. The fact that $F^{-1}$ preserves order then follows since $F$ has finite order under composition with itself.

For example, consider the dihedral group $I_2(m)$ with Coxeter generating set $\{a,b\}$, and consider the standard Coxeter element $c=ab$. (We think of $a$ and $b$ as adjacent reflections of the regular $m$-gon, and $c$ as a rotation.) The Hasse diagram of $NC(I_2(m))$ with respect to the Coxeter element $c$ looks like
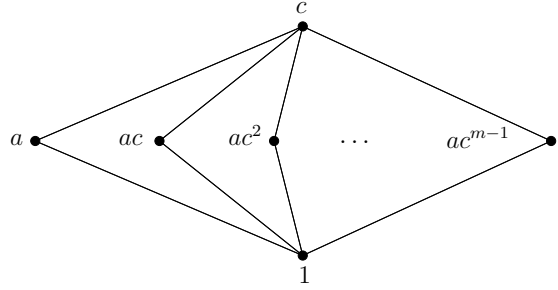
\begin{center}
\scalebox{.8}{\input{i_2_m.pstex_t}}
\end{center}
It is clear that the group of automorphisms of the abstract poset $NC(I_2(m))$ consists of {\em all} permutations on the set of $m$ reflections. However, it may happen that some of these automorphisms are not relevant to the structure theory of noncrossing partitions. In general, we will isolate a certain subgroup of automorphisms that is algebraically and combinatorially significant.
 
Because the Coxeter diagram of a finite irreducible Coxeter group $W$ is a tree, there exists a unique partition of the generators into two sets $S=S_\ell\sqcup S_r$ such that the elements of $S_\ell$ commute pairwise, as do the elements of $S_r$. (Think of a proper $2$-coloring of the Coxeter graph.) Let $\ell$ denote the product of the generators $S_\ell$ (in any order) and let $r$ denote the product of the generators $S_r$, so that (in particular) $\ell^2=1$ and $r^2=1$. Clearly the product $c=\ell r$ is a Coxeter element of $W$, and in this case we call the triple $(c,\ell,r)$ a {\sf bipartite Coxeter element}. Recall that the Coxeter elements of $W$ form a conjugacy class. Then, since all simple generating sets are conjugate, we see that every Coxeter element of $W$ can be expressed as a bipartite Coxeter element for some choice of $S$.

Now define two maps on $NC(W)$ with respect to the Coxeter element $c=\ell r$ by setting
\begin{equation}
\label{eq:LandR}
R(\pi):=r\pi^{-1}r\quad\text{ and }\quad L(\pi):=\ell\pi^{-1}\ell
\end{equation}
for all $\pi\in NC(W)$. It turns out that these maps are automorphisms.
 
\begin{lemma}
\label{lemma:LandR}
The maps $L$ and $R$ are poset automorphisms of $NC(W)$.
\end{lemma}
 
\begin{proof}
We will prove the result for $R$. The proof for $L$ is similar.
 
First notice that the map $R:W\to W$ is invertible with $R^{-1}=R$. We must show that $R$ maps $NC(W)$ to itself, and that it preserves order.

So consider $\pi$ and $\mu$ in $NC(W)$ with $\pi\leq_T\mu\leq_T c$. By reversing the reduced $T$-words for these elements, the subword property of the absolute order \ref{subword} implies that $\pi^{-1}\leq_T\mu^{-1}\leq_T c^{-1}=r\ell$. Then since conjugation by any fixed element --- in particular, by $r$ --- is an automorphism of the absolute order, we have
\begin{equation*}
r\pi^{-1}r\leq_T r\mu^{-1}r\leq_T r(r\ell)r=\ell r =c,
\end{equation*}
which proves the result.
\end{proof}
 
Observe that the composition $L\circ R$ is the same as conjugation by $c$,
\begin{equation*}
L\circ R\,(\pi)=\ell r \pi r\ell= c\pi c^{-1},
\end{equation*}
which is an automorphism of order $h$, so the maps $L$ and $R$ generate a dihedral group of automorphisms of $NC(W)$ of order $2h$. Returning to the example $W=I_2(m)$, notice that the triple $(c,a,b)$ is a bipartite Coxeter element. Since $L$ and $R$ are automorphisms, they fix the minimum element $1$ and the maximum element $c$ of $NC(I_2(m))$, so we need only investigate the action of $L$ and $R$ on the set of reflections
\begin{equation*}
T=\{a,ac,ac^2,\ldots,ac^{m-1}\}.
\end{equation*}
Here we have $L(ac^i)=ac^{-i}aa=ac^{-i}$ and $R(ac^i)=bc^{-i}ab=ac^{-i+2}$ for all $1\leq i\leq m$. Thus the maps $L$ and $R$ generate a dihedral group of automorphisms (acting on the $m$-gon with vertices labelled $a,ac,ac^2,\ldots,ac^{m-1}$, clockwise; this is {\em not} the same $m$-gon on which $I_2(m)$ acts). In general, this group is much smaller than the full automorphism group of $NC(I_2(m))$.

Our motivation for considering the group $\langle L,R\rangle$ comes from the case of the symmetric group --- type $A_{n-1}$. There is a nice way to represent the elements of $NC(A_{n-1})$ pictorially (see Figure \ref{fig:noncrossing}), and in this case the group $\langle L,R\rangle$ is just the dihedral group of motions acting on the picture. We will discuss this in Section \ref{sec:classicalautomorphisms}.
 
Now we define a generalization of this dihedral action for all positive integers $k$. Recall that $NC_{(k)}(W)$ is the set of $k$-delta sequences under the componentwise partial order.

\begin{definition}
For all $(\delta)_k=(\delta_1,\ldots,\delta_k)\in NC_{(k)}(W)$ define
\begin{align}
\label{eq:R*}
R^*(\delta)_k & :=\left(R(\delta_k),R(\delta_{k-1}),\ldots,R(\delta_2),R(\delta_1)\right)\quad\text{and} \\
\label{eq:L*}
L^*(\delta)_k &:=\left(L(\delta_1),R(\delta_k),R(\delta_{k-1}),\ldots,R(\delta_2)\right).
\end{align}
\end{definition}
When $k=1$, notice that $L^*$ and $R^*$ coincide with $L$ and $R$, respectively. We claim that $L^*$ and $R^*$ are automorphisms in general.
 
\begin{lemma}
\label{lemma:L*andR*}
The maps $L^*$ and $R^*$ are poset automorphisms of $NC_{(k)}(W)$.
\end{lemma}
 
\begin{proof}
Note that $L^*$ and $R^*$ are invertible on $W^k$ with $(L^*)^{-1}=L^*$ and $(R^*)^{-1}=R^*$, and that they preserve componentwise order by Lemma \ref{lemma:LandR}. We must show that $R^*(\delta)_k$ and $L^*(\delta)_k$ are delta sequences.

First we will show that $R^*(\delta)_k$ is a delta sequence. If $K:NC(W)\to NC(W)$ denotes the Kreweras complement $K_1^c$, then by Lemma \ref{lemma:delta} $(2)$, we must show that $R(\delta_i)\leq_T K\left( R(\delta_j)\right)$ for all $1\leq i<j\leq k$. So fix $1\leq i<j\leq k$ and note that the following sequence of relations are equivalent:
\begin{eqnarray*}
R(\delta_i) & \leq_T & K\left( R(\delta_j)\right), \\
K^{-1}\left( R(\delta_i)\right) & \geq_T & R(\delta_j), \\
R\left( K^{-1}\left( R(\delta_i)\right) \right) &\geq_T &\delta_j, \\
r\left( c(r\delta_i^{-1} r)^{-1}\right)^{-1} r &\geq_T & \delta_j, \\ 
\delta_i^{-1} c &\geq_T & \delta_j,\\
K(\delta_i) &\geq_T &\delta_j.
\end{eqnarray*}
The final statement is true since $(\delta)_k$ is a delta sequence. Then to see that $L^*(\delta)_k$ is a delta sequence, fix $1<i\leq k$ and note further that the following sequence of relations are equivalent:
\begin{eqnarray*}
R(\delta_i) &\leq_T & K\left( L(\delta_1)\right), \\
\delta_i &\leq_T & R\left( K\left( L(\delta_1)\right) \right), \\
\delta_i &\leq_T & r\left( (\ell\delta_1^{-1}\ell)^{-1}c \right)^{-1} r, \\
\delta_i &\leq_T & \delta_1^{-1} c,\\
\delta_i &\leq_T & K(\delta_1).
\end{eqnarray*}
Again, the last statement is true since $(\delta)_k$ is a delta sequence.
\end{proof}

Finally, consider the composition of $L^*$ and $R^*$.

\begin{definition}
\label{def:C*}
Let $C^*$ denote the composition $L^*\circ R^*$ given by
\begin{eqnarray}
\label{eq:L*R*}
C^* (\delta)_k & =& \left( L\circ R(\delta_k),R^2(\delta_1),R^2(\delta_2),\ldots,R^2(\delta_{k-1})\right)\\
&=& (c\delta_k c^{-1},\delta_1,\delta_2,\ldots,\delta_{k-1}), \nonumber
\end{eqnarray}
for all $k$-delta sequences $(\delta)_k\in NC_{(k)}(W)$.
\end{definition}

This $C^*$ is an automorphism of $NC_{(k)}(W)$ (by Lemma \ref{lemma:L*andR*}), and it evidently has order $kh$. Hence $L^*$ and $R^*$ generate a dihedral group of poset automorphisms of $NC_{(k)}(W)$ with order $2kh$. We will see in Section \ref{sec:classicalautomorphisms} that this dihedral action also has a natural pictorial interpretation in the case $W=A_{n-1}$.

\section{Fuss-Catalan and Fuss-Narayana Numbers}
\label{sec:fusscat}
Now we return to a discussion of numerology. In the last section, we examined the basic structural properties of $NC^{(k)}(W)$; here we will consider enumerative questions. We will present a series of results, all of them proved in a case-by-case manner, that suggest hidden depth in the definition of $NC^{(k)}(W)$.

Recall from Section \ref{sec:invariant} the definition of the degrees $d_1\leq d_2\leq\cdots\leq d_n$, the exponents $m_1\leq m_2\leq\cdots\leq m_n$, and the Coxeter number $h=d_n$ of the rank $n$ finite Coxeter group $W$. Since the $k$-divisible noncrossing partitions $NC^{(k)}(W)$ are a generalization of $NC(W)$, their cardinality is a generalization of the Coxeter-Catalan number.

\begin{definition}
Define the {\sf Fuss-Catalan polynomial} associated to $W$ by
\begin{equation}
\label{eq:fusscat}
\Cat^{(k)}(W):=\prod_{i=1}^n \frac{kh+d_i}{d_i}=\frac{1}{\abs{W}}\prod_{i=1}^n (kh+d_i).
\end{equation}
If $k$ is a positive integer, we call $\Cat^{(k)}(W)$ the {\sf Fuss-Catalan number}.
\end{definition}

It is evident that this formula is a very natural generalization of the Coxeter-Catalan number $\Cat(W)=\Cat^{(1)}(W)$. In order to distinguish these among the myriad generalizations of the Catalan numbers, we have chosen the name ``Fuss-Catalan'', which has the advantages of memorability and historical accuracy. When Leonhard Euler became nearly blind after eye surgery in 1772, he wrote to Daniel Bernoulli asking him to send a mathematical assistant. It was the young Niklaus Fuss who accepted the job, and worked with Euler for ten years until the elder's death in 1783. The two had such a close relationship that Fuss even married Euler's granddaughter, and he was the author of the most famous eulogy to Euler (see \cite{ranjan}). Fuss was also a successful mathematician in his own right, and in 1791 he published a paper \cite{fuss} proving that the number of ways to dissect a convex $(kn+2)$-gon into $(k+2)$-gons is given by
\begin{equation}
\label{eq:typeAfusscat}
\Cat^{(k)}(A_{n-1})= \frac{1}{n}\binom{(k+1)n}{n-1}.
\end{equation}
Interestingly, this publication predated the birth of Eug\`ene Charles Catalan by 23 years! In the modern literature, there is no standard notation for \eqref{eq:typeAfusscat}, but the names ``Fuss numbers'' and ``Fuss-Catalan numbers'' have both commonly been used. For details, see the historical references in Przytycki and Sikora~\cite{p-s:fuss} or Fomin and Reading~\cite{fomin-reading}.

We will refer to the combinatorics that surrounds the numbers $\Cat^{(k)}(W)$ as the {\sf Fuss-Catalan combinatorics} of $W$. Chapter \ref{sec:fusscatcomb} is devoted exclusively to this topic.

To prove that the Fuss-Catalan numbers count $k$-divisible noncrossing partitions, we need to count multichains. If $P$ is a poset, let $\Z(P,k)$ denote the number of $k$-multichains $x_1\leq x_2\leq\cdots\leq x_k$ in $P$. It is well-known that this quantity is polynomial in $k$, and the expression $\Z(P,k)$ is known as the {\sf zeta polynomial} of the poset. Actually, most authors refer to $\Z(P,k+1)$ as the ``zeta polynomial'' (see \cite[Chapter 3.11]{stanley:ec1}) since this counts multichains of {\em length $k$} (our $(k+1)$-multichains). We hope that our use will be clear in context and that no confusion will result from the change of index.

The following result appeared for the first time as Proposition 9 in Chapoton~\cite{chapoton:one}, where it was verified case-by-case. In the classical types, the result is equivalent to formulas of Edelman \cite[Theorem 4.2]{edelman:kdiv}, Reiner \cite[Proposition 7]{reiner}, and Athanasiadis and Reiner \cite[Theorem 1.2 $(iii)$]{athanasiadis-reiner}. Reiner verified the exceptional types by computer.

\begin{theorem}
\label{th:NCzeta}
The zeta polynomial of $NC(W)$ is equal to the Fuss-Catalan polynomial of $W$,
\begin{equation*}
\Z(NC(W),k)=\Cat^{(k)}(W).
\end{equation*}
\end{theorem}

This immediately implies
\begin{theorem}
\label{th:NCkcount}
When $k$ is a finite positive integer, we have
\begin{equation*}
\abs{NC^{(k)}(W)}=\Cat^{(k)}(W).
\end{equation*}
\end{theorem}

We emphasize that, at this writing, {\em a uniform proof of this fact is not known}, despite the elegance of the formula for $\Cat^{(k)}(W)$. This is an important open problem.

The Fuss-Catalan polynomial is just one enumerative invariant associated to the $k$-divisible noncrossing partitions. Because the poset $NC^{(k)}(W)$ is graded (Theorem \ref{th:semilattice}), we should also consider its rank numbers. The rank numbers of $NC(W)$ are commonly known as the ``Narayana numbers'', so we suggest the following notation.

\begin{definition}
\label{def:fussnar}
For each $0\leq i\leq n$, define the {\sf Fuss-Narayana polynomial}
\begin{equation}
\label{eq:fussnar}
\Nar^{(k)}(W,i):=\#\left\{ (\pi)_k\in NC^{(k)}(W): \rk(\pi)_k=\ell_T(\pi_1)=i\right\}.
\end{equation}
When $k\in\posint$, we call $\Nar^{(k)}(W,i)$ the {\sf Fuss-Narayana number}.
\end{definition}

We will shortly verify that $\Nar^{(k)}(W,i)$ is indeed a polynomial in $k$, but first let us consider some basic cases. When $k=1$, the poset $NC^{(1)}(W)=NC(W)$ is self-dual, hence $\Nar^{(1)}(W,i)=\Nar^{(1)}(W,n-i)$ for all $1\leq i\leq n$. And since $NC(W)$ has a maximum and minimum element, we have $\Nar^{(1)}(W,n)=\Nar^{(1)}(W,0)=1$. In general, the poset $NC^{(k+1)}(W)$ has a maximum element and its minimal elements are in bijection with $NC^{(k)}(W)$ (Corollary \ref{cor:minimal}), giving
\begin{equation*}
\Nar^{(k+1)}(W,n)=1\qquad\text{ and }\qquad \Nar^{(k+1)}(W,0)=\Cat^{(k)}(W).
\end{equation*}
We will see, however, that the Fuss-Narayana numbers are quite a bit less straightforward for $1<i<n$.

Because the Fuss-Narayana polynomials count certain $k$-multichains in $NC(W)$, it is easy to show that they are polynomials in $k$. Moreover, we can show that they have rational coefficients that alternate in sign.

\begin{theorem}
\label{th:fussnarpoly}
For $k\in\posint$ and $1\leq i\leq n$, we have
\begin{enumerate}
\item $\Nar^{(k)}(W,i)$ is a polynomial in $k$ of degree $n-i$.
\label{naritem:one}
\item The associated polynomial
\begin{equation*}
(-1)^{n-i}\cdot\abs{W}\cdot\Nar^{(-k)}(W,i)
\end{equation*}
has positive integer coefficients in all degrees from $0$ to $n-i$.
\label{naritem:two}
\end{enumerate}
\end{theorem}

\begin{proof}
If $\pi$ is in $NC(W)$, let $W_\pi$ denote a parabolic subgroup of $W$ with Coxeter element $\pi$ (Theorem \ref{th:coxelement}), and note that the interval $[\pi,c]$ in $NC(W)$ is isomorphic to $NC(W_{\pi^{-1}c})$ (Proposition \ref{prop:subintervals}). Since $\Nar^{(k)}(W,i)$ counts the number of $k$-multichains in $NC(W)$ whose bottom element has rank $i$, it is equal to the sum of zeta polynomials
\begin{equation*}
\Nar^{(k)}(W,i)=\sum_{\{\pi\in NC(W):\ell_T(\pi)=i\}} \Z([\pi,c],k-1),
\end{equation*}
where $\Z([\pi,c],k-1)=\Cat^{(k-1)}(W_{\pi^{-1}c})$ by Theorem \ref{th:NCzeta}. Then \eqref{naritem:one} follows since, for all $\pi\in NC(W)$ with $\ell_T(\pi)=i$, $\Cat^{(k-1)}(W_{\pi^{-1}c})$ is evidently a polynomial in $k$ with degree $n-i$ and positive leading coefficient.

Now recall that the exponents $m_1,m_2,\ldots,m_n$ of $W$ are a permutation of the numbers $h-m_1,h-m_2,\ldots,h-m_n$ when $W$ is a real reflection group. Applying this to the formula for the Fuss-Catalan number, we get the polynomial
\begin{equation*}
(-1)^n\cdot\Cat^{(-k-1)}(W)=\frac{1}{\abs{W}}\prod_{j=1}^n (kh+d_j-2)
\end{equation*}
which has positive coefficients in all degrees since $d_j\geq 2$ for all $1\leq j\leq n$. Finally, we conclude that the polynomial
\begin{equation*}
(-1)^{n-i}\cdot\abs{W}\cdot\Nar^{(-k)}(W,i)=\sum_{\{\pi\in NC(W):\ell_T(\pi)=i\}} \frac{\abs{W}}{\abs{W_{\pi^{-1}c}}} \Cat^{(-k-1)}(W_{\pi^{-1}c})
\end{equation*}
has positive integer coefficients in all degrees from $0$ to $n-i$, proving \eqref{naritem:two}.
\end{proof}

At the risk of lessening the suspense, we now present:

\begin{figure}
\vspace{.1in}
\begin{center}
\scalebox{0.82}{
\begin{tabular}{|c|c|l|}
\hline
$W$ & $i$ & $\Nar^{(k)}(W,i)$\\
\hline
\hline
$A_{n-1}$ & $i$ & $\textstyle{\frac{1}{n}\binom{n}{i}\binom{kn}{n-i-1}}$\\
\hline
$B_n$ & $i$ & $\textstyle{\binom{n}{i}\binom{kn}{n-i}}$\\
\hline
$D_n$ & $i$ & $\textstyle{\binom{n}{i}\binom{k(n-1)}{n-i}+\binom{n-2}{i}\binom{k(n-1)+1}{n-i}}$\\
\hline
\multirow{3}{*}{$I_2(m)$}
& $2$ & $1$\\
& $1$ & $mk$\\
& $0$ & $k(mk-m+2)/2$\\
\hline
\multirow{4}{*}{$H_3$}
& $3$ & $1$\\
& $2$ & $15k$\\
& $1$ & $5k(5k-2)$\\
& $0$ & $k(5k-2)(5k-4)/3$\\
\hline
\multirow{5}{*}{$H_4$}
& $4$ & $1$\\
& $3$ & $60k$\\
& $2$ & $k(465k-149)/2$\\
& $1$ & $15k(3k-1)(5k-3)$\\
& $0$ & $k(3k-1)(5k-3)(15k-14)/4$\\
\hline
\multirow{5}{*}{$F_4$}
& $4$ & $1$\\
& $3$ & $24k$\\
& $2$ & $k(78k-23)$\\
& $1$ & $12k(3k-1)(2k-1)$\\
& $0$ & $k(3k-1)(2k-1)(6k-5)/2$\\
\hline
\multirow{7}{*}{$E_6$}
& $6$ & $1$\\
& $5$ & $36k$\\
& $4$ & $12k(21k-4)$\\
& $3$ & $9k(4k-1)(18k-5)$\\
& $2$ & $2k(4k-1)(3k-1)(30k-13)$\\
& $1$ & $6k(4k-1)(3k-1)(2k-1)(12k-7)/5$\\
& $0$ & $k(4k-1)(3k-1)(2k-1)(12k-7)(6k-5)/30$\\
\hline
\multirow{8}{*}{$E_7$}
& $7$ & $1$\\
& $6$ & $63k$\\
& $5$ & $21k(63k-11)/2$\\
& $4$ & $21k(9k-2)(27k-7)/2$\\
& $3$ & $21k(9k-2)(3k-1)(63k-23)/8$\\
& $2$ & $3k(9k-2)(3k-1)(9k-4)(207k-103)/40$\\
& $1$ & $9k(9k-2)(3k-1)9k-4)(9k-5)(3k-2)/40$\\
& $0$ & $k(9k-2)(3k-1)(9k-4)(9k-5)(3k-2)(9k-8)/280$\\
\hline
\multirow{9}{*}{$E_8$}
& $8$ & $1$\\
& $7$ & $120k$\\
& $6$ & $35k(105k-17)/2$\\
& $5$ & $45k(5k-1)(45k-11)$\\
& $4$ & $k(5k-1)(10350k^2-6675k+1084)/2$\\
& $3$ & $15k(5k-1)(3k-1)(5k-2)(30k-13)$\\
& $2$ & $5k(5k-1)(3k-1)(5k-2)(15k-8)(195k-107)/48$\\
& $1$ & $5k(5k-1)(3k-1)(5k-2)(15k-8)(5k-3)(15k-11)/56$\\
& $0$ & $k(5k-1)(3k-1)(5k-2)(15k-8)(5k-3)(15k-11)(15k-14)/1344$\\
\hline
\end{tabular}}
\end{center}
\caption{Fuss-Narayana polynomials for the finite irreducible Coxeter groups}
\label{fig:fussnar}
\end{figure}

\begin{theorem}
\label{th:fussnar}
Figure \ref{fig:fussnar} contains the complete list of Fuss-Narayana polynomials for the finite irreducible Coxeter groups.
\end{theorem}

\begin{proof}
The case of the dihedral group $I_2(m)$ is trivial, since we already know that $\Nar^{(k)}(I_2(m),2)=1$ and $\Nar^{(k)}(I_2(m),0)=\Cat^{(k-1)}(I_2(m))=k(mk-m+2)/2$. We then use the fact that
\begin{equation*}
\Cat^{(k)}(I_2(m))=\sum_{i=0}^2 \Nar^{(k)}(I_2(m),i).
\end{equation*}
Each of the formulas for the groups $A_{n-1}$, $B_n$ and $D_n$ represents a theorem, and these appear in Chapter \ref{sec:classical} as Theorem \ref{th:fussnarA}, Theorem \ref{th:fussnarB} and Theorem \ref{th:fussnarD}, respectively. The Fuss-Narayana polynomials for the exceptional groups have been computed in \verb|Maple|, using John Stembridge's \verb|posets| and \verb|coxeter| packages \cite{stembridge}. The procedures are available from the author upon request.
\end{proof}

Perhaps the most remarkable thing about this chart is the fact that the same polynomials have been observed independently by Fomin and Reading~\cite{fomin-reading} and Athanasiadis\footnote{Athanasiadis' results apply only in the case of Weyl groups, and it is not clear in this case that the numbers are polynomial in $k$.}~\cite{athanasiadis:nar}, in very different circumstances (see Chapter \ref{sec:fusscatcomb}).

Given a complete classification like this, we can now make several observations that immediately achieve the status of ``theorems'', although they are more like ``true conjectures''.

\begin{theorem}
If $i\in\{0,1,n-1,n\}$ and the degrees $d_1\leq d_2\leq\cdots\leq d_n$ are arranged in increasing order, then we have
\begin{equation}
\label{eq:fussnarguess}
\Nar^{(k)}(W,i)=\prod_{j=0}^{n-i} \frac{kh-d_j+2}{d_j}.
\end{equation}
\end{theorem}
That is, we have
\begin{equation*}
\begin{array}{lcl}
\Nar^{(k)}(W,n) & = & 1,\\
\Nar^{(k)}(W,n-1) & = & nk,\\
\Nar^{(k)}(W,1) & = & \frac{h}{\abs{W}}\prod_{j=0}^{n-1} (kh-d_j+2),\quad\text{and}\\
\Nar^{(k)}(W,0) & = & \Cat^{(k-1)}(W)
\end{array}
\end{equation*}
for all finite Coxeter groups $W$. The last two formulas follow from the facts that $\abs{W}$ is equal to the product of the degrees, and that the set $\{m_1,\ldots,m_n\}$ of exponents coincides with the set $\{h-m_1,\ldots,h-m_n\}$, where $h$ is the Coxeter number.

However, observe that something very strange happens for $2\leq i\leq n-2$. It seems that formula \eqref{eq:fussnarguess} is ``almost correct'' in the sense that $\Nar^{(k)}(W,i)$ is ``almost divisible'' over $\rational[k]$ by the factor
\begin{equation*}
\prod_{j=0}^{n-i} (kh-d_j+2),
\end{equation*}
but the precise way in which this fails is difficult to understand. In general, all but one of the irreducible factors of $\Nar^{(k)}(W,i)$ is explained by this formula. In fact, it would be a theorem that ``$\Nar^{(k)}(W,i)$ splits in $\rational[k]$'', except for the unique exception of $\Nar^{(k)}(E_8,4)$, which has a single irreducible quadratic factor! Fomin and Reading~\cite{fomin-reading} (who independently observed this phenomenon) made a table of correction factors to quantify by how much formula \eqref{eq:fussnarguess} fails, but they were unable to explain these factors.

\begin{problem}
Find a uniform formula for $\Nar^{(k)}(W,i)$.
\end{problem}

Another ``true conjecture'' is the following.

\begin{theorem}
\label{th:averagerank}
If an element of $NC^{(k)}(W)$ is chosen uniformly at random, its expected rank is $n/(k+1)$.
\end{theorem}

\begin{proof}
This is equivalent to the formula
\begin{equation*}
\frac{1}{\Cat^{(k)}(W)} \sum_{i=0}^n i\cdot\Nar^{(k)}(W,i) = \frac{n}{k+1},
\end{equation*}
which can be observed from the data in Figure \ref{fig:fussnar}.
\end{proof}

When $k=1$, this theorem tells us that the expected rank is $n/2$, which also follows, for example, from the fact that $NC(W)$ is self-dual. However, in the general case, it is surprising that the average rank has such a simple expression. In particular, this theorem says that the elements of $NC^{(k)}(W)$ become more concentrated near the bottom as $k$ grows.

We can imagine a uniform proof of this fact if we think of the Fuss-Narayana numbers $\Nar^{(k)}(W,i)$ as the $h$-vector of Fomin and Reading's {\sf generalized cluster complex} $\Delta^{(k)}(W)$ (see Chapter \ref{sec:fusscatcomb}). They showed that each codimension $1$ face of $\Delta^{(k)}(W)$ is contained in exactly $(k+1)$ top-dimensional faces \cite[Proposition 3.10]{fomin-reading}. Assuming that the complex $\Delta^{(k)}(W)$ is {\sf shellable}, Theorem \ref{th:averagerank} would follow. This problem is open, but Tzanaki has constructed a shelling in the classical types $A$ and $B$~\cite{tzanaki}, and she is working on the general problem\footnote{ Athanasiadis and Tzanaki \cite{athanasiadis-tzanaki} have recently costructed a shelling in general, and shown that the complex $\Delta^{(k)}(W)$ is $(k+1)$-Cohen-Macaulay.}. Of course, such a proof would give no insight about the $k$-divisible noncrossing partitions without some understood relationship between $NC^{(k)}(W)$ and $\Delta^{(k)}(W)$, which we currently do not have.

It is also interesting to consider the unimodality of rank sequences. A sequence of nonzero integers $\{\alpha_i\}_{i=0}^n$ is said to be {\sf unimodal} if there exists $0\leq i\leq n$ such that
\begin{equation*}
\alpha_1\leq\alpha_2\leq\cdots\leq \alpha_i\geq\cdots\geq\alpha_n,
\end{equation*}
and in this case we say that $i$ is a {\sf mode} of the sequence.

\begin{theorem}
\label{th:unimodal}
The sequence $\{ \Nar^{(k)}(W,i)\}_{i=0}^n$ is unimodal for all $k\in\posint$.
\end{theorem}

\begin{proof}
One can observe case-by-case from Figure \ref{fig:fussnar} that, for each $1\leq i\leq n$, the polynomial
\begin{equation*}
(\Nar^{(k)}(W,i))^2-\Nar^{(k)}(W,i-1)\cdot\Nar^{(k)}(W,i+1),
\end{equation*}
is positive for all real numbers $k\geq 1$ (consider the derivative). This shows, in particular, that the sequence is log-concave, and it is well-known that this implies unimodality.
\end{proof}

There are two interesting questions related to unimodality.

\begin{problem}
Where is the mode of the sequence $\{\Nar^{(k)}(W,i)\}_{i=0}^n$?
\end{problem}

We suspect that the mode $j$ satisfies $j\geq n/(k+1)$, but we have not investigated this. Reiner has conjectured that such a result might follow if the Fomin-Reading complex $\Delta^{(k)}(W)$ is shown to be $(k+1)$-Cohen Macaulay in the sense of Baclawski~\cite{baclawski} (personal communication)\footnotemark[\value{footnote}]. In general, the study of inequalities in $h$-vectors is an important area of research.

An {\sf antichain} in a poset $P$ is a set of elements $\{x_1,x_2,\ldots,x_k\}\in P$ that are pairwise incomparable. If the greatest size of an antichain in $P$ is equal to the largest rank number, we say that $P$ has the {\sf Sperner property}, or $P$ is a {\sf Sperner poset}. If, in addition, the cardinality of the disjoint union of any $r$ antichains is less than or equal to the sum of the $r$ largest rank numbers, we say that $P$ is {\sf strongly Sperner}. Either property implies that $P$ is rank-unimodal.

It is known that $NC(A_{n-1})$ and $NC(B_n)$ are both strongly Sperner, as proved by Simion and Ullman \cite[Theorem 2]{simion-ullman} and Reiner \cite[Theorem 13]{reiner}, respectively; both proofs depended on a standard construction called a {\sf symmetric chain decomposition}. This suggests the following problem.

\begin{problem}
Is the poset $NC^{(k)}(W)$ strongly Sperner?
\end{problem}

Finally, we present a curious equidistribution property. When $k=1$ and $W$ is a Weyl group, Stanley's $g$-theorem implies that $NC^{(k)}(W)$ is rank-unimodal because its rank numbers are the $h$-vector of a convex polytope $\Delta^{(1)}(W)$ (see \cite{chapoton-fomin-zelevinsky}). For general $k$, we might hope to deduce unimodality from the $k=1$ case. Recall that the poset $NC^{(k)}(W)$ contains $\Cat^{(k-1)}(W)$ many maximal intervals, and by the above remark each of these is unimodal. Then one could possibly prove unimodality by understanding how these maximal intervals ``zip'' together. Consideration of this ``zipping'' leads to the following conjecture.

\begin{conjecture}
\label{conj:hofh}
If we choose an $\ell$-multichain uniformly at random from the set
\begin{equation*}
\left\{ ((\pi^1)_k,(\pi^2)_k,\ldots,(\pi^\ell)_k)\in (NC^{(k)}(W))^\ell: \ell_T(\pi_1^1)=i\right\},
\end{equation*}
then the expected number of maximal intervals in $NC^{(k)}(W)$ containing this multichain is
\begin{equation}
\label{eq:mysterynumbers}
\frac{\Nar^{(k)}(W,n-i)}{\Nar^{(1)}(W,n-i)},
\end{equation}
and this is independent of the integer $\ell$.
\end{conjecture}

The mysterious numbers \eqref{eq:mysterynumbers} describe the amount of ``overlap'' among maximal intervals at each rank. When $i=n$ the intervals overlap completely (since every maximal interval contains the Coxeter element $(c,c,\ldots,c)\in W^k$), and when $n=0$ the intervals overlap not at all (since each minimal element of $NC^{(k)}(W)$ is contained in precisely one maximal interval). It is worth mentioning that the numbers \eqref{eq:mysterynumbers} are not, in general, integers. These numbers will reappear later in Conjecture \ref{conj:dualF}. Finally, Conjecture \ref{conj:hofh} suggests that the covering of $NC^{(k)}(W)$ by its maximal intervals is a structure worthy of more study.

One may take from this section the following question.

\begin{problem}
What is the true nature of the Fuss-Narayana polynomials related to a finite Coxeter group $W$? One way to approach these polynomials is via the poset $NC^{(k)}(W)$, but there are other definitions (see Chapter \ref{sec:fusscatcomb}). Give a uniform explanation for the case-by-case observations of this section.
\end{problem}

In the next section we will say much more about $\ell$-multichains in $NC^{(k)}(W)$.

\section{The Iterated Construction and Chain Enumeration}
\label{sec:iterated}
In this section, we will study further ``homological'' properties of delta sequences. We will prove a general meta-structural property of $k$-divisible noncrossing partitions (Theorem \ref{th:nckl=nclk}), and as a by-product, we will obtain the zeta polynomimal of $NC^{(k)}(W)$. Afterwards, in Section \ref{sec:shelling}, we use the zeta polynomial to obtain topological information about $NC^{(k)}(W)$.

When defining the poset of $k$-divisible noncrossing partitions, we were faced with a choice: to use the language of multichains or that of delta sequences. Since the clearest definition of $NC^{(k)}(W)$ is by the componentwise partial order on delta sequences, why discuss multichains at all? It turns out that the full richness of the subject requires both perspectives. In particular, the problem of multichain enumeration is well-understood. Also, using multichains allows us to generalize the definition of $NC^{(k)}(W)$.

\begin{definition}
\label{def:P^k}
Given an induced subposet $P$ of $NC(W)$, define $P^{(k)}$ as the subposet of $NC^{(k)}(W)$ consisting of $k$-multichains in $P$,
\begin{equation*}
P^{(k)}:=\left\{ (p_1,p_2,\ldots,p_k)\in NC^{(k)}(W): p_i\in P\text{ for all } 1\leq i\leq k \right\}.
\end{equation*}
\end{definition}

Our goal is to iterate the definition of $k$-divisible noncrossing partitions, to define a poset $(NC^{(k)})^{(\ell)}$ of ``$\ell$-divisible, $k$-divisible noncrossing partitions'' for all pairs $k,\ell$ of positive integers. (Here we will assume $k$ and $\ell$ are finite.) In view of Definition \ref{def:P^k}, we need only express $NC^{(k)}(W)$ as an induced subposet of $NC(W^k)$ in some ``natural'' way.

Recall that $NC^{(k)}(W)$ is {\em isomorphic} to a subposet of $NC(W^k)$: the map $\partial$ is an anti-isomorphism from $NC^{(k)}(W)$ to $NC_{(k)}(W)$, and $NC_{(k)}(W)$ is identified with an order ideal in $NC(W^k)$ (Lemma \ref{lemma:orderideal}). This gives an embedding
\begin{equation*}
\partial: NC^{(k)}(W)\hookrightarrow (NC(W^k))^*
\end{equation*}
of $NC^{(k)}(W)$ into the {\em dual poset} of $NC(W^k)$. Since $NC(W^k)$ is self-dual, $NC^{(k)}(W)$ is certainly isomorphic to some order filter in $NC(W^k)$, but this is not canonical since it depends on a choice of anti-automorphism. If we fix some anti-automorphism $\Psi$ of $NC(W^k)$, then $\Psi\circ\partial$ is an embedding of $NC^{(k)}(W)$ as an order filter in $NC(W^k)$,
\begin{equation*}
\Psi\circ\partial:NC^{(k)}(W)\hookrightarrow NC(W^k).
\end{equation*}
We fix a notation for this embedding.

\begin{definition}
Given an anti-automorphism $\Psi$ of $NC(W^k)$, let $NC_\Psi^{(k)}(W)$ denote the image of the embedding $\Psi\circ\partial:NC^{(k)}(W)\hookrightarrow NC(W^k)$.
\end{definition}

\begin{figure}
\vspace{.1in}
\begin{center}
\input{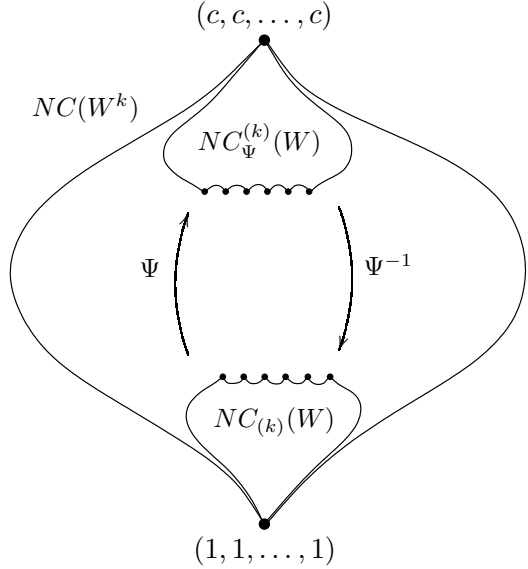}
\end{center}
\caption{$NC^{(k)}(W)$ as an order filter in $NC(W^k)$}
\label{fig:embed}
\end{figure}

Figure \ref{fig:embed} shows $NC^{(k)}(W)$ embedded as an order filter in $NC^{(k)}(W)$. Now it seems reasonable to consider the poset
\begin{equation*}
(NC_\Psi^{(k)}(W))^{(\ell)}
\end{equation*}
as some sort of ``$\ell$-divisible, $k$-divisible noncrossing partitions''. Our only concern is that this definition depends on the choice of $\Psi$. Luckily, it turns out that any ``reasonable'' choice for $\Psi$ yields an isomorphic poset.

If $K$ is the Kreweras complement on $NC(W)$, note that every even power of $K$ is an automorphism of $NC(W)$,
\begin{equation}
\label{eq:krew-type1}
K^{2i}(\pi)=c^{-i}\pi c^i,
\end{equation}
corresopnding to conjugation by $c^i$, and every odd power of $K$,
\begin{equation}
\label{eq:krew-type2}
K^{2i+1}(\pi)=c^{-i}\pi^{-1} c^{i+1},
\end{equation}
is an anti-automorphism of $NC(W)$. Since the Coxeter element $c$ has order $h$, there are $h$ distinct anti-automorphisms of this type, and we give them a special name.

\begin{definition}
If $K$ is the Kreweras complement on $NC(W)$, the odd powers of $K$ are called {\sf Kreweras-type} anti-automorphisms of $NC(W)$.
\end{definition}

At last, we are ready to define the poset $(NC^{(k)}(W))^{(\ell)}$.

\begin{definition}
\label{def:NCkl}
Fix a Kreweras-type anti-automorphism $\Psi$ of $NC(W^k)$. Then, for finite integers $k$ and $\ell$, we define
\begin{equation*}
(NC^{(k)}(W))^{(\ell)} := (NC_\Psi^{(k)}(W))^{(\ell)}
\end{equation*}
following Definition \ref{def:P^k}.
\end{definition}

By construction, this is well-defined.

\begin{lemma}
$(NC^{(k)})^{(\ell)}$ is well-defined up to isomorphism.
\end{lemma}

\begin{proof}
Let $\Psi$ be a Kreweras-type anti-automorphism of $NC(W^k)$ and let $(c)_k$ denote the Coxeter element $(c,c,\ldots,c)\in W^k$. If $K$ is the Kreweras complement on $NC(W^k)$, then notice from equations \eqref{eq:krew-type1} and \eqref{eq:krew-type2} that $\Psi((\pi)_k)$ is equal to $(c)_k^{-i} K((\pi)_k) (c)_k^i$ for some $i$. Since conjugation by $(c)_k^i$ is an automorphism of $NC(W^k)$ it follows that
\begin{equation*}
(NC_\Psi^{(k)})^{(\ell)}\cong (NC_K^{(k)})^{(\ell)},
\end{equation*}
and since $\Psi$ was arbitrary this proves the result.
\end{proof}

Now we examine the structure of $(NC^{(k)}(W))^{(\ell)}$ and give some evidence that this is the ``correct'' generalization of $NC^{(k)}(W)$. One immediate problem is the accumulating mountain of notation. Because the Kreweras-type complement $\Psi$ in Definition \ref{def:NCkl} is arbitrary, we will fix 
\begin{equation*}
\Psi((\pi)_k):=K^{-1}((\pi)_k)=(c)_k(\pi)_k^{-1}
\end{equation*}
from now on, since this yields the cleanest notation. An element of $(NC^{(k)}(W))^{(\ell)}$ is an $\ell$-multichain of $k$-multichains in $NC(W)$. To extend the vector notation $(\pi)_k=(\pi_1,\pi_2,\ldots,\pi_k)$, we will abbreviate an element of $(NC^{(k)}(W))^{(\ell)}$ as a matrix
\begin{equation}
\label{eq:matrix}
((\pi)_k)^l:=\begin{pmatrix}
\pi_1^1 & \pi_1^2 & \cdots & \pi_1^l \\
\pi_2^1 & \pi_2^2 & \cdots & \pi_2^l \\
\vdots  & \vdots  & \ddots & \vdots  \\
\pi_k^1 & \pi_k^2 & \cdots & \pi_k^l
\end{pmatrix},
\end{equation}
where each column $(\pi^j)_k := (\pi_1^j,\pi_2^j,\ldots,\pi_k^j)$ is a $k$-multichain in $NC(W)$, and the sequence of columns $((\pi^1)_k,(\pi^2)_k,\ldots,(\pi^\ell)_k)$ is an $\ell$-multichain in $NC^{(k)}(W)$. Notice that each row of the matrix \eqref{eq:matrix} is also an $\ell$-multichain in $NC(W)$ (Lemma \ref{lemma:horizontal}), but that the sequence of rows is {\em not} a $k$-multichain in $NC^{(\ell)}(W)$. It is worth mentioning that, if we swap $K$ for $K^{-1}$, then the rows will form a multichain in $NC^{(\ell)}(W)$, whereas the columns will no longer be a multichain in $NC^{(k)}(W)$.

Independent of these considerations, we might simply define $(NC^{(k)}(W))^{(\ell)}$ as the set of $k\times\ell$ matrices $((\pi)_k)^\ell$ of elements from $NC(W)$ with weakly increasing rows and columns, and the additional property that
\begin{equation}
\label{eq:nc^k^l}
(\pi_i^{j+1})^{-1}\pi_{i+1}^{j+1}\leq_T (\pi_i^j)^{-1}\pi_{i+1}^j
\end{equation}
for all $1\leq i\leq k$ and $1\leq j\leq \ell$. It is helpful to extend Notation \ref{notation:extend}, setting
\begin{eqnarray*}
\pi_i^j&:=&1\quad\text{ whenever } i<1 \text{ or } j<1 \quad\text{and}\\
\pi_i^j&:=&c \quad\text{ in all other undefined cases, }
\end{eqnarray*}
so that formulas such as \eqref{eq:nc^k^l} make sense for arbitrary integers $i$ and $j$. By analogy with Section \ref{sec:meta}, we might think of $((\pi)_k)^\ell$ as an infinite matrix in which only the entries in the upper-left $k\times\ell$ submatrix are possibly not equal to $c$. This perspective allows direct comparison of matrices with different $k$ and $\ell$ values.

In terms of the matrix notation, the partial order $(NC^{(k)}(W))^{(\ell)}$ (with respect to $\Psi=K^{-1}$) has the following characterization.

\begin{lemma}
\label{lemma:ncklchar}
Given $((\pi)_k)^\ell$ and $((\mu)_k)^\ell$ in $(NC^{(k)}(W))^{(\ell)}$, we have $((\pi)_k)^\ell\leq ((\mu)_k)^\ell$ if and only if
\begin{equation*}
(\mu_i^j)^{-1}\mu_{i+1}^j(\mu_{i+1}^{j+1})^{-1}\mu_i^{j+1}\leq_T (\pi_i^j)^{-1}\pi_{i+1}^j(\pi_{i+1}^{j+1})^{-1}\pi_i^{j+1}
\end{equation*}
for all $i,j\in\integers$.
\end{lemma}

\begin{proof}
This follows directly from Definition \ref{def:NCkl}, taking $\Psi=K^{-1}$.
\end{proof}

Now we present the main result of this section. This theorem simultaneously justifies the definition of $(NC^{(k)}(W))^{(\ell)}$ and the use of the superscript notation. Notice that the proof is very similar to that of Theorem \ref{th:orderideal}, and it uses the same sort of ``diagram chasing''.

\begin{theorem}
\label{th:nckl=nclk}
For all finite positive integers $k$ and $\ell$, we have
\begin{equation*}
(NC^{(k)}(W))^{(\ell)}\cong NC^{(k\ell)}(W) \cong(NC^{(\ell)}(W))^{(k)}.
\end{equation*}
\end{theorem}

\begin{proof}
It will suffice to show the first isomorphism. The second isomorphism follows by reversing the roles of $k$ and $\ell$.

So consider $((\pi)_k)^\ell$ in $(NC^{(k)}(W))^{(\ell)}$ and define a map $\Delta$ taking $((\pi)_k)^\ell$ to the $k\ell$-tuple $\Delta\left(((\pi)_k)^\ell\right):=$
\begin{equation}
\label{eq:klmchain}
\begin{array}{cccccccccccc}
( & \!\!\!\!\!\!\pi_1^1 &\!\!\!\!\!,& \!\!\!\!\!\pi_2^1(\pi_2^2)^{-1}\pi_1^2 &\!\!\!\!\!,& \!\!\!\!\!\pi_2^1(\pi_2^3)^{-1}\pi_1^3 &\!\!\!\!\!,& \!\!\!\!\!\ldots &\!\!\!\!\!,& \!\!\!\!\!\pi_2^1(\pi_2^\ell)^{-1}\pi_1^\ell\,, & \\
& \!\!\!\!\!\!\pi_2^1 &\!\!\!\!\!,& \!\!\!\!\!\pi_3^1(\pi_3^2)^{-1}\pi_2^2 &\!\!\!\!\!,& \!\!\!\!\!\pi_3^1(\pi_3^3)^{-1}\pi_2^3 &\!\!\!\!\!,& \!\!\!\!\!\ldots &\!\!\!\!\!,& \!\!\!\!\!\pi_3^1(\pi_3^\ell)^{-1}\pi_2^\ell\,, & \\
& \!\!\!\!\!\vdots & & & & & & & & & \\
& \!\!\!\!\!\!\pi_{k-1}^1 &\!\!\!\!\!,& \!\!\!\!\pi_k^1(\pi_k^2)^{-1}\pi_{k-1}^2 &\!\!\!\!\!,& \!\!\!\!\!\pi_k^1(\pi_k^3)^{-1}\pi_{k-1}^3 &\!\!\!\!\!,& \!\!\!\!\!\ldots &\!\!\!\!\!,& \!\!\!\!\!\pi_k^1(\pi_k^\ell)^{-1}\pi_{k-1}^\ell\,, & \\
& \!\!\!\!\!\!\pi_k^1 &\!\!\!\!\!,& \!\!\!\!\!\pi_k^2 &\!\!\!\!\!,& \!\!\!\!\!\pi_k^3 &\!\!\!\!\!,& \!\!\!\!\!\ldots &\!\!\!\!\!,& \!\!\!\!\!\pi_k^\ell & \!\!\!\!\!\!).
\end{array}
\end{equation}
We claim that $\Delta$ is the desired isomorphism $(NC^{(k)}(W))^{(\ell)}\cong NC^{(k\ell)}(W)$. To prove this, we must show (1) that $\Delta\left( ((\pi)_k)^\ell \right)$ is in $NC^{(k\ell)}(W)$; (2) that $\Delta$ is a bijection between $(NC^{(k)}(W))^{(\ell)}$ and $NC^{(k\ell)}(W)$; and (3) that $\Delta$ and $\Delta^{-1}$ preserve order.

To see (1), first note that the fact that $((\pi^1)_k,(\pi^2)_k,\ldots,(\pi^\ell)_k)$ is an $\ell$-multichain in $NC^{(k)}(W)$ is equivalent to the conditions
\begin{equation}
\label{eq:pimchain}
(\pi_i^\ell)^{-1}\pi_{i+1}^\ell\leq_T(\pi_i^{\ell-1})^{-1}\pi_{i+1}^{\ell-1}\leq_T\cdots\leq_T (\pi_i^2)^{-1}\pi_{i+1}^2\leq_T(\pi_i^1)^{-1}\pi_{i+1}^1,
\end{equation}
for all $1\leq i\leq k$. But since  $\pi^1_{i+1}\geq_T K_1^{\pi^1_{i+1}}(\pi^1_i)=(\pi_i^1)^{-1}\pi_{i+1}^1$, we see that $\pi_{i+1}^1$ is above every element in the multichain \eqref{eq:pimchain}. Applying the order-reversing map $(K^{\pi_{i+1}^1})^{-1}$ to equation \eqref{eq:pimchain}, we get
\begin{equation}
\label{eq:pimchain2}
\pi_i^1\leq_T\pi_{i+1}^1(\pi_{i+1}^2)^{-1}\pi_i^2\leq_T\cdots\leq_T\pi_{i+1}^1(\pi_{i+1}^\ell)^{-1}\pi_i^\ell\leq_T\pi_{i+1}^1.
\end{equation}
Concatenating the multichains \eqref{eq:pimchain2} for $1\leq i\leq k$, we conclude that $\Delta\left( ((\pi)_k)^\ell\right)$ is a $k\ell$-multichain, proving (1).

To show (2), we give an algorithm to compute $\Delta^{-1}\left( (\mu)_{kl}\right)$ for an arbitrary multichain $(\mu)_{k\ell}\in NC^{(kl)}(W)$. Then we will show that $\Delta^{-1}$ maps $NC^{(k\ell)}(W)$ into $(NC^{(k)}(W))^{(\ell)}$. First, consider the case when $(\mu)_{k\ell}=\Delta\left( ((\pi)_k)^\ell\right)$ for some $((\pi)_k)^\ell\in (NC^{(k)}(W))^{(\ell)}$.
Observing equation \eqref{eq:klmchain}, we see that $(\mu)_{k\ell}$ already contains the entries in the left column and bottom row of $((\pi)_k)^\ell$. This allows us to retrieve the rest of the entries of $((\pi)_k)^\ell$ by ``unzipping'' each of the columns in \eqref{eq:klmchain} as follows. To unzip the $j$-th column of $((\pi)_k)^\ell$, for $1\leq j\leq \ell$, we first set $\pi_k^j:=\mu_{(k-1)\ell+j}$ and then recursively define
\begin{equation*}
\pi_{k-i}^j:=\pi_{k-i+1}^j(\pi_{k-i+1}^1)^{-1}\mu_{(k-i)\ell+j},
\end{equation*}
where $i$ runs from $1$ to $k$. Do this for each $j$, beginning with $j=1$. This process inverts the map $\Delta$. Next, note that this algorithm can be applied equally well to an arbitrary multichain $(\mu)_{k\ell}\in NC^{(k\ell)}(W)$, and the resulting matrix, which we denote by $\Delta^{-1}\left( (\mu)_{k\ell}\right)$, is uniquely determined. The fact that $\Delta^{-1}\left( (\mu)_{k\ell}\right)$ is in $(NC^{(k)}(W))^{(\ell)}$ then follows by reversing the steps in the proof of (1).

Finally, it is easy to see that both $\Delta$ and $\Delta^{-1}$ preserve order, since the relations $((\pi)_k)^\ell\leq ((\mu)_k)^\ell$ in $(NC^{(k)}(W))^{(\ell)}$ and $\Delta\left( ((\pi)_k)^\ell\right)\leq\Delta\left( ((\mu)_k)^\ell \right)$ in $NC^{(k\ell)}(W)$ are both equivalent to the set of conditions
\begin{equation*}
(\mu_i^j)^{-1}\mu_{i+1}^j(\mu_{i+1}^{j+1})^{-1}\mu_i^{j+1}\leq_T (\pi_i^j)^{-1}\pi_{i+1}^j(\pi_{i+1}^{j+1})^{-1}\pi_i^{j+1}
\end{equation*}
for all integers $i$ and $j$, by Lemma \ref{lemma:ncklchar}. Hence they are equivalent to each other, proving (3).
\end{proof}

Like Theorem \ref{th:orderideal}, this result has a homological feel. For instance, the map $\Delta$ could be seen as forming a ``long exact sequence'' from an exact sequence of ``chain complexes'' (multichains). Or, $\Delta$ might be thought of as computing the ``total homology'' of the ``double chain complex'' $((\pi)_k)^\ell$, where the horizontal and vertical homology are the delta sequences induced by the multichain of rows and multichain of columns, respectively. Then Theorem \ref{th:nckl=nclk} proves the ``commutativity of taking double homology''. We are curious how far this analogy can go.

The key in this proof was to deform the matrix $((\pi)_k)^\ell)$ so that every element of the first row is below $\pi_2^1$, every element of the second row is below $\pi_3^1$, and so on, creating a $k\ell$-multichain. Thinking of $NC(W)$ as a lattice, there is an obvious way to accomplish this using meets, and it turns out to be equivalent to the map $\Delta$ \eqref{eq:klmchain}.

\begin{lemma}
The isomorphism $\Delta:(NC^{(k)}(W))^{(\ell)}\to NC^{(k\ell)}(W)$ sends the matrix
\begin{equation*}
((\pi)_k)^\ell=\begin{pmatrix}
\pi_1^1 & \pi_1^2 & \cdots & \pi_1^l \\
\pi_2^1 & \pi_2^2 & \cdots & \pi_2^l \\
\vdots  & \vdots  & \ddots & \vdots  \\
\pi_k^1 & \pi_k^2 & \cdots & \pi_k^l
\end{pmatrix}
\end{equation*}
to the $kl$-multichain
\begin{equation*}
\begin{array}{ccccccccccc}
( & \!\!\!\!\!\!\pi_1^1\wedge\pi_2^1 &\!\!\!\!\!,& \!\!\!\!\! \pi_1^2\wedge\pi_2^1 &\!\!\!\!\!,& \!\!\!\!\!\pi_1^3\wedge \pi_2^1 &\!\!\!\!\!,& \!\!\!\!\!\ldots &\!\!\!\!\!,& \!\!\!\!\!\pi_1^l\wedge\pi_2^1\,, & \\
  & \!\!\!\!\!\!\pi_2^1\wedge\pi_3^1 &\!\!\!\!\!,& \!\!\!\!\!\pi_2^2\wedge\pi_3^1 &\!\!\!\!\!,& \!\!\!\!\!\pi_2^3\wedge\pi_3^1 &\!\!\!\!\!,& \!\!\!\!\!\ldots &\!\!\!\!\!,& \!\!\!\!\!\pi_2^l\wedge\pi_3^1\,, & \\
  & \!\!\!\!\!\vdots & & & & & & & & & \\
  & \!\!\!\!\!\!\pi_{k-1}^1\wedge\pi_k^1 &\!\!\!\!\!,& \!\!\!\!\pi_{k-1}^2\wedge\pi_k^1 &\!\!\!\!\!,& \!\!\!\!\!\pi_{k-1}^3\wedge\pi_k^1 &\!\!\!\!\!,& \!\!\!\!\!\ldots &\!\!\!\!\!,& \!\!\!\!\!\pi_{k-1}^l\wedge\pi_k^1\,, & \\
  & \!\!\!\!\!\!\pi_k^1\wedge c &\!\!\!\!\!,& \!\!\!\!\!\pi_k^2\wedge c &\!\!\!\!\!,& \!\!\!\!\!\pi_k^3\wedge c &\!\!\!\!\!,& \!\!\!\!\!\ldots &\!\!\!\!\!,& \!\!\!\!\!\pi_k^l\wedge c & \!\!\!\!\!\!).
\end{array}
\end{equation*}
\end{lemma}

\begin{proof}
To prove this, we note that Lemma \ref{lemma:joins} $(2)$ can be restated in terms of meets. 

Given $\mu\leq_T\nu$ in $NC(W)$, let $K$ denote the Kreweras complement $K_\mu^\nu$ and recall that $[\mu,\nu]\cong NC(W')$, where $W'=W_{\mu^{-1}\nu}$ (Proposition \ref{prop:subintervals}). Now suppose that $\pi$, $\sigma$, and $K(K^{-1}(\pi)K^{-1}(\sigma))$ are in $NC(W')$ with $K^{-1}(\pi)\leq_T K^{-1}(\pi) K^{-1}(\mu)$. Since $K^{-1}(\pi)$, $K^{-1}(\sigma)$ and $K^{-1}(\pi)K^{-1}(\mu)$ are all in $NC(W')$, Lemma \ref{lemma:joins} $(2)$ (restricted to the parabolic subgroup $W'$) implies that
$$
K^{-1}(\pi)K^{-1}(\sigma)=K^{-1}(\pi)\vee K^{-1}(\sigma),
$$
or
\begin{eqnarray*}
K(K^{-1}(\pi)K^{-1}(\sigma)) &=& K(K^{-1}(\pi)\vee K^{-1}(\sigma)) \\
&=& K(K^{-1}(\pi))\wedge K(K^{-1}(\sigma))\\
&=& \pi\wedge\sigma.
\end{eqnarray*}

Now, our goal is to show that
\begin{equation}
\label{eq:goal}
\pi_{i+1}^1(\pi_{i+1}^j)^{-1}\pi_i^j =\pi_i^j\wedge \pi_{i+1}^1,
\end{equation}
for all integers $i$ and $j$. Since both the rows and the columns of $((\pi)_k)^\ell$ are weakly increasing, we have a square of the form
\begin{equation*}
\begin{array}{ccc}
\pi_i^1 & \leq & \pi_i^j \\
\turndown{\leq} & & \turndown{\leq}\\
\pi_{i+1}^1 & \leq & \pi_{i+1}^j
\end{array}
\end{equation*}
and we may restrict our attention to the interval $[1,\pi_{i+1}^j]$. Recall from Theorem \ref{th:nckl=nclk} that we have
\begin{equation*}
\pi_{i+1}^1(\pi_{i+1}^j)^{-1}\pi_i^j\leq_T \pi_{i+1}^1\leq_T \pi_{i+1}^j,
\end{equation*}
so that
\begin{equation*}
(K_1^{\pi_{i+1}^j})^{-1}(\pi_{i+1}^1)\leq_T (K_1^{\pi_{i+1}^j})^{-1}(\pi_{i+1}^1(\pi_{i+1}^j)^{-1}\pi_i^j)=(K_1^{\pi_{i+1}^j})^{-1}(\pi_i^j) (K^{\pi_{i+1}^j})^{-1}(\pi_{i+1}^1),
\end{equation*}
and both elements of this inequality are in $[1,\pi_{i+1}^j]$. Finally, notice that
\begin{equation*}
K^{\pi_{i+1}^j}\left( (K_1^{\pi_{i+1}^j})^{-1}(\pi_i^j) (K^{\pi_{i+1}^j})^{-1}(\pi_{i+1}^1)  \right) = \pi_{i+1}^1(\pi_{i+1}^j)^{-1}\pi_i^j,
\end{equation*}
hence equation \eqref{eq:goal} follows from the modified Lemma \ref{lemma:joins} above.
\end{proof}

The most startling thing about this lattice characterization of the map $\Delta$ is the fact that it is invertible; under usual conditions, the map $(\pi,\mu)\mapsto \pi\wedge\mu$ ``forgets'' information about the elements $\pi$ and $\mu$.

Finally, we can count the $\ell$-multichains in $NC^{(k)}(W)$, refined by the rank of the bottom element.

\begin{theorem}
\label{th:zetapoly}
For all positive integers $k$ and $\ell$, we have
\begin{enumerate}
\item The number of $\ell$-multichains of $k$-divisible noncrossing partitions is equal to the Fuss-Catalan numer $\Cat^{(k\ell)}(W)$:
\begin{equation*}
\Z(NC^{(k\ell)}(W),1)=\Z(NC^{(k)}(W),\ell)=\Z(NC(W),k\ell)=\Cat^{(k\ell)}(W).
\end{equation*}
\item The number of $\ell$-multichains $((\pi)_k)^\ell=((\pi^1)_k,\ldots,(\pi^\ell)_k)$ of $k$-divisible noncrossing partitions whose bottom element has rank $i$ is given by the Fuss-Narayana number $\Nar^{(k\ell)}(W,i)$:
\begin{equation*}
\Nar^{(k\ell)}(W,i)=\#\left\{ ((\pi)_k)^\ell\in (NC^{(k)}(W))^{(\ell)}:\rk(\pi^1)_k=\ell_T(\pi_1^1)=i\right\}.
\end{equation*}
\end{enumerate}
\end{theorem}

\begin{proof}
This follows immediately from Theorems \ref{th:NCzeta} and \ref{th:nckl=nclk}, and the definition of the Fuss-Narayana polynomials \ref{def:fussnar}.
\end{proof}

That is, the zeta polynomial of $NC^{(k)}(W)$ is $\Cat^{(k\ell)}(W)$, regarded as a polynomial in $\ell$. This encodes a lot of enumerative information. In particular, the zeta polynomial implicitly counts the number of chains (with no repetition) of all sizes.

If $P$ is a finite poset, let $b_i$ denote the number of {\sf $i$-chains} $x_1< x_2<\cdots <x_i$ in $P$. To describe a $k$-multichain that contains precisely the elements $\{x_1,x_2,\ldots,x_i\}$, we must specify which elements are repeated. Since there will be $k-i$ total repetitions, this amounts to choosing a $(k-i)$-multiset (set with possible repetition) from the $i$-set $\{1,2,\ldots,i\}$, and there are\
\begin{equation*}
\left(\binom{i}{i-k}\right)=\binom{k-1}{i-1}
\end{equation*}
ways to do this (see \cite[Chapter 1.2]{stanley:ec1}). Thus, we have the following expression for the zeta polynomial of $P$:
\begin{equation}
\label{eq:zetachains}
\Z(P,k)=\sum_{i\geq 1} b_i\binom{k-1}{i-1}.
\end{equation}
Since $P$ is finite and chains do not contain repetition, there exists some maximum $d$ such that $b_d\neq 0$, and we see that $\Z(P,k)$ is a polynomial in $k$ of degree $d-1$, with leading coefficient $b_d/(d-1)!$. This leads to the following corollary.

\begin{corollary}
\label{cor:maxchains}
The number of maximal chains in $NC^{(k)}(W)$ is equal to
\begin{equation*}
n!(kh)^n/\abs{W},
\end{equation*}
where $h$ is the Coxeter number and $n$ is the rank of $W$.
\end{corollary}

\begin{proof}
Since every maximal chain in $NC^{(k)}(W)$ has $n+1$ elements, the leading coefficient of the zeta polynomial $\Z(NC^{(k)}(W),\ell)$ is equal to $b_{n+1}/ n!$ and is given by the limit
\begin{equation*}
\lim_{\ell\to\infty} \frac{\Z(NC^{(k)}(W),\ell)}{\ell^n}=\frac{1}{\abs{W}}\lim_{\ell\to\infty} \frac{1}{\ell^n}\prod_{i=1}^n (\ell k h+d_i).
\end{equation*}
\end{proof}

For example, consider the symmetric group $A_{n-1}$ with rank $n-1$ and Coxeter number $h=n$. In this case, Corollary \ref{cor:maxchains} says that the lattice of classical noncrossing partitions $NC(A_{n-1})$ ($k=1$) contains $n^{n-2}$ maximal chains. Algebraically, this means that a fixed Coxeter element has $n^{n-2}$ different reduced $T$-words. This number is familiar in combinatorics, since it is the number of labelled trees on $n$ vertices, and it counts a family of objects called {\sf parking functions}. In~\cite{stanley:parking}, Stanley constructed a bijection between maximal chains in $NC(A_{n-1})$ and parking functions, using certain edge-labellings of the poset $NC(A_{n-1})$. In this sense, we might consider the maximal chains in $NC^{(k)}(W)$ as {\sf generalized parking functions}.

We have seen that the zeta polynomial counts chains and multichains. But the zeta polynomials also encodes topological information, as we describe in the next section.

\section{Shellability and Euler Characteristics}
\label{sec:shelling}
The centerpiece of this section is a joint result with Hugh Thomas (Theorem \ref{th:ELNCk}) in which we prove that the order complex of the poset $NC^{(k)}(W)$ is shellable (and hence Cohen-Macaulay). Combining this with the zeta polynomial, we are able to compute the Euler characteristics and homotopy types of some related complexes. First, we define the relevant terminology.

Recall that an {\sf abstract simplicial complex} on a set $X$ is a collection $\Delta$ of subsets of $X$ such that $\{x\}\in\Delta$ for all $x\in X$, and such that $\Delta$ is closed under taking subsets. The elements of $\Delta$ with cardinality $i$ are called {\sf $i$-faces}, or {\sf $(i-1)$-dimensional faces}; a maximal face of $\Delta$ is called a {\sf facet}; and we say $\Delta$ is {\sf pure of dimension $d$} if all of its facets are $d$-dimensional. The {\sf order complex} $\Delta(P)$ of a poset $P$ is the simplicial complex with a face $\{x_1,x_2,\ldots,x_i\}$ for each chain $x_1<\cdots <x_i$ in $P$. Note that $\Delta(P)$ is pure of dimension $d$ if and only if $P$ is graded of height $d$.

If $\Delta$ is a finite, pure, simplicial complex, a {\sf shelling} of $\Delta$ is a total ordering of the facets $F_1,F_2,\ldots,F_r$ such that $F_j\cap(\cup_{i=1}^{j-1} F_i)$ is a nonempty union of maximal faces of $F_j$ for all $2\leq j\leq r$. The paper \cite{bjorner} by Bj\"orner is a standard reference on the application of shellings to posets, and this theory now forms a major branch in the tree of algebraic combinatorics. We will use just a few facts.

If $P$ is a finite graded poset and the order complex $\Delta(P)$ is shellable, we say that $P$ is a {\sf shellable poset}. The utility of this concept is that it tells us a great deal about the topology of the complex. In particular, shellability implies that the complex $\Delta(P)$ is homotopy equivalent to a wedge of top-dimensional spheres. It also implies that the Stanley-Reisner ring of $\Delta(P)$ is Cohen-Macaulay.

A common way to establish the shellability of a poset is via an edge-labelling of its Hasse diagram. If $P$ is a poset let $C(P)$ denote its set of cover relations,
\begin{equation*}
C(P):=\left\{ (x,y)\in P\times P: x\prec y\right\}.
\end{equation*}
We say that a function $\lambda$ from $C(P)$ to some poset $\Lambda$ is an {\sf edge-labelling} of $P$ by $\Lambda$, and an unrefinable chain $x_1\prec x_2\prec\cdots x_r$ in $P$ is called {\sf rising} if
\begin{equation*}
\lambda(x_1,x_2)\leq\lambda(x_2,x_3)\leq\cdots \leq \lambda(x_{r-1},x_r).
\end{equation*}
We call $\lambda$ an {\sf edge-lexicographic labelling} (or an {\sf EL-labelling}) of $P$ if each interval $[x,y]$ in $P$ contains a {\em unique} rising maximal chain, and the labelling of this chain is {\em lexicographically first} among all maximal chains in $[x,y]$. Finally, note that the facets of $\Delta(P)$ correspond to the maximal chains of $P$. It is proven in \cite{bjorner} that an EL-labelling on a {\em bounded} poset $P$ (with maximum element $\hat{1}$ and minimum element $\hat{0}$) induces a shelling of $\Delta(P)$ by taking the lexicographic order on facets. (The bounded assumption is necessary so that all maximal chains may be compared.) In this case, we say that $P$ is {\sf EL-shellable}. For more details about shellable posets, see Bj\"orner~\cite{bjorner} or the more modern survey~\cite{bjorner:topological}, which gives a general introduction to the use of topological methods in combinatorics.

The following theorem was proved recently by Athanasiadis, Brady and Watt \cite[Theorem 1.1]{athanasiadis-brady-watt}.

\begin{theorem}
\label{th:ELNC}
For each finite Coxeter group $W$, the lattice of noncrossing partitions $NC(W)$ is EL-shellable.
\end{theorem}

Recall that the edges in the Hasse diagram of $NC(W)$ are naturally labelled by reflections $T$: We have $\pi\prec\mu$ in $NC(W)$ if and only if $\mu=\pi t$ for some $t\in T$ with $\ell_T(\mu)=\ell_T(\pi)+1$. In this case, the reflection $t$ is unique. To prove their result, Athanasiadis, Brady and Watt defined a total order on the set $T$ (relative to the fixed Coxeter element $c$) such that the natural edge-labelling by $T$ becomes EL. We will call this the {\sf ABW order} on $T$. Since we need only its existence, we will not describe the details here.

Now we will show that the poset $NC^{(k)}(W)$ is shellable. Note that this is more difficult than it may seem at first glance. Since $NC^{(k)}(W)$ is not bounded (it has no minimum element), an EL-labelling of its Hasse diagram {\em does not imply shellability}.  (It is possible to define an unbounded, EL-labelled poset whose order complex is {\em not} shellable.) Instead, we will construct an EL-labelling of the associated bounded poset $NC^{(k)}(W)\cup\{\hat{0}\}$ with a minimum element $\hat{0}$ formally adjoined. Recall from Section \ref{sec:semilattice} that this poset is a lattice.

\begin{theorem}[joint with Hugh Thomas]
\label{th:ELNCk}
For each finite Coxeter group $W$ and each positive integer $k$, the lattice $NC^{(k)}(W)\cup\{\hat{0}\}$ is EL-shellable.
\end{theorem}

\begin{proof}
Consider the bounded poset $NC_{(k)}(W)\cup\{\hat{1}\}$ of $k$-delta sequences with a maximum element $\hat{1}$ formally adjoined. We will prove that $NC_{(k)}(W)\cup\{\hat{1}\}$ is EL-shellable, and the result follows by duality.

First, denumerate the reflections $T=\{t_1,t_2,\ldots,t_N\}$ by the ABW order. It is well-known that this induces an EL-labelling of the lattice $NC(W^k)$ in the following way. Recall that $NC(W^k)$ is edge-labelled by the set of reflections
\begin{equation*}
T^k=\{ t_{i,j}=(1,1,\ldots,t_j,\ldots,1) : 1\leq i,j\leq N\}
\end{equation*}
(where $t_j$ occurs in the $i$-th entry of $t_{i,j}$) of the Coxeter group $W^k$. In this case, the cover relations of index $i$ are labelled by reflections $t_{i,j}$ for some $1\leq j\leq N$. Then the {\sf lex ABW order}
\begin{equation*}
t_{1,1},t_{1,2},\ldots,t_{1,N},t_{2,1},\ldots,t_{2,N},\ldots,t_{k,1},\ldots,t_{k,N}
\end{equation*}
induces an EL-shelling of $NC(W^k)$. It is a general phenomenon that the direct product of EL-labellings induces an EL-labelling \cite[Theorem 4.3]{bjorner}.

Now recall that $NC_{(k)}(W)$ is an order ideal in $NC(W^k)$ (Lemma \ref{lemma:orderideal}), so clearly the lex ABW order restricts to an EL-labelling of the Hasse diagram of $NC_{(k)}(W)$. The Hasse diagram of $NC_{(k)}(W)\cup\{\hat{1}\}$ is obtained from $NC_{(k)}(W)$ by adding an edge of the form $(\delta)_k\prec\hat{1}$ for each maximal element $(\delta)_k\in NC_{(k)}(W)$. Our goal is to label these new edges in such a way that the EL property extends to the new labelling. We do this by fixing a symbol $\lambda$ and labelling each new edge $(\delta)_k\prec\hat{1}$ by this symbol. Then we order the new label set $T^k\cup\{\lambda\}$ by placing $\lambda$ between $t_{1,N}$ and $t_{2,1}$ in the lex ABW order. That is, $\lambda$ comes {\em after} the index $1$ reflections, and {\em before} all others. We claim that this is an EL-labelling of $NC_{(k)}(W)\cup\{\hat{1}\}$.

The proof uses induction on the rank of $W$. To verify that the edge-labelling of $NC_{(k)}(W)\cup\{\hat{1}\}$ is EL, we must show that every interval in $NC_{(k)}(W)\cup\{\hat{1}\}$ has a unique rising maximal chain that is lexicographically first. Intervals of the form $[(\delta)_k,(\varepsilon)_k]$ with $(\delta)_k$ and $(\varepsilon)_k$ in $NC_{(k)}(W)$ trivially satisfy this property since the lex ABW order is EL on $NC_{(k)}(W)$. Now consider the unique maximal interval $[(1)_k,\hat{1}]$. If a maximal chain in this interval is rising, then its labels other than $\lambda$ must consist entirely of index $1$ reflections, hence this chain must pass through the element $(c,1,\ldots,1)$. Since the interval $[(1)_k,(c,1,\ldots,1)]$ contains a unique rising maximal chain, so does $[(1)_k,\hat{1}]$. Furthermore, this chain is lexicographically first among maximal chains in $[(1)_k,(c,1,\ldots,1)]\cup\{\hat{1}\}$, and it lexicographically precedes all other maximal chains in $[(1)_k,\hat{1}]$ since each of these must contain a label of index $\geq 2$. Finally, consider an interval of the form $[(\delta)_k,\hat{1}]$ where $(\delta)_k$ is not equal to $(1)_k$. By Theorem \ref{th:orderideal}, this interval is isomorphic to $NC_{(k)}(W')\cup\{\hat{1}\}$ for some proper parabolic subgroup $W'$ of $W$. Since the isomorphism in Theorem \ref{th:orderideal} preserves the index of edge labels, we may use exactly the same argument for the interval $[(\delta)_k,\hat{1}]$ as we did for $[(1)_k,\hat{1}]$ above.
\end{proof}

A few remarks: The fact that $NC_{(k)}(W)$ has an EL-labelling follows easily from Athanasiadis, Brady and Watt~\cite{athanasiadis-brady-watt} and the fact that $NC_{(k)}(W)$ is an order ideal in $NC(W^k)$. The real difficulty in Theorem \ref{th:ELNCk} is to extend this labelling to the bounded poset $NC_{(k)}(W)\cup\{\hat{1}\}$. The essential insight here was provided by Hugh Thomas. Once defined, it is relatively straightforward to verify the EL property. We remarked earlier that an EL-labelled, unbounded poset need not have any nice topological properties. The fact that the ABW shelling extends to $NC^{(k)}(W)$ in a straightforward way is one more reason to believe that $NC^{(k)}(W)$ is the ``correct'' generalization of the noncrossing partitions.

Suppose that $P$ is a finite, graded poset with rank function $\rk:P\to\integers$ and rank set $\{0,1,\ldots,n\}$. Then for each rank subset $R\subseteq \{0,1,\ldots,n\}$, we define the {\sf rank-selected subposet}
\begin{equation*}
P_R:=\left\{ x\in P: \rk(x)\in R\right\}
\end{equation*}
as the induced subposet of $P$ consisting of the elements with rank in $R$. It is well-known  \cite[Theorem 4.1]{bjorner} that if $P$ is shellable, then all of its rank-selected subposets are shellable. In particular, we have the following.

\begin{corollary}
The poset $NC^{(k)}(W)$ of $k$-divisible noncrossing partitions is shellable.
\end{corollary}

\begin{proof}
Restrict to the rank set $R=\{1,2,\ldots,n\}$ in $NC^{(k)}(W)\cup\{\hat{0}\}$.
\end{proof}

Now let us see what topological information we can squeeze from the zeta polynomial. The zeta polynomial $\Z(P,k)$ was {\em defined} as the number of $k$-multichains $x_1\leq x_2\leq\cdots\leq x_k$ in $P$; so, in particular, $\Z(P,1)=\abs{P}$. We might wonder if the number $\Z(P,k)$ has any significance when evaluated at non-positive integers $k$. Recall (equation \eqref{eq:zetachains}) that $\Z(P,k)$ is given by
\begin{equation*}
\Z(P,k)=\sum_{i\geq 0} b_i\binom{k-1}{i-1},
\end{equation*}
where $b_i$ is the number of $i$-chains $x_1<x_2<\cdots<x_i$ in $P$. Evaluating at $k=0$, and using the fact that $\binom{-p}{q}=(-1)^q\binom{p+q-1}{q}$ when $q$ is a nonnegative integer, we find that
\begin{equation*}
\Z(P,0)=\sum_{i\geq 1} b_i\binom{-1}{i-1}=\sum_{i\geq 1} (-1)^{i-1} b_i=\chi(\Delta(P)),
\end{equation*}
which is the usual (non-reduced) Euler characteristic of the order complex $\Delta(P)$.  When the poset $P$ possesses a maximum element $\hat{1}$ or a minimum element $\hat{0}$, then the Euler characteristic is $\chi(\Delta(P))=1$ since the complex $\Delta(P)$ is contractible (it has a cone point).\footnote{We might think of this as counting the unique ``$0$-multichain'' in $P$ (which is the unique empty face of $\Delta(P)$).}

Evaluating at $k=-1$, we get
\begin{equation}
\label{eq:k=-1}
\Z(P,-1)=\sum_{i\geq 1} b_i\binom{-2}{i-1}=\sum_{i\geq 1} (-1)^{i-1}\,i\cdot b_i.
\end{equation}
But the most interesting case for us occurs when $k=-2$:

When $P$ is bounded, the order complex $\Delta(P)$ is not very interesting since it is contractible. In this case, it is common to consider instead the order complex with the cone points $\hat{0}$ and $\hat{1}$ deleted. There is a standard result \cite[Propositions 3.8.5 and 3.11.1]{stanley:ec1} that relates the topology of this ``deleted'' order complex to the zeta polynomial and the M\"obius function of the poset $P$. (For information on M\"obius functions and M\"obius inversion, see \cite[Chapter 3]{stanley:ec1}.)

\begin{theorem}
If $P$ is a finite poset with $\hat{0}$ and $\hat{1}$, then
\begin{equation*}
\Z(P,-2)=\mu_P(\hat{0},\hat{1})=\tilde\chi\left( \Delta(P\setminus\{\hat{0},\hat{1}\})\right),
\end{equation*}
where $\mu_P$ is the M\"obius function of $P$ and $\tilde\chi$ is the reduced Euler characteristic.
\end{theorem}

If $\Delta$ is a simplicial complex and $f_i$ counts the $i$-dimensional faces of $\Delta$ for $i\geq 0$, recall that the {\sf Euler characteristic} $\chi(\Delta)$ and {\sf reduced Euler characteristic} $\tilde{\chi}(\Delta)$ of $\Delta$ are defined by
\begin{equation*}
\chi(\Delta):=\sum_{i\geq 0} (-1)^i f_i\qquad\text{ and }\qquad \tilde\chi(\Delta):=\chi(\Delta)-1.
\end{equation*}
The {\sf wedge sum} of simplicial complexes $\Delta'\vee\Delta'':=\Delta'\sqcup\Delta''/\{x_0\sim y_0\}$ is the quotient of the disjoint union by identifying the point $\{x_0\}\in\Delta'$ with the point $\{y_0\}\in\Delta''$. The main advantages of the {\em reduced} Euler characteristic over the usual Euler characteristic are the following:
if $\Delta$ is a $d$-dimensional sphere with $d\geq 1$ then $\tilde\chi(\Delta)=(-1)^d$; the reduced Euler characteristic adds over wedge sums, $\tilde\chi(\Delta'\vee\Delta'')=\tilde\chi(\Delta')+\tilde\chi(\Delta'')$; and the reduced Euler characteristic respects quotients, $\tilde\chi(\Delta'/\Delta'')=\tilde\chi(\Delta')-\tilde\chi(\Delta'')$.

Setting $P=NC(W)$ and combining with Theorem \ref{th:NCzeta} we obtain the following well-known\footnote{See Chapoton \cite{chapoton:one}. The type A version of this formula was known to Kreweras \cite{kreweras}.} formula for the reduced Euler characteristic of the deleted order complex of $NC(W)$:
\begin{equation*}
\tilde\chi\left( \Delta(NC(W)\setminus\{1,c\})\right) =\Z(NC(W),-2)=\Cat^{(-2)}(W).
\end{equation*}
The number $\Cat^{(-2)}(W)$ will recur frequently, so we set down some notation. Since the Fuss-Catalan number $\Cat^{(k)}(W)$ is a polynomial, we can formally evaluate it at $-k-1$. Using again the fact that the numbers $h-m_1,\ldots,h-m_n$ are a permutation of the exponents $m_1,\ldots,m_n$, we get
\begin{eqnarray*}
\Cat^{(-k-1)}(W) &=& \frac{1}{\abs{W}}\prod_{i=1}^n ((-k-1)h+d_i)\\
&=& \frac{1}{\abs{W}}\prod_{i=1}^n (-kh+(-h+m_i)+1)\\
&=& \frac{1}{\abs{W}}\prod_{i=1}^n (-kh-m_i+1)\\
&=& \frac{(-1)^n}{\abs{W}}\prod_{i=1}^n (kh+d_i-2),
\end{eqnarray*}
and we give this formula a special name.
\begin{definition}
The {\sf positive Fuss-Catalan polynomial} associated to $W$ is
\begin{equation}
\label{eq:posfusscat}
\Cat_+^{(k)}(W):=(-1)^n\Cat^{(-k-1)}(W)=\frac{1}{\abs{W}}\prod_{i=1}^n (kh+d_i-2).
\end{equation}
In the case $k=1$ we write $\Cat_+(W)=\Cat_+^{(1)}(W)=(-1)^n\Cat^{(-2)}(W)$.
\end{definition}

Restating the above observations, we have the following result from \cite{athanasiadis-brady-watt}.
\begin{theorem}
\label{th:NChomotopy}
The deleted order complex $\Delta(NC(W)\setminus\{1,c\})$ is homotopic to a wedge of $\Cat_+(W)$ many $(n-2)$-dimensional spheres.
\end{theorem}

\begin{proof}
First, note that $\Delta(NC(W)\setminus\{1,c\})$ is a pure $(n-2)$-dimensional simplicial complex, since each maximal chain in $NC(W)\setminus\{1,c\}$ has $n-1$ elements. Because $NC(W)\setminus\{1,c\}$ is a rank-selected subposet of $NC(W)$, Theorem \ref{th:ELNC} implies that its order complex is shellable, and hence it is homotopic to a wedge of $(n-2)$-dimensional spheres. Finally, recall that
\begin{equation*}
\tilde\chi\left(\Delta(NC(W)\setminus\{1,c\})\right)=\Cat^{(-2)}(W)=(-1)^n\Cat_+(W).
\end{equation*}
Since the reduced Euler characteristic of an $(n-2)$-dimensional sphere is $(-1)^{n-2}=(-1)^n$, and this adds over wedge sums, we get the result.
\end{proof}

Thus the positive Coxeter-Catalan number $\Cat_+(W)$ has a topological interpretation in terms of the noncrossing partitions, and we might well call it the ``{\em topological} Coxeter-Catalan number''. The notation ``positive'' for these numbers is motivated by a connection with cluster theory (see Section \ref{sec:cluster}).

Considering the above result, one might guess that the positive Fuss-Catalan numbers $\Cat_+^{(k)}(W)$ also have a topological interpretation. However, there is an immediate problem in trying to generalize Theorem \ref{th:NChomotopy} to the $k$-divisible noncrossing partitions; that is, the poset $NC^{(k)}(W)$ does not have a minimum element to delete! There are two possible ways to proceed.

First, we consider the order complex of $NC_{(k)}(W)$ with only the top element $(c)_k$ deleted.

\begin{theorem}
\label{th:NCkhomotopy}
For all positive integers $k$, the complex $\Delta(NC^{(k)}(W)\setminus\{(c)_k\})$ has reduced Euler characteristic
\begin{equation}
\label{eq:kec}
\tilde\chi\left( \Delta(NC^{(k)}(W)\setminus\{(c)_k\})\right)= (-1)^{n-1}\Cat_+^{(k-1)}(W),
\end{equation}
and hence it is homotopic to a wedge of $\Cat_+^{(k-1)}(W)$ many $(n-1)$-dimensional spheres.
\end{theorem}

\begin{proof}
Because $NC^{(k)}(W)\setminus\{(c)_k\}$ is a rank-selected subposet of $NC^{(k)}(W)$, Theorem \ref{th:ELNCk} tells us that its order complex is shellable, and hence homotopic to a wedge of top-dimensional spheres (in this case, $(n-1)$-dimensional spheres). We will be done if we can compute the reduced Euler characteristic.

So let $b_i$ denote the number of $i$-chains $(\pi^1)_k< (\pi^2)_k<\cdots <(\pi^i)_k$ in $NC^{(k)}(W)$, let $c_i$ denote the number of $i$-chains in $NC^{(k)}(W)\setminus\{(c)_k\}$, and set $b_0=c_0=1$. In this case it is easy to see that $b_i=c_{i-1}+c_i$ for all $i\geq 1$, since each $(i-1)$-chain in $NC^{(k)}(W)\setminus\{(c)_k\}$ extends to a unique $i$-chain in $NC^{(k)}(W)$. Applying formula \eqref{eq:k=-1}, we have a telescoping sum
\begin{eqnarray*}
\Z(NC^{(k)}(W),-1) &=& \sum_{i\geq 1} (-1)^{i-1}\, i\cdot b_i \\ 
&=& \sum_{i\geq 1} (-1)^{i-1}c_{i-1} \\
&=& -\tilde\chi\left( \Delta( NC^{(k)}(W)\setminus\{(c)_k\} ) \right).
\end{eqnarray*}
On the other hand, Theorem \ref{th:zetapoly} implies that
\begin{equation*}
\Z(NC^{(k)}(W),-1)=\Z(NC(W),-k)=\Cat^{(-k)}(W)=(-1)^n\Cat_+^{(k-1)},
\end{equation*}
which completes the proof.
\end{proof}

Following this, the material from Section \ref{sec:meta} suggests that we should also look at the quotient complexes $\Delta(NC^{(\ell)}(W))/\Delta(NC^{(k)}(W))$ for $1\leq k\leq\ell$, where $NC^{(k)}(W)$ is embedded as an order filter in $NC^{(\ell)}(W)$. We can immediately compute the reduced Euler characteristic of this complex.

\begin{corollary}
Let $NC^{(k)}(W)$ be any isomorphic copy of the $k$-divisible noncrossing partitions embedded within $NC^{(\ell)}(W)$, for $1\leq k\leq \ell$. Then we have
\begin{equation}
\label{eq:formula1}
\tilde\chi\left( \frac{\Delta(NC^{(\ell)}(W))}{\Delta(NC^{(k)}(W))}\right)=(-1)^{n-1}\left( \Cat_+^{(\ell-1)}(W)-\Cat_+^{(k-1)}(W)\right).
\end{equation}
\end{corollary}

\begin{proof}
This follows directly from Theorem \ref{th:NCkhomotopy} and properties of the reduced Euler characteristic.
\end{proof}

So the positive Fuss-Catalan numbers definitely play a role in the topology of the $k$-divisible noncrossing partitions. However, these results describe a new type of behavior that was not observed in the $k=1$ case. Indeed, if we set $k=1$ in equation \eqref{eq:kec}, we merely recover the fact that
\begin{equation*}
\tilde\chi\left(\Delta(NC(W)\setminus\{c\})\right)=(-1)^{n-1}\Cat_+^{(0)}=0,
\end{equation*}
which is obvious because the complex $\Delta(NC(W)\setminus\{c\})$ is contractible (the identity $1\in W$ is a cone point). That is, Theorem \ref{th:NCkhomotopy} is {\em not} a generalization of Theorem \ref{th:NChomotopy}.

To truly generalize Theorem \ref{th:NChomotopy}, we should instead consider the order complex of $NC^{(k)}(W)$ with its maximum element $(c)_k$ and {\em all of its minimal elements} deleted. In this case, we have a conjecture and not a theorem. Let $\mins$ denote the set of minimal elements in $NC^{(k)}(W)$.

\begin{conjecture}
\label{conj:homotopy}
For all positive integers $k$, the order complex of $NC^{(k)}(W)\setminus(\{(c)_k\}\cup\mins)$ has reduced Euler characteristic
\begin{equation}
\label{eq:formula2}
(-1)^n\left( \Cat_+^{(k)}(W)-\Cat_+^{(k-1)}(W)\right),
\end{equation}
and it is homotopic to a wedge of $(n-2)$-dimensional spheres.
\end{conjecture}

This statement {\em is} a generalization of Theorem \ref{th:NChomotopy}, because in the case $k=1$, this restricts to
\begin{equation*}
\tilde\chi(NC(W)\setminus\{1,c\})=(-1)^n\left( \Cat_+^{(1)}(W)-\Cat_+^{(0)}(W)\right)=(-1)^n\Cat_+(W),
\end{equation*}
as desired. Notice, also, the coincidence between formulas \eqref{eq:formula1} and \eqref{eq:formula2}. This suggests a possible method of proof for the above conjecture.

\begin{problem}
Prove Conjecture \ref{conj:homotopy}, perhaps by exhibiting a homotopy equivalence between the complexes
\begin{equation*}
\frac{\Delta(NC^{(k+1)}(W))}{\Delta(NC^{(k)}(W))}\qquad \text{ and }\qquad \Delta(NC^{(k)}(W)\setminus(\{(c)_k\}\cup\mins)).
\end{equation*}
\end{problem}

Finally, we suggest a problem for further study. In the paper \cite{stanley:posets} Stanley initiated the representation theory of finite posets. If $P$ is a finite, graded poset with $\hat{0}$ and $\hat{1}$, then any group $G$ of automorphisms of $P$ also acts on the order complex $\Delta(P)$. If the $d$-dimensional complex $\Delta(P)$ is shellable, and hence Cohen-Macaulay, then only its top reduced homology group $\tilde{H}_d(\Delta(P),\complex)$ is nonzero, and the representation of $G$ acting on this group may be interesting. In particular, the character of $g\in G$ acting on $\tilde{H}_d(\Delta(P),\complex)$ is given by the M\"obius invariant $\mu_{P^g}(\hat{0},\hat{1})$, computed in the sublattice $P^g\subseteq P$ fixed pointwise by $g$.

Recall from Section \ref{sec:automorphisms} that there is a natural dihedral group of poset automorphisms of $NC(W)$, generated by the maps $L$ and $R$. In an unpublished manuscript \cite{montenegro}, Montenegro investigated the action of $\langle L,R\rangle$ on the lattice $NC(A_{n-1})$ and computed the character of its top homology representation (see also \cite{reiner}). More recently, Callan and Smiley have considered enumerative questions related to the type $A$ case \cite{callan-smiley}.

For general $k$, we have defined a dihedral group of automorphisms of $NC_{(k)}(W)$ generated by maps $L^*$ and $R^*$ in  Section \ref{sec:automorphisms}. Since we now know  that the lattice $NC_{(k)}(W)\cup\{\hat{1}\}$ is shellable, it may be interesting to study the action of $\langle L^*,R^*\rangle$ on its order complex.

\begin{problem}
\label{prob:homologyrep}
Investigate the action of the dihedral group $\langle L^*,R^*\rangle$ on the lattice $NC_{(k)}(W)\cup\{\hat{1}\}$. What are the properties of the fixed point lattices? What is the character of this action on the top homology of the order complex? Does this representation have some significance?
\end{problem}

This completes our uniform study of the ``type $W$'' $k$-divisible noncrossing partitions $NC^{(k)}(W)$. We hope that the contents of this chapter will inspire further investigations.

In the next chapter, we turn to a concrete study of the poset $NC^{(k)}(W)$ for the classical finite Coxeter groups.


\chapter{The Classical Types}
\label{sec:classical}

In the classification of finite irreducible Coxeter groups, the three infinite families $A_{n-1}$ (the {\sf symmetric group}: symmetries of the regular simplex with $n$ vertices), $B_n=C_n$ (symmetries of the hypercube/hyperoctahedron) and $D_n$ (an index $2$ subgroup of $B_n$) are known as the {\sf classical groups}. In addition to this, there is the family of {\sf dihedral groups} $I_2(m)$ (symmetries of a regular $m$-gon) as well as six {\sf exceptional groups}: $H_3$, $H_4$, $F_4$, $E_6$, $E_7$ and $E_8$. To say that one has proved a result about all finite Coxeter groups in a {\sf case-by-case} manner means that one has proved a separate theorem for each of the three classical families (the dihedral groups are usually very easy to deal with since they all have rank $2$), and verified the result for the remaining exceptional groups by hand, or using a computer. Typically, a result about all finite Coxeter groups can be reduced to the irreducible case.

We are fortunate to have this complete classification, since it allows us to make a lot of conjectures, and to ``prove'' these conjectures quickly, usually before they are really understood. In this sense, the subject of finite Coxeter groups has a very experimental flavor.

This chapter is very concrete. In the first two sections, we will recall the theory of type $A$ (``classical'') noncrossing partitions. After this, we will explore our results from Chapter \ref{sec:kdiv} in the context of the classical groups, and prove some case-by-case results.
Traditionally, much of the combinatorics of Coxeter groups was understood first in the type $A$ case, before it was satisfactorily generalized to other types. This was certainly the case with the type $A$ noncrossing partitions, which were studied as early as 1972~\cite{kreweras}, and much of our type $W$ terminology is inspired by this context. For instance, the term ``noncrossing partition'' comes from the fact that $NC(W)$ can be realized as a partial order on certain ``noncrossing'' set partitions when $W$ is of classical type. We will see that this context also inspires the term ``$k$-divisible''.

\section{Classical Noncrossing Partitions}
\label{sec:classicalNC}
First, we introduce the idea of a classical noncrossing partition. The term ``noncrossing'' originated in the paper {\em Sur les partitions non crois\'ees d'un cycle}~\cite{kreweras}, published by Kreweras in 1972, in the inaugural volume of Discrete Mathematics. This was the first study of these objects using modern algebraic-combinatorial methods such as the M\"obius function and poset theory. However, the notion of a noncrossing partition is very elementary, and had probably appeared many times before. One of the earliest\footnote{See the Introduction (Stream \ref{sec:stream3}) for more details.} appearances is in Becker \cite{becker3}, where they were called ``planar rhyme schemes''. The survey~\cite{simion} by Simion gives an excellent account of the history of classical noncrossing partitions. The type $B$ classical noncrossing partitions are due to Reiner~\cite{reiner} and the type $D$ noncrossing partitions are due to Athanasiadis and Reiner~\cite{athanasiadis-reiner}. We will discuss these in sections \ref{sec:typeB} and \ref{sec:typeD}, respectively.

To begin, fix $[n]:=\{1,2,\ldots,n\}$ and let $\P=\{P_1,P_2,\ldots,P_m\}$ denote a partition of the set $[n]$, where we call $P_i$ a {\sf block} of $\P$ for $1\leq i\leq m$. We say that two blocks $P_i\neq P_j$ {\sf cross} if there exist $1\leq a<b<c<d\leq n$ with $\{a,c\}\subseteq P_i$ and $\{b,d\}\subseteq P_j$. If $P_i$ and $P_j$ do not cross for all $1\leq i<j\leq m$, we say that $\P$ is a {\sf noncrossing partition} of $[n]$.

\begin{figure}
\vspace{.1in}
\begin{center}
\input{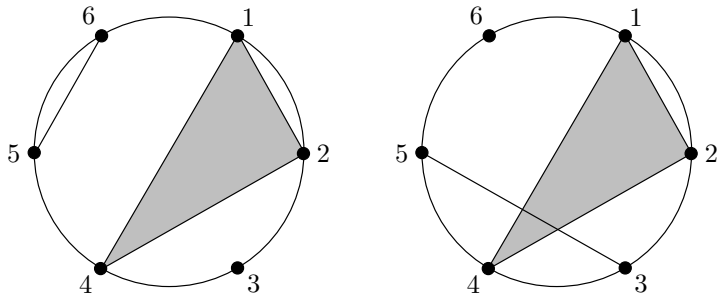}
\end{center}
\caption{A noncrossing and a crossing partition of the set $[6]$}
\label{fig:ncexample}
\end{figure}

This property is made clear if we think of $[n]$ as labelling the vertices of a regular $n$-gon, clockwise. If we identify each block of $\P$ with the convex hull of its corresponding vertices, then we see that $\P$ is noncrossing precisely when its blocks are pairwise disjoint (that is, they don't ``cross''). In this sense, we see that the property of ``noncrossing'' really depends on the cyclic order and not the linear order on $[n]$. Figure \ref{fig:ncexample} shows that $\{\{1,2,4\},\{3\},\{5,6\}\}$ is a noncrossing partition of the set $[6]$, whereas $\{\{1,2,4\},\{3,5\},\{6\}\}$ is {\sf crossing}. We will call this the {\sf circular representation} of the partition. When describing a specific partition of $[n]$, we will usually list the blocks in increasing order of their minimum elements, although this will not be important.

The set of noncrossing partitions of $[n]$ forms a poset under refinement of partitions, with maximum element $\hat{1}_n=\{\{1,2,\ldots,n\}\}$ and minimum element $\hat{0}_n=\{\{1\},\{2\},\ldots,\{n\}\}$. Unlike the general algebraic case, it is easy to verify that this poset is a lattice, since the meet of two partitions $\P$ and $\Q$ is just their coarsest common refinement (whose blocks are obtained by intersecting the blocks of $\P$ with the blocks of $\Q$). The lattice is also graded, with rank function given by $n$ minus the number of blocks,
\begin{equation}
\label{eq:rankblocks}
\rk(\P):= n-\abs{\P}.
\end{equation}

\begin{definition}
Let $NC(n)$ denote the lattice of noncrossing partitions of the set $[n]=\{1,2,\ldots,n\}$, under the refinement partial order.
\end{definition}

The following fundamental results were proved by Kreweras.

\begin{theorem}[\cite{kreweras}]\hspace{.1in}
\label{th:kreweras}
\begin{enumerate}
\item $NC(n)$ is counted by the classical Catalan number,
\begin{equation*}
\abs{NC(n)}=\Cat(n):=\frac{1}{n}\binom{2n}{n-1}.
\end{equation*}
\item $NC(n)$ is ranked by the classical Narayana numbers,
\begin{equation*}
\#\big\{\P\in NC(n):\abs{\P}=i\big\}=\Nar(n,i):=\frac{1}{n}\binom{n}{i}\binom{n}{i-1}.
\end{equation*}
\item The zeta polynomial of $NC(n)$ (which counts the number of $k$-multichains $\P_1\leq\cdots\leq\P_k$ in $NC(n)$) is given by the classical Fuss-Catalan number,
\begin{equation*}
\Z(NC(n),k)=\Cat^{(k)}(n):=\frac{1}{n}\binom{(k+1)n}{n-1}.
\end{equation*}
\end{enumerate}
\end{theorem}

\begin{proof}
All of these follow from the stronger result \cite[Theorem 4]{kreweras} which we will present later as Theorem \ref{th:krewerastype}.
\end{proof}

The relationship between the classical noncrossing partitions $NC(n)$ and the symmetric group was discussed briefly by Kreweras \cite{kreweras} and it was studied more explicitly by Biane in 1997 (his interest in the subject was motived by the theory of free probability; see \cite{biane:probability}). Recall that the reflection generating set of $A_{n-1}$ is the set $T=\{(ij):1\leq i<j\leq n\}$ of transpositions, and that the Coxeter elements of $A_{n-1}$ are precisely the $n$-cycles. In our notation, Biane proved \cite[Theorem 1]{biane:crossings} that $NC(n)$ and $NC(A_{n-1})$ are isomorphic posets.

\begin{theorem}
\label{th:biane}
Given a permutation $\pi$ of $[n]$, let $\set{\pi}$ denote the partition of $[n]$ consisting of the orbits of $\pi$. The map $\pi\mapsto\set{\pi}$ is a poset isomorphism
\begin{equation*}
NC(A_{n-1})\longrightarrow NC(n).
\end{equation*}
\end{theorem}
The proof will follow from two fundamental lemmas.

\begin{lemma}
\label{lemma:cycles}
Consider transposition $t=(ij)$ and permutation $\pi$ in $A_{n-1}$.
\begin{itemize}
\item If symbols $i$ and $j$ occur in the same cycle of $\pi$, then in both $t\pi$ and $\pi t$ this cycle is split into two, each containing one of $i$ and $j$.
\item If the symbols $i$ and $j$ occur in two different cycles of $\pi$, then in $t\pi$ and $\pi t$ these two cycles are joined into one.
\end{itemize}
In both cases, the orientations of the two smaller cycles agree with the orientation of the larger cycle.
\end{lemma}

\begin{proof}
We will prove the statements for $t\pi$. The proof for $\pi t$ is analogous. Suppose first that $i$ and $j$ occur in the same cycle of $\pi$, shown below. If we perform $\pi$ and then $t$, observe how the cycle breaks into two:
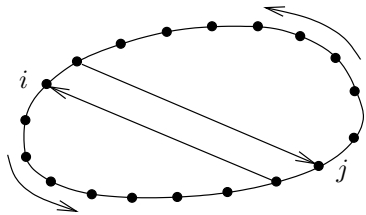
\begin{equation}
\label{fig:break}
\input{cyclebreak.pstex_t}
\end{equation}
Now suppose that $i$ and $j$ occur in two different cycles of $\pi$, shown below. If we perform $\pi$ and then $t$, observe how the two cycles are stitched into one:
\begin{center}
\input{cyclestitch.pstex_t}
\end{center}
\end{proof}

\begin{lemma}
\label{lemma:numbercycles}
Let $\ell_T:A_{n-1}\to\integers$ be the word length on $A_{n-1}$ with respect to the generating set $T$ of transpositions. For a permutation $\pi\in A_{n-1}$ we have
\begin{equation*}
\ell_T(\pi)= n - \text{the number of cycles in } \pi.
\end{equation*}
\end{lemma}
\begin{proof}
Suppose $\pi$ has $n-k$ cycles and note that we can interpret $\ell_T(\pi)$ as the length of a geodesic from $\pi$ to the identity $1$ in the Cayley graph. By Lemma \ref{lemma:cycles}, each step we take in the Cayley graph either joins two cycles or breaks a cycle in two. Since we must break $k$ cycles to reach $1$ from $\pi$, we must take at least $k$ steps, or $\ell_T(\pi)\geq k$. Finally, it is always possible to travel from $\pi$ to $1$ in exactly $k$ steps: recursively, if symbols $i$ and $j$ occur in the same cycle we multiply by the transposition $t=(ij)$. Thus $\ell_T(\pi)=k$.
\end{proof}

\begin{proof}[Proof of Theorem \ref{th:biane}]
Consider $NC(A_{n-1})$ with respect to an $n$-cycle $c\in A_{n-1}$ and note that $NC(A_{n-1})$ consists of those permutations $\pi$ that lie on a geodesic from $c$ to $1$ in the Cayley graph $(A_{n-1},T)$. By Lemma \ref{lemma:numbercycles}, all such $\pi$ are obtained from $c$ by repeated splitting of cycles, as pictured in \eqref{fig:break}. Observe that such splittings preserve the property of ``noncrossing'', and since every noncrossing partition may be obtained in this way we conclude that $\pi\mapsto\set{\pi}$ is a surjection $NC(A_{n-1})\to NC(n)$. Furthermore, if we consider the cycle $c$ as oriented ``counterclockwise'' --- as in \eqref{fig:break} --- all cycles split off from this will also be counterclockwise. Thus every noncrossing partition is achieved only once, and $\pi\mapsto\set{\pi}$ is an injection. Finally, since the partial order on $NC(A_{n-1})$ corresponds to ``splitting cycles'' and the partial order on $NC(n)$ corresponds to ``splitting blocks'' we have an isomorphism.
\end{proof}

Note that $\pi\mapsto\set{\pi}$ gives an isomorphism $NC(A_{n-1})\to NC(n)$ for {\em any} choice of $n$-cycle $c$, where the notion of ``crossing'' in $NC(n)$ is now interpreted with respect to the cyclic order on $[n]$ induced by $c$. (That is, in the circular representation we use $c$ to label the vertices of the $n$-gon.) For a given $c$, we will denote the inverse isomorphism $NC(n)\to NC(A_{n-1})$ by
\begin{equation}
\label{eq:parttoperm}
\P\mapsto\overset{\rightarrow_c}{\P}.
\end{equation}
This map (which Biane \cite{biane:crossings} calls the ``trace'' map) assigns a cyclic permutation to each block of $\P$, oriented according to $c$, and then takes the product of these cycles. When the choice of $n$-cycle (Coxeter element) $c$ is unambiguous, we will denote the map \eqref{eq:parttoperm} simply by $\P\mapsto\overset{\rightarrow}{\P}$.

\begin{figure}
\vspace{.1in}
\begin{center}
\input{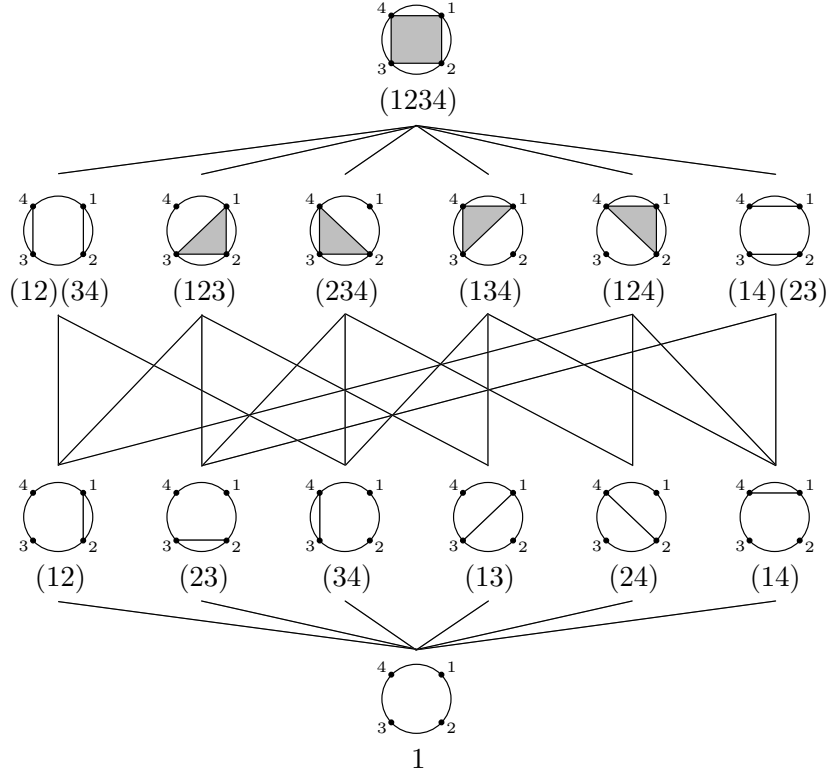}
\end{center}
\caption{$NC(4)$ is isomorphic to $NC(A_3)$}
\label{fig:nc4}
\end{figure}

Figure \ref{fig:nc4} (reproduced from Chapter \ref{sec:introduction}) displays the isomorphism between $NC(4)$ and $NC(A_3)$. Compare this to Figure \ref{fig:nc(a3)} which displays the entire Cayley graph of $(A_3,T)$.

Working simultaneously to Biane, Reiner~\cite{reiner} also generalized the classical noncrossing partitions in the context of reflection groups. His starting point was the fact that the lattice $\Pi(n)$ of {\em unrestricted} partitions of the set $[n]$ is isomorphic to the intersection lattice of the Coxeter arrangement of type $A$ (known as the {\sf braid arrangement}). Recall that the Coxeter arrangement $\A(W)$ of the finite Coxeter system $(W,S)$ is the set of reflecting hyperplanes
\begin{equation*}
\A(W)=\left\{ \alpha^\perp: t_\alpha\in T \right\}
\end{equation*}
for the geometric representation $\sigma:W\hookrightarrow GL(V)$ (see Section \ref{sec:coxsystems}). If $\A(W)=\{\alpha_1^\perp,\ldots,\alpha_N^\perp\}$, then we define the {\sf partition lattice} of $W$ as the set of intersections of reflecting hyperplanes
\begin{equation*}
\Pi(W):=\left\{ \cap_{i\in I} \alpha_i^\perp:I\subseteq\{1,\ldots,N\}\right\}
\end{equation*}
(where we understand the empty intersection to be the whole space, $\cap_{i\in\emptyset} \alpha_i^\perp=V$), and partially order these by {\em reverse}-inclusion of subspaces. If $\reals^n$ has standard basis vectors $\{e_1,\ldots,e_n\}$ then the braid arrangement $\A(A_{n-1})$ can be described concretely as the set of hyperplanes $\{H_{ij}: 1\leq i<j\leq n\}$, where $H_{ij}$ is the zero set of the linear form $e_i-e_j$. In this case it is clear that $\Pi(A_{n-1})$ is isomorphic to the classical partition lattice $\Pi(n)$: the partition $\P\in\Pi(n)$ corresponds to the subspace of $\reals^n$ in which the $i$-th and $j$-th coordinates are equal whenever $i$ and $j$ are in the same block of $\P$. Reiner's idea was to look for a subposet of $\Pi(W)$ whose elements can be called ``noncrossing'' in some natural sense. He did this for the classical types $B$ and $D$, using geometric reasoning, and he extended many of the known combinatorial results. It was later found that his type $B$ noncrossing partitions are isomorphic to $NC(B_n)$, but his type $D$ generalization turned out not to agree with the Coxeter element definition of $NC(D_n)$ (see Sections \ref{sec:typeB} and \ref{sec:typeD}).

One year later, in 1998, Birman, Ko and Lee released an influential paper about the word and conjugacy problems in braid groups, in which they defined a new monoid presentation of the braid group~\cite{birman-ko-lee}. In the study of this presentation, they noticed exactly the same ``noncrossing'' property of the Cayley graph $(A_{n-1},T)$ that Biane had observed, but they used the term ``obstructing'' instead of ``crossing'' and they were not aware of the notion of a noncrossing partition. Following Birman-Ko-Lee, there was a flurry of activity in the field of combinatorial group theory. Two independent streams emerged.

In the first, Brady isolated the importance of the Coxeter elements~\cite{brady2}, and he used the Birman-Ko-Lee monoid presentation to construct new $K(\pi,1)$'s for the braid groups \cite{brady}. In \cite{brady}, Brady considered the poset $NC(A_{n-1})$ as part of a Garside structure, and so he needed to investigate the lattice property. This led him to discover the work of Kreweras~\cite{kreweras} and Reiner \cite{reiner} and the combinatorial literature on noncrossing partitions. At this point, he collaborated with Watt to study the poset $NC(W)$ for a general finite Coxeter group. In~\cite{brady-watt:partialorder} they defined a partial order on the orthogonal group and proved the uniqueness property of moved spaces (Theorem \ref{theorem:hardMov}). Then in \cite{brady-watt:kpione} they gave the first published definition of $NC(W)$ and extended the work \cite{brady} to all finite type Artin groups. One of the motivations for \cite{brady-watt:kpione} was to demonstrate that their $NC(D_n)$ differed from Reiner's type $D$ noncrossing partitions.

In the second stream, Bessis, Digne and Michel also recognized the importance of the Coxeter elements in the Birman-Ko-Lee monoid, which they interpreted from the perspective of Springer theory \cite{bessis-digne-michel}. They realized how to define $NC(W)$ for other Coxeter groups, but they followed Birman-Ko-Lee in calling these elements ``non-obstructing'' and so they did not discover the combinatorial literature on noncrossing partitions. Bessis went on to write an extensive study \cite{bessis:dual} of the poset $NC(W)$, in which he generalized the Birman-Ko-Lee monoid to all finite type Artin groups. This work was itself quite influential; it contains the first mention of a ``dual Coxeter system'' (see Definition \ref{def:dualCoxeter}).

After Bessis' paper \cite{bessis:dual} appeared on the arXiv in 2001, Brady contacted him and shared his knowledge of classical noncrossing partitions and Reiner's work \cite{reiner}. Bessis, Digne and Michel were able to switch the terminology from ``non-obstructing'' to ``non-crossing'' in their paper \cite{bessis-digne-michel} before publication. Meanwhile, Biane, Goodman and Nica~\cite{biane:typeB} independently discovered the Cayley graph interpretation of Reiner's type $B$ noncrossing partitions and applied this to free probability \cite{biane:typeB}. They also became aware of Brady-Watt and Bessis before publication. By 2003, everyone was on the same page.

Due to the fact that these researchers were working in different fields --- Reiner in combinatorics; Brady, Watt and Bessis in combinatorial group theory; and Biane in free probability --- it took some years for them to realize that they were working with the same objects. Once this coincidence became apparent, a workshop was held at the American Institute of Mathematics in January 2005 \cite{armstrong:braids} at which all of the different perspectives came together for the first time.

\section{The Classical Kreweras Complement}
\label{sec:classicalkreweras}
In Section \ref{sec:shifting} we encountered a family of anti-automorphisms $K_\mu^\nu:[\mu,\nu]\to[\mu,\nu]$, defined by $K_\mu^\nu(\pi)=\mu\pi^{-1}\nu$, which exhibit the local self-duality of the poset $NC(W)$. Now we will examine the classical (type $A$) analogues of these maps. We use the letter $K$ in honor of Germain Kreweras, who defined and used the type $A$ version of the map $K_1^c:NC(W)\to NC(W)$ in the seminal paper \cite{kreweras}. Nica and Speicher~\cite{nica-speicher:ntuples}  later considered the type $A$ version of $K_1^\nu:[1,\nu]\to[1,\nu]$, which they called the ``relative Kreweras complement''. In this section, we will recall both of these constructions, and take the generalization to its natural conclusion by defining the classical version of $K_\mu^\nu$ for all $\mu\leq_T \nu$ in $NC(A_{n-1})$. The classical Kreweras complement is essential to the main results of this chapter.

First, we need some arithmetic on set partitions. Given $a\in\integers\setminus\{0\}$ and a set of integers $X\subset\integers$, we define two new sets with the same cardinality as $X$: the {\sf translation}, $X+a:=\{x+a:x\in X\}$, and the {\sf dilation}, $aX:=\{ax:x\in X\}$. These operations extend to partitions in an obvious way.

\begin{definition}
Let $\P=\{P_1,P_2,\ldots,P_m\}$ be a partition of $[n]$ and consider $a\in\integers\setminus\{0\}$. Then
\begin{equation*}
\P+a:=\{ P_1+a,P_2+a,\ldots,P_m+a\}
\end{equation*}
is a partition of $[n]+a$, called the {\sf translation} of $\P$ by $a$, and
\begin{equation*}
a\P:=\{aP_1,aP_2,\ldots,aP_m\}
\end{equation*}
is a partition of $a[n]$, called the {\sf dilation} of $\P$ by $a$.
\end{definition}

This notation allows us to express the ``interleaving'' of partitions.

\begin{definition}
\label{def:interleave}
Let $\P$ and $\Q$ be partitions of $[n]$. The partition
\begin{equation*}
\langle\P,\Q\rangle:=(2\P-1)\cup(2\Q)
\end{equation*}
of $[2n]$ is called the {\sf interleaving} of $\P$ and $\Q$.
\end{definition}
This is an intuitive idea if we think of the integers $[2n]$ labelling the vertices of a regular $2n$-gon: $\P$ defines a partition on the odd vertices $\{1,3,\ldots,2n-1\}$ and $\Q$ defines a partition on the even vertices $\{2,4,6,\ldots,2n\}$. For example, we have
\begin{equation*}
\Big\langle \left\{\{1\},\{2,3,4\}\right\},\left\{\{1,2\},\{3,4\}\right\}\Big\rangle=\{\{1\},\{2,4\},\{3,5,7\},\{6,8\}\}.
\end{equation*}
Notice here that $\P\in NC(4)$ and $\Q\in NC(4)$ are noncrossing, but that the interleaving $\langle\P,\Q\rangle$ is {\em crossing}. If we begin with noncrossing partitions $\P$ and $\Q$, it is natural to ask when the interleaving $\langle\P,\Q\rangle$ will also be noncrossing. The answer leads to the definition of the classical Kreweras complement. We follow \cite[Section 3]{kreweras}.

\begin{definition}
\label{def:kreweras}
Given $\P\in NC(n)$, the classical {\sf Kreweras complement} $K(\P)$ of $\P$ is defined to be the coarsest partition $\Q$ of $[n]$ such that the interleaving $\langle\P,\Q\rangle$ is noncrossing.
\end{definition}

\begin{figure}
\vspace{.1in}
\begin{center}
\input{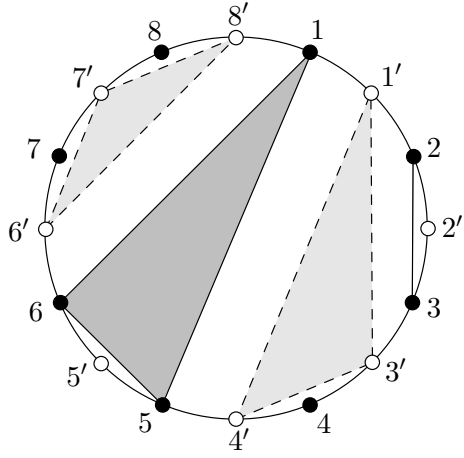}
\end{center}
\caption{Illustration of the Kreweras complement}
\label{fig:kreweras}
\end{figure}

This definition is best understood pictorially. For example, Figure \ref{fig:kreweras} demonstrates that
\begin{equation*}
K\Big(\{\{1,5,6\},\{2,3\},\{4\},\{7\},\{8\}\} \Big)=\{\{1,3,4\},\{2\},\{5\},\{6,7,8\}\}.
\end{equation*}

Note that the composition $K\circ K$ is just a counterclockwise rotation of the circular diagram, hence $K$ defines a self-bijection on $NC(n)$. It is easy to see that $K$ reverses order, since if two blocks of $\P$ are joined into one, then some block of $K(\P)$ is split into two blocks. Moreover, if the partition $\P$ has $m$ blocks, then the Kreweras complement $K(\P)$ has $n-m+1$ blocks, since each block added to $\P$ reduces the number of blocks in $K(\P)$ by $1$. 

Now consider the isomorphism $NC(n)\cong NC(A_{n-1})$ (Theorem \ref{th:biane}) with respect to the Coxeter element $c=(12\cdots n)$. In particular, one may observe from Figure \ref{fig:kreweras} that we have
\begin{equation*}
(156)(23)\cdot (134)(678) = (12345678),
\end{equation*}
where the product is taken in the symmetric group $A_7$. That is, the permutation $(134)(678)=\left((156)(23)\right)^{-1}\cdot (12345678)$ is the group-theoretical Kreweras complement of the permutation $(156)(23)$. Recall that $\P\mapsto\overset{\rightarrow_c}{\P}$ is the map that takes a partition to its corresponding permutation with respect to the cyclic order $c$. In general, Kreweras noticed \cite[Section 3]{kreweras} that
\begin{equation}
\label{eq:krewworks}
\overset{\rightarrow_c}{\P}\,\cdot\,\overset{\longrightarrow_c}{K(\P)}=c,
\end{equation}
where the product is taken in $A_{n-1}$. It follows immediately that
\begin{equation*}
K(\{\pi\})=\{\pi^{-1}c\}=\{K_1^c(\pi)\}
\end{equation*}
for all $\pi\in A_{n-1}$, hence the classical Kreweras complement agrees with the group-theoretical Kreweras complement on $NC(A_{n-1})$ and there should be no ambiguity in using the letter $K$ for both purposes.

Notice, it is the planarity of the circular representation of a noncrossing partition that allows the existence of this sort of complementation map. Indeed, the lattice $\Pi(n)$ of unrestricted set partitions has no such planarity property, and it is {\em not} self-dual. The relationship between the posets $NC(n)$ and $\Pi(n)$ is very interesting, and we will return to this later.

Now we wish to generalize the Kreweras complement to each interval in $NC(n)$ --- to construct the classical version of the map $K_\mu^\nu$ for all $\mu\leq_T\nu$ in $NC(A_{n-1})$. We will follow the examples of Nica and Speicher in~\cite{nica-speicher:ntuples} and~\cite{speicher:multiplicative}.

We have seen that every interval in $NC(W)$ is isomorphic to $NC(W')$ for some parabolic subgroup $W'\subseteq W$, and hence is self-dual. Speicher explained the local self-duality of the lattice $NC(n)$ by noting that every interval $[\M,\N]$ in $NC(n)$ decomposes as a direct product of posets $\prod_i NC(n_i)$ for some integers $n_i\leq n$. Fixing such a decomposition, there is an obvious anti-automorphism on the interval $[\M,\N]$: namely, the direct product of the Kreweras complements from each factor. After establishing some notation, we will show that this idea leads to the correct definition.

First, for any finite set of integers $X\subseteq\integers$ it is natural to consider the lattice of partitions of $X$ that are noncrossing with respect to the usual linear order on $X$. We will denote this lattice by $NC(X)$. Then we have $NC(X)\cong NC(\abs{X})$ with respect to the order-preserving bijection $X\leftrightarrow \{1,2,\ldots,\abs{X}\}$. Let $\hat{1}_X=\{X\}$ and $\hat{0}_X=\{\{x\}:x\in X\}$ denote the maximum and minimum elements of this lattice, respectively, and let $K_X$ denote the classical Kreweras complement on $NC(X)$. Further, for any partition $\P\in NC(X)$ and any subset $U\subseteq X$, let $\P| U\in NC(U)$ denote the noncrossing partition of $U$ that is the restriction of $\P$. Following Speicher \cite[Proposition 1]{speicher:multiplicative} we obtain a decomposition exhibiting the local self-duality of $NC(n)$.

\begin{theorem}[\cite{speicher:multiplicative}]
\label{th:interval}
Each interval $[\M,\N]$ in $NC(n)$ decomposes as a direct product of posets
\begin{equation}
\label{eq:interval}
[\M,\N]\cong \prod_{i,j} NC(X_{i,j}),
\end{equation}
where $X_{i,j}$ is the $j$-th block of the partition $K_{N_i}^{-1}(\M| N_i)\in NC(N_i)$, and $N_i$ is the $i$-th block of $\N$, ordered arbitrarily.
\end{theorem}

\begin{proof}
Note that each interval $[\M,\N]$ in $NC(n)$ decomposes according to the blocks of $\N$ in an obvious way,
\begin{equation}
\label{eq:NCdecomp}
[\M,\N]\cong \prod_i \left[\M| N_i,\hat{1}_{N_i}\right],
\end{equation}
where $N_i$ is the $i$-th block of $\N$. This isomorphism is canonical, given by sending $\P\in[\M,\N]$ to the sequence of restrictions $\left(\P|N_i\right)_i\in\prod_i \left[ \M|\N_i,\hat{1}_{N_i}\right]$.

By self-duality of $NC(N_i)$, each of these upper intervals is isomorphic to a lower interval in $NC(N_i)$,
\begin{equation}
\label{eq:NCdecomp2}
\left[\M|N_i,\hat{1}_{N_i}\right]\cong \left[\hat{0}_{N_i},K_{N_i}^{-1}\left(\M|N_i\right)\right],
\end{equation}
but this isomorphism is not canonical since it depends on a choice of anti-\linebreak[4]automorphism (here we have chosen $K_{N_i}^{-1}$). Applying the decomposition \eqref{eq:NCdecomp} to the interval on the right side of \eqref{eq:NCdecomp2} yields
\begin{equation}
\label{eq:NCdecomp3}
\left[\hat{0}_{N_i},K_{N_i}^{-1}\left(\M|N_i\right)\right]\cong\prod_j \left[\hat{0}_{X_{i,j}},\hat{1}_{X_{i,j}}\right]=\prod_j NC(X_{i,j}),
\end{equation}
where $X_{i,j}$ is the $j$-th block of $K_{N_i}^{-1}\left(\M|N_i\right)$. Composing \eqref{eq:NCdecomp}, \eqref{eq:NCdecomp2} and \eqref{eq:NCdecomp3} yields the desired isomorphism.
\end{proof}

This decomposition is not canonical but it is unique up to poset isomorphism. Thus, for each interval $[\M,\N]$ in $NC(n)$ there is a corresponding multiset of integers $\left\{ \abs{X_{i,j}}\right\}_{i,j}$ describing its isomorphism type. This observation has a nice algebraic interpretation. If we consider the isomorphism $NC(n)\cong NC(A_{n-1})$, we know that every interval $[\M,\N]\subseteq NC(n)$ is isomorphic to $NC(W')$, where $W'$ is some parabolic subgroup of $A_{n-1}$. The parabolic subgroups of a given finite Coxeter group $W$ have isomorphism types given by the induced subgraphs of its Coxeter diagram. Since the Coxeter diagram of $A_{n-1}$ is a chain (Figure \ref{fig:coxdiagrams}), all of its induced subgraphs are disjoint unions of chains. That is, every parabolic subgroup $W'$ of $A_{n-1}$ has isomorphism type
\begin{equation*}
W'\cong A_{x_1}\times A_{x_2}\times\cdots\times A_{x_r},
\end{equation*}
and hence
\begin{equation*}
[\M,\N]\cong NC(W')\cong \prod_{i} NC(A_{x_i})\cong \prod_i NC(x_i+1)
\end{equation*}
for the canonical multiset of integers $\{x_i+1\}_i=\left\{\abs{X_{i,j}}\right\}_{i,j}$. Thus, we might describe this multiset of integers as the {\sf type} of the interval. For other families of finite Coxeter groups $W$ the situation is not as straightforward.

Consider again \eqref{eq:interval}. Explicitly, the isomorphism $[\M,\N]\to \prod_{i,j} NC(X_{i,j})$ is given by
\begin{equation}
\label{eq:intisom}
\P\mapsto \left( K_{N_i}^{-1}\left(\P|N_i\right)|X_{i,j}\right)_{i,j},
\end{equation}
and the inverse isomorphism is given by
\begin{equation}
\label{eq:intinverse}
\left(\Q_{i,j}\right)_{i,j}\mapsto \cup_i K_{N_i}\left( \cup_j\Q_{i,j}\right),
\end{equation}
where each $\Q_{i,j}$ is in $NC(X_{i,j})$. It is now clear how to define an anti-automorphism of the interval $[\M,\N]\subseteq NC(n)$. We follow \cite[Definition 2.4]{nica-speicher:ntuples}.

\begin{definition}
\label{def:relkreweras}
Given $\M\leq\N$ in $NC(n)$, the {\sf relative Kreweras complement} $K_\M^\N$ is the unique map completing the square
\[
\xymatrix{
[\M,\N] & \displaystyle{\prod} NC(X_{i,j})\ar[l]\\
[\M,\N] \ar [u]^{K_\M^\N}\ar[r] & \displaystyle{\prod} NC(X_{i,j})\ar[u]_{\prod K_{X_{i,j}}}\\
}
\]
where the bottom arrow is given by \eqref{eq:intisom}, the top arrow by \eqref{eq:intinverse}, and the right arrow is the direct product of the Kreweras complements $K_{X_{i,j}}$ on the factors $NC(X_{i,j})$.
\end{definition}

Since this notation is a bit opaque, we will examine some important special cases. If $\M=\hat{0}_n$ and $\N=\hat{1}_n$ then the decomposition \eqref{eq:interval} reduces to the obvious
\begin{equation*}
\left[\hat{0}_n,\hat{1}_n\right]\cong NC(n),
\end{equation*}
and in this case the map $K_{\hat{0}_n}^{\hat{1}_n}$ is equal to the classical Kreweras complement $K$, as desired. When possible, we will drop the subscript $\hat{0}_n$ or the superscript $\hat{1}_n$.

\begin{figure}
\vspace{.1in}
\begin{center}
\input{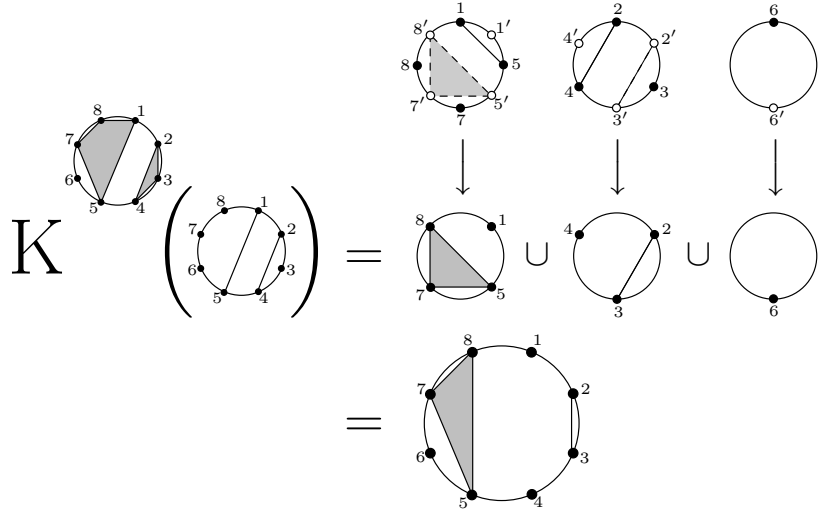}
\end{center}
\caption{Illustration of the relative Kreweras complement}
\label{fig:relKexample}
\end{figure}

The most important special case from our perspective is the ``half-general'' case, when the bottom is fixed $\M=\hat{0}_n$ and the top $\N$ is free. In this case, the interval $\left[\hat{0}_n,\N\right]$ decomposes naturally according to the blocks of $\N$,
\begin{equation*}
\left[ \hat{0}_n,\N\right]\cong \prod_{i=1}^t NC(N_i),
\end{equation*}
where $\N=\{N_1,\ldots,N_t\}$. Then given any $\P\leq\N$ in $NC(n)$, it is easy to check that the relative Kreweras complement $K^\N(\P)$ is given by
\begin{equation*}
K^\N(\P)=\bigcup_{i=1}^t K_{N_i}\left(\P|N_i\right).
\end{equation*}
That is, we divide $\P$ according to the blocks of $\N$, and take the usual Kreweras complement within each block. Then we put the pieces back together. Figure \ref{fig:relKexample} shows an example of this type of calculation. Notice that if $\P$ has $m$ blocks and $\N$ has $t$ blocks, then $K^\N(\P)$ has $n-m+t$ blocks. Indeed, we saw before that $K_{N_i}\left(\P|N_i\right)$ has $\abs{N_i}-\abs{\P|N_i}+1$ blocks. Summing over $i$, we get
\begin{equation}
\label{eq:relKblocks}
\sum_{i=1}^t\left( \abs{N_i}-\abs{\P|N_i}+1\right) = \sum_{i=1}^t \abs{N_i}-\sum_{i=1}^t \abs{\P|N_i}+t=n-m+t.
\end{equation}
Of course, this also follows from the fact that $K^\N$ is an anti-automorphism, since it sends $\P$ to an element of complementary rank within $\left[\hat{0}_n,\N\right]$.

Finally, we will show that this classical relative Kreweras complement agrees with the group-theoretical Kreweras complement. Having done this, we will proceed to a discussion of the classical $k$-divisible noncrossing partitions.

\begin{lemma}
\label{lemma:mess}
For all $\M\leq\P\leq\N$ in $NC(n)$ we have
\begin{equation*}
K_\M^\N(\P)=\M\vee K^\N(\P).
\end{equation*}
\end{lemma}

\begin{proof}
This follows from two elementary facts, which are clear after considering circular representations. First, for any sets $U\subseteq X\subseteq[n]$ and partition $\P\in NC(X)$ we have
\begin{equation}
\label{eq:restrict}
K_X(\P)| U= K_U\left(\P| U\right).
\end{equation}
Second, for any $\P$ and $\Q$ in $NC(X)$ with $\Q=\{Q_1,\ldots,Q_m\}$ we have
\begin{equation}
\label{eq:meet}
\bigcup_{i=1}^m \P|Q_i=\Q\wedge \P.
\end{equation}
Now, to compute $K_\M^\N(\P)$ we follow Definition \ref{def:relkreweras} and the following sequence of notationally daunting but easy steps:
\begin{align*}
\P \mapsto & \left( K_{N_i}^{-1}(\P|N_i)|X_{i,j} \right)_{i,j} & \text{from \eqref{eq:intisom}}\\
  \mapsto & \left[ K_{X_{i,j}} \left( K_{N_i}^{-1}(\P|N_i)|X_{i,j} \right)\right]_{i,j} & \text{by definition}\\
  = & \left[ (\P|N_i)|X_{i,j} \right]_{i,j} & \text{from \eqref{eq:restrict}}\\
  \mapsto & \cup_i K_{N_i}\left( \cup_j (\P|N_i)|X_{i,j}\right) & \text{from \eqref{eq:intinverse}}\\
  = & \cup_i K_{N_i}\left( K_{N_i}^{-1}(\M|N_i)\wedge (\P|N_i)\right) & \text{from \eqref{eq:meet}}\\
  = & \cup_i \left[ (\M|N_i)\vee K_{N_i}(\P|N_i)\right] \\
  = & \left(\cup_i (\M|N_i)\right) \vee \left(\cup_i K_{N_i}(\P|N_i)\right) \\
  = & (\N\wedge\M) \vee K^\N(\P) = \M \vee K^\N(\P).
\end{align*}
\end{proof}

\begin{theorem}
For all $\mu\leq_T\pi\leq_T\nu$ in $NC(A_{n-1})$ we have
\begin{equation*}
\left\{ K_\mu^\nu(\pi)\right\} = \left\{ \mu\pi^{-1}\nu\right\}= K_{\{\pi\}}^{\{\nu\}}(\{\pi\}),
\end{equation*}
where $\pi\mapsto\{\pi\}$ is the isomorphism $NC(A_{n-1})\to NC(n)$ in Theorem \ref{th:biane}.
\end{theorem}

\begin{proof}
The fact that $\left\{ K_1^\nu(\pi)\right\}=K^{\{\nu\}}(\{\pi\})$ was proved by Nica and Speicher \cite[Section 2.5]{nica-speicher:ntuples}, and it follows from an argument exactly analogous to \eqref{eq:krewworks}. Applying Lemmas \ref{lemma:kreweras} and \ref{lemma:mess}, we have
\begin{align*}
\set{K_\mu^\nu(\pi)}= & \set{\mu\vee K^\nu(\pi)}=\set{\mu}\vee\set{K^\nu(\pi)} \\
  = & \set{\mu}\vee K^{\set{\nu}}(\set{\pi}) = K_{\set{\mu}}^{\set{\nu}}(\set{\pi}).
\end{align*}
\end{proof}

Now that we have established the precise relationship between the Kreweras complements on $NC(A_{n-1})$ and $NC(n)$, we will freely move between the group-theoretical and the classical point of view.

\section{Classical $k$-Divisible Noncrossing Partitions}
\label{sec:classicalNCk}
We have come to the motivating example of this memoir, as we outlined earlier in Section \ref{sec:motivation}. Here, we will describe the classical $k$-divisible noncrossing partitions that have motivated our algebraic work in Chapter \ref{sec:kdiv}. For all positive integers $n$ and $k$, it turns out that the poset $NC^{(k)}(A_{n-1})$ (Definition \ref{def:NCk}) is isomorphic to the subposet of $NC(kn)$ consisting of ``$k$-divisible partitions'' --- that is, each of whose blocks has cardinality divisible by $k$. This result (Theorem \ref{th:mainisom}) together with a type $B$ analogue (Theorem \ref{th:mainisomB}) motivates our general use of the term ``$k$-divisible''.

We begin with the definition.

\begin{definition}
\label{def:classicalNCk}
Let $NC^{(k)}(n)$ denote the induced subposet of $NC(kn)$ consisting of partitions in which each block has cardinality divisible by $k$.
\end{definition}

\begin{figure}
\vspace{.1in}
\begin{center}
\input{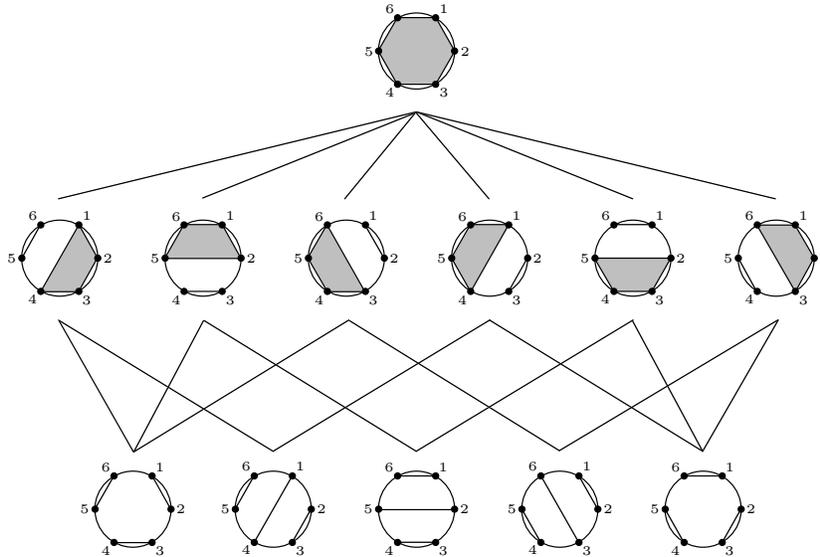}
\end{center}
\caption{The poset $NC^{(2)}(3)$ of $2$-divisible noncrossing partitions of the set $[6]$}
\label{fig:nc^2(3)}
\end{figure}

Figure \ref{fig:nc^2(3)} displays the poset $NC^{(2)}(3)$, consisting of $2$-divisible noncrossing partitions of the set $[6]$. We notice immediately that this poset is isomorphic to $NC^{(2)}(A_2)$ as shown in Figure \ref{fig:nc^2(a2)}. However, it is not transparent how to formalize this isomorphism.

Since the coarsening of partitions preserves the property of $k$-divisibility, $NC^{(k)}(n)$ is an order filter in $NC(kn)$. Hence the $k$-divisible noncrossing partitions form a graded join-semilattice. As in the $k=1$ case \eqref{eq:rankblocks}, the rank function is given by $n$ minus the number of blocks, $\rk(\P)=n-\abs{\P}$. This poset was introduced by Edelman~\cite{edelman:kdiv}, who calculated many of its enumerative invariants, including the zeta polynomial. It was later considered by Stanley~\cite{stanley:parking} in connection with parking functions and the Ehrenborg quasisymmetric function of a poset.

\subsection{Shuffle Partitions}

Recall that the poset $NC^{(k)}(W)$ was defined as the componentwise order on delta sequences. To really understand the isomorphism between Figures \ref{fig:nc^2(a2)} and \ref{fig:nc^2(3)}, we must confront the idea of a {\em classical} delta sequence. However, since our notion of a delta sequence (Definition \ref{def:mchaindsequence}) is essentially algebraic, it is not clear what the definition should be in the classical case. Following Theorem \ref{th:biane}, notice that $NC^k(A_{n-1})=NC(A_{n-1}^k)$ is isomorphic to the poset $NC^k(n):=(NC(n))^k$ via the map
\begin{equation}
\label{eq:kmap}
(\pi)_k=(\pi_1,\pi_2,\ldots,\pi_k)\mapsto (\{\pi\})_k=(\{\pi_1\},\{\pi_2\},\ldots,\{\pi_k\}).
\end{equation}
Thus, we make the following definition.

\begin{definition}
\label{def:classicaldelta}
We say that a sequence $(\Q)_k\in NC^k(n)$ is a {\sf (classical) delta sequence} if it is the image of some delta sequence $(\pi)_k\in NC(A_{n-1}^k)$ under the isomorphism \eqref{eq:kmap}.
\end{definition}

Unlike a delta sequence, the notion of a multichain makes perfect sense in $NC(n)$, and since the reciprocal bijections $\partial$ \eqref{partial} and $\sint$  \eqref{sint} are also transferred to $NC^k(n)$ via \eqref{eq:kmap}, one could alternatively define classical delta sequences in terms of multichains.

\begin{lemma}
The sequence $(\Q)_k$ in $NC^k(n)$ is a delta sequence if and only if there exists some multichain $(\P)_k=(\P_1,\ldots,\P_k)$ in $NC^k(n)$ such that
\begin{equation*}
(\Q)_k=\partial(\P)_k=\left( K^{\P_2}(\P_1),\ldots,K^{\P_k}(\P_{k-1}),K(\P_k)\right).
\end{equation*}
\end{lemma}

This definition still seems rather arbitrary.  We will show that the nature of classical delta sequences is best expressed by the idea of interleaving partitions. Generalizing Definition \ref{def:interleave}, we can interleave a sequence of partitions as follows.

\begin{definition}
\label{def:shuffle}
Given a sequence $(\Q)_k$ of partitions of $[n]$, the partition
\begin{equation*}
\langle\Q\rangle_k=\langle\Q_1,\Q_2,\ldots,\Q_k\rangle:=\bigcup_{i=1}^k \left(k\Q_i-(k-i)\right)
\end{equation*}
of $[kn]$ is called the {\sf shuffle partition} of $(\Q)_k$.
\end{definition}

Again, this is an intuitive idea if we consider the circular representation of the shuffle partition $\langle\Q\rangle_k$ on the $kn$-gon: $\Q_1$ defines a partition on the vertices $\{1,k+1,2k+1,\ldots,k(n-1)+1\}$, $\Q_2$ defines a partition on $\{2,k+2,2k+2,\ldots,k(n-1)+2\}$, and so on. Then, if the sequence $(\Q)_k$ consists of noncrossing partitions, it is natural to ask when the shuffle $\langle\Q\rangle_k$ of these partitions will be noncrossing. In the case $k=2$, the answer led to the definition of the Kreweras complement (Definition \ref{def:kreweras}). The solution to the general problem explains the significance of the classical delta sequences.

\begin{theorem}
\label{th:shuffle}
Given $(\Q)_k\in NC^k(n)$, the shuffle partition $\langle\Q\rangle_k$ is noncrossing if and only if $(\Q)_k$ is a delta sequence.
\end{theorem}

\begin{proof}
If $\langle\Q\rangle_k$ is noncrossing then $(\Q)_k$ must be in $NC^k(n)$. Recalling Definition \ref{def:kreweras}, it is easy  to see that $\langle\Q\rangle_k$ is noncrossing if and only if we have $\Q_j\leq K(\Q_i)$ for all $1\leq i<j\leq k$. By Lemma \ref{lemma:delta}, this property is equivalent to the fact that $(\Q)_k$ is a delta sequence.
\end{proof}

In this sense, the notion of a delta sequence generalizes the definition of the Kreweras complement. In general, we will use the term ``shuffle partition'' to describe any noncrossing partition that arises in this way.

\begin{definition}
\label{def:shufflepartitions}
The partitions in $NC(kn)$ of the form
\begin{equation*}
NC_{(k)}(n):=\left\{\langle\Q\rangle_k:(\Q)_k\in NC^k(n) \text{ is a delta sequence }\right\}
\end{equation*}
are called {\sf $k$-shuffle partitions}.
\end{definition}

\begin{figure}
\vspace{.1in}
\begin{center}
\input{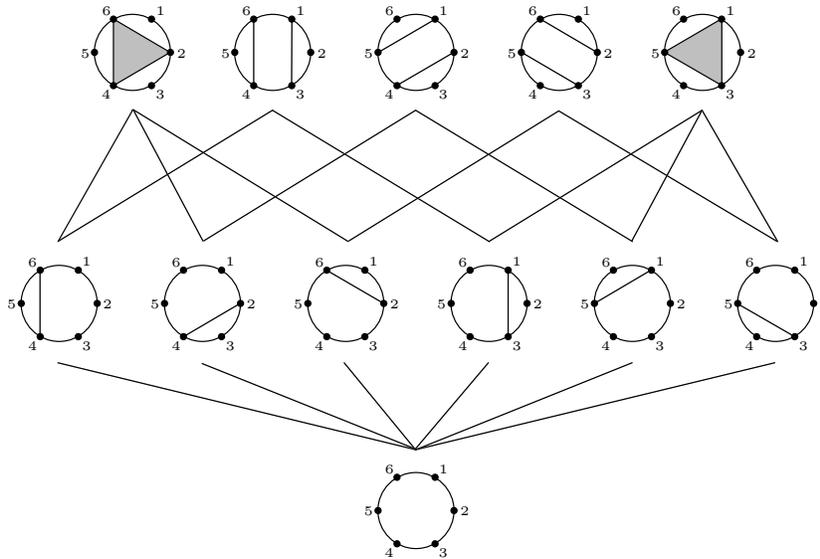}
\end{center}
\caption{The poset $NC_{(2)}(3)$ of $2$-shuffle noncrossing partitions of the set $[6]$}
\label{fig:2-shuffle}
\end{figure}

It is easy to see that a $k$-shuffle partition $\Q$ is characterized by the fact that integers in the same block of $\Q$ are congruent modulo $k$, and since this property is preserved by refinement it follows that $NC_{(k)}(n)$ is an order ideal in $NC(kn)$. For example, Figure \ref{fig:2-shuffle} displays the $2$-shuffle noncrossing partitions of the set $[6]$, in which each block is contained within $\{1,3,5\}$ or $\{2,4,6\}$. Notice that this is an order ideal in $NC(6)$.

Now we have an order ideal of $k$-shuffle partitions and an order filter of $k$-divisible partitions in $NC(kn)$. Considering Figures \ref{fig:nc^2(3)} and \ref{fig:2-shuffle}, it seems plausible that these two posets are dual to each other. Indeed they are.

\begin{lemma}
\label{lemma:kdivkshuffle}
The Kreweras complement $K$ (or any odd power of $K$) on $NC(kn)$ is an anti-isomorphism between the order filter $NC^{(k)}(n)$ of $k$-divisible partitions and the order ideal $NC_{(k)}(n)$ of $k$-shuffle partitions.
\end{lemma}

\begin{proof}
We must show that $K$ takes $k$-divisible partitions to $k$-shuffle partitions and that $K^{-1}$ takes $k$-shuffle partitions to $k$-divisible partitions. To see this, we consider the pictorial representation of the Kreweras complement (Figure \ref{fig:kreweras}). When $k=1$, the result is trivial, so suppose that $k>1$.

Now let $\Q\in NC(kn)$ be $k$-divisible, and consider $1\leq i<j\leq kn$ with $i$ and $j$ in the same block of $K(\Q)$. Then the set $\{i+1,\ldots,j\}$ must be equal to a union of blocks of $\Q$. Since $\Q$ is $k$-divisible, we conclude that $j-i$ is divisible by $k$, and hence $K(\Q)$ is $k$-shuffle.

Conversely, suppose that $\Q\in NC(kn)$ is {\em not} $k$-divisible. In this case, there exist $1\leq i\leq j<kn$ such that $\{i,\ldots,j\}$ is a union of blocks of $\Q$, exactly one of which has size not divisible by $k$, and consequently $j-i+1$ is not divisible by $k$. Since $i$ and $j+1$ are contained in the same block of $K^{-1}(\Q)$, we conclude that $K^{-1}(\Q)$ is not $k$-shuffle.
\end{proof}

Before moving on, we note that the concept of a shuffle is common in combinatorics. Given two ``alphabets'' $\{a_1,a_2,\ldots,a_n\}$ and $\{b_1,b_2,\ldots,b_n\}$, a {\sf shuffle} is any word containing all $2n$ symbols in which the $a$'s and $b$'s occur in their natural order. Greene defined a partial order on shuffles~\cite{greene} which was later studied by Simion and Stanley~\cite{simion-stanley} and shares many features in common with $NC(n)$ (see~\cite{simion-stanley,speicher:multiplicative}). The $2$-shuffle partitions we have defined correspond to the {\sf regular shuffle}
\begin{equation*}
a_1b_1a_2b_2\cdots a_nb_n,
\end{equation*}
and in general one may define a subposet of $NC(kn)$ corresponding to any shuffle on $k$ alphabets of size $n$. Of course, the isomorphism type of this poset will depend only on the cyclic order on the shuffle word. This family of shuffle subposets of $NC(kn)$ interpolates somehow between the extreme cases $NC^k(n)$ and $NC^{(k)}(n)$. Perhaps this is an interesting idea to pursue.

\subsection{The Main Isomorphism}
The characterization of the classical delta sequences as shuffle partitions is the key to understanding the classical $k$-divisible noncrossing partitions. Bringing together our observations so far, it is clear how to write down an isomorphism between $NC^{(k)}(A_{n-1})$ and $NC^{(k)}(n)$. The following result is the prime motivation for everything else in this memoir.
\begin{theorem}
\label{th:mainisom}
The map $(\pi)_k\mapsto K^{-1}\Big\langle\partial\big(\{\pi\}\big)_k\Big\rangle$, given explicitly by
\begin{equation}
\label{eq:mainisom}
(\pi_1,\ldots,\pi_k)\mapsto K^{-1}\bigg\langle\Big( K^{\{\pi_2\}}(\{\pi_1\}),\ldots,K^{\{\pi_k\}}(\{\pi_{k-1}\}),K(\{\pi_k\})\Big)\bigg\rangle,
\end{equation}
is an isomorphism from the $k$-divisible noncrossing partitions of the symmetric group $NC^{(k)}(A_{n-1})$ to the classical $k$-divisible noncrossing partitions $NC^{(k)}(n)$.
\end{theorem}

\begin{proof}
Given two delta sequences $(\Q)_k$ and $(\R)_k$ in $NC^k(n)$, note that the shuffles are related by $\langle\Q\rangle_k\leq \langle\R\rangle_k$ in $NC(kn)$ if and only if $(\Q)_k\leq (\R)_k$ componentwise in $NC^k(n)$, and by Definition \ref{def:classicaldelta} this happens if and only if $(\delta)_k\leq (\varepsilon)_k$, where $(\Q)_k=(\{\delta\})_k$ and $(\R)_k=(\{\varepsilon\})_k$, which is equivalent to $\sint(\varepsilon)_k\leq\sint(\delta)_k$ in $NC^{(k)}(W)$. Thus, the map $(\pi)_k\mapsto\left\langle\partial\left(\{\pi\}\right)_k\right\rangle$ is an anti-isomorphism from $NC^{(k)}(A_{n-1})$ to the order ideal of $k$-shuffle partitions in $NC(kn)$.

On the other hand, Lemma \ref{lemma:kdivkshuffle} says that $K^{-1}$ is an anti-isomorphism from $k$-shuffle partitions to $k$-divisible partitions in $NC(kn)$. The composition of these two anti-isomorphisms is the desired isomorphism.
\end{proof}

\begin{figure}
\vspace{.1in}
\begin{center}
\input{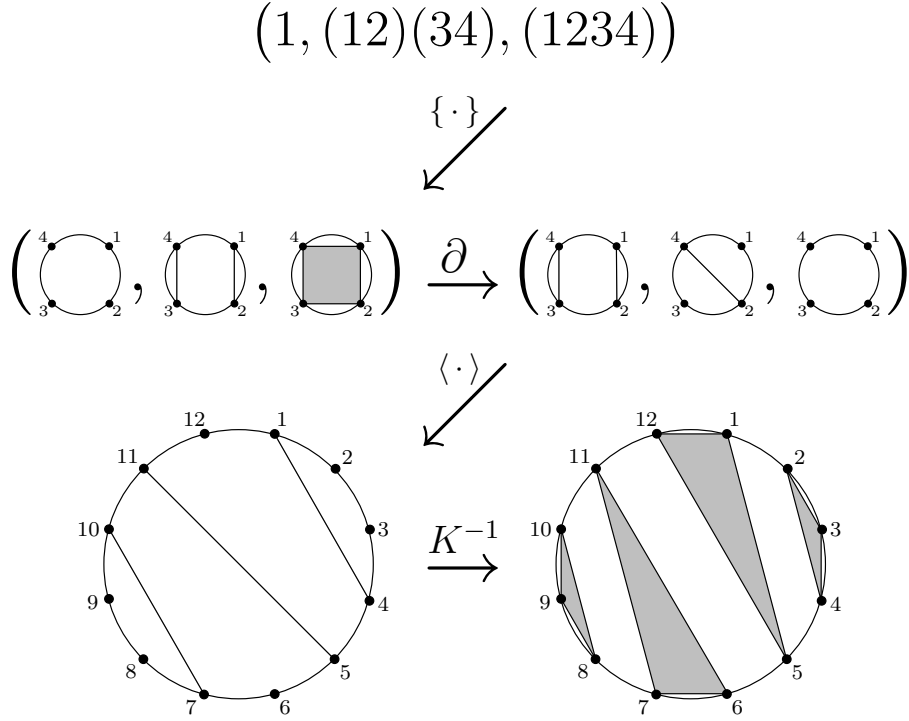}
\end{center}
\caption{An example of the isomorphism \eqref{eq:mainisom}}
\label{fig:mainisomexample}
\end{figure}

This proof is a trivial verification because we have set up the correct machinery. Notice that we could equally well use any odd power of $K$ to define this isomorphism; we have chosen $K^{-1}$ simply for notational convenience. Since the map \eqref{eq:mainisom} is notationally dense it is more illuminating to see an example. Figure \ref{fig:mainisomexample} shows how to compute each step of the isomorphism, sending the $3$-multichain $(1,(12)(34),(1234))$ in $NC(A_3)$ (with respect to the Coxeter element $(1234)$) to the $3$-divisible partition
\begin{equation*}
\{\{1,5,12\},\{2,3,4\},\{6,7,11\},\{8,9,10\}\}.
\end{equation*}

\subsection{Combinatorial Properties}

Before we apply Theorem \ref{th:mainisom} to the classical types in the next few sections, it is worth examining some combinatorial properties of the map \eqref{eq:mainisom}. In particular, notice that the final three steps in Figure \ref{fig:mainisomexample} involve only classical noncrossing partitions. This map has independent interest. Given a $k$-multichain of noncrossing partitions $(\P)_k\in NC^k(n)$, we set
\begin{equation}
\nabla(\P)_k:= K^{-1}\big\langle\partial(\P)_k\big\rangle.
\end{equation}
The next result follows immediately from Theorem \ref{th:mainisom}.

\begin{corollary}
\label{cor:kmchaintokdiv}
The map $(\P_1,\P_2,\ldots,\P_k)\mapsto \nabla(\P_1,\P_2,\ldots,\P_k)$ is a bijection from $k$-multichains in $NC(n)$ to $k$-divisible noncrossing partitions of $[kn]$.
\end{corollary}

Even though it is well-known that both of these objects are counted by the Fuss-Catalan number $\Cat^{(k)}(n)=\frac{1}{n}\binom{(k+1)n}{n-1}$, an explicit bijection has not appeared before in the literature. Combining the work of Dershowitz~\cite{dershowitz} and Edelman~\cite{edelman:trees}, one obtains a bijection implicitly, but this involves passing to ordered trees. The bijection $\nabla$ is direct and involves only the relative Kreweras complement. (In the case $k=2$, this bijection also appears in Nica and Speicher~\cite[Proposition 2.6]{nica-speicher:ntuples}.) Furthermore, $\nabla$ preserves some important information about the multichain.

\begin{corollary}
For all multichains $(\P)_k\in NC^k(n)$, the rank of $\P_1$ in $NC(n)$ is equal to the rank of $\nabla(\P)_k$ in $NC^{(k)}(n)$.
\end{corollary}

This follows immediately from Theorem \ref{th:mainisom}, since the rank of a multichain in $NC^{(k)}(A_{n-1})$ is $\rk(\pi)_k=\ell_T(\pi_1)$ (see Theorem \ref{th:semilattice}), but it is also easy to show directly. Recall that the rank function of $NC^{(k)}(n)$ is given by $n$ minus the number of blocks, $\rk(\P)=n-\abs{\P}$. Using property \eqref{eq:relKblocks} of the relative Kreweras complement, we see that the $k$-divisible partition
\begin{equation*}
\nabla(\P)_k=K^{-1} \Big\langle\big( K^{\P_2}(\P_1),\ldots,K^{\P_k}(\P_{k-1}),K(\P_k)\big)\Big\rangle,
\end{equation*}
where we set $\P_{k+1}=\hat{1}_n$, contains
\begin{equation*}
kn-\sum_{i=1}^k \left(n-\abs{\P_{i+1}}+\abs{\P_i}\right) +1 = \abs{\P_1}
\end{equation*}
blocks. In particular, when the bottom element contains $n$ blocks (that is, when $\P_1=\hat{0}_n$) then the $k$-divisible partition $\nabla(\P)_k$ also contains $n$ blocks and they must all have cardinality $k$. These partitions are precisely the minimal elements of $NC^{(k)}(n)$ and we will call them {\sf $k$-equal}. Thus $\nabla$ defines a bijection between $(k-1)$-multichains in $NC(n)$ and the $k$-equal noncrossing partitions of $[kn]$.

\begin{corollary}
\label{cor:kequal}
The map $(\P_1,\P_2,\ldots,\P_{k-1})\mapsto \nabla(\hat{0}_n,\P_1,\P_2,\ldots,\P_{k-1})$ 
is a bijection from $(k-1)$-multichains in $NC(n)$ to noncrossing partitions of $[kn]$ in which each block has size $k$ ($k$-equal partitions).
\end{corollary}

Note that this result is the combinatorial version of Corollary \ref{cor:minimal}. In the case $k=2$, this defines a bijection between noncrossing partitions of $[n]$ and noncrossing ``pairings'' of $[2n]$. This operation is well-known, and it is sometimes referred to as ``thickening'' (see \cite{hall}).

The next result is a stronger combinatorial property of $\nabla$ that includes the above facts as special cases. Given a partition $\P$ of the set $[n]$, its multiset of block sizes is called the {\sf type} of the partition. More generally,

\begin{definition}
Given a partition $\P=(P_1,\ldots,P_m)$ of a finite set $X\subseteq\integers$ and a positive integer $k$, define the {\sf $k$-type} of $\P$ as the multiset
\begin{equation*}
\lambda^{(k)}(\P):=\left\{ \frac{\abs{\P_i}}{k}: 1\leq i\leq m\right\}
\end{equation*}
of block sizes of $\P$, each divided by $k$.
\end{definition}

Notice that $\lambda^{(k)}(\P)$ is an {\sf integer partition} of the integer $\abs{X}/k$ precisely when $\P$ is a $k$-divisible partition. In particular, when $\P\in NC^{(k)}(W)$ we write $\lambda^{(k)}(\P)\vdash n$ to denote the fact that $\lambda^{(k)}(\P)$ is an integer partition of $n$. Since order does not matter in a multiset, it is conventional to write the elements in weakly decreasing order, say $\lambda^{(k)}(\P)=\{\lambda_1\geq\lambda_2\geq\cdots\geq\lambda_r\}$, or shortened to
\begin{equation*}
\lambda^{(k)}(\P)=(\lambda_1,\lambda_2,\ldots,\lambda_r).
\end{equation*}
For example, among the $2$-divisible noncrossing partitions of the set $[6]$ (Figure \ref{fig:nc^2(3)}) there are five partitions of $2$-type $(1,1,1)$, six partitions of $2$-type $(2,1)$ and one partition of $2$-type $(3)$. Notice that the $1$-type of a partition $\P$ is simply the usual type, and we write
\begin{equation*}
\lambda(\P):=\lambda^{(1)}(\P)
\end{equation*}
in this case.

The map $\nabla$ preserves type in the following sense.

\begin{theorem}
\label{th:typepreserving}
We have
\begin{equation*}
\lambda^{(k)}\left(\nabla(\P)_k\right)=\lambda(\P_1)
\end{equation*}
for all multichains $(\P)_k\in NC^k(n)$.
\end{theorem}

\begin{proof}
We will actually prove a stronger result. We will show that the restriction of the $k$-divisible partition $\nabla(\P)_k$ to the set
\begin{equation*}
k[n]-(k-1)=\{ 1,k+1,2k+1,\ldots,k(n-1)+1\}
\end{equation*}
is equal to the partition $k\P_1-(k-1)$ (that is, the partition induced by $\P_1$ and the order-preserving bijection $[n]\leftrightarrow (k[n]-(k-1))$). Then since $\nabla(\P)_k$ is $k$-divisible, the result follows. Observing the pictorial representation of the Kreweras complement (Figure \ref{fig:kreweras}), we can see that the restriction of $\nabla(\P)_k$ to the set $k[n]-(k-1)$ is equal to the partition $k\Q-(k-1)$, where
\begin{equation*}
\Q=\wedge_{i=1}^k K^{-1}\left( K^{\P_{i+1}}(\P_i)\right).
\end{equation*}
For example, in Figure \ref{fig:mainisomexample}, consider the delta sequence 
\begin{align*}
&(K^{\P_2}(\P_1),K^{\P_3}(\P_2),K(\P_3))= \\
&\qquad\Big( \big\{\{1,2\},\{3,4\}\}\,\,,\,\,\{\{1\},\{2,4\},\{3\}\}\,\,,\,\,\{\{1\},\{2\},\{3\},\{4\}\big\}\Big)
\end{align*}
and the corresponding shuffle partition
\begin{equation}
\label{eq:shuffleof12}
\S=\big\{\{1,4\},\{2\},\{3\},\{5,11\},\{6\},\{7,10\},\{8\},\{9\},\{12\}\big\}.
\end{equation}
Notice that the restriction of $K^{-1}(\S)$ to the set $\{1,4,7,10\}$ is the partition 
\begin{equation*}
\{\{1\},\{4\},\{7\},\{10\}\}.
\end{equation*}
If we observe how $K^{-1}(\S)$ is computed, it is determined completely by restrictions imposed by the elements of the delta sequence. In particular, it is not difficult to see that the connections among the set $\{1,4,7,10\}$ in $K^{-1}(\S)$ are given precisely by the intersection
\begin{equation*}
K^{-1}(K^{\P_2}(\P_1))\wedge K^{-1}(K^{\P_3}(\P_2))\wedge K^{-1}(K(\P_3)),
\end{equation*}
as claimed.

Now given $(\P)_k\in NC^{(k)}(n)$, suppose that $(\pi)_k$ is the multichain in $NC(A_{n-1})$ satisfying $(\P)_k=(\{\pi\})_k$. Then Lemma \ref{lemma:joins} shows that
\begin{eqnarray*}
\pi_1^{-1}c&=&\textstyle{\prod}_{i=1}^k \pi_i^{-1}\pi_{i+1}\\
K_1^c(\pi_1) &=&\textstyle{\prod}_{i=1}^k K^{\pi_{i+1}}(\pi_i)\\
K_1^c(\pi_1) &=&\vee_{i=1}^k K_1^{\pi_{i+1}}(\pi_i)\\
\pi_1 &=& (K_1^c)^{-1}( \vee_{i=1}^k K_1^{\pi_{i+1}}(\pi_i))\\
\pi_1 &=&\wedge_{i=1}^k (K_1^c)^{-1}(K_1^{\pi_{i+1}}(\pi_i)).
\end{eqnarray*}
Applying the reverse map $(\pi)_k\mapsto (\P)_k$ to the last equation, we conclude that $\P_1=\Q$ as desired.
\end{proof}

So the isomorphism \eqref{eq:mainisom} preserves not only the rank, but also the structure of the bottom element of the multichain. Following this, it makes sense to define the ``type'' $\lambda(\P)_k$ of a multichain $(\P)_k$ to be equal to the type $\lambda(\P_1)$ of its bottom element $\P_1$. It is no accident that we also use the word ``type'' to refer to the isomorphism class of a finite Coxeter system $(W,S)$. In Chapter \ref{sec:fusscatcomb}, we will use this idea to define the type of an element in $NC^{(k)}(W)$.

Finally, observe how the map $\nabla$ acts on $\ell$-multichains of $k$-divisible noncrossing partitions. If $(\P_1,\P_2,\ldots,\P_\ell)\in (NC^{(k)}(n))^\ell$ is a multichain of $k$-divisible partitions, then Corollary \ref{cor:kmchaintokdiv} says that $\nabla(\P_1,\P_2,\ldots,\P_\ell)$ is an $\ell$-divisible noncrossing partition of $[k\ell n]$; but more is true. This partition is actually $k\ell$-divisible.

\begin{corollary}
\label{cor:lmchainofkdiv}
The map $\nabla$ is a bijection from $\ell$-multichains of $k$-divisible noncrossing partitions of $[kn]$ to $k\ell$-divisible noncrossing partitions of $[k\ell n ]$.
\end{corollary}

\begin{proof}
Consider an $\ell$-multichain $(\P_1,\P_2,\ldots,\P_\ell)$ in $(NC^{(k)}(n))^\ell$. By Theorem \ref{th:typepreserving} we have $\lambda^{(\ell)}\left(\nabla(\P)_\ell\right)=\lambda(\P_1)$. That is, the cardinality of each block of $\nabla(\P)_\ell$, when divided by $\ell$, is still divisible by $k$.
\end{proof}

This Corollary can be thought of as the classical version of Theorem \ref{th:nckl=nclk}.

\subsection{Automorphisms}
\label{sec:classicalautomorphisms}

Finally, we describe a classical interpretation of the dihedral group of automorphisms defined in Section \ref{sec:automorphisms}. To do this, we need to interpret the $n$-cycle $c=(12\cdots n)$ as a ``bipartite Coxeter element''. This cannot be done using the generating set $S$ of adjacent transpositions, so we pass to a different Coxeter generating set.

\begin{definition}
Define the {\sf snake generating set} $\S$ of $A_{n-1}$ as the disjoint union of the two sets of transpositions
\begin{align*}
S_\ell & :=\left\{ (i,n-i+2): 2\leq i\leq \lceil n/2\rceil\right\}, \text{ and}\\
S_r &:=\left\{ (i,n-i+1): 1\leq i\leq \lfloor n/2\rfloor\right\},
\end{align*}
and set $\ell=\prod_{t\in S_\ell} t$ and $r=\prod_{t\in S_r} t$.
\end{definition}

\begin{figure}
\vspace{.1in}
\begin{center}
\input{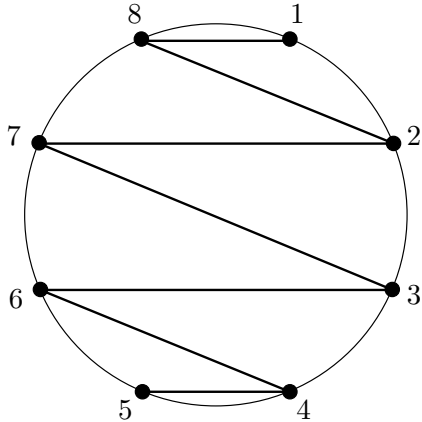}
\end{center}
\caption{The snake generating set for $A_7$}
\label{fig:snake}
\end{figure}

For instance, the snake generating set for $A_7$ consists of $S_\ell=\{(28),(37),(46)\}$ and $S_r=\{(18),(27),(36),(45)\}$. Figure \ref{fig:snake} displays this generating set, with horizontal lines representing transpositions in $S_r$ and diagonal lines representing transpositions in $S_\ell$. The union $S_\ell\sqcup S_r$ is evidently a Coxeter generating set for $A_7$ since the corresponding Coxeter diagram is a chain. Notice also that the corresponding bipartite Coxeter element is
\begin{equation*}
\ell\cdot r= (28)(37)(46)\cdot(18)(27)(36)(45) = (12345678),
\end{equation*}
as desired.

Now consider the noncrossing partition $\pi=(156)(23)\in NC(A_7)$ as displayed in Figure \ref{fig:kreweras}. Conjugating by the element $r=(18)(27)(36)(45)$ gives $r\pi r=(384)(67)$, which is {\em not} in $NC(A_7)$ because the $3$-cycle $(384)$ is oriented counterclockwise. However, if we invert the element $r\pi r$, we get $r\pi^{-1} r=(348)(67)$, which {\em is} in $NC(A_7)$ since its orbits on the set $[8]$ are noncrossing and all of its cycles are oriented clockwise. That is, the map $R:NC(A_7)\to NC(A_7)$ defined by $R(\pi)=r\pi^{-1}r$ gives a reflection of the circular representation across the vertical bisector, the line through $4'$ and $8'$ in Figure \ref{fig:kreweras}. Similarly, the map $L(\pi)=\ell\pi^{-1}\ell$ gives a reflection of the circular representation across the line through vertices $1$ and $5$ in the figure. In general, we have the following characterization.

\begin{lemma}
The automorphisms $L$ and $R$ on $NC(A_{n-1})$, with respect to the snake generating set, generate the dihedral group of motions on the circular representation.
\end{lemma}

In fact, this property motivates the definition of the maps $L$ and $R$ in the general case (equation \eqref{eq:LandR}). Furthermore, the main isomorphism \eqref{eq:mainisom} gives a similar pictorial characterization of the automorphisms $R^*$ \eqref{eq:R*} and $L^*$ \eqref{eq:L*} acting on $NC_{(k)}(A_{n-1})$.

\begin{lemma}
If we transfer the maps $L^*$ and $R^*$, with respect to the snake generating set, from $NC_{(k)}(A_{n-1})$ to $NC^{(k)}(n)$ via the composition of $\sint$ and \eqref{eq:mainisom}, then these generate the dihedral group of motions on the circular representation of the $k$-divisible noncrossing partitions.
\end{lemma}

Since the notation here is so dense, it is more helpful to describe the main idea of the proof. Based on Theorem \ref{th:shuffle}, we may think of $NC_{(k)}(A_{n-1})$ as the poset of $k$-shuffle noncrossing partitions of $[kn]$. Now consider what happens if we reflect a shuffle partition across the vertical bisector: the partition on the congruence class $k[n]+i$ switches with the partition on the congruence class $k[n]+(k-i+1)$, and both of these partitions are ``reflected'' from left to right by the map $R$. For example, if we consider the $3$-shuffle partition \eqref{eq:shuffleof12} of the set $[12]$ in the bottom left corner of Figure \ref{fig:mainisomexample}, this reflection corresponds to:

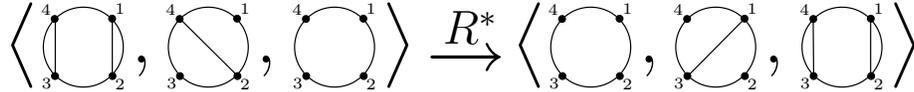
\begin{equation*}
\input{shuffleR.pstex_t}
\end{equation*}
This explains the form of $R^*$. To understand the map $L^*$, consider the reflection of a $k$-shuffle partition across the diameter through the symbol $1$. In this case, the partition on congruence class $k[n]+1$ remains in that class and is ``reflected'' by $L$, whereas the partitions in congruence classes $k[n]+i$ and $k[n]+(k-i+2)$ are switched for all $2\leq i\leq k$, and each is ``reflected'' by $R$. For the shuffle partition in Figure \ref{fig:mainisomexample}, this corresponds to reflection in the line through vertices $1$ and $7$:

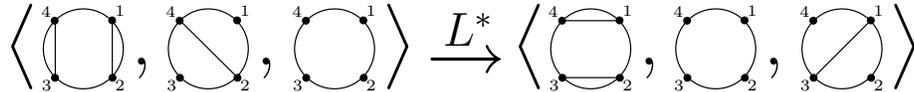
\begin{equation*}
\input{shuffleL.pstex_t}
\end{equation*}
Finally, it is perhaps easiest to see that the composition $L^*\circ R^*$ \eqref{eq:L*R*} has the correct form. This map shifts each partition into the next congruence class, clockwise, and rotates the indices on the final partition, since its ``starting'' index has changed. Clearly this corresponds to $1/kn$ of a full clockwise rotation of the circular representation.

The precise details here are not important. We only wish to convey the fact that the algebraically defined automorphisms $L^*$ and $R^*$ on $NC_{(k)}(W)$ from Section \ref{sec:automorphisms} are geometrically motivated.

\section{Type $A$}
\label{sec:typeA}

Theorem \ref{th:mainisom} reduces the study of the algebraically defined $NC^{(k)}(A_{n-1})$ to the study of the classical $k$-divisible noncrossing partitions $NC^{(k)}(n)$. This connection is a two-way street.

In one direction, all of the results we have established in Chapter \ref{sec:kdiv} now have classical analogues. Some of the enumerative results were known before (see~\cite{edelman:kdiv}), but many of the structural results are new. In particular, the results of Section \ref{sec:meta} were not observed before in the classical case. Transferring Lemma \ref{lemma:contained}, we get specific embeddings of $NC^{(k)}(n)$ into $NC^{(\ell)}(n)$ for all $k\leq\ell$. Theorem \ref{th:orderideal} can also be transferred, but this result is perhaps not surprising, since given $\P=(P_1,\ldots,P_m)\in NC^{(k)}(W)$, it is intuitively clear that the principal order ideal is isomorphic to a direct product
\begin{equation*}
\Lambda(\P)\cong \prod_{i=1}^m NC^{(k)}(P_i)
\end{equation*}
over the blocks of $\P$. The results about shellability and Euler characteristics are also new in the classical setting.

In the other direction, the relationship between $NC^{(k)}(A_{n-1})$ and $NC^{(k)}(n)$ makes it much easier to compute enumerative properties of $NC^{(k)}(A_{n-1})$. For example, we recall the following result of Edelman \cite[Theorem 4.2]{edelman:kdiv}. Given an $\ell$-multichain $\P_1\leq\P_2\leq\cdots\leq\P_\ell$ in $NC^{(k)}(n)$, define its {\sf rank-jump vector} by
\begin{equation*}
(j_1,j_2,\ldots,j_{\ell+1}):=(r_1-0,r_2-r_1,\ldots,r_\ell-r_{\ell+1},(n-1)-r_\ell),
\end{equation*}
where we set $r_i=\rk(\P_i)=n-\abs{\P_i}$ for all $1\leq i\leq \ell$. Note, in particular, that the entries of the rank-jump vector sum to $n-1$, the height of the poset. Edelman proved the following fact about $NC^{(k)}(n)$, which now also applies to $NC^{(k)}(A_{n-1})$.

\begin{theorem}[\cite{edelman:kdiv}]
\label{th:edelman}
The number of $\ell$-multichains in $NC^{(k)}(A_{n-1})$ with rank-jump vector $(j_1,j_2,\ldots,j_{\ell+1})$ is equal to
\begin{equation*}
\frac{1}{n}\binom{n}{j_1}\binom{kn}{j_2}\binom{kn}{j_3}\cdots\binom{kn}{j_{\ell+1}}.
\end{equation*}
\end{theorem}

This remarkable formula contains a lot of information. For instance, summing over all rank-jump vectors we obtain the zeta polynomial,
\begin{equation}
\label{eq:edelman}
\Z(NC^{(k)}(A_{n-1}),\ell)=\frac{1}{n}\binom{(k\ell+1)n}{n-1},
\end{equation}
and evaluating at
\begin{equation*}
(j_1,j_2,\ldots,j_{n+1})=(0,1,1,\ldots,1,1,0)
\end{equation*}
tells us that the number of maximal chains in $NC^{(k)}(A_{n-1})$ is equal to $(kn)^{n-1}/n$, both of which we had computed before by different means (Theorem \ref{th:zetapoly}, Corollaries \ref{cor:maxchains} and \ref{cor:lmchainofkdiv}). However, Theorem \ref{th:edelman} also gives us the type $A$ Fuss-Narayana numbers, which we were unable to compute before.

\begin{theorem}
\label{th:fussnarA}
The type $A$ Fuss-Narayana polynomial is
\begin{equation*}
\Nar^{(k)}(A_{n-1},i)=\frac{1}{n}\binom{n}{i}\binom{kn}{n-1-i}.
\end{equation*}
\end{theorem}

\begin{proof}
Evalute formula \eqref{eq:edelman} at $(j_1,j_2)=(i,n-1-i)$.
\end{proof}

Another classical result due to Kreweras counts the noncrossing partitions by type. Given an integer partition $\lambda$, let $m_i$ be the number of entries of $\lambda$ equal to $i$, and define
\begin{equation*}
m_\lambda:= m_1!m_2!m_3!\cdots.
\end{equation*}
Kreweras proved the following \cite[Theorem 4]{kreweras}. This is the ``stronger result'' that we mentioned in the proof of Theorem \ref{th:kreweras}. As is common, we will denote the number of parts of $\lambda$ by $l(\lambda)$.

\begin{theorem}[\cite{kreweras}]
\label{th:krewerastype}
The number of noncrossing partitions of type $\lambda\vdash n$ with $l(\lambda)=i$ is equal to
\begin{equation*}
\frac{n!}{m_\lambda(n-i+1)!}.
\end{equation*}
\end{theorem}

For example, the number of noncrossing partitions with type $\lambda=(2,1,1)$ is 
\begin{equation*}
\frac{4!}{1!2!(4-3+1)!}=6,
\end{equation*}
as can be observed in Figure \ref{fig:nc4}. Applying the results of the last section, we get a stronger version of Kreweras' theorem.

\begin{theorem}
\label{th:kltype}
Consider positive integers $k$, $\ell$, $n$ and integer partition $\lambda\vdash n$ with number of parts $l(\lambda)=i$. The number of $\ell$-multichains $\P_1\leq\cdots\leq\P_\ell$ in $NC^{(k)}(n)$ whose bottom element has $k$-type $\lambda^{(k)}(\P_1)=\lambda$ is equal to
\begin{equation}
\label{eq:kltype}
\frac{(k\ell n )!}{m_\lambda(k\ell n-i+1)!}.
\end{equation}
\end{theorem}

\begin{proof}
First note that we have $m_{\lambda^{(k)}(\P_1)}=m_\lambda$ for all $k$, since $m_\lambda$ encodes only multiplicities. Then, by Theorem \ref{th:krewerastype}, formula \eqref{eq:kltype} counts the number of noncrossing partitions of $[k\ell n]$ with $k\ell$-type $\lambda$. Apply Theorem \ref{th:typepreserving} and Corollary \ref{cor:lmchainofkdiv}.
\end{proof}

Thus, for example, the number of $2$-multichains in $NC^{(2)}(3)$ whose bottom element has $2$-type $(1,1,1)$ is equal to
\begin{equation*}
\frac{(12)!}{3!(12-3+1)!}=22.
\end{equation*}
(see Figure \ref{fig:nc^2(3)}). Notice also that by setting $\ell=1$ and summing formula \eqref{eq:kltype} over all $\lambda\vdash n$ with $l(\lambda)=i$ gives another verification of the type $A$ Fuss-Narayana number.

Theorems \ref{th:edelman} and \ref{th:kltype} are similar in that they both count $\ell$-multichains of $k$-divisible noncrossing partitions with respect to different statistics. It is natural to ask whether there exists a closed formula that simultaneously generalizes formulas \eqref{eq:edelman} and \eqref{eq:kltype}. Previously, we have identified the ``type'' of an $\ell$-multichain $\P_1\leq\cdots\leq\P_\ell$ in $NC^{(k)}(n)$ as the $k$-type of its bottom element $\lambda^{(k)}(\P_1)$, but this is really only a partial accounting of the type of the whole multichain. It would be more precise to specify the $k$-type of each element of the multichain.

\begin{definition}
\label{def:totaltype}
Given an $\ell$-multichain $\P_1\leq\cdots\leq\P_\ell$ in $NC^{(k)}(n)$, define the {\sf total $k$-type} of this multichain as
\begin{equation*}
\Lambda^{(k)}(\P_1,\P_2,\ldots,\P_\ell):=\left(\lambda^{(k)}(\P_1),\ldots,\lambda^{(k)}(\P_\ell)\right).
\end{equation*}
\end{definition}

Notice that the total $k$-type $\Lambda^{(k)}(\P)_\ell$ generalizes both the usual type and the rank-jump vector of the multichain. Thus, an answer to the following question will simultaneously generalize \eqref{eq:edelman} and \eqref{eq:kltype}.

\begin{problem}
\label{prob:totaltype}
Does there exist a closed formula counting multichains in $NC^{(k)}(n)$ by total type?\footnote{Krattenthaler and M\"uller have answered this question with several remarkable formulas \cite{kratt-muller}.}
\end{problem}

It is interesting to note that the total $k$-type of an $\ell$-multichain in $NC^{(k)}(W)$ is an $\ell$-multichain in the {\sf dominance order} on integer partitions of $n$. Given two partitions $\lambda=(\lambda_1,\ldots,\lambda_m)\vdash n$ and $\mu=(\mu_1,\ldots,\mu_t)\vdash n$, we say that {\sf $\mu$ dominates $\lambda$}, and we write $\lambda\unlhd\mu$, if
\begin{equation*}
\lambda_1+\lambda_2+\cdots+\lambda_i\leq \mu_1+\mu_2+\cdots+\mu_i
\end{equation*}
for all $i$, where we take $\lambda_i=0$ for $i>m$ and $\mu_i=0$ for $i>t$. (This is an important concept in the representation theory of symmetric groups, see \cite[Chapter 2.2]{sagan}.) It is easy to see that $\P\leq\Q$ in $NC^{(k)}(n)$ implies $\lambda^{(k)}(\P)\unlhd\lambda^{(k)}(\Q)$.

\section{Type $B$}
\label{sec:typeB}

In \cite{reiner}, Reiner introduced the notion of a classical type $B$ noncrossing partition and he generalized many of the known properties of type $A$ partitions. Here, we generalize Reiner's construction to $k$-divisible partitions of type $B$.

If $\{e_1,e_2,\ldots,e_n\}$ is the standard basis of $\reals^n$, then the type $B_n$ crystallographic root system has $n^2$ positive roots, consisting of
\begin{equation*}
\{e_i\}_{1\leq i\leq n}\quad\text{and}\quad\{e_i\pm e_j\}_{1\leq i<j\leq n},
\end{equation*}
and the type $C_n$ crystallographic root system has positive roots
\begin{equation*}
\{2e_i\}_{1\leq i\leq n}\quad\text{and}\quad \{e_i\pm e_j\}_{1\leq i<j\leq n}.
\end{equation*}
When $n\geq 3$, these root systems are not congruent, but their corresponding Weyl groups (generated by reflections orthogonal to the roots) are clearly the same. Hence, we will only speak of the finite Coxeter group $B_n$.

As well as being the group of symmetries of the hypercube and hyperoctahedron in $n$-dimensions, it is sometimes convenient to think of $B_n$ as the group of {\sf signed permutations}. If we set
\begin{equation*}
[\pm n]:=\{1,2,\ldots,n,-1,-2,\ldots,-n\},
\end{equation*}
then a signed permutation of the set $[\pm n]$ is one that commutes with the {\sf antipodal permutation} $(1,-1)(2,-2)\cdots (n,-n)$. (When discussing signed permutations, we will use commas in the cycle notation.)

\begin{definition}
\label{def:typeB}
$B_n$ is the subgroup of permutations of $[\pm n]$ that centralizes the antipodal map $(1,-1)(2,-2)\cdots (n,-n)$.
\end{definition}

It is straightforward to verify that the Weyl group $B_n$ is isomorphic to the group of signed permutations by sending reflections to permutations:
\begin{equation*}
\begin{array}{rcll}
t_{e_i} &\mapsto& (i,-i)&\quad\text{ for all }1\leq i\leq n,\\
t_{e_i-e_j} &\mapsto& (i,j)(-i,-j)&\quad\text{ for all }1\leq i<j\leq n,\\
t_{e_i+e_j} &\mapsto& (i,-j)(j,-i)&\quad\text{ for all }1\leq i<j\leq n,
\end{array}
\end{equation*}
where $t_\alpha$ denotes the reflection in the hyperplane $\alpha^\perp$, for $\alpha\in\reals^n$. Thus, $B_n$ is isomorphic to a subgroup of the symmetric group $A_{2n-1}$ of permutations of the set $[2n]$. Henceforth, we will fix the specific inclusion
\begin{equation*}
B_n\hookrightarrow A_{2n-1}
\end{equation*}
by identifying the sets $[\pm n]$ and $[2n]$ in the obvious way: sending $i\mapsto i$ for $i\in\{1,2,\ldots,n\}$ and $i\mapsto n-i$ for $i\in\{-1,-2,\ldots,-n\}$. Not surprisingly, 
we can describe the lattice of type $B$ noncrossing partitions as a sublattice of $NC(A_{2n-1})\cong NC(2n)$.

First notice that the square of the Kreweras complement $K^2:NC(2n)\to NC(2n)$ is equivalent to $1/2n$ of a full counterclockwise rotation of the circular representation (Figure \ref{fig:kreweras}). This corresponds algebraically to conjugation by the Coxeter element $(K_1^c)^2(\pi)=c^{-1}\pi c$. Hence $K^{2n}:NC(2n)\to NC(2n)$ is the antipodal map, or one half of a full rotation.

\begin{definition}
For all positive integers $n$, let $\widetilde{NC}(2n)$ denote the sublattice of $NC(2n)$ fixed by the antipodal map $K^{2n}:NC(2n)\to NC(2n)$. This is called the lattice of {\sf (classical) type $B$ noncrossing partitions}.
\end{definition}

\begin{figure}
\vspace{.1in}
\begin{center}
\input{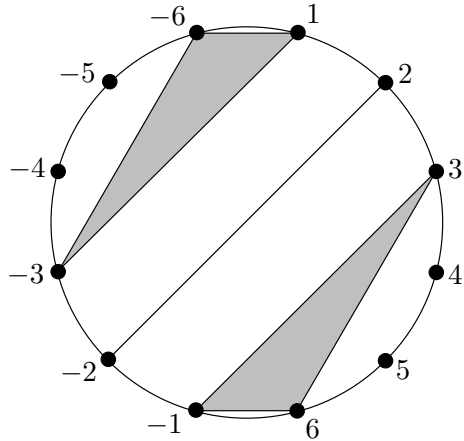}
\end{center}
\caption{An example of a type $B$ noncrossing partition}
\label{fig:typeB}
\end{figure}

In the spirit of signed permutations, we will represent a type $B$ noncrossing partition by labelling the vertices of the $2n$-gon by the integers
\begin{equation*}
1,2,\ldots,n,-1,-2,\ldots,-n,
\end{equation*}
instead of the usual $1,2,\ldots,2n$. For example, Figure \ref{fig:typeB} displays the classical type $B$ noncrossing partition
\begin{equation}
\label{eq:Bpartition}
\{\{1,-3,-6\},\{2,-2\},\{3,6,-1\},\{4\},\{5\},\{-4\},\{-5\}\}.
\end{equation}
Since type $B$ partitions are invariant under the antipodal map, another name sometimes used is {\sf centrally symmetric noncrossing partitions}.

Now, it is not difficult to show that $\widetilde{NC}(2n)$ is isomorphic to the lattice $NC(B_n)$. After Reiner defined these partitions in \cite{reiner}, there were three independent proofs of the following fact, given by Biane, Goodman and Nica \cite{biane:typeB},  Bessis \cite{bessis:dual} and Brady and Watt \cite{brady-watt:kpione}. Considering the action of $B_n$ on $[\pm n]$ by signed permutations, let
\begin{equation*}
\pi\mapsto\{\pi\}
\end{equation*}
denote the map that sends a permutation $\pi$ to the partition of $[\pm n]$ by orbits of $\pi$.

\begin{theorem}[\cite{biane:typeB,bessis:dual,brady-watt:kpione}]
The map $\pi\mapsto\{\pi\}$ is a poset isomorphism from $NC(B_n)$ to $\widetilde{NC}(2n)$.
\end{theorem}

Here, again, we have used the identification $[\pm n]\leftrightarrow [2n]$. The above theorem also says that $NC(B_n)$ with respect to the Coxeter element
\begin{equation*}
(1,2,\ldots,n,-1,-2,\ldots,-n)
\end{equation*}
is isomorphic to a sublattice of $NC(A_{2n-1})$ with respect to the Coxeter element $c=(12\cdots(2n))$: $NC(B_n)$ is precisely the sublattice fixed by conjugation by $c^n$.

\begin{figure}
\vspace{.1in}
\begin{center}
\input{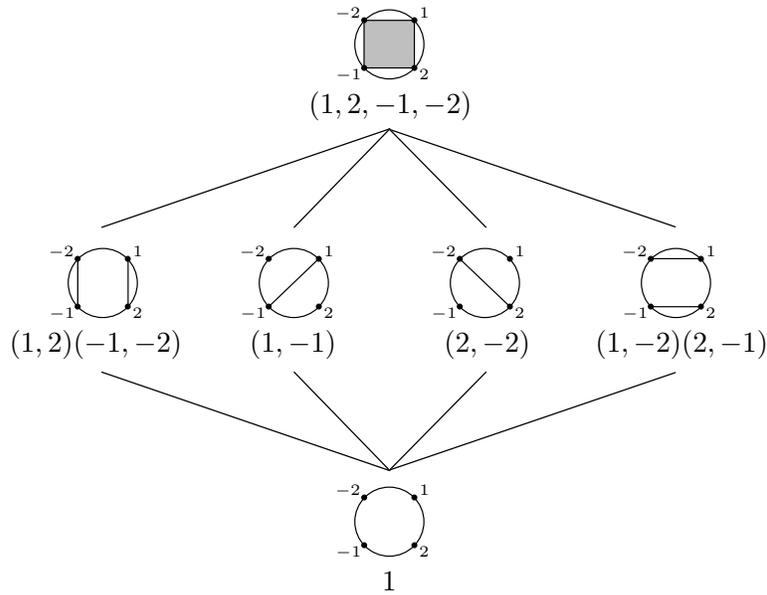}
\end{center}
\caption{$\widetilde{NC}(4)$ is isomorphic to $NC(B_2)$}
\label{fig:ncb4}
\end{figure}

For example, Figure \ref{fig:ncb4} displays the Hasse diagram of $\widetilde{NC}(4)$ and illustrates the isomorphism with $NC(B_2)$. Notice that $NC(B_2)$ is a sublattice of $NC(A_3)$ (Figure \ref{fig:nc4}), but that the rank function is not preserved. This must be the case since the lattice $NC(B_n)$ has height $n$, whereas $NC(A_{2n-1})$ has height $2n$. This discrepancy between rank functions was explained by Reiner \cite{reiner}: Given a partition $\P\in\widetilde{NC}(2n)$, any block that contains a pair of antipodal points $\{i,-i\}$ is called a {\sf zero block} of $\P$. Clearly, $\P$ can contain at most one zero block and the {\sf nonzero blocks} come in pairs.

\begin{definition}
Given a partition $\P\in\widetilde{NC}(2n)$, let $\nz(\P)$ denote the number of {\em pairs} of nonzero blocks in $P$.
\end{definition}

The rank function of $\widetilde{NC}(2n)$ is given by $n$ minus the number of pairs of nonzero blocks,
\begin{equation*}
\rk(\P)=n-\nz(\P).
\end{equation*}

This clearly differs from the type $A$ rank function. For example, the type $B$ partition \eqref{eq:Bpartition} shown in Figure \ref{fig:typeB} has seven blocks, consisting of one zero block $\{2,-2\}$ and three pairs of nonzero blocks. Hence this partition has rank $6-3=3$ in $\widetilde{NC}(12)$; but it has rank $12-7=5$ as an element of $NC(12)$.

The Kreweras complement, however, {\em is} preserved from type $A$ to type $B$. Since $B_n$ is a subgroup of $A_{2n-1}$, we have $K_\mu^\nu(\pi)=\mu\pi^{-1}\nu\in B_n$ whenever $\mu$, $\pi$ and $\nu$ are in $B_n$. This means that the map $K_\mu^\nu$ restricts to an anti-automorphism of the interval $[\mu,\nu]\subseteq NC(B_n)$ and it can be computed pictorially (Figures \ref{fig:kreweras} and \ref{fig:relKexample}).

Now we introduce the classical $k$-divisible type $B$ noncrossing partitions.

\begin{definition}
\label{def:typeBkdiv}
Let $\widetilde{NC}\hspace{0in}^{(k)}(2n)$ denote the induced subposet of $\widetilde{NC}(2kn)$ consisting of partitions in which each block has size divisible by $k$.
\end{definition}

\begin{figure}
\vspace{.1in}
\begin{center}
\input{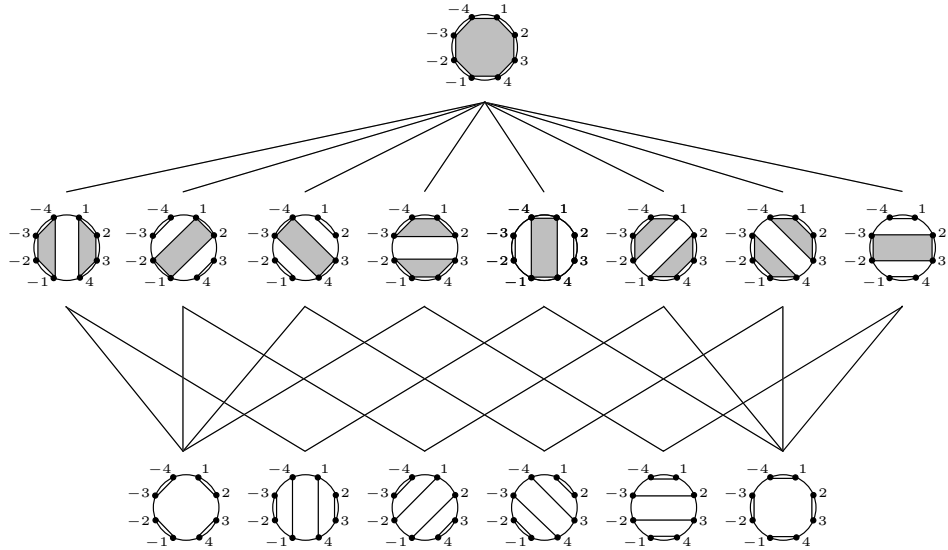}
\end{center}
\caption{The poset $\widetilde{NC}\hspace{0in}^{(2)}(4)$ of $2$-divisible centrally symmetric noncrossing partitions of the set $[\pm 4]$}
\label{fig:nc^2(b2)}
\end{figure}

For example, Figure \ref{fig:nc^2(b2)} displays the Hasse diagram of $\widetilde{NC}\hspace{0in}^{(2)}(4)$ consisting of $2$-divisible centrally symmetric noncrossing partitions of $[\pm 4]$. Notice that this is isomorphic to an order filter in $\widetilde{NC}(8)$. In general, because the Kreweras complement restricts to subgroups, it is easy to see that $NC^{(k)}(B_n)$ is an induced subposet of $NC^{(k)}(A_{2n-1})$. Thus we might hope to have a type $B$ version of Theorem \ref{th:mainisom}. The following is the main result of this section.

\begin{theorem}
\label{th:mainisomB}
Consider $NC(B_n)$ as the sublattice of $NC(A_{2n-1})$ fixed under conjugation by $c^n$. Then the map $(\pi)_k\mapsto K^{-1}\left\langle\partial\left(\{\pi\}\right)_k\right\rangle$ from Theorem \ref{th:mainisom}, when restricted to $NC^{(k)}(B_n)$, is an isomorphism $NC^{(k)}(B_n)\cong \widetilde{NC}\hspace{0in}^{(k)}(2n)$.
\end{theorem}

\begin{proof}
This is a straightforward consequence of Theorem \ref{th:mainisom}. We need only show that this map, as illustrated in Figure \ref{fig:mainisomexample}, preserves the property of central symmetry.

So let $(\pi)_k$ be an element of $NC^{(k)}(B_n)\subseteq NC^{(k)}(A_{n-1})$ and let $(\P)_k=(\{\pi\})_k$ be the corresponding multichain of partitions. Note that the property of central symmetry is preserved under the Kreweras complement, and by taking meets and joins. Applying Lemma \ref{lemma:joins} and Theorem \ref{th:kreweras}, we see that the reciprocal bijections $\partial$ \eqref{partial} and $\sint$ \eqref{sint} have the classical forms
\begin{eqnarray*}
\partial(\P)_k &=& \left( K^{\P_2}(\P_1),\ldots,K^{\P_k}(\P_{k-1}),K(\P_k)\right)\\
 &=& \left(\P_2\wedge K(\P_1),\ldots,\P_k\wedge K(\P_{k-1}),K(\P_k)\right)
\end{eqnarray*}
and
\begin{equation*}
\sint(\Q)_k=\left( \Q_0,\Q_0\vee\Q_1,\ldots,\Q_0\vee\Q_1\vee\cdots\vee\Q_{k-1}\right),
\end{equation*}
where we set $\Q_0=K^{-1}(\Q_1\vee\Q_2\vee\cdots\vee\Q_k)$. It follows that all of the entries of the multichain $(\P)_k$ are centrally symmetric if and only if all of the entries of the delta sequence $\partial(\P)_k$ are centrally symmetric. And from the circular representation (Figure \ref{fig:mainisomexample}), it is clear that the elements of $\partial(\P)_k$ are centrally symmetric if and only if the $k$-divisible partition $K^{-1}\left\langle\partial(\P)_k\right\rangle$ is centrally symmetric.
\end{proof}

In fact, we have chosen the example in Figure \ref{fig:mainisomexample} so that it also represents the type $B$ isomorphism; note that the elements of the input multichain and the output $k$-divisible noncrossing partition are all centrally symmetric. As we did for type $A$, we will now obtain information about the algebraically defined $NC^{(k)}(B_n)$ by studying the classical $k$-divisible type $B$ noncrossing partitions. Many results can be obtained by straightforward generalization of Edelman \cite{edelman:kdiv} and Reiner \cite{reiner}.

Define the rank-jump vector of a multichain $\P_1\leq\P_2\leq\cdots\leq\P_\ell$ in $\widetilde{NC}\hspace{0in}^{(k)}(2n)$ to be
\begin{equation*}
(j_1,j_2,\ldots,j_{\ell+1}):=(r_1-0,r_2-r_1,\ldots,r_\ell-r_{\ell-1},n-r_\ell),
\end{equation*}
where we set $r_i=\rk(\P_i)=n-\nz(\P_i)$ for all $1\leq i\leq \ell$. Again, the entries of the rank-jump vector sum to $n$, the height of the poset $\widetilde{NC}\hspace{0in}^{(k)}(2n)$. We have the following type $B$ analogue of Theorem \ref{th:edelman}.

\begin{theorem}
\label{th:rankjumpB}
The number of $\ell$-multichains in $NC^{(k)}(B_n)$ with rank-jump vector $(j_1,j_2,\ldots,j_{\ell+1})$ is equal to
\begin{equation}
\label{eq:rankjumpB}
\binom{n}{j_1}\binom{kn}{j_2}\binom{kn}{j_3}\cdots\binom{kn}{j_{\ell+1}}.
\end{equation}
\end{theorem}

\begin{proof}
Using a straightforward generalization of Theorem 4.2 in Edelman \cite{edelman:kdiv} and Proposition 7 in Reiner \cite{reiner}, we get a bijection from the set
\begin{equation*}
\left\{ (L,R_1,R_2,\ldots,R_\ell): L\subseteq [n], R_i\subseteq [kn], \sum_{i=1}^\ell \abs{R_i}=\abs{L}\right\},
\end{equation*}
to the set of multichains $\P_1\leq\cdots\leq\P_\ell$ in $\widetilde{NC}\hspace{0in}^{(k)}(2n)$, in which $\abs{R_i}=\rk(\P_{i+1})-\rk(\P_i)$ for all $1\leq i\leq \ell$ (setting $\P_{\ell+1}=\{[\pm kn]\}$) and $\abs{L}=\rk(\P_1)$. The result follows.
\end{proof}

Again, this gives the zeta polynomial
\begin{equation}
\label{eq:zetaB}
\Z(NC^{(k)}(B_n),\ell)=\binom{(k\ell+1)n}{n}
\end{equation}
and tells us that the number of maximal chains in $NC^{(k)}(B_n)$ equals $(kn)^n$, which we knew before. It also gives us the type $B$ Fuss-Narayana numbers.

\begin{theorem}
\label{th:fussnarB}
The type $B$ Fuss-Narayana polynomial is
\begin{equation*}
\Nar^{(k)}(B_n,i)=\binom{n}{i}\binom{kn}{n-i}.
\end{equation*}
\end{theorem}

\begin{proof}
Evaluate formula \eqref{eq:rankjumpB} at $(j_1,j_2)=(i,n-i)$.
\end{proof}

There is also a corresponding notion of ``type'' for type $B$ noncrossing partitions, introduced by Athanasiadis \cite{athanasiadis:classical}. Here is the $k$-divisible generalization.

\begin{definition}
Given a $k$-divisible partition $\P\in\widetilde{NC}\hspace{0in}^{(k)}(2n)$, the {\sf $k$-type} of $\P$, denoted by $\lambda^{(k)}(\P)$, is the integer partition that has one entry equal to $\abs{P}/k$ for each pair $\{P,-P\}$ of nonzero blocks in $\P$. Thus $\lambda^{(k)}(\P)$ is a partition of the integer $(n-m)/2$, where $m$ is the size of the zero block in $\P$.
\end{definition}

Be careful not to confuse this with the earlier (type $A$) notion of type. For example, the partition in Figure \ref{fig:typeB} has type $(3,1,1)\vdash 5$ as an element of $\widetilde{NC}(12)$, but it has type $(3,3,2,1,1,1,1)\vdash 12$ as an element of $NC(12)$. Athanasiadis \cite[Theorem 2.3]{athanasiadis:classical} proved the following $B$-analogue of Theorem \ref{th:krewerastype}.

\begin{theorem}[\cite{athanasiadis:classical}]
\label{th:Btype}
The number of type $B$ noncrossing partitions of type $\lambda$ with number of parts $l(\lambda)=i$ is equal to
\begin{equation*}
\frac{n!}{m_\lambda (n-i)!}.
\end{equation*}
\end{theorem}

Following this, we have a $B$-analogue of Theorem \ref{th:kltype}.

\begin{theorem}
\label{th:klBtype}
Consider positive integers $k$, $\ell$, $n$ and an integer partition $\lambda\vdash n'\leq n$ with number of parts $l(\lambda)=i$. The number of $\ell$-multichains in $\widetilde{NC}\hspace{0in}^{(k)}(2n)$ whose bottom element has $k$-type $\lambda^{(k)}(\P_1)=\lambda$ is equal to
\begin{equation}
\label{eq:klBtype}
\frac{(k\ell n)!}{m_\lambda (k\ell n-i)!}.
\end{equation}
\end{theorem}

\begin{proof}
By Theorem \ref{th:Btype}, this number counts the $k\ell$-divisible type $B$ noncrossing partitions with $k\ell$-type $\lambda$. Now apply Theorem \ref{th:typepreserving} and Corollary \ref{cor:lmchainofkdiv}. (Strictly speaking, we need the stronger result in the proof of Theorem \ref{th:typepreserving} to control information about the zero block.)
\end{proof}

For example, the number of $2$-multichains in $\widetilde{NC}\hspace{0in}^{(2)}(4)$ whose bottom element has $2$-type $(1,1)$ is equal to
\begin{equation*}
\frac{8!}{2!(2-2)!}=28,
\end{equation*}
which can be verified in Figure \ref{fig:nc^2(b2)}. Generalizing Definition \ref{def:totaltype}, we can define a notion of ``total type'' for type $B$ noncrossing partitions, and ask the following question.

\begin{problem}
Find a formula counting $\ell$-multichains in $\widetilde{NC}\hspace{0in}^{(k)}(2n)$ by total type.\footnote{Krattenthaler and M\"uller have now answered this question \cite{kratt-muller}.}
\end{problem}

\begin{figure}
\vspace{.1in}
\begin{center}
\input{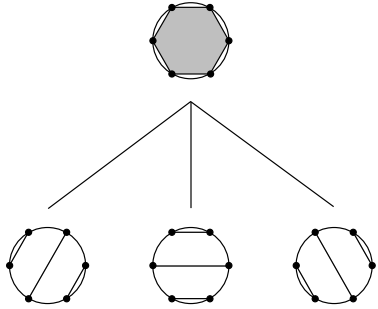}
\end{center}
\caption{The mysterious poset $\widetilde{NC}\hspace{0in}^{(2)}(3)$}
\label{fig:mystery}
\end{figure}

Finally, we suggest an interesting problem for further study. Consider the poset of centrally symmetric $2$-divisible noncrossing partitions of the set $[6]$, as shown in Figure \ref{fig:mystery}. This is an induced subposet of $NC^{(2)}(3)$ (Figure \ref{fig:nc^2(3)}) and we might assume that it falls under our algebraic theory of $k$-divisible noncrossing partitions. However, recall that $NC^{(k)}(B_n)$ is isomorphic to the poset of $k$-divisible centrally symmetric partitions of $[2kn]$. Since $3$ is an {\em odd} integer, the poset in Figure \ref{fig:mystery} does not correspond to anything algebraic in our current framework. This is a mystery.

We expand Definition \ref{def:typeBkdiv} to include this case.

\begin{definition}
For all positive integers $k$ and $n$, let $\widetilde{NC}\hspace{0in}^{(k)}(n)$ denote the subposet of $NC^{(k)}(n)$ that is fixed under the antipodal map $K^{kn}:NC(kn)\to NC(kn)$.
\end{definition}

When $n$ is even, this reduces to Definition  \ref{def:typeBkdiv}, and we have seen (Theorem \ref{th:mainisomB}) that in this case
\begin{equation*}
\widetilde{NC}\hspace{0in}^{(k)}(n)\cong NC^{(k)}(B_{n/2}).
\end{equation*}
If $k$ and $n$ are both odd, then $kn$ is also odd and hence $K^{kn}:NC(kn)\to NC(kn)$ is an {\em anti-}automorphism. In this case, the poset $\widetilde{NC}\hspace{0in}^{(k)}(n)$ is empty.

However, when $k$ is even and $n$ is odd, something strange happens. The product $kn$ is even, so that $K^{kn}:NC(kn)\to NC(kn)$ still represents the antipodal automorphism. In this case $\widetilde{NC}\hspace{0in}^{(k)}(n)$ is an interesting subposet of $NC^{(k)}(n)$ that does not correspond to anything we recognize. What is going on here? The following conjecture generalizing \eqref{eq:zetaB} deepens the mystery.

\begin{conjecture}
\label{conj:mysteryzeta}
Consider integers $k$ and $n$ where $n$ is even and $k$ is arbitrary, or where $n$ is odd and $k$ is even. The zeta polynomial of $\widetilde{NC}\hspace{0in}^{(k)}(n)$ is given by
\begin{equation}
\label{eq:mysteryzeta}
\Z(\widetilde{NC}\hspace{0in}^{(k)}(n),\ell)=\binom{\lfloor (k\ell+1)n/2\rfloor}{\lfloor n/2\rfloor}.
\end{equation}
\end{conjecture}

\begin{problem}
\label{prob:typeBmystery}
Describe the poset $\widetilde{NC}\hspace{0in}^{(k)}(n)$ with $k$ even and $n$ odd. Does it have rank-selection and type-selection formulas similar to \eqref{eq:rankjumpB} and \eqref{eq:klBtype}? Is it related in any way to the algebraic theory of the posets $NC^{(k)}(W)$, where $W$ is a finite Coxeter group? When $k$ and $n$ are both odd, does there exist a nonempty poset to take the place of $\widetilde{NC}\hspace{0in}^{(k)}(n)$ that has zeta polynomial given by \eqref{eq:mysteryzeta}?
\end{problem}

\section{Type $D$}
\label{sec:typeD}

The crystallographic root system of type $D_n$ has $n(n-1)$ positive roots, consisting of
\begin{equation*}
\{ e_i\pm e_j\}_{1\leq i<j\leq n}
\end{equation*}
in $\reals^n$. Hence $D_n$ is a sub-root system of $B_n$ and $C_n$, and the Weyl group $D_n$ is an (index $2$) subgroup of the signed permutation group $B_n$, generated by the elements $(i,j)(-i,-j)$ and $(i,-j)(j,-i)$ for all $1\leq i<j\leq n$. However, $D_n$ is not the group of symmetries of a regular polytope, since its Coxeter diagram is branched (Figure \ref{fig:coxdiagrams}). One can also see from the Coxeter diagrams that the root system $D_3$ is congruent to the root system $A_3$, which has positive roots
\begin{equation*}
\{e_i-e_j\}_{1\leq i<j\leq 4}.
\end{equation*}

The groups $D_n$ are the most idiosyncratic of the classical reflection groups, and when proving a case-by-case result the type $D$ case is often the hardest to establish. The typical progression of a theory is as follows: The type $A$ combinatorics is classical and has been studied for some time. The type $B$ combinatorics is a straightforward conceptual generalization of the type $A$ case, and these together suggest a uniform definition. Finally, the type $D$ combinatorics is tailored to fit the algebraic definition. Rarely is the type $D$ combinatorics intuitive enough to be discovered before its algebraic version.

This was precisely the progression of understanding for the noncrossing partitions. The notion of a type $D$ noncrossing partition was first considered by Reiner \cite{reiner}, where he also introduced the type $B$ noncrossing partitions. His type $D$ lattice was geometrically motivated and it had many nice properties, but it turned out not to agree with the eventual algebraic definition of $NC(D_n)$ (see remarks at the end of Section \ref{sec:classicalNC}). After Brady and Watt demonstrated this discrepancy \cite{brady-watt:kpione}, Athanasiadis and Reiner  \cite{athanasiadis-reiner} modified the type $D$ combinatorics to agree with the algebraic version. However, the ``correct'' type $D$ combinatorics is less intuitive than Reiner's original formulation.

\begin{figure}
\vspace{.1in}
\begin{center}
\input{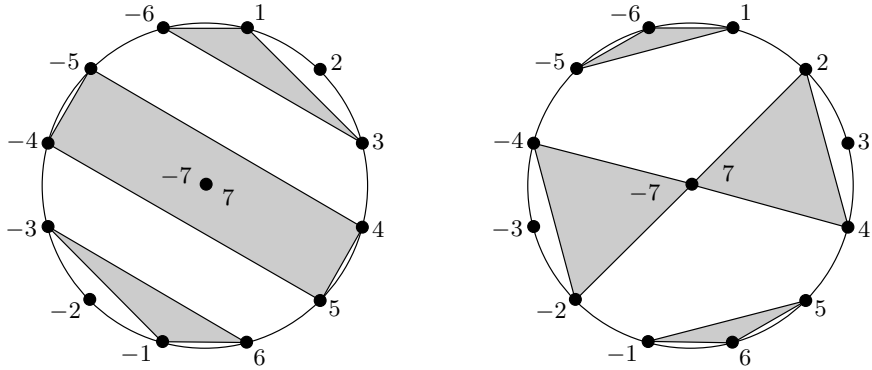}
\end{center}
\caption{Two type $D$ noncrossing partitions}
\label{fig:typeD}
\end{figure}

Athanasiadis and Reiner defined a classical type $D$ noncrossing partition as follows. Label the vertices of a regular $(2n-2)$-gon clockwise with the integers
\begin{equation*}
1,2,\ldots,n-1,-1,-2,\ldots,-(n-1),
\end{equation*}
and label the centroid by both $n$ and $-n$. Given a partition $\P=(P_1,P_2,\ldots,P_m)$ of the set $[\pm n]$ and a block $P\in\P$, identify $P$ with the convex hull $\rho(P)$ of the corresponding vertices. Two blocks $P_i$ and $P_j$ are said to {\sf cross} if one of $\rho(P_i)$ and $\rho(P_j)$ contains a point of the other in its relative interior. Notice that we can have $\rho(P_i)=\rho(P_j)$ if and only if $\{P_i,P_j\}=\{\{n\},\{-n\}\}$. We say $\P$ is a {\sf type $D$ noncrossing partition} if it is centrally symmetric and if $P_i$ and $P_j$ do not cross for all $1\leq i<j\leq m$. Figure \ref{fig:typeD} displays the type $D$ noncrossing partitions
\begin{equation*}
\{\{1,3,-6\},\{2\},\{4,5,7,-4,-5,-7\},\{6,-1,-3\},\{-2\}\}
\end{equation*}
and
\begin{equation*}
\{\{1,-5,-6\},\{2,4,7\},\{3\},\{5,6,-1\},\{-2,-4,-7\},\{-3\}\}
\end{equation*}
of $[\pm 7]$, the first of which has a zero block and the second of which does not. Athanasiadis and Reiner showed that the poset of type $D$ noncrossing partitions of $[\pm n]$ (which is evidently a lattice) is isomorphic to $NC(D_n)$ \cite[Theorem 1.1]{athanasiadis-reiner}.

Thus, we would like to define a poset of classical ``$k$-divisible'' type $D$ noncrossing partitions that generalizes Theorem 1.1 \cite{athanasiadis-reiner} to all positive integers $k$. The obvious idea is to define a $k$-divisible type $D$ partition as a type $D$ partition of $[\pm kn]$ in which each block has size divisible by $k$. However, this definition does not coincide with $NC^{(k)}(D_n)$. There are two difficulties. First, the Kreweras complement does not have an interpretation in this case that is as purely geometric as in the type $A$ case. Second, it is not clear how to label the centroid for higher $k$. Perhaps we should label the vertices of a regular $k(2n-2)$-gon clockwise by the integers
\begin{equation*}
1,2,\ldots,n-1,n+1,n+2,\ldots,2n-1,\ldots,(k-1)n+1,(k-1)n+2,\ldots,kn-1,
\end{equation*}
followed by their negatives, and label the centroid by the set
\begin{equation*}
\{n,2n,\ldots,kn,-n,-2n,\ldots,-kn\}.
\end{equation*}
We do not know the answer.

\begin{problem}
\label{prob:typeD}
Find a combinatorial realization of the poset $NC^{(k)}(D_n)$ generalizing Athanasiadis and Reiner \cite{athanasiadis-reiner}.\footnote{Krattenthaler has now given an elegant solution to this problem using annular noncrossing partitions (see \cite[Section 7]{kratt-muller} and \cite{kratt:typeD}). See the Acknowledgements.}
\end{problem}

However, we are still able to compute the type $D$ Fuss-Narayana numbers, using the following result of Athanasiadis and Reiner \cite[Theorem 1.2 $(ii)$]{athanasiadis-reiner}.

\begin{theorem}[\cite{athanasiadis-reiner}]
The number of multichains $\pi_1\leq\cdots\leq\pi_k$ in $NC(D_n)$ with rank jump vector $(j_1,\ldots,j_{k+1})$ is equal to
\begin{equation}
\label{eq:ath-reiner}
2\binom{n-1}{j_1}\cdots\binom{n-1}{j_{k+1}}+\sum_{m=1}^{k+1} \binom{n-1}{j_1}\cdots\binom{n-2}{j_m-2}\cdots\binom{n-1}{j_{k+1}}.
\end{equation}
\end{theorem}

\begin{theorem}
\label{th:fussnarD}
The type $D$ Fuss-Narayana polynomial is
\begin{equation}
\label{eq:fussnarD}
\Nar^{(k)}(D_n,i)=\binom{n}{i}\binom{k(n-1)}{n-i}+\binom{n-2}{i}\binom{k(n-1)+1}{n-i}.
\end{equation}
\end{theorem}

\begin{proof}
The number $\Nar^{(k)}(D_n,i)$ counts the number of multichains $\pi_1\leq\cdots\leq\pi_k$ in $NC(D_n)$ in which $\pi_1$ has rank $i$. So we should set $j_1=i$ in \eqref{eq:ath-reiner} and then sum over all compositions $(j_2,\ldots,j_{k+1})$ of $n-i$. Doing this for the left summand of \eqref{eq:ath-reiner} yields
\begin{equation*}
2\binom{n-1}{i}\binom{k(n-1)}{n-i}.
\end{equation*}
Now consider the $k+1$ terms in the right summand of \eqref{eq:ath-reiner}. If we perform our summation on each of them, we get
\begin{equation*}
\binom{n-2}{i-2}\binom{k(n-1)}{n-i}
\end{equation*}
when $m=1$ and
\begin{equation}
\label{eq:blim}
\binom{n-1}{i}\binom{k(n-1)-1}{n-i-2}
\end{equation}
when $m\in\{2,3,\ldots,k+1\}$. Multiplication of \eqref{eq:blim} by $k$ yields
\begin{equation*}
\binom{n-2}{i}\binom{k(n-1)}{n-i-1}.
\end{equation*}
Putting all of this together, we conclude that $\Nar^{(k)}(D_n,i)$ is equal to
\begin{equation*}
2\binom{n-1}{i}\binom{k(n-1)}{n-i}+\binom{n-2}{i-2}\binom{k(n-1)}{n-i}+\binom{n-2}{i}\binom{k(n-1)}{n-i-1},
\end{equation*}
which simplifies to \eqref{eq:fussnarD}.
\end{proof}

The lack of a combinatorial realization of the poset $NC^{(k)}(D_n)$ is the biggest hole in the basic theory described in this memoir. It will be important to find such a combinatorial realization so that further results about $NC^{(k)}(W)$ can be proven case-by-case. For example, Krattenthaler has recently studied the $M$-triangle of the poset $NC^{(k)}(W)$ \cite{krattenthaler1,krattenthaler2}, based on an earlier circulated version of this memoir (see Section \ref{sec:triangles}). He was able to prove the theorem in question (related to Conjecture \ref{conj:triangles}) in all cases except type $D$, and here the proof was incomplete since it required a type $D$ version of  Theorems \ref{th:edelman} and \ref{th:rankjumpB}. Hopefully this hole will be filled soon.\footnotemark[\value{footnote}] 


\chapter{Fuss-Catalan Combinatorics}
\label{sec:fusscatcomb}

In this chapter we continue the discussion of Fuss-Catalan numerology begun in Section \ref{sec:fusscat}. It turns out that many of the enumerative formulas described in Section \ref{sec:fusscat} have been observed independently in two other (a priori, unrelated) contexts. The first is the work of Athanasiadis on ``nonnesting partitions'' \cite{athanasiadis:cat,athanasiadis:nar} and the second is the work of Fomin and Reading on ``generalized cluster complexes'' \cite{fomin-reading}.

We will introduce both of these topics and describe some of the exciting coincidences between them. It seems that the three families of objects --- noncrossing partitions, nonnesting partitions and cluster complexes --- are related in deep ways that we do not yet understand. In particular, we make several conjectures below. The existence of mysterious relationships between the three main families of ``Fuss-Catalan objects'' provides a strong motivation to study each of these families individually. We denote the combinatorics surrounding these objects as the {\sf Fuss-Catalan combinatorics of finite Coxeter groups}.

In the final section, we will suggest two directions for future research that are related to the Fuss-Catalan combinatorics but may include a broader range of ideas.

\section{Nonnesting Partitions}
\label{sec:nonnesting}

\subsection{Classical Nonnesting Partitions}
\label{sec:classicalNN}

As before, let $$\P=\{P_1,P_2,\ldots,P_m\}$$ denote a partition of the set $[n]$. We say that two blocks $P_i\neq P_j$ {\sf nest} if there exist $1\leq a<b<c<d\leq n$ with $\{a,d\}\subseteq P_i$ and $\{b,c\}\subseteq P_j$ and there does not exist $b<e<c$ with $e\in P_i$. If $P_i$ and $P_j$ do not nest for all $1\leq i<j\leq m$, we say that $\P$ is a {\sf (classical) nonnesting partition} of $[n]$. Let $NN(n)$ denote the set of nonnesting partitions of $[n]$.

It is important to notice that the property of ``nonnesting'' really depends on the {\em linear} order on $[n]$ and not the {\em cyclic} order. For instance, $\{\{1,2\},\{3,4\}\}$ is a nonnesting partition of the set $[4]$, but if we shift each index by $1$ modulo $4$ we obtain $\{\{1,4\},\{2,3\}\}$, which is {\em nesting}. For this reason, it is not easy to verify the property of ``nonnesting'' in the circular representation, and we need some other way to depict a nonnesting partition.

\begin{figure}
\vspace{.1in}
\begin{center}
\input{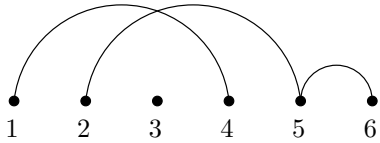}
\end{center}
\caption{The bump diagram of a nonnesting partition}
\label{fig:bumps}
\end{figure}

We define the {\sf bump diagram} of a partition $\P$ of $[n]$ as follows. Place $n$ dots in a row, labelled $1,2,\ldots,n$ from left to right. Then for each block $P=\{i_1,i_2,\ldots,i_r\}$ of $\P$ with $1\leq i_1<i_2<\cdots <i_r\leq n$ we draw an upper semicircular arc joining vertices $i_j$ and $i_{j+1}$ for all $1\leq j<r$. We call each of these semicircular arcs a {\sf bump}. Thus, a partition is nonnesting precisely when its bumps do not ``nest''. For example, Figure \ref{fig:bumps} displays the bump diagram of the nonnesting partition $\{\{1,4\},\{2,5,6\},\{3\}\}$ of the set $[6]$. A nice feature of the bump diagram is that it also detects the property ``noncrossing'': a partition is noncrossing precisely when its bumps do not ``cross''. That is, we can use bump diagrams as a common language for both classical nonnesting and classical noncrossing partitions, although this may obscure some of the structure of the noncrossing partitions.

For more on classical nonnesting partitions, see Athanasiadis \cite{athanasiadis:classical}.

\subsection{Antichains in the Root Poset}
\label{sec:antichains}

What do nonnesting partitions have to do with Coxeter groups? Again, a nonnesting partition of $[n]$ can be defined in terms of the symmetric group $A_{n-1}$ in such a way that the definition generalizes to other reflection groups. In this case, however, we will see that {\em the definition makes sense only for Weyl groups}, since it requires a crystallographic root system. Recall from Section \ref{sec:rootsystems} that a finite spanning set of vectors $\Phi$ in $\reals^n$ is called a {\sf crystallographic root system} if the intersection of the line $\reals\alpha$ with $\Phi$ is equal to $\{\alpha,-\alpha\}$ for all $\alpha\in\Phi$; if the orthogonal reflection $t_\alpha$ in the hyperplane $\alpha^\perp$ satisfies $t_\alpha\Phi=\Phi$ for all $\alpha\in\Phi$; and if the number $2(\alpha,\beta)/(\alpha,\alpha)$ is an integer for all $\alpha,\beta\in\Phi$. Now fix a choice of positive system $\Phi=\Phi^+\sqcup \Phi^-$ and simple system $\Pi$. The defining property of the simple system $\Pi$ is the fact that every root can be uniquely expressed as an integer linear combination of simple roots in which the coefficients are all nonnegative, or all nonpositive. This allows us to define a partial order on the positive roots.

\begin{definition}
\label{def:rootposet}
For all positive roots $\alpha,\beta\in\Phi^+$, we say that $\alpha\leq\beta$ if and only if $\beta-\alpha$ is in the {\em nonnegative} integer span of the simple roots $\Pi$. We denote the set $\Phi^+$ together with this partial order as $(\Phi^+,\leq)$, and call this the {\sf root poset} of the corresponding Weyl group $W$.
\end{definition}

\begin{figure}
\vspace{.1in}
\begin{center}
\input{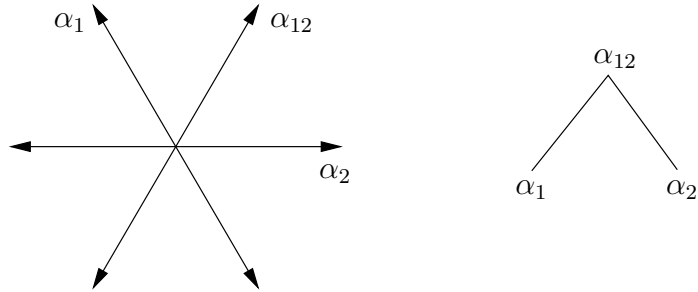}
\end{center}
\caption{The root system and root poset of type $A_2$}
\label{fig:a2roots}
\end{figure}

For example, Figure \ref{fig:a2roots} displays the crystallographic root system of type $A_2$ and the corresponding root poset. Here we have chosen the simple roots $\Pi=\{\alpha_1,\alpha_2\}$ with corresponding positive roots $\Phi^+=\{\alpha_1,\alpha_2,\alpha_{12}\}$, where $\alpha_{12}=\alpha_1+\alpha_2$. Notice, in particular, that the simple roots are the minimal elements of the root poset, and there is a unique maximum, or {\sf highest root}. The root poset plays an important role in the representation theory of Lie algebras and Lie groups.

We use the term ``root poset'' to distinguish from the distinct notion of a ``root order'' (which is used, for instance, in Bj\"orner and Brenti \cite{bjorner-brenti}). We should note that the root posets of type $B_n$ and $C_n$ are isomorphic, even though the root systems themselves are not. Thus, we will not make a distinction between them.

The following definition is due to Postnikov (see \cite[Remark 2]{reiner}). Recall that an {\sf antichain} in a poset is a set of pairwise-incomparable elements.

\begin{definition}
\label{def:postnikov}
Given a Weyl group $W$, an antichain in the root poset $(\Phi^+,\leq)$ is called a {\sf nonnesting partition}. Let $NN(W)$ denote the set of nonnesting partitions of $W$.
\end{definition}

Observe from Figure \ref{fig:a2roots} that the nonnesting partitions of type $A_2$ consist of
\begin{equation*}
NN(W)=\left\{ \emptyset,\{\alpha_1\},\{\alpha_2\},\{\alpha_{12}\},\{\alpha_1,\alpha_2\}\right\}.
\end{equation*}
This set contains five elements, which is, not coincidentally, the Coxeter-Catalan number $\Cat(A_2)=5$.

Unlike the noncrossing partitions, the classical nonnesting partitions do not have a long history, and they were defined {\em simultaneously} with their algebraic generalization (Definition \ref{def:postnikov}). To see how a classical nonnesting partition corresponds to an antichain in the root poset $(\Phi^+,\leq)$ of type $A_{n-1}$, consider the group $A_{n-1}$ with respect to the Coxeter generating set $S$ of adjacent transpositions. If $\{e_1,\ldots,e_n\}$ is the standard basis for $\reals^n$, it is usual to define the positive roots as
\begin{equation*}
\Phi^+(A_{n-1})=\{ e_i-e_j :1\leq i<j\leq n\},
\end{equation*}
even though this root system is not essential (it has a fixed line spanned by $e_1+e_2+\cdots +e_n$). The correspondence with classical nonnesting partitions is straightforward.

\begin{lemma}
The map sending the positive root $e_i-e_j$ to the bump $(i,j)$ is a bijection from $NN(A_{n-1})$ to $NN(n)$.
\end{lemma}

\begin{figure}
\vspace{.1in}
\begin{center}
\input{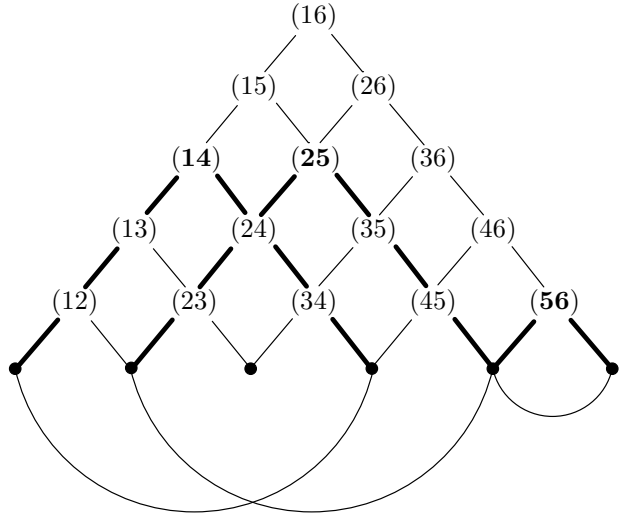}
\end{center}
\caption{The antichain $\{(14),(25),(56)\}$ in the root poset of type $A_5$ corresponds to the nonnesting partition $\{\{1,4\},\{2,5,6\},\{3\}\}$ of the set $[6]$}
\label{fig:bumps-antichains}
\end{figure}

This correspondence is very visual, and we can understand it by observing Figure \ref{fig:bumps-antichains}. Note that the positive roots are in bijection with reflections $T$ since each positive root is the positive-pointing normal to one of the reflecting hyperplanes. In this case, the positive root $e_i-e_j$ corresponds to the transposition $(ij)\in T$. So we think of $(\Phi^+,\leq)$ as a partial order on the set $T$ of transpositions. Observe that this poset has a very regular form, and that the antichain $\{(14),(25),(56)\}$ in the type $A_5$ root poset corresponds to the nonnesting partition $\{\{1,4\},\{2,5,6\},\{3\}\}$ of the set $[6]$ (compare with Figure \ref{fig:bumps}). It is clear from the picture that two bumps will nest precisely when their corresponding roots are comparable in the root order.

Observe also from Figure \ref{fig:bumps-antichains} that there is a bijection from nonnesting partitions of type $A_{n-1}$ to Dyck paths. (A {\sf Dyck path} is a path in the integer lattice $\integers^2$ from the point $(0,0)$ to the point $(n,n)$ using only steps of the form $(1,0)$ and $(0,1)$, and such that every point $(i,j)$ of the path satisfies $i\leq j$.) Indeed, an antichain in the type $A_{n-1}$ root poset corresponds to a Dyck path from $(0,0)$ to $(n,n)$ with a ``peak'' at each element of the antichain. It is well-known that the number of Dyck paths with $i$ peaks is equal to the Narayana number $\Nar(n,i)=\frac{1}{n}\binom{n}{i}\binom{n}{i-1}$. Also, an antichain with $i$ elements corresponds to a nonnesting partition of $[n]$ with $n-i$ blocks, since adding a bump reduces the number of blocks by one. Thus,

\begin{theorem}\hspace{.1in}
\begin{enumerate}
\item The number of nonnesting partitions of $[n]$ is equal to
\begin{equation*}
\abs{NN(n)}=\Cat(A_{n-1})=\frac{1}{n}\binom{2n}{n-1}.
\end{equation*}
\item The  number of nonnesting partitions of $[n]$ with $i$ blocks is equal to
\begin{equation*}
\Nar(A_{n-1},i)=\frac{1}{n}\binom{n}{i}\binom{n}{i-1}.
\end{equation*}
\end{enumerate}
\end{theorem}

Comparing with Theorem \ref{th:kreweras}, we see that the nonnesting partitions and noncrossing partitions are equidistributed by the number of blocks. Athanasiadis proved a strengthening of this \cite[Theorem 2.5]{athanasiadis:classical}. Recall that the {\sf type} of a partition $\P$ of $[n]$ is the integer partition $\lambda\vdash n$ encoding the block sizes of $\P$.

\begin{theorem}[\cite{athanasiadis:classical}]
\label{th:athtypeA}
The number of nonnesting partitions of $[n]$ of type $\lambda\vdash n$ with number of parts $l(\lambda)=i$ is equal to
\begin{equation*}
\frac{n!}{m_\lambda (n-i+1)!}.
\end{equation*}
\end{theorem}

Thus, comparing with Theorem \ref{th:krewerastype},  the classical nonnesting and classical noncrossing partitions are equidistributed by type. This is the beginning of a mystery.

We can naturally think of the set $NN(n)$ as a poset under refinement of partitions. It turns out that $NN(4)$ is isomorphic to $NC(4)$ (see Figure \ref{fig:nc4}), but that $NN(n)$ and $NC(n)$ differ as posets for all $n>4$. In general, the poset $NN(n)$ is {\em not} a lattice, it is {\em not} graded, it is {\em not} self-dual, and indeed we do not know of any nice poset-theoretical property that is satisfies.\footnote{One gets the impression that $NN(n)$ arises from $NC(n)$ by a hammer blow.} Yet it is clear that the nonnesting and noncrossing partitions are closely related somehow. It is worth mentioning that the intuitive map that locally converts each crossing to a nesting
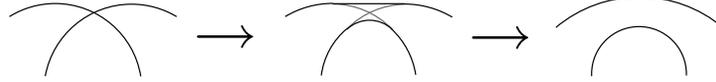
\begin{center}
\input{uncross.pstex_t}
\end{center}
defines a bijection from $NN(n)$ to $NC(n)$ that preserves the number of blocks (see \cite{kasraoui-zeng}). However, this bijection does not preserve type and it has been resistant to generalization.

The type $B$ and type $D$ nonnesting partitions were also studied by Athanasiadis \cite{athanasiadis:classical}, and they were shown to be equidistributed with the type $B$ and $D$ noncrossing partitions according to type \cite{athanasiadis-reiner}. We would like to compare $NN(W)$ and $NC(W)$ for a general Weyl group.

\subsection{Partitions and the Lattice of Parabolic Subgroups}
\label{sec:parabolic}

Let $W$ be a Weyl group. There is a natural order on the set of antichains $NN(W)$ that generalizes the refinement order on classical nonnesting partitions.

\begin{definition}
\label{def:NNrefinement}
Define the {\sf refinement order} on $NN(W)$ by setting
\begin{equation*}
A\leq B\quad\Longleftrightarrow\quad \bigcap_{\alpha\in A} \alpha^\perp \supseteq \bigcap_{\beta\in B} \beta^\perp.
\end{equation*}
\end{definition}

If $\Phi$ is the root system associated with finite Coxeter group $W$, recall that the {\sf Coxeter arrangement} $\A(W)$ is the collection of reflecting hyperplanes
\begin{equation*}
\A(W)=\left\{ \alpha^\perp : \alpha\in\Phi^+\right\},
\end{equation*}
and the {\sf partition lattice} $\Pi(W)$ is the collection of intersections of reflecting hyperplanes
\begin{equation*}
\Pi(W)=\left\{ \cap_{\alpha\in A} \,\alpha^\perp : A\subseteq \Phi^+\right\},
\end{equation*}
ordered by reverse-inclusion of subspaces. The minimum element of this lattice is the whole space $V$, which we understand to be the empty intersection, and the maximum element is the trivial subspace $\{0\}$. Thus, the function $A\mapsto \cap_{\alpha\in A} \alpha^\perp$ is an order preserving map from the refinement order on nonnesting partitions $NN(W)$ to the  partition lattice $\Pi(W)$. Athanasiadis and Reiner \cite[Corollary 6.2]{athanasiadis-reiner} proved that this map is actually a poset embedding.

\begin{theorem}[\cite{athanasiadis-reiner}]
\label{th:NNembed}
If $W$ is a Weyl group, then the map $A\mapsto \cap_{\alpha\in A}\alpha^\perp$ defines a poset embedding $NN(W)\hookrightarrow \Pi(W)$.
\end{theorem}

Also, since the elements of an antichain are linearly independent (see \cite{sommers}), this embedding sends an antichain of cardinality $i$ to a subspace of codimension $i$. In Proposition \ref{prop:Movembedding}, we saw that the lattice $NC(W)$ embeds into the lattice of subspaces of $V$ under inclusion via the map $\pi\mapsto\Mov(\pi)$. This implies the following.

\begin{proposition}
\label{prop:NCembed}
If $W$ is a finite Coxeter group, then the map $\pi\mapsto\Fix(\pi)$ defines a poset embedding $NC(W)\hookrightarrow \Pi(W)$.
\end{proposition}

\begin{proof}
We have seen in Proposition \ref{prop:Movembedding} that the map $\pi\mapsto\Mov(\pi)$ is an embedding of $NC(W)$ into the lattice of subspaces of $V$ under inclusion. Since $\Fix(\pi)=\Mov(\pi)^\perp$, we need only show that $\Fix(\pi)$ is in $\Pi(W)$ for all $\pi\in NC(W)$. If wet let $T_\pi=\{t_\alpha:t_\alpha\leq_T\pi , t_\alpha\in T\}$ denote the set of reflections below $\pi$ in absolute order, we will be done if we can show that
\begin{equation*}
\Fix(\pi)=\bigcap_{t_\alpha\in T_\pi} \Fix(t_\alpha),
\end{equation*}
since $\Fix(t_\alpha)=\alpha^\perp$ for all vectors $\alpha$. But note that $\Fix(\pi)$ is contained in $\Fix(t_\alpha)$ for all $t_\alpha\in T_\pi$ (Theorem \ref{th:basicMov} $(2)$), hence $\Fix(\pi)$ is contained in the intersection $\cap_{t_\alpha\in T_\pi} \Fix(t_\alpha)$; and since $\pi$ has a reduced $T$-word of length $r=n-\dim\Fix(\pi)$ (Theorem \ref{th:basicMov}), Carter's Lemma \ref{lemma:carter} implies that $\Fix(\pi)$ and $\cap_{t_\alpha\in T_\pi} \Fix(t_\alpha)$ have the same dimension, hence they are equal.
\end{proof}

Thus, both $NC(W)$ and $NN(W)$ have natural embeddings into the partition lattice $\Pi(W)$ and may be directly compared.

To further this comparison, we consider an additional important interpretation of the partition lattice. If $(W,S)$ is a finite Coxeter system, recall (Section \ref{sec:NC}) that a {\sf parabolic subgroup} of $W$ is any subgroup of the form $wW_Iw^{-1}$, where $w\in W$ and $W_I=\langle I\rangle$ is the subgroup generated by some subset $I\subseteq S$ of simple generators. Let $\L(W)$ denote the set of parabolic subgroups of $W$ partially ordered by inclusion. The following natural result can best be descibed as ``folklore''. The earliest reference known to us --- Barcelo and Ihrig \cite[Theorem 3.1]{barcelo-ihrig}, in 1999 --- is surprisingly recent.

\begin{theorem}[\cite{barcelo-ihrig}]
Consider a finite Coxeter group $W$ together with its geometric representation $\sigma:W\hookrightarrow GL(V)$. The partition lattice $\Pi(W)$ is isomorphic to the lattice $\L(W)$ of parabolic subgroups of $W$ via the following reciprocal maps: send each subspace $U\in\Pi(W)$ to its isotropy subgroup
\begin{equation*}
U\mapsto \Iso(U):=\left\{ w\in W : \sigma(w)(u)=u \text{ for all } u\in U\right\},
\end{equation*}
and send each parabolic subgroup $W'\in\L(W)$ to its fixed subspace
\begin{equation*}
W'\mapsto \Fix(W'):=\left\{ u\in V: \sigma(w)(u)=u \text{ for all } w\in W'\right\}.
\end{equation*}
\end{theorem}

Now we can think of $NC(W)$ and $NN(W)$ as embedded subposets of the partition lattice $\Pi(W)$ or the lattice of parabolic subgroups $\L(W)$. Following Theorem \ref{th:athtypeA}, we might wonder if there is some general notion of ``type equidistribution'' for these posets. Thinking of the embedding of $NC(A_{n-1})$ into the Cayley graph $(A_{n-1},T)$ (Theorem \ref{th:biane}), we note that the ``type'' of a partition is the same as its ``cycle type'' as a permutation. Since cycle types and conjugacy classes coincide in the symmetric group, we might guess that ``conjugacy class'' is the correct general notion of ``type''. Athanasiadis and Reiner proved that --- indeed --- $NC(W)$ and $NN(W)$ are equidistributed in $\L(W)$ by $W$-conjugacy class \cite[Theorem 6.3]{athanasiadis-reiner}.

Given $\pi\in NC(W)$, let $W_\pi=\langle T_\pi\rangle$ denote the subgroup of $W$ generated by the reflections $T_\pi$ below $\pi$ in absolute order. Since $\pi$ is a Coxeter element in some parabolic subgroup $W'\subseteq W$ (Theorem \ref{th:coxelement}), and since absolute order restricts well to parabolic subgroups, Theorem \ref{th:basicMov} \eqref{basicMov3} and Lemma \ref{lemma:basicCox} \eqref{basicCox3} imply that $T_\pi$ is precisely the reflection generating set of $W'$. That is, we have $W_\pi=W'$, and we conclude that $W_\pi$ is itself a parabolic subgroup. Furthermore, given any set of positive roots $A\subseteq\Phi^+$, let $W_A=\Iso\left(\cap_{\alpha\in A}\alpha^\perp\right)$ denote the isotropy group of the intersection of the corresponding hyperplanes. Following Athanasiadis and Reiner \cite[Theorem 6.3]{athanasiadis-reiner}, we have

\begin{theorem}
\label{th:ath-reiner}
Let $W$ be a Weyl group, and fix choices of Coxeter element $c$ and positive system $\Phi^+$ for $W$. Let
\begin{equation*}
f:NC(W)\hookrightarrow\L(W)\quad\text{ and }\quad g:NN(W)\hookrightarrow \L(W)
\end{equation*}
denote the embeddings of $NC(W)$ and $NN(W)$ into the lattice of parabolic subgroups given by $\pi\mapsto W_\pi$ and $A\mapsto W_A$, respectively. Then if $\mathscr{O}$ is any orbit in $\L(W)$ under $W$-conjugation we have
\begin{equation*}
\abs{f^{-1}(\mathscr{O})}=\abs{g^{-1}(\mathscr{O})}.
\end{equation*}
\end{theorem}

\begin{proof}
Athanasiadis and Reiner actually proved this for the embeddings into the partition lattice $\Pi(W)$ given in Theorem \ref{th:NNembed} and Proposition \ref{prop:NCembed}, and for orbits in $\Pi(W)$ under the geometric action of $W$. We must show that composition with the isomorphism $U\mapsto\Iso(U)$ gives the desired result.

Given any $t\in T$ with $t\leq_T \pi$, notice that $t\in\Iso(\Fix(\pi))$, since $\Fix(t)\supseteq\Fix(\pi)$ by Theorem \ref{th:basicMov}. Hence $W_\pi\subseteq\Iso(\Fix(\pi))$. But since the isomorphism $U\mapsto\Iso(U)$ is rank-preserving, we must have $W_\pi=\Iso(\Fix(\pi))$. By definition, $A\mapsto \Iso\left(\cap_{\alpha\in A}\alpha^\perp\right)$ is the composition of the embedding $NN(W)\hookrightarrow \Pi(W)$ given by  $A\mapsto \cap_{\alpha\in A}\alpha^\perp$ and the isomorphism $U\mapsto \Iso(U)$ from $\Pi(W)$ to $\L(W)$.

Finally, given $U\in\Pi(W)$ and $w\in W$, notice that $\Iso(\sigma(w)U)=w\,\Iso(U)w^{-1}$, where $\sigma$ is the geometric representation of $W$. To see this, consider $w'\in\Iso(U)$ and observe that
\begin{equation*}
\sigma(ww'w^{-1})\sigma(w)u=\sigma(w)(\sigma(w')u)=\sigma(w)u,
\end{equation*}
for all $u\in U$, so that $ww'w^{-1}\in\Iso(\sigma(w)U)$. This implies that $w\,\Iso(U)w^{-1}$ is a subgroup of $\Iso(\sigma(w)U)$, and again, since they have the same rank, they are equal. Hence $W$-orbits in $\Pi(W)$ correspond to orbits in $\L(W)$ under $W$-conjugation.
\end{proof}

This remarkable theorem deserves some comment: First, notice that the equidistribution property {\em does not depend on the choices of Coxeter element and positive root system}. Athanasiadis and Reiner proved this result in a case-by-case way, using earlier results of Athanasiadis for types $A$ and $B$ \cite{athanasiadis:classical}; type $D$ results proved in the same paper \cite[Corollaries 5.2 and 5.4]{athanasiadis-reiner}; and computer verification for the exceptional types. In type $A$, Athanasiadis also gave a bijective proof \cite[Section 4]{athanasiadis:classical}, but the bijection did not extend to other types. As of this writing, {\em no uniform proof is known}.

As we have said, the posets $NC(W)$ and $NN(W)$ seem not to have much in common, yet Theorem \ref{th:ath-reiner} shows that the noncrossing partitions and nonnesting partitions are related enumeratively in a very deep way. Following Athanasiadis and Reiner, we suggest the following problem.

\begin{problem}
\label{prob:NC=NN}
What is the nature of the relationship between $NC(W)$ and $NN(W)$? Find a ``type-preserving'' bijection that explains Theorem \ref{th:ath-reiner} in a uniform way.
\end{problem}

Indeed, we do not currently know of {\em any} bijection between $NC(W)$ and $NN(W)$. This problem is particularly important since we also do not have a uniform proof that the noncrossing partitions are counted by the Coxeter-Catalan number. For the nonnesting partitions, such a uniform proof {\em does} exist (see the Introduction), and it would be nice to connect the noncrossing partitions to this uniform theory.

\subsection{Geometric Multichains of Filters}
\label{sec:geometricmchains}

Instead of resolving this mystery, we will now generalize it. In this memoir, we have defined a generalization of the noncrossing partitions for each positive integer $k$. Recently, Athanasiadis has developed a parallel generalization of the nonnesting partitions. To describe this, we must consider a different (perhaps more natural) partial order on the set of antichains $NN(W)$.

If $W$ is a Weyl group and $A$ is an antichain in the root poset $(\Phi^+,\leq)$, let
\begin{equation*}
\filter(A):=\bigcup_{\alpha\in A}\filter(\alpha)=\left\{ \beta\in\Phi^+ :\alpha\leq\beta \text{ for some }\alpha\in A\right\}
\end{equation*}
denote the {\sf order filter} in $(\Phi^+,\leq)$ generated by $A$. It is well-known that the set of order filters in a poset is in bijection with the set of antichains: the antichain corresponding to a filter is its set of minimal elements.

\begin{definition}
Define the {\sf filter order} on $NN(W)$ by setting
\begin{equation*}
A\leq B\quad\Longleftrightarrow\quad\filter(A)\subseteq\filter(B).
\end{equation*}
\end{definition}

The filter order on $NN(W)$ seems to be more interesting than the refinement order (Definition \ref{def:NNrefinement}) from a poset-theoretic point of view. Indeed, since the filter order on $NN(W)$ is the same as the lattice of order ideals of the dual root poset $NN(W)\cong J\left( (\Phi^+,\leq)^*\right)$, it is a {\sf distributive lattice} \cite[Theorem 3.4.1]{stanley:ec1}. There is also a nice geometric interpretation of this poset.

Given a Weyl group $W$, consider its crystallographic root system $\Phi=\Phi^+\sqcup\Phi^-$ with respect to the inner product $(\cdot,\cdot)$, and define the {\sf Shi arrangement} as the collection of hyperplanes
\begin{equation*}
\Shi(W):=\left\{ H_\alpha^0:\alpha\in\Phi^+\right\}\cup\left\{H_\alpha^1:\alpha\in\Phi^+\right\}
\end{equation*}
where $H_\alpha^i=\{x\in V: (x,\alpha)=i\}$ for all integers $i$. That is, the Shi arrangement consists of the Coxeter arrangement, together with the affine hyperplane extended by one unit in the direction of each of the positive roots. The region $\left\{ x\in V: (x,\alpha)>0 \text{ for all } \alpha\in\Phi^+\right\}\subseteq V$ is called the {\sf positive cone}, and we call any chamber of the Shi arrangement inside the positive cone a {\sf positive chamber}.\footnote{In Chapter \ref{sec:introduction}, we frequently referred to the Catalan arrangement. It should be noted that, within the positive cone, the Catalan and Shi arrangements coincide. We have chosen to emphasize the Shi arrangement at this point because of its relation to ``parking functions'' and ``diagonal harmonics'' (Section \ref{sec:future}).} Cellini and Papi proved the following.

We say that a chamber $\C$ is {\sf on the positive side} of hyperplane $H_\alpha^i$ if we have $(x,\alpha)>i$ for all $x\in\C$. 

\begin{theorem}[\cite{cellini-papi}]
\label{th:cellini-papi}
Given a chamber $\C$ in the Shi arrangement, let $\filter(\C)\subseteq\Phi^+$ denote the set of positive roots $\alpha$ such that $\C$ is on the positive side of $H_\alpha^1$. Then the map $\C\mapsto\filter(\C)$ is a bijection from the positive chambers to the order filters in the root poset.
\end{theorem}

\begin{figure}
\vspace{.1in}
\begin{center}
\input{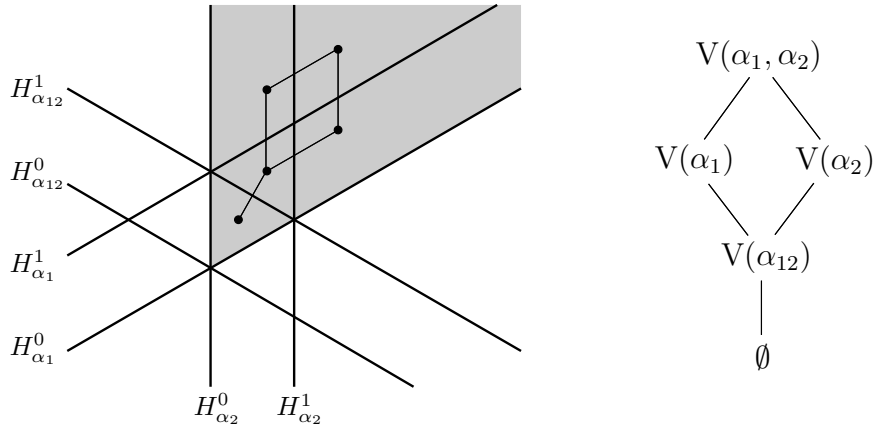}
\end{center}
\caption{The Shi arrangement $\Shi(A_2)$ and filter order on $NN(A_2)$}
\label{fig:shi(a2)}
\end{figure}

Again, this is easiest to understand pictorially. Figure \ref{fig:shi(a2)} displays the Shi arrangement of type $A_2$ with positive cone shaded, and indicates the bijection with order filters in the root poset (compare this with Figure \ref{fig:a2roots}). In this case, the chamber corresponding to the filter $\filter(\alpha_1)=\{\alpha_1,\alpha_{12}\}$ is on the positive sides of the hyperplanes $H_{\alpha_1}^1$ and $H_{\alpha_{12}}^1$ and is on the negative side of the hyperplane $H_{\alpha_2}^1$. 

Furthermore, we can define a partial order on the set of positive chambers in the Shi arrangement by declaring a cover relation $\C\prec\mathcal{C'}$ whenever the chambers $\C$ and $\mathcal{C'}$ share a wall and this wall separates $\mathcal{C'}$ from the origin. Considering Figure \ref{fig:shi(a2)}, we can see that this partial order on chambers is isomorphic to the partial order on filters $NN(W)$.

The key to understanding Athanasiadis' generalization of the nonnesting partitions is to consider a certain generalization of the Shi hyperplane arrangement.

\begin{definition}
For each positive integer $k$, define the {\sf extended Shi arrangement} as the collection of hyperplanes
\begin{equation*}
\Shi^{(k)}(W):=\bigcup_{i=-k+1}^k \left\{ H_\alpha^i:\alpha\in\Phi^+\right\};
\end{equation*}
that is, the Coxeter arrangement together with the first $k$ positive affine extensions and the first $k-1$ negative affine extensions.
\end{definition}

By studying characteristic polynomials, Athanasiadis was able to count the positive chambers in the extended Shi arrangement \cite[Corollary 1.3]{athanasiadis:cat} and found that they are counted by a ``generalized Catalan number''. In our notation,

\begin{theorem}
\label{th:ath-shi}
Given a Weyl group $W$ and positive integer $k$, the number of positive chambers in the extended Shi arrangement is equal to the Fuss-Catalan number $\Cat^{(k)}(W)$ \eqref{eq:fusscat} and the number of bounded positive chambers is equal to the positive Fuss-Catalan number $\Cat_+^{(k)}(W)$ \eqref{eq:posfusscat}. 
\end{theorem}

The following result had been conjectured by Edelman and Reiner \cite[Conjecture 3.3]{edelman-reiner} and Athanasiadis \cite[Question 6.2]{athanasiadis:survey} before it was proved uniformly by Yoshinaga \cite[Theorem 1.2]{yoshinaga}. The characteristic polynomial of a hyperplane arrangement will be discussed briefly in Section \ref{sec:triangles}.

\begin{theorem}[\cite{yoshinaga}]
Given a Weyl group $W$ and positive integer $k$, the characteristic polynomial of the extended Shi arrangement $\Shi^{(k)}(W)$ is given by
\begin{equation*}
\chi(\Shi^{(k)}(W),t)=(t-kh)^n,
\end{equation*}
where $h$ is the Coxeter number and $n$ is the rank of $W$. Hence the arrangement has a total of $(kh+1)^n$ chambers and $(kh-1)^n$ bounded chambers.
\end{theorem}

\begin{figure}
\vspace{.1in}
\begin{center}
\input{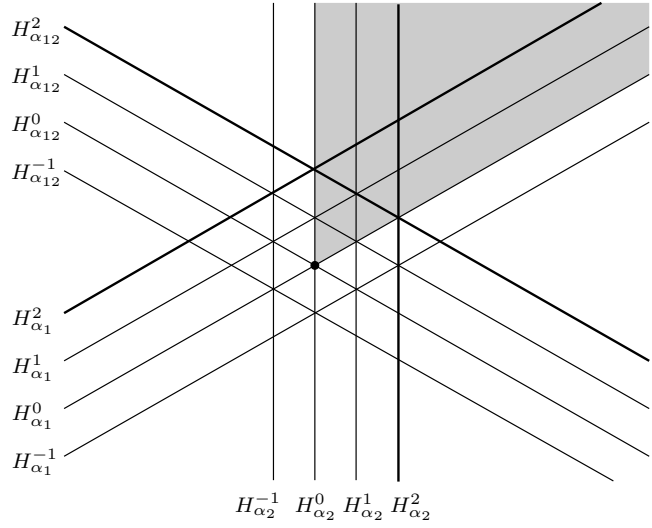}
\end{center}
\caption{The extended Shi arrangement $\Shi^{(2)}(A_2)$}
\label{fig:shi^2(a2)}
\end{figure}

For example, consider the Weyl group $A_2$ with rank $n=2$ and Coxeter number $h=3$. Figure \ref{fig:shi^2(a2)} displays the extended Shi arrangement $\Shi^{(2)}(A_2)$ with the positive cone shaded. Notice that there are $\Cat^{(2)}(A_2)=12$ positive chambers and $\Cat_+^{(2)}(A_2)=7$ bounded positive chambers. Furthermore, there are a total of $(2\cdot3+1)^2=49$ chambers and $(2\cdot3-1)^2=25$ bounded chambers.

Following Cellini and Papi (Theorem \ref{th:cellini-papi}), we see that each positive chamber of the extended Shi arrangement corresponds to a collection of order filters in the root poset.

\begin{definition}
Given a positive chamber $\C$ of $\Shi^{(k)}(W)$, define
\begin{equation*}
\filter_i(\C):=\left\{ \alpha\in\Phi^+: (v,\alpha)>i \text{ for all } v\in\C\right\}
\end{equation*}
for each integer $1\leq i\leq k$. By Theorem \ref{th:cellini-papi}, each of these is an order filter in the root poset.
\end{definition}

\begin{figure}
\vspace{.1in}
\begin{center}
\input{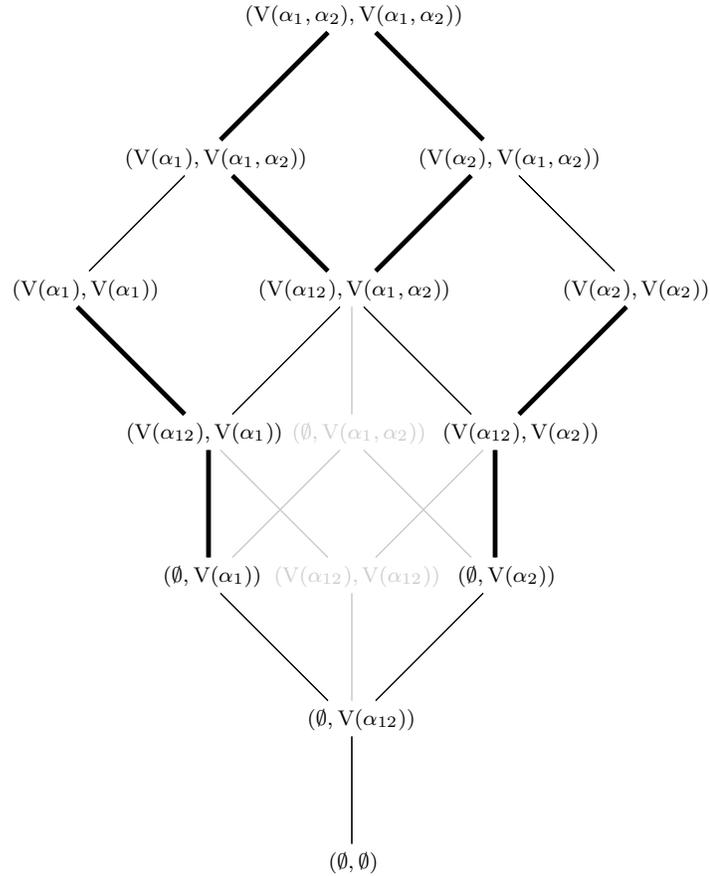}
\end{center}
\caption{The filter order on $NN^{(2)}(A_2)$}
\label{fig:nn^2(a2)}
\end{figure}

If $\C$ is a positive chamber of $\Shi^{(k)}(W)$ then it is also easy to see that the filters $\filter_i(\C)$, $1\leq i\leq k$, form a multichain under inclusion:
\begin{equation*}
\filter_k(\C)\subseteq\filter_{k-1}(\C)\subseteq\cdots\subseteq\filter_2(\C)\subseteq\filter_1(\C).
\end{equation*}
Thus, the function
\begin{equation*}
\C\mapsto \left(\filter_k(\C),\ldots,\filter_2(\C),\filter_1(\C)\right)
\end{equation*}
defines an injective map from the positive chambers of $\Shi^{(k)}(W)$ to $k$-multichains in the filter order on $NN(W)$. However, this map is {\em not} surjective. For example, the multichains $\left(\filter(\alpha_{12}),\filter(\alpha_{12})\right)$ and $\left(\emptyset,\filter(\alpha_1,\alpha_2)\right)$ do not correspond to positive chambers of $\Shi^{(2)}(A_2)$ (see Figures \ref{fig:shi^2(a2)} and \ref{fig:nn^2(a2)}). It is natural then to ask for a characterization of the multichains in the filter order on $NN(W)$ that do correspond to positive Shi chambers. The answer to this question was explained by Athanasiadis in \cite{athanasiadis:cat} and \cite{athanasiadis:nar}.

Given two subsets $A,B\subseteq\Phi^+$ of positive roots, define
\begin{equation*}
A+B:=\left\{ \alpha+\beta: \alpha\in A, \beta\in \,B\right\}.
\end{equation*}
The next definition is from \cite[page 180]{athanasiadis:nar}.

\begin{definition}
Suppose that $\filter_k \subseteq \filter_{k-1}\subseteq\cdots\subseteq \filter_1\subseteq \filter_0=\Phi^+$ is a $k$-multichain of filters in the root poset and set $\filter_i=\filter_k$ for $i>k$. We say that this multichain is {\sf geometric} if
\begin{equation*}
\left( \filter_i+\filter_j\right)\cap\Phi^+ \subseteq \filter_{i+j}
\end{equation*}
holds for all indices $i,j$ and
\begin{equation*}
\left( \Lambda_i+\Lambda_j\right)\cap\Phi^+\subseteq \Lambda_{i+j}
\end{equation*}
holds for all indices $i,j\geq 1$ with $i+j\leq k$, where $\Lambda_i=\Phi^+\setminus \filter_i$.
\end{definition}

These geometric multichains are precisely those that correspond to positive Shi chambers.

\begin{theorem}[\cite{athanasiadis:nar}]
Given a Weyl group $W$ and positive integer $k$, the map
\begin{equation*}
\C\mapsto \left(\filter_k(\C),\ldots,\filter_2(\C),\filter_1(\C)\right)
\end{equation*}
is a bijection from positive chambers of $\Shi^{(k)}(W)$ to geometric $k$-multichains of filters.
\end{theorem}

Letting $NN^{(k)}(W)$ denote the set of geometric $k$-multichains in the filter order on $NN(W)$\footnote{Notice the notational harmony with $NC^{(k)}(W)$, which is defined as the set of $k$-multichains in $NC(W)$.}, we have a generalization of the filter order.
\begin{definition}
Define the {\sf filter order} on $NN^{(k)}(W)$ by setting
\begin{equation*}
(A_k,\ldots,A_1)\leq (B_k,\ldots,B_1)
\end{equation*}
whenever $A_i\subseteq B_i$ for all $1\leq i\leq k$.
\end{definition}

Figure \ref{fig:nn^2(a2)} displays the Hasse diagram of the filter order on $NN^{(2)}(A_2)$. Again, notice how this relates to the positive chambers in Figure \ref{fig:shi^2(a2)}. We partially order the positive chambers by setting $\C\prec\C'$ whenever $\C$ and $\C'$ share a wall that separates $\C'$ from the origin. Then this partial order on chambers is isomorphic to the filter order on $NN^{(2)}(A_2)$. There is also a finer structure to observe in Figure \ref{fig:nn^2(a2)}.

If $\C$ is a chamber in a hyperplane arrangement, we call the hyperplane $H$ a {\sf wall} of $\C$ if it supports one of the facets of $\C$, and denote the set of walls of $\C$ by $\WL(\C)$. Furthermore, if the wall $H\in\WL(\C)$ separates $\C$ from the origin (including the case when $H$ contains the origin), we call $H$ a {\sf floor} of $\C$, and otherwise we call it a {\sf ceiling} of $\C$ (gravity points toward the origin). Let $\FL(\C)$ denote the set of floors of $\C$ and let $\CL(\C)$ denote the set of ceilings, so that $\WL(\C)=\FL(\C)\sqcup\CL(\C)$. In the extended Shi arrangement, floors and ceilings come in several ``colors''.

\begin{definition}
Given a positive chamber $\C$ in $\Shi^{(k)}(W)$ with corresponding multichain $(\filter)_k\in NN^{(k)}(W)$, define the set of {\sf $i$-colored floor roots}
\begin{equation*}
\FL_i(\C)=\FL_i(\filter)_k :=\left\{ \alpha\in\Phi^+ : H_\alpha^i\in\FL(\C)\right\}
\end{equation*}
and the set of {\sf $i$-colored ceiling roots}
\begin{equation*}
\CL_i(\C)=\CL_i(\filter)_k:=\left\{ \alpha\in\Phi^+ :H_\alpha^i\in\CL(\C)\right\}
\end{equation*}
for all $1\leq i\leq k$.
\end{definition}

Notice that floors and ceilings correspond to cover relations in $NN^{(k)}(W)$, and we can read these quite easily from Figure \ref{fig:nn^2(a2)}. The edges in the Hasse diagram connected to a multichain correspond to its walls, with the edges below being floors and the edges above being ceilings. Here bold edges have color $2$ and regular edges have color $1$. The slope of the edge indicates which root it corresponds to: slope $-1$ for $\alpha_1$, slope $1$ for $\alpha_2$ and vertical edges for $\alpha_{12}$. (The lighter elements and edges indicate the two $2$-multichains in $NN(A_2)$ that are {\em not} geometric, and hence do not belong to $NN^{(2)}(A_2)$.) For example, the multichain $(\filter(\alpha_{12}),\filter(\alpha_1))$ has walls $\{ H_{\alpha_2}^1,H_{\alpha_1}^2,H_{\alpha_{12}}^2\}$, with floors $\{H_{\alpha_{12}}^2\}$ and ceilings $\{H_{\alpha_2}^1,H_{\alpha_1}^2\}$.

The walls and ceilings have interesting enumerative properties. Athanasiadis referred to the set $\FL_k(\filter)_k$ as the set of {\sf rank $k$ indecomposables} of the multichain $(\filter)_k$, and he was able to count the elements of $NN^{(k)}(W)$ by their number of rank $k$ indecomposables when $W$ is of classical type. In Figure \ref{fig:nn^2(a2)}, if we count multichains by their number of $2$-colored floors (rank $2$ indecomposables) notice that there are $5$ multichains with zero, $6$ multichains with one, and $1$ multichain with two rank $2$ indecomposables. Compare this with the rank numbers of the poset $NC^{(2)}(A_2)$ in Figure \ref{fig:nc^2(a2)}. In general, we believe the following is true.

\begin{conjecture}
\label{conj:NCnar=NNnar}
Given a Weyl group $W$ and positive integer $k$, the number of elements of $(\filter)_k\in NN^{(k)}(W)$ with $\abs{\FL_k(\filter)_k}=i$ is equal to the number of elements of $NC^{(k)}(W)$ with rank $i$; that is, the Fuss-Narayana number $\Nar^{(k)}(W,i)$ \eqref{eq:fussnar}.
\end{conjecture}

In the classical types, this follows from Proposition 5.1 in Athanasiadis \cite{athanasiadis:nar}, where he computed these numbers. We expect that it will also hold for the exceptional types, but a computer verification has not been carried out. In principle, one would need a ``trick'' to accomplish this, since it involves verification for an infinite number of values of $k$. (In our computations for Theorem \ref{th:fussnar}, we were able to use the ``trick'' of zeta polynomials.) This enumerative correspondence is surprising since Athanasiadis' work and the work in this memoir were done independently and were motivated by different subjects.

It will be interesting to study the floors and ceilings statistics in more depth. We believe that the sets $\FL_i(V)_k$ and $\CL_i(V)_k$ always have cardinality $\leq n$ (this is true for $\FL_i(V)_k$ since Athanasiadis showed it is an antichain in the root poset). Hence we make the following definition.

\begin{definition}
Given a Weyl group $W$ and positive integer $k$, define the {\sf $i$-th floor vector} as $(fl_i(0),fl_i(1),\ldots,fl_i(n))$, where
\begin{equation*}
fl_i(j):=\#\left\{ (\filter)_k\in NN^{(k)}(W) : \abs{\FL_i(\filter)_k}=j\right\},
\end{equation*}
and define the {\sf $i$-th ceiling vector} as $(cl_i(0),cl_i(1),\ldots,cl_i(n))$, where
\begin{equation*}
cl_i(j):=\#\left\{ (\filter)_k\in NN^{(k)}(W) : \abs{\CL_i(\filter)_k}=j\right\}.
\end{equation*}
\end{definition}

Notice that the entries of the $i$-th floor vector sum to $\Cat^{(k)}(W)$ for all $0\leq i\leq n$, as do the entries of the $i$-th ceiling vector, and the $k$-th floor vector contains the Fuss-Narayana numbers (Conjecture \ref{conj:NCnar=NNnar}). We do not understand these vectors in general, but we have a conjecture.

\begin{conjecture}
\label{conj:floors=ceilings}
We have $fl_i(j)=cl_i(j)$ for all $1\leq i\leq k$ and $0\leq j\leq n$.
\end{conjecture}

That is, the number of chambers with $j$ $i$-colored floors is equal to the number of chambers with $j$ $i$-colored ceilings. This seems to indicate the existence of some sort of duality. However, there does {\em not} exist any map on chambers that swaps $i$-colored floors and $i$-colored ceilings.

Conjecture \ref{conj:NCnar=NNnar} states that the sets $NC^{(k)}(W)$ and $NN^{(k)}(W)$ are equidistributed by ``rank''. Following Theorem \ref{th:ath-reiner}, perhaps there is also an equidistribution by ``type''.

\begin{conjecture}
\label{conj:NCk=NNk}
Let $W$ be a Weyl group, and fix a choice of Coxeter element $c$ and positive system $\Phi^+$ for $W$. Let
\begin{equation*}
f_k:NC^{(k)}(W)\to\L(W)
\end{equation*}
denote the map $(\pi)_k\mapsto W_{\pi_1}$ that sends a multichain to the parabolic subgroup determined by its bottom element $\pi_1$, and let
\begin{equation*}
g_k:NN^{(k)}\to \L(W)
\end{equation*}
denote the map $(\filter)_k\mapsto W_{\FL_k(\filter)_k}$ that sends a multichain $(\filter)_k$ to the parabolic subgroup generated by its set $\FL_k(\filter)_k$ of $k$-colored floors. Then if $\mathscr{O}$ is any orbit in $\L(W)$ under $W$-conjugation, we have
\begin{equation*}
\abs{f_k^{-1}\left(\mathscr{O}\right)}=\abs{g_k^{-1}\left(\mathscr{O}\right)}.
\end{equation*}
\end{conjecture}

Note: In contrast with Theorem \ref{th:ath-reiner}, the maps $f_k$ and $g_k$ are not injective in general. Also, we expect that the same result holds when $g_k$ is defined in terms of the map $(\filter)_k\mapsto W_{\CL_k(\filter)_k}$ that sends a filter to the parabolic generated by its $k$-colored ceilings. In order to prove this conjecture, one would only need to show that it holds for one particular choice of Coxeter element $c$ and positive system $\Phi^+$. Then, since all Coxeter elements are $W$-conjugate and any positive system can be sent into any other by an element of $W$, the result follows.

The next problem generalizes Problem \ref{prob:NC=NN}.

\begin{problem}
What is the true nature of the relationship between $NC^{(k)}(W)$ and $NN^{(k)}(W)$? Find a ``type-preserving'' bijection between them. Which elements of $NC^{(k)}(W)$ should the bounded Shi chambers correspond to? Is there some natural partial order on $NN^{(k)}(W)$ that generalizes the refinement order? Can $NN^{(k)}(A_{n-1})$ be realized as a poset of partitions of the set $[kn]$? Is there also some equidistribution by total type (see Definition \ref{def:totaltype})? Is there a natural generalization of $NN^{(k)}(W)$ for noncrystallographic Coxeter groups?
\end{problem}

We have more questions than answers right now.

\section{Cluster Complexes}
\label{sec:cluster}

\subsection{Polygon Triangulations and the Associahedron}
\label{sec:associahedron}
The study of triangulations of a convex polygon goes back at least to Leonhard Euler. In a 1751 letter to Christian Goldbach, he suggested a method for counting these, but he could not prove that it worked. Later, Euler communicated the problem to Hungarian mathematician Jan Andrej Segner, who then provided a correct method in 1756. In the late 1830's, Joseph Liouville posed as an open challenge to prove whether Euler's original method was correct. He received many solutions, including one from a mathematician named Eug\`ene Charles Catalan (for more information, see \cite{ranjan}). They had all proven the following.

\begin{theorem}
The number of triangulations of a convex $(n+2)$-gon is equal to the Catalan number $\Cat(A_{n-1})=\frac{1}{n}\binom{2n}{n-1}$.
\end{theorem}

The modern history of polygon triangulations began in the early 1960's with James Stasheff's PhD thesis \cite{stasheff}, in which he set down the foundations of a homotopy-invariant notion of associativity. The motivating question seems at first unrelated to combinatorics: to what extent can one define a projective space over the octonions? Milnor had provided a construction of projective spaces for an arbitrary topological group, but this construction depended in an essential way on the associativity of the group operation. Stasheff considered the problem of constructing projective spaces for a topological space with a continuous multiplication that {\em might not be associative}. (When the multiplication is associative {\em up to homotopy}, such a space is called an {\sf $H$-space}, or {\sf Hopf space}.) An essential part of Stasheff's construction was a cell complex $K_i$ that serves to replace the standard cube $I^i$ as a model for singular homology. In the 1970's this work led Hilton and Roitberg to a negative answer to the homotopy version of Hilbert's Fifth Problem: does an $H$-space that is also a manifold have to be a Lie group?

The complex $K_i$ defined by Stasheff has since become important in several areas of mathematics, including combinatorics. It turns out that $K_i$ can be realized as a {\sf simple convex polytope} in $\reals^i$ (here ``simple'' means that every vertex is contained in the same number of edges), hence it makes sense to speak of the dual {\sf simplicial polytope} (in which every maximal face is a simplex). For ease of presentation, we prefer to discuss the simplicial version first.

By a {\sf diagonal} of a convex polygon, we mean any edge between vertices, excluding the pre-existing edges of the polygon. Two diagonals are said to {\sf cross} if they intersect in their interiors (i.e. not at an endpoint). We say that a set of diagonals is a {\sf partial triangulation} if it may be completed to a full triangulation; that is, if the diagonals are pairwise noncrossing.

\begin{definition}
Let $\Delta(n)$ denote the abstract (flag) simplicial complex whose vertices are the diagonals of a convex $(n+2)$-gon and whose faces are partial triangulations. We call $\Delta(n)$ the {\sf simplicial associahedron}.
\end{definition}

Recall that a simplex with $d+1$ vertices has {\sf dimension} $d$, and note that a full triangulation of a convex $(n+2)$-gon contains $n-1$ diagonals. Thus $\Delta(n)$ is a pure $(n-2)$-dimensional complex with $\Cat(A_{n-1})=\frac{1}{n}\binom{2n}{n-1}$ maximal faces (``facets'') --- one corresponding to each triangulation --- and
\begin{equation*}
\binom{n+2}{2}-(n+2)=\frac{n(n+1)}{2}-1
\end{equation*}
vertices --- one corresponding to each diagonal. Notice also that two facets of $\Delta(n)$ share a codimension $1$ face whenever they are related by a {\sf diagonal flip}; that is, when one diagonal is replaced by the other in some quadrilateral. In general, the number of faces of $\Delta(n)$ containing $i$ diagonals is given by
\begin{equation}
\label{eq:KC}
\frac{1}{i+1}\binom{n-1}{i}\binom{n+i+1}{i}.
\end{equation}
The solution to this counting problem was written down by Kirkman \cite{kirkman} in 1857 but the first proof was given by Cayley \cite{cayley} in 1891; the numbers \eqref{eq:KC} are today known as the {\sf Kirkman-Cayley numbers}.\footnote{Kirkman and Cayley actually solved a more general problem; see Section \ref{sec:generalcluster} below.}

Now we may easily define Stasheff's complex.

\begin{definition}
Let $\nabla(n)$ denote the polytopal complex that is dual to $\Delta(n)$. That is, the vertices of $\nabla(n)$ correspond to triangulations of a convex $(n+2)$-gon and for each partial triangulation $D$ there is a corresponding face of $\nabla(n)$, consisting of triangulations that contain $D$. We call $\nabla(n)$ the {\sf simple associahedron}.
\end{definition}

Stasheff did not use the language of polygon triangulations. Instead, since he was interested in associativity, each vertex of his complex $K_i$ corresponded to a different way to parenthesize a string of length $i+2$. The fact that $\nabla(n)$ is combinatorially equivalent to $K_{n-1}$ then follows from the well-known relationship between polygon triangulations and parenthesizations of a string.

The so called ``Stasheff polytope'' has an interesting history, and there have been many different proofs that this complex is a polytope. According to Stasheff, the first proof was given by Milnor (unpublished, unrecorded). Later, the same complex arose independently in combinatorics, where it was called the ``associahedron''. It was defined in this context by Haiman (see \cite{lee}), who also provided a proof of polytopality (unpublished manuscript). The first published proof was given by Lee \cite{lee}. More recently, Loday \cite{loday} and Postnikov \cite{postnikov} have constructed geometric realizations. The associahedron also occurs in the general theory of fiber polytopes by Billera and Sturmfels \cite{billera-sturmfels1,billera-sturmfels2} and the theory of secondary polytopes by Gelfand, Kapranov and Zelevinsky \cite{GKZ}. A type $B$ associahedron was defined by Simion \cite{simion:associahedron} (see also \cite{hohlweg-lange}). For many beautiful pictures of associahedra, see the survey \cite{fomin-reading:survey} by Fomin and Reading.

\begin{figure}
\vspace{.1in}
\begin{center}
\input{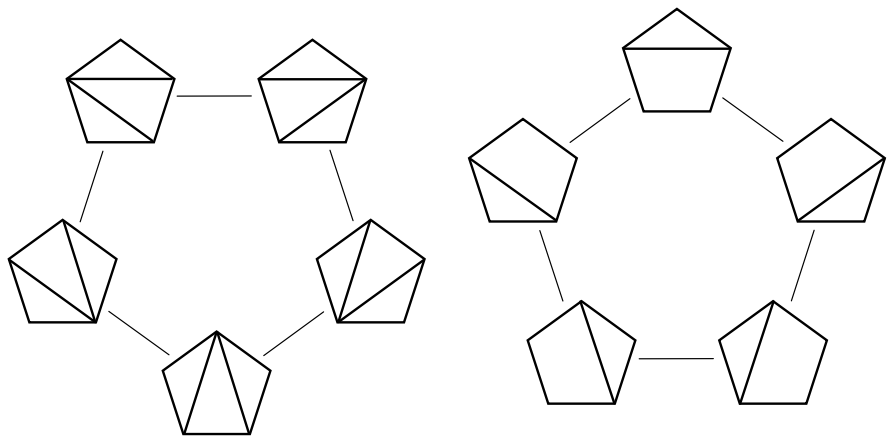}
\end{center}
\caption{The simple associahedron $\nabla(3)$ and simplicial associahedron $\Delta(3)$}
\label{fig:associahedra}
\end{figure}

Although the simple associahedron $\nabla(n)$ is historically more important, the simplicial associahedron $\Delta(n)$ is a more natural object and it is easier to work with. For instance, it is a {\sf flag complex} (or a {\sf clique complex}); that is, a set of vertices forms a face if and only if the vertices are connected pairwise by edges. Thus it is defined entirely by its set of vertices and the binary relation called ``crossing''. Figure \ref{fig:associahedra} displays the simple associahedron $\nabla(3)$ and the simplicial associahedron $\Delta(3)$. Note that there are $\Cat(A_2)=\frac{1}{3}\binom{6}{2}=5$ triangulations and $\frac{3\cdot 4}{2}-1=5$ diagonals of a convex pentagon.

Now we return to our home base, the subject of finite reflection groups. Following the pattern of the noncrossing and the nonnesting partitions, we should look for a definition of the associahedron that explicitly recognizes a connection with the symmetric group, and that is suitable for generalization. In this case, the complexes in Figure \ref{fig:associahedra} should correspond to the group $A_2$.

\subsection{The Cluster Complex}

Such a generalization has recently been found by Fomin and Zelevinsky, in connection with their theory of {\sf cluster algebras}. They introduced the idea of a cluster algebra in \cite{fomin-zelevinsky:foundations} as an algebraic and combinatorial framework for the theory of dual canonical bases and total positivity in Lie groups, and cluster algebras have recently found unexpected applications in Teichm\"uller theory and the theory of quiver representations. Many of these connections are described in the survey \cite{fomin-zelevinsky:survey}.

We will not describe the algebraic framework here, and instead focus on the combinatorics. The starting point is the finite type classification \cite{fomin-zelevinsky:finitetype}, in which the finite type cluster algebras are shown to correspond to crystallographic root systems. The key combinatorial idea is an algebraic generalization of ``polygon triangulations''. Although there do not exist cluster algebras for the noncrystallographic finite Coxeter groups, much of the combinatorics {\em can} be generalized. An excellent introduction to the combinatorics of cluster algebras is given in the survey by Fomin and Reading \cite{fomin-reading:survey}. We follow their approach here.

Let $W$ be a finite Coxeter group with (possibly noncrystallographic) root system $\Phi$ as described in Section \ref{sec:rootsystems}. Given a positive system $\Phi^+$ with corresponding simple system $\Pi$, we define a special set of roots.

\begin{definition}
The set of positive roots together with the negative simple roots
\begin{equation*}
\Pge:=\Phi^+\sqcup \left(-\Pi\right)
\end{equation*}
is called the set of {\sf almost-positive roots}.
\end{definition}

The elements of the set $\Pge$ will play the role of ``diagonals''. If $\Phi$ is the root system of type $A_{n-1}$, recall that the simple roots $\Pi$ are in bijection with adjacent transpositions and the positive roots $\Phi^+$ are in bijection with all transpositions, so the set $\Pge$ has cardinality
\begin{equation*}
\abs{\Pge}=\binom{n}{2}+(n-1)=\frac{n(n+1)}{2}-1,
\end{equation*}
which is equal to the number of diagonals of a convex $(n+2)$-gon. 

\begin{figure}
\vspace{.1in}
\begin{center}
\input{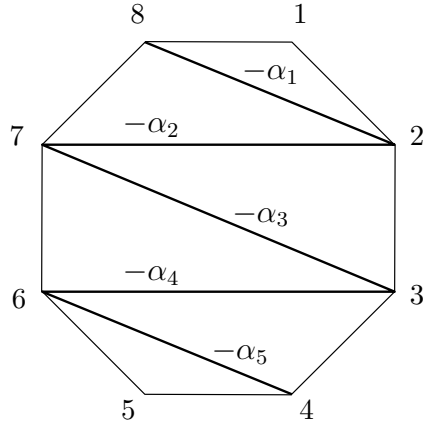}
\end{center}
\caption{The snake of negative simple roots of type $A_5$}
\label{fig:snake2}
\end{figure}

\begin{figure}
\vspace{.1in}
\begin{center}
\input{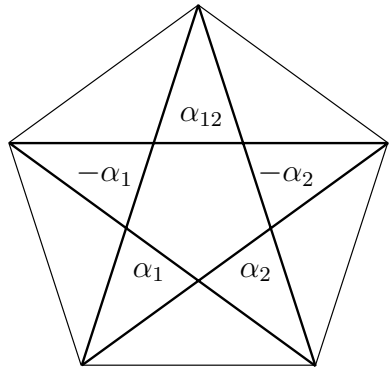}
\end{center}
\caption{Labelling of the diagonals of a pentagon by the almost-positive roots of type $A_2$}
\label{fig:pentagon}
\end{figure}

Fomin and Zelevinsky described a natural correspondence between almost-positive roots and diagonals that is closely related to the ``snake generating set'' shown before in Figure \ref{fig:snake}. Consider the regular $(n+2)$-gon with vertices labelled clockwise by $1,2,\ldots,n+2$. If the simple roots are denumerated
\begin{equation*}
\Pi=\left\{ \alpha_1,\alpha_2,\ldots,\alpha_{n-1}\right\},
\end{equation*}
then we label the diagonal connecting vertices $\lceil\frac{i}{2}\rceil+1$ and $n-\lfloor\frac{i}{2}\rfloor+2$ by the negative simple root $-\alpha_i$ for all $1\leq i\leq n-1$. This is called the ``snake'' of negative simple roots. Figure \ref{fig:snake2} displays these diagonals for the group $A_5$. Note that any other diagonal will necessarily cross at least one of these negative simple diagonals. If a given diagonal crosses with the set of negative diagonals $\{-\alpha_i:i\in I\subseteq[n-1]\}$ then we label it by the positive root $\sum_{i\in I}\alpha_i$. It can easily be seen that this determines a bijection between diagonals of the convex $(n+2)$-gon and almost-positive roots of type $A_{n-1}$. For example, consider the root system of type $A_2$, as shown in Figure \ref{fig:a2roots}. Figure \ref{fig:pentagon} displays the labelling of the diagonals of a pentagon by the almost-positive roots $\Pge=\{-\alpha_1,-\alpha_2,\alpha_1,\alpha_2,\alpha_{12}\}$.

Now, in order to define an associahedron, we must be able to detect when two almost-positive roots ``cross''. When one of the roots is negative this is easy: the root $\alpha\in\Pge$ crosses with the negative simple root $-\alpha_i$ if and only if $\alpha_i$ occurs in the simple root expansion of $\alpha$. In particular, the set of negative simple roots are mutually noncrossing, which is clear from the diagram (Figure \ref{fig:snake2}).
But how can we detect whether two positive roots cross? Considering rotations of the polygon in Figure \ref{fig:snake2}, we notice that every diagonal may be rotated so that it coincides with an element of the snake, and clearly two diagonals will cross if and only if their rotated images cross. Thus, we wish to define an algebraic ``rotation'' on the set $\Pge$ that corresponds to rotation of the polygon in type $A$.

This was accomplished by Fomin and Zelevinsky \cite{fomin-zelevinsky:ysystems}, using the idea of a ``tropical Coxeter element''. Recall the definition of a bipartite Coxeter element from Section \ref{sec:automorphisms}. If we set $t_i=t_{\alpha_i}$ for all $1\leq i\leq n-1$, then $S=\{t_1,t_2,\ldots,t_{n-1}\}$ is the simple generating set for $A_{n-1}$ corresponding to the simple roots $\Pi=\{\alpha_1,\alpha_2,\ldots,\alpha_{n-1}\}$. Since the Coxeter diagram of a finite Coxeter group is a forest, the simple generators can be partitioned into disjoint sets $S=S_\ell\sqcup S_r$ such that the elements of $S_\ell$ commute pairwise, as do the elements of $S_r$. Set $\ell=\prod_{t\in S_\ell} t$ and $r=\prod_{t\in S_r} t$. Then we regard the triple $(c,\ell,r)$ with $c=\ell r$ as a bipartite Coxeter element. Amazingly, this construction can be ``tropicalized'' to create a ``rotation'' map on the set $\Pge$.

Let $\Pi_\ell$ and $\Pi_r$ denote the simple roots corresponding to the sets $S_\ell$ and $S_r$, respectively. If $\alpha\in\Pi$, it is well-known that the simple reflection $t_\alpha$ sends the set $\Phi^+\cup\{-\alpha\}$ to itself \cite[Proposition 4.1]{humphreys}. Furthermore, if $\alpha\in\Pi_\ell$, then $t_\alpha$ acts as the identity on the set $\Pi_\ell\setminus\{\alpha\}$ since the elements of $\Pi_\ell$ are mutually orthogonal. That is, the map $\ell$ is an involution on $\Phi^+\cup(-\Pi_\ell)$, and similarly $r$ is an involution on $\Phi^+\cup(-\Pi_r)$. We have shown that the following maps are well-defined.

\begin{definition}
\label{def:tropical}
Given a bipartite Coxeter element $(c,\ell,r)$ for finite Coxeter group $W$ with almost-positive roots $\Pge$, define involutions $\tau_\ell$ and $\tau_r$ from the set $\Pge$ to itself by
\begin{equation*}
\tau_\ell(\alpha) :=\begin{cases} \alpha &\text{ if } \alpha\in -\Pi_r \\ \ell(\alpha) &\text{ otherwise} \end{cases} \quad\text{ and }\quad \tau_r(\alpha) :=\begin{cases} \alpha &\text{ if } \alpha\in -\Pi_\ell \\ r(\alpha) &\text{ otherwise} \end{cases}.
\end{equation*}
The composition $\tau:=\tau_\ell\circ\tau_r$ is called a {\sf tropical Coxeter element}.
\end{definition}

What kind of map is this tropical Coxeter element? We have defined it only on the set $\Pge$, and clearly it does {\em not} extend linearly to the whole space $V$. However, it does extend to $V$ as a {\em piecewise-linear map}, and there is a precise way in which this can be considered a ``tropicalization'' of the Coxeter element $c$ (see \cite{fomin-zelevinsky:ysystems}). Here we will consider only the action of $\tau$ on the set $\Pge$; not on the whole space.

\begin{figure}
\vspace{.1in}
\begin{center}
\hspace{-.5in}\input{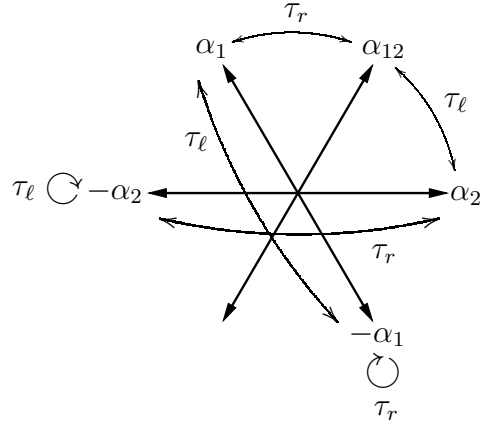}
\end{center}
\caption{The action of involutions $\tau_\ell$ and $\tau_r$ on the almost-positive roots of type $A_2$}
\label{fig:taus}
\end{figure}

Let us investigate the action of $\tau_\ell$ and $\tau_r$ on the root system of type $A_2$. Without loss of generality, we set $\Pi_\ell=\{\alpha_1\}$ and $\Pi_r=\{\alpha_2\}$. Then one can verify that the involutions $\tau_\ell$ and $\tau_r$ act on the set of almost-positive roots $\Pge=\{-\alpha_1,-\alpha_2,\alpha_1,\alpha_2,\alpha_{12}\}$ as in Figure \ref{fig:taus}. In particular, notice that $\tau_\ell$ and $\tau_r$ act as reflections on the pentagon in Figure \ref{fig:pentagon}, such that the composition $\tau=\tau_\ell\circ\tau_r$ gives a {\em counterclockwise rotation} ($\tau_\ell$ and $\tau_r$ almost act as reflections on the root system in Figure \ref{fig:taus}, but they are ``bent''). We conclude that the maps $\tau_\ell$ and $\tau_r$ generate the dihedral group of motions on the diagonals in Figure \ref{fig:pentagon}.

This is exactly what was desired, and we can now apply this to determine when two positive roots cross. For example, to determine whether $\alpha_1$ and $\alpha_2$ cross in Figure \ref{fig:pentagon}, we apply $\tau$ until one of them becomes negative; in this case, $\tau(\alpha_1)=\alpha_2$ and $\tau(\alpha_2)=-\alpha_2$. Since $\alpha_2$ and $-\alpha_2$ cross ($\alpha_2$ contains $\alpha_2$ in its simple root expansion), so do $\alpha_1$ and $\alpha_2$.

This inspires the following general definition.

\begin{definition}
\label{def:pgecross}
Given a finite Coxeter group $W$ with root system $\Phi=\Phi^+\sqcup\Phi^-$ and tropical Coxeter element $\tau$, we say that the almost-positive roots $\alpha,\beta\in\Pge$ {\sf cross} if there exists some $r$ such that $\tau^r(\alpha)$ is a negative root and $-\tau^r(\alpha)$ occurs in the simple root expansion of $\tau^r(\beta)$.
\end{definition}

Based on the above observations, this definition seems to work in type $A$, but for a general root system there is something to prove. Three things must be shown, and these were proved by Fomin and Zelevinsky \cite{fomin-zelevinsky:ysystems}. Although their paper discusses only crystallographic root systems, the result we state holds in general.

\begin{theorem}[\cite{fomin-zelevinsky:ysystems}]
Let $W$ be a finite Coxeter group with (possibly noncrystallographic) root system $\Phi=\Phi^+\sqcup\Phi^-$ and tropical Coxeter element $\tau$.
\begin{enumerate}
\item Any $\tau$-orbit on $\Pge$ has nonempty intersection with $-\Pi$.
\item Given $\alpha,\beta\in\Pge$, $\alpha$ crosses with $\beta$ if and only if $\beta$ crosses with $\alpha$.
\item Given $\alpha,\beta\in\Pge$, $\alpha$ crosses with $\beta$ if and only if $\tau(\alpha)$ crosses with $\tau(\beta)$. 
\end{enumerate}
\end{theorem}

The first part of the theorem guarantees that any two almost-positive roots may be compared; the second and third parts show that Definition \ref{def:pgecross} is well-defined. Of course, the definition of crossing depends on the particular choice of tropical Coxeter element $\tau$, but different choices turn out to give equivalent results. Finally, this allows us to define the Fomin-Zelevinsky generalization of the associahedron. Its formal name is the cluster complex.

\begin{definition}
Given a finite Coxeter group $W$, define $\Delta(W)$ as the flag simplicial complex whose faces are the mutually noncrossing subsets of the almost-positive roots $\Pge$. The maximal faces of $\Delta(W)$ are called {\sf clusters} and $\Delta(W)$ is called the {\sf cluster complex} corresponding to $W$.
\end{definition}

Thus, the cluster complex $\Delta(W)$ is a generalization of the simplicial associahedron to all finite Coxeter groups, where the elements of $\Pge$ correspond to ``diagonals'' and the clusters correspond to ``triangulations''. The complex also preserves the idea that  all triangulations of a given type should contain the same number of diagonals: in general, if $W$ has rank $n$ then each cluster of type $W$ contains $n$ almost-positive roots \cite[Theorem 1.8]{fomin-zelevinsky:ysystems}, hence the complex $\Delta(W)$ is {\sf pure} of dimension $n-1$. Fomin and Zelevinsky also showed that $\Delta(W)$ is the dual complex of a complete simplicial fan \cite[Corollary 1.11]{fomin-zelevinsky:ysystems}, and so it has the homotopy type of a sphere. One could define the dual sphere $\nabla(W)$ as the ``simple associahedron'' of type $W$, but since it adds nothing we will not do this.

The numerology of the cluster complex is quite remarkable. If we let $f_i$ denote the number of faces of $\Delta(W)$ with cardinality $i$ (and dimension $i-1$), then the sequence of face numbers $(f_0,f_1,\ldots,f_n)$ is called the {\sf $f$-vector} of the complex. Associated to this $f$-vector is another vector that is closely related to the topology of the complex. Define the {\sf $h$-vector} to be the sequence $(h_0,h_1,\ldots,h_n)$, where
\begin{equation*}
\sum_{i=0}^n f_i\,(x-1)^{n-i} = \sum_{i=0}^n h_i \,x^{n-i}.
\end{equation*}
In our case, we will always have $h_n=1$, which follows from the fact that $\Delta(W)$ is a sphere. For example, the $f$-vector of the type $A_3$ cluster complex is $(1,9,21,14)$, corresponding to the numbers of partial triangulations of a hexagon (counted by the number of diagonals), and the $h$-vector is $(1,6,6,1)$. Notice that the entries of the $h$-vector are the rank-numbers of the lattice of noncrossing partitions $NC(A_3)$ (Figure \ref{fig:nc4}). We set down one more notation.

\begin{definition}
We say that a cluster in $\Delta(W)$ is {\sf positive} if it contains only positive roots.
\end{definition}

The following results have all been proven case-by-case.

\begin{theorem}
Consider the cluster complex $\Delta(W)$ where $W$ is a finite Coxeter group with rank $n$.
\begin{enumerate}
\item The number of clusters in $\Delta(W)$ (maximal faces) is equal to the Coxeter-Catalan number $\Cat^{(1)}(W)$ \eqref{eq:fusscat}.
\label{fz1}
\item The number of positive clusters in $\Delta(W)$ is equal to the positive Coxeter-Catalan number $\Cat_+^{(1)}(W)$ \eqref{eq:posfusscat}. 
\label{fz2}
\item The entries of the $h$-vector of $\Delta(W)$ are the Coxeter-Narayana numbers. That is, $h_i=\Nar^{(1)}(W,n-i)$ \eqref{eq:fussnar}.
\label{fz3}
\end{enumerate}
\end{theorem}

\begin{proof}
\eqref{fz1} and \eqref{fz2} are Propositions  3.8 and 3.9 in \cite{fomin-zelevinsky:ysystems}, respectively. Fomin and Zelevinsky did not mention Narayana numbers in the paper \cite{fomin-zelevinsky:ysystems}, but \eqref{fz3} was probably verified soon after the paper appeared. A more general version can be found in Fomin and Reading \cite{fomin-reading} (see Theorem \ref{th:kcluster} below).
\end{proof}

By now, we have encountered these numbers several times. Recall that the Coxeter-Catalan number $\Cat^{(1)}(W)$ counts the elements of the noncrossing partition lattice $NC(W)$ and the antichains $NN(W)$ in the root poset (when $W$ is a Weyl group). The Narayana numbers are the rank numbers of the lattice $NC(W)$ and they count antichains $NN(W)$ by cardinality. The appearance of the positive Coxeter-Catalan number $\Cat_+^{(1)}(W)$ here explains the notation, since this number counts the number of ``positive'' clusters in $\Delta(W)$. Recall that $\Cat_+^{(1)}(W)$ is also the number of bounded positive regions in the Shi arrangement (see Theorem \ref{th:ath-shi}), or equivalently, the number of antichains in $NN(W)$ consisting of non-simple roots. $\Cat_+^{(1)}(W)$ also counts the number of homotopy spheres in the deleted order complex of $NC(W)$ (Theorem \ref{th:NChomotopy}).

Furthermore, the emergence of these numbers in a new context allows a new chance for understanding. Since $\Delta(W)$ is a simplicial homotopy sphere, the Dehn-Sommerville relations (or Poincar\'e duality) tell us that the $h$-vector is symmetric: $h_i=h_{n-i}$ for all $0\leq i\leq n$. For the case of $NC(W)$, this is also clear since the poset is self-dual and hence rank-symmetric. However, the symmetry of the Narayana numbers is {\em not at all clear} from the antichains $NN(W)$. In this case, symmetry implies that the number of antichains with $i$ elements is equal to the number of antichains with $n-i$ elements, and {\em no bijective proof of this fact is currently known} (see Conjecture \ref{conj:floors=ceilings} in the case $k=1$).

Finally, in the case that $W$ is a Weyl group, Chapoton, Fomin and Zelevinsky \cite{chapoton-fomin-zelevinsky} have constructed a polytopal realization of the cluster complex $\Delta(W)$. Combined with the $g$-theorem, this then implies that the $h$-vector is unimodal. This is the {\em only known method} that can be used to prove unimodality in a uniform way (see Theorem \ref{th:unimodal}). The following expansion of Problem \ref{prob:NC=NN} is of central importance to this rapidly growing field.

\begin{problem}
What is the nature of the relationships between the objects $NC(W)$, $NN(W)$ and $\Delta(W)$? What is an appropriate notion of ``type'' for the clusters in $\Delta(W)$? Find bijections between these objects that preserve type and that respect the positive Catalan numbers.
\end{problem}

\begin{update}[July 2007]
We should mention some recent progress towards this goal. While the general problem remains open, there have been two successful studies linking the noncrossing partitions $NC(W)$ with the cluster complex $\Delta(W)$. In both of these studies, a notion of parabolic ``type'' for clusters is evident.

$\bullet$
In the same paper in which they first proved the lattice property for $NC(W)$, Brady and Watt  \cite{brady-watt:lattice} showed how to construct the complex $\Delta(W)$ purely in terms of noncrossing partitions. Based on this work, Athanasiadis, Brady and Watt then constructed a shelling \cite{athanasiadis-brady-watt} for the lattice $NC(W)$. Although \cite{brady-watt:lattice} did not immediately explain the numerological coincidences between $NC(W)$ and $\Delta(W)$, Athanasiadis, Brady, McCammond and Watt have followed \cite{brady-watt:lattice} with a study \cite{ABMW} in which they give a uniform bijective correspondence between $NC(W)$ and $\Delta(W)$ that preserves the main Catalan statistics.

$\bullet$
Reading has independently constructed a uniform bijection \cite{reading:sortable} between noncrossing partitions and clusters. To do this, he defined, for each Coxeter element $c\in W$, the notion of {\sf $c$-sortable elements} in the group $W$ (these generalize the classical  {\sf stack-sortable permutations}); he  then provided bijections from these elements to the elements of $NC(W)$ and to the facets of $\Delta(W)$ (the ``clusters''). This work grew out of Reading's thesis work on {\sf Cambrian lattices} \cite{reading} (these generalize the classical {\sf Tamari lattices}). Reading and Speyer have recently given a complete account \cite{reading-speyer} of the relationships between the Cambrian lattices, $c$-sortable elements and cluster complexes.

The relationship between noncrossing partitions $NC(W)$ and nonnesting partitions $NN(W)$ remains a complete mystery.

\end{update}

\subsection{The Generalized Cluster Complex}
\label{sec:generalcluster}
Finally, we consider a Fuss-Catalan generalization of the cluster complex. Recently, Fomin and Reading defined a simplicial complex $\Delta^{(k)}(W)$ for each finite Coxeter group $W$ and positive integer $k$ that generalizes the Fomin-Zelevinsky cluster complex $\Delta(W)$. We will see that this complex is a combinatorial analogue of the geometric $k$-multichains of filters $NN^{(k)}(W)$ defined by Athanasiadis (Section \ref{sec:nonnesting}) and our poset of $k$-divisible noncrossing partitions $NC^{(k)}(W)$. Unlike the original work on $NN^{(k)}(W)$ and $NC^{(k)}(W)$, Fomin and Reading had in mind the specific goal of constructing a complex with the ``correct'' Fuss-Catalan numerology.

The type $A$ version of the generalized cluster complex has a natural interpretation in terms of higher polygon dissections. Note that a convex polygon with $N$ vertices can be dissected into $(k+2)$-gons if and only if $N=kn+2$ for some positive integer $n$. In this case, the dissection will consist of $n$ $(k+2)$-gons. Recall from Section \ref{sec:fusscat} that Euler's mathematical assistant and grandson-in-law Niklaus Fuss proved the following result in 1791.

\begin{theorem}[\cite{fuss}]
For all positive integers $n$ and $k$, the number of dissections of a convex $(kn+2)$-gon into $(k+2)$-gons is equal to the Fuss-Catalan number
\begin{equation}
\label{eq:fuss}
\Cat^{(k)}(n)=\frac{1}{n}\binom{(k+1)n}{n-1}.
\end{equation}
\end{theorem}

For example, this says that there are $\Cat^{(2)}(3)=12$ ways to dissect a convex octagon into quadrilaterals (``quadrangulations of an octagon''). In general, we will say that a collection of diagonals in a convex $(kn+2)$-gon is a {\sf $(k+2)$-angulation} if it dissects the polygon into $(k+2)$-gons. Hence, a $(k+2)$-angulation of a $(kn+2)$-gon must contain $n-1$ diagonals. This suggests an interesting idea. Perhaps the definition of the simplicial associahedron $\Delta(n)$ can be mimicked in order to construct a flag complex on the set of diagonals of a convex $(kn+2)$-gon, each of whose maximal faces is a $(k+2)$-angulation. This idea is on the right track, but we run into a problem since not all diagonals of the $(kn+2)$-gon occur in a $(k+2)$-angulation. In order to make this work, we must restrict our notion of a ``diagonal''.

\begin{definition}
A diagonal of a convex $(kn+2)$-gon is called a {\sf $k$-divisible diagonal} if it may participate in a $(k+2)$-angulation. That is, a $k$-divisible diagonal dissects the larger polygon into two polygons each of which has a number of vertices congruent to $2$ modulo $k$.
\end{definition}

Thus, a polygon will possess a $k$-divisible diagonal if and only if it has $kn+2$ vertices for some $n$. We say that a collection of mutually noncrossing $k$-divisible diagonals is a {\sf partial $(k+2)$-angulation}. (This has also been called a {\sf $k$-divisible dissection}; see Tzanaki \cite{tzanaki}).

Part of the motivation for the recent study of these objects was the paper of Przytycki and Sikora \cite{p-s:fuss}, in which they gave a simple bijective proof \cite[Corollary 0.2]{p-s:fuss} of the classical result of Kirkman \cite{kirkman} and Cayley \cite{cayley}: the number of partial $(k+2)$-angulations of a convex $(kn+2)$-gon containing $i$ diagonals is equal to
\begin{equation*}
\frac{1}{i+1}\binom{n-1}{i}\binom{kn+i+1}{i}.\footnote{This is the ``more general problem'' mentioned in Section \ref{sec:associahedron}.}
\end{equation*}
Przytycki and Sikora's paper also contains several notes on the history of polygon dissections.

This leaves us with a suitable definition of a ``generalized'' associahedron. The following definition was suggested by Reiner (see \cite{tzanaki}).

\begin{definition}
Let $\Delta^{(k)}(n)$ denote the abstract (flag) simplicial complex whose vertices are the $k$-divisible diagonals of a convex $(kn+2)$-gon, and whose faces are partial $(k+2)$-angulations. We will call this the {\sf simplicial $k$-divisible associahedron}.
\end{definition}

\begin{figure}
\vspace{.1in}
\begin{center}
\scalebox{0.8}{\input{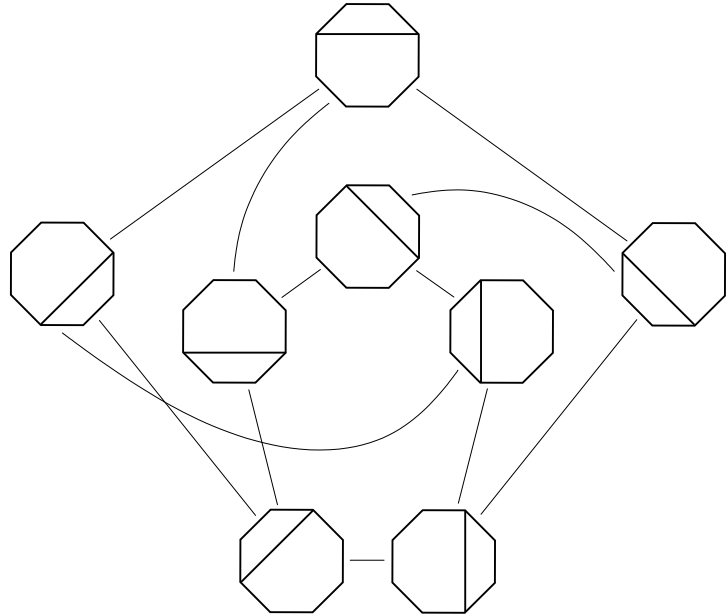}}
\end{center}
\caption{The $2$-divisible associahedron $\Delta^{(2)}(3)$ of $2$-divisible dissections of a convex octagon}
\label{fig:2-associahedron}
\end{figure}

Here we use the term ``associahedron'' to suggest an analogy with the complex $\Delta(n)$, and {\em not} to suggest that $\Delta^{(k)}(n)$ is a polytope. In general, it is not even a sphere\footnote{For this reason, it is not at all clear how to define a corresponding {\em simple} $k$-divisible associahedron $\nabla^{(k)}(n)$.}.
For example, Figure \ref{fig:2-associahedron} displays the $2$-divisible associahedron $\Delta^{(2)}(3)$ of partial quadrangulations of a convex octagon. Notice that the complex has $\frac{1}{2}\binom{2}{1}\binom{8}{1}=8$ vertices, corresponding to the $2$-divisible diagonals, and $\Cat^{(2)}(3)=12$ edges, corresponding to the quadrangulations; thus its $f$-vector is $(1,8,12)$. This implies that the $h$-vector of this complex is $(1,6,5)$, which is exciting since these are the Fuss-Narayana numbers corresponding to $k=2$ and $W=A_2$. Furthermore, one can see (for instance, from the Euler characteristic) that the complex $\Delta^{(2)}(3)$ is homotopy equivalent to a wedge of $\Cat_+^{(1)}(3)=5$ circles. We will see that these suggestive observations are not coincidental.

At the suggestion of Reiner, Tzanaki was the first to study this complex, as well as a type $B$ analogue. In the paper \cite{tzanaki}, she proved that the complex $\Delta^{(k)}(W)$, where $W=A_{n-1}$ or $B_n$, is shellable; its $h$-vector is given by the Fuss-Narayana numbers $\Nar^{(k)}(W,i)$ in Figure \ref{fig:fussnar}; and it is homotopic to a wedge of $\Cat_+^{(k-1)}(W)$ many spheres. Soon after, Fomin and Reading \cite{fomin-reading} demonstrated how to extend the definition of $\Delta^{(k)}(W)$ to all finite Coxeter groups and positive integers $k$.

The key to the type $A$ version --- as with the case of the cluster complex $\Delta(A_{n-1})$ --- is to label the
\begin{equation*}
\frac{1}{2}\binom{n-1}{1}\binom{kn+2}{1}=(n-1)+k\binom{n}{2}
\end{equation*}
$k$-divisible diagonals of a convex $(kn+2)$-gon by a suitable set of almost-positive roots; in this case the set of {\em colored} almost-positive roots.

\begin{definition}
Given a noncrystallographic root system $\Phi=\Phi^+\sqcup\Phi^-$, define the set of {\sf colored almost-positive roots} $\Pge^{(k)}$ to consist of $k$ distinct copies of the positive roots $\Phi^+$ (labelled with superscripts $1,2,\ldots,k$), together with a unique copy of the negative simple roots $-\Pi$ (labelled with superscript $1$):
\begin{equation*}
\Pge^{(k)}:=\left\{ \alpha^i :\alpha\in\Phi^+, i\in\{1,\ldots,k\}\right\}\cup\left\{ \alpha^1:\alpha\in-\Pi\right\}.
\end{equation*}
\end{definition}

\begin{figure}
\vspace{.1in}
\begin{center}
\input{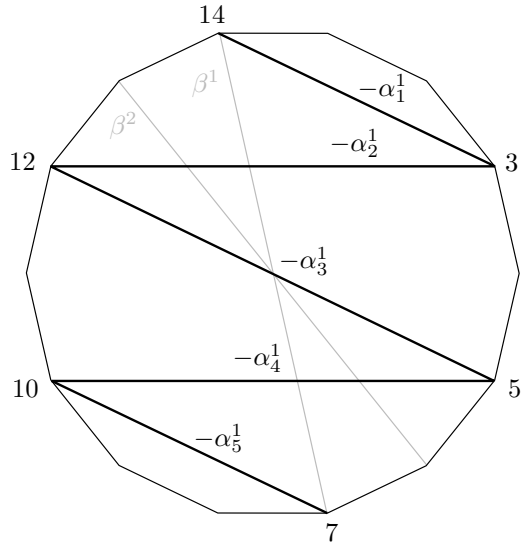}
\end{center}
\caption{The $2$-snake of type $A_5$}
\label{fig:2-snake}
\end{figure}

Notice that the set $\Pge^{(k)}$ of type $A_{n-1}$ contains $(n-1)+k\binom{n}{2}$ elements, as desired. Fomin and Reading defined a labelling of the $k$-divisible diagonals of a convex $(kn+2)$-gon by the colored almost-positive roots $\Pge^{(k)}$ of type $A_{n-1}$:

Again, let $\Pi=\{\alpha_1,\alpha_2,\ldots,\alpha_{n-1}\}$ denote the simple roots of type $A_{n-1}$, and consider a convex $(kn+2)$-gon with vertices labelled clockwise by $1,2,\ldots,kn+2$. Generalizing the snake of Figure \ref{fig:snake2}, we label the diagonal connecting vertices $\lceil\frac{i}{2}\rceil k+1$ and $k(n-\lfloor\frac{i}{2}\rfloor)+2$ by the negative simple root $-\alpha_i^1$, for all $1\leq i\leq n-1$. This collection of $k$-divisible diagonals is called the {\sf $k$-snake}; notice that it forms a $(k+2)$-angulation. Figure \ref{fig:2-snake} shows the $2$-snake of type $A_5$, where $k=2$ and $n=6$.

To label the remaining $k\binom{n}{2}$ $k$-divisible diagonals, notice that if $\beta=\sum_{i\in I}\alpha_i$ is a positive root for some $I\subseteq[n-1]$, then there are exactly $k$ distinct $k$-divisible diagonals that cross with the set of snake diagonals labelled by $\{-\alpha_i^1:i\in I\}$, and these form a contiguous sequence under rotation of the $(kn+2)$-gon. We label these diagonals as $\beta^1,\beta^2,\ldots,\beta^k$, {\em counterclockwise}. For example, in Figure \ref{fig:2-snake}, we have displayed the diagonals $\beta^1$ and $\beta^2$, where $\beta$ is the positive root $\alpha_2+\alpha_3+\alpha_4$.

Now, in order to construct a generalized cluster complex, it remains to find a ``rotation map'' on the set $\Pge^{(k)}$, with which we might define the notion of ``noncrossing''. This is accomplished by the following \cite[Definition 3.3]{fomin-reading}.

\begin{definition}[\cite{fomin-reading}]
\label{def:tropical*}
Let $\tau$ be a tropical Coxeter element as in Definition \ref{def:tropical}. Then for all $\alpha^i\in\Pge^{(k)}$, we set
\begin{equation*}
\tau^*(\alpha^i):=\begin{cases} \alpha^{i+1} & \text{ if } \alpha\in\Phi^+ \text{ and } i<k \\ (\tau(\alpha))^i & \text{ otherwise }\end{cases}
\end{equation*}
\end{definition}

\begin{figure}
\vspace{.1in}
\begin{center}
\scalebox{0.8}{\input{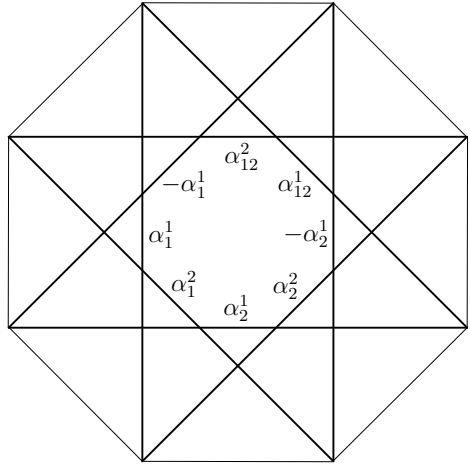}}
\end{center}
\caption{The collection of $2$-divisible diagonals of an octagon, labelled by the colored almost-positive roots $\Pge^{(2)}$ of type $A_2$}
\label{fig:2-diagonals}
\end{figure}

To see how the map $\tau^*$ works, consider the colored almost-positive roots $\Pge^{(2)}$ of type $A_2$. These label the $2$-divisible diagonals of an octagon, as shown in Figure \ref{fig:2-diagonals}. Here, we suppose again that $\Pi_\ell=\{\alpha_1\}$ and $\Pi_r=\{\alpha_2\}$. One can see that $\tau^*$ acts precisely as $1/8$ of a counterclockwise rotation on the octagon. For example, recall that $\tau$ is the counterclockwise rotation on the pentagon in Figure \ref{fig:pentagon}, so starting with the diagonal $-\alpha_1^1$, we get $\tau^*(-\alpha_1^1)=(\tau(-\alpha_1))^1=\alpha_1^1$. Then, since $1<k=2$ we have $\tau^*(\alpha_1^1)=\alpha_1^2$, and then $\tau^*(\alpha_1^2)=(\tau(\alpha_1))^1=\alpha_2^1$. In general, considering Figure \ref{fig:2-snake}, one observes that the map $\tau^*$ is designed to act as a counterclockwise rotation in general.

It is interesting that Fomin and Reading did not study the full dihedral group of motions on the $(kn+2)$-gon. This could easily be done by defining ``reflection'' maps $\tau_\ell^*$ and $\tau_r^*$ generalizing $\tau_\ell$ and $\tau_r$ in Definition \ref{def:tropical}, such that $\tau_\ell^*\circ\tau_r^*=\tau^*$. However, it turns out that these maps are complicated to write down, and it is probably not worth the trouble unless one is going to use them for something. Notice that the dihedral group $\langle\tau_\ell^*,\tau_r^*\rangle$ acting on the set $\Pge^{(k)}$ is strongly analogous to the the dihedral group $\langle L^*,R^*\rangle$ that acts on the $k$-divisible noncrossing partitions $NC^{(k)}(W)$, developed in Section \ref{sec:automorphisms}.

It is clear from Figure \ref{fig:2-snake} that a diagonal $\beta^j$ crosses with a diagonal of the form $-\alpha_i^1$ if and only if $\alpha_i$ occurs in the simple root expansion of $\beta$. Generalizing Definition \ref{def:pgecross}, we define a binary relation --- called``crossing'' --- on the set $\Pge^{(k)}$.

\begin{definition}
\label{def:pgekcross}
Given a finite Coxeter group $W$ with root system $\Phi=\Phi^+\sqcup\Phi^-$, a positive integer $k$, and tropical Coxeter element $\tau$, we say that the colored almost-positive roots $\alpha^i,\beta^j\in\Pge^{(k)}$ {\sf cross} if there exists some integer $r$ such that $(\tau^*)^r(\alpha^i)=(\alpha')^1$ and $(\tau^*)^r(\beta^j)=(\beta')^{j'}$, where $\alpha'$ is a negative root and $-\alpha'$ occurs in the simple root expansion of $\beta'$.
\end{definition}

Again, this is well-defined, as shown by Fomin and Reading.

\begin{theorem}[\cite{fomin-reading}]
Given the setup in Definition \ref{def:pgekcross}, we have
\begin{enumerate}
\item Any $\tau^*$-orbit on $\Pge^{(k)}$ has nonempty intersection with $-\Pi^1$.
\item Given $\alpha^i,\beta^j\in\Pge^{(k)}$, $\alpha^i$ crosses with $\beta^j$ if and only if $\beta^j$ crosses with $\alpha^i$.
\item Given $\alpha^i,\beta^j\in\Pge^{(k)}$, $\alpha^i$ crosses with $\beta^j$ if and only if $\tau^*(\alpha^i)$ crosses with $\tau^*(\beta^j)$.
\end{enumerate}
\end{theorem}

Thus, finally, we have the ``generalized cluster complex'' \cite[Definition 3.8]{fomin-reading}. Again, the isomorphism type is independent of the tropical Coxeter element used.

\begin{definition}[\cite{fomin-reading}]
For a finite Coxeter group $W$, let $\Delta^{(k)}(W)$ denote the flag simplicial complex whose faces are the mutually noncrossing subsets of colored almost-positive roots $\Pge^{(k)}$. The maximal faces of $\Delta^{(k)}(W)$ are called {\sf colored clusters}, and we call $\Delta^{(k)}(W)$ the {\sf $k$-divisible cluster complex} corresponding to $W$.
\end{definition}

Fomin and Reading did not name the complex $\Delta^{(k)}(W)$; we use the term ``$k$-divisible'' to connect with the notions of $k$-divisible polygon dissections and $k$-divisible noncrossing partitions. We see that $\Delta^{(k)}(A_{n-1})$ is a pure complex of dimension $n-2$, since each $(k+2)$-angulation contains $n-1$ diagonals. In general, Fomin and Reading showed that $\Delta^{(k)}(W)$ is a pure complex of dimension $n-1$ when $W$ has rank $n$ \cite[Theorem 3.9]{fomin-reading}.

Let us now examine the complex $\Delta^{(2)}(A_2)$, displayed in Figure \ref{fig:2-associahedron}. Comparing with Figure \ref{fig:2-diagonals}, we see how the vertices of this complex correspond to the set $\Pge^{(2)}$. The two pentagons indicated correspond to two embedded copies of the cluster complex $\Delta(A_2)$, with vertex sets
\begin{equation*}
\{-\alpha_1^1,-\alpha_2^1,\alpha_1^1,\alpha_2^1,\alpha_{12}^1\}\quad\text{ and }\quad\{-\alpha_1^1,-\alpha_2^1,\alpha_1^2,\alpha_2^2,\alpha_{12}^2\}.
\end{equation*}
In general, the set $-\Pi^1\cup(\Phi^+)^i$ consisting of the negative simple roots and the positive $i$-colored roots is an embedded copy of $\Delta(W)$ inside $\Delta^{(k)}(W)$, for all $1\leq i\leq k$; and, moreover, if $I\subseteq[k]$ is a set of colors, then the subcomplex of $\Delta^{(k)}(W)$ induced by the vertices $-\Pi^1\cup\left( \cup_{i\in I}(\Phi^+)^i\right)$ is isomorphic to $\Delta^{(\abs{I})}(W)$ \cite[Corollary 3.6]{fomin-reading}. Also, many more subcomplexes of $\Delta^{(k)}(W)$ can be obtained by applying the symmetry group $\langle \tau_\ell^*,\tau_r^*\rangle$ to those above.
It is interesting how the complex $\Delta^{(k)}(W)$ is built out of complexes $\Delta^{(i)}(W)$ for $1\leq i\leq k$, which overlap in a very symmetric way. Compare this phenomenon to Lemma \ref{lemma:contained}, which describes embeddings of $k$-divisible noncrossing partition posets.

Finally, and perhaps most importantly, Fomin and Reading observed that $\Delta^{(k)}(W)$ is truly a Fuss-Catalan object, in our sense of the term.

\begin{definition}
We say a colored cluster in $\Delta^{(k)}(W)$ is {\sf positive} if it contains only positive colored roots.
\end{definition}

\begin{theorem}
\label{th:kcluster}
Consider the $k$-divisible cluster complex $\Delta^{(k)}(W)$, where $W$ is a finite Coxeter group and $k$ is a positive integer.
\begin{enumerate}
\item The number of maximal faces (colored clusters) in $\Delta^{(k)}(W)$ is equal to the Fuss-Catalan number $\Cat^{(k)}(W)$ \eqref{eq:fusscat}.
\label{item:kcluster1}
\item The number of positive colored clusters in $\Delta^{(k)}(W)$ is equal to the positive Fuss-Catalan number $\Cat_+^{(k)}(W)$ \eqref{eq:posfusscat}.
\label{item:kcluster2}
\item The entries of the $h$-vector of $\Delta^{(k)}(W)$ are the Fuss-Narayana numbers. That is, $h_i=\Nar^{(k)}(W,n-i)$ \eqref{eq:fussnar}.
\label{item:kcluster3}
\end{enumerate}
\end{theorem}

\begin{proof}
\eqref{item:kcluster1} and \eqref{item:kcluster2} appear as Propositions 8.4 and 12.4, respectively, in \cite{fomin-reading}, and the $h$-vector \eqref{item:kcluster3} was computed in Theorem 10.2 and displayed in Table 3 of \cite{fomin-reading}. The numbers can be observed case-by-case to coincide with the Fuss-Narayana numbers in our Figure \ref{fig:fussnar}.
\end{proof}

Our running example of $k=2$ and $W=A_2$ is now complete. We invite the reader at this point to compare Figures \ref{fig:nc^2(a2)}/\ref{fig:nc^2(3)}, \ref{fig:shi^2(a2)}/\ref{fig:nn^2(a2)} and \ref{fig:2-associahedron}/\ref{fig:2-diagonals}, and to observe the striking relationships between the objects $NC^{(2)}(A_2)$, $NN^{(2)}(A_2)$ and $\Delta^{(2)}(A_2)$. The following problem is the climax of this memoir. When we speak of the subject of ``Fuss-Catalan combinatorics'', we are referring to this:

\begin{problem}
\label{prob:fusscat}
Explain the relationships between
\begin{equation*}
NC^{(k)}(W),\quad NN^{(k)}(W)\quad\text{ and }\quad \Delta^{(k)}(W).
\end{equation*}
\end{problem}

This completes our exposition of the Fuss-Catalan combinatorics. To end this section, and the memoir, we give a few hints as to where Problem \ref{prob:fusscat} may lead.

\subsection{Addendum --- Relating $NC^{(k)}(W)$ and $\Delta^{(k)}(W)$}

The notion of a $k$-Cohen-Macaulay complex was defined by Baclawski \cite{baclawski}: a simplicial complex $\Delta$ is said to be {\sf $k$-Cohen-Macaulay} if the subcomplex obtained from $\Delta$ by deleting any $k-1$ of its vertices remains Cohen-Macaulay in the usual sense. Athanasiadis and Tzanaki \cite{athanasiadis-tzanaki} showed that the $k$-divisible cluster complex and its subcomplex of positive clusters are both $(k+1)$-Cohen-Macaulay.

\begin{definition}
Let $\Delta_+^{(k)}(W)$ denote the subcomplex of $\Delta^{(k)}(W)$ induced by the set $\cup_{i=1}^k (\Phi^+)^i$ of colored positive roots. This is called the {\sf positive $k$-divisible cluster complex}.
\end{definition}

\begin{theorem}[\cite{athanasiadis-tzanaki}]\hspace{.1in}
\begin{enumerate}
\item The complex $\Delta^{(k)}(W)$ with any $k$ vertices deleted is shellable and hence Cohen-Macaulay.
\item The complex $\Delta_+^{(k)}(W)$ with any $k$ vertices deleted is shellable and hence Cohen-Macaulay.
\end{enumerate}
\end{theorem}

In particular, the complex $\Delta_+^{(k)}(W)$ is homotopy equivalent to a wedge of $(n-1)$-dimensional spheres. Considering Theorem \ref{th:kcluster} \eqref{item:kcluster3}, perhaps these spheres can be enumerated. Athanasiadis and Tzanaki did this, and also related the result to our deleted order complex $\Delta\left(NC_{(k)}(W)\setminus\{(1)_k\}\right)$.

\begin{theorem}[\cite{athanasiadis-tzanaki}]
\label{th:homotopy}
\hspace{.1in}
\begin{enumerate}
\item The complex $\Delta_+^{(k)}(W)$ is homotopy equivalent to a wedge of $\Cat_+^{(k-1)}(W)$ many $(n-1)$-dimensional spheres.
\item The complex $\Delta_+^{(k)}(W)$ is homotopy equivalent to the order complex \break $\Delta\left(NC_{(k)}(W)\setminus\{(1)_k\}\right)$.
\end{enumerate}
\end{theorem}

This gives an alternate proof of our Theorem \ref{th:NCkhomotopy}. However, our result that the complex $\Delta\left(NC_{(k)}(W)\cup\{\hat{1}\}\right)$ is shellable (Theorem \ref{th:ELNCk}) is stronger than Theorem \ref{th:homotopy}, which computes only the homotopy type.

We expect that this work will continue, and that a relationship between the $k$-divisible noncrossing partitions $NC^{(k)}(W)$ and the $k$-divisible cluster complex $\Delta^{(k)}(W)$ analogous to that described in \cite{ABMW} (for $k=1$) will be attained.

\section{Chapoton Triangles}
\label{sec:triangles}

Some of the most remarkable enumerative conjectures in this field have been made by Fr\'ed\'eric Chapoton. In the two papers {\em Enumerative properties of generalized associahedra} \cite{chapoton:one} and {\em Sur le nombre de r\'eflexions pleines dans les groupes de Coxeter finis} \cite{chapoton:two}, he observed several interesting numerological features of the Coxeter-Catalan objects $NC(W)$, $NN(W)$ and $\Delta(W)$. In particular, he defined a family of two-variable generating functions, which he called ``triangles''. Here, we will extend Chapoton's definitions and conjectures to the case of the Fuss-Catalan combinatorics.

For now, we suppose that $W$ is a finite Weyl group with rank $n$ and bipartite Coxeter element $(c,\ell,r)$ (see Section \ref{sec:automorphisms}). Recall that $NC_{(k)}(W)$ is the graded meet-semilattice of $k$-delta sequences. This poset has height $n$ with a unique minimum element $(1)_k\in W$, $\Cat^{(k-1)}(W)$ maximal elements, and $\Cat^{(k)}(W)$ elements in total. Also, the deleted order complex of this poset is homotopy equivalent to a wedge of $\Cat_+^{(k-1)}(W)$ many $(n-1)$-dimensional spheres.

In general, if $P$ is a poset, there is an important function $\mu$ from the set of intervals $\Int(P)$ to the integers, called the {\sf M\"obius function} of the poset. It is defined recursively by setting
\begin{align*}
\mu(x,x) & := 1 \quad\text{ for all } x\in P, \text{ and }\\
\mu(x,y) & := -\sum_{x<z<y} \mu(x,z) \quad\text{ for all } x<y \text{ in } P.
\end{align*}
The most important feature of this function is the {\sf M\"obius inversion formula} (see \cite[Chapter 3.7]{stanley:ec1}): if $f$ and $g$ are functions $P\to\integers$, then the two formulas
\begin{equation*}
g(x)=\sum_{y\leq x} f(y) \quad\text{ and }\quad f(x)=\sum_{y\leq x} g(y)\,\mu(x,y)
\end{equation*}
are equivalent. In the case that $P$ is a graded poset with minimum element $\hat{0}$ and height $n$, the {\sf characteristic polynomial} of $P$ is defined as
\begin{equation*}
\chi(P,q):=\sum_{x\in P} \mu(\hat{0},x)\,q^{n-\rk(x)}.
\end{equation*}
In particular, if $P$ is the intersection lattice of a hyperplane arrangement, then $\chi(P,q)$ is called the characteristic polynomial of the arrangement.

Recall that $NN^{(k)}(W)$ is the set of geometric $k$-multichains of filters in the root poset $(\Phi^+,\leq)$, or equivalently the set of positive chambers in the extended Shi arrangement $\Shi^{(k)}(W)$. If $\filter$ is an element of $NN^{(k)}(W)$, it has a corresponding set $\FL_i(\filter)\in\Phi^+$ of $i$-colored floors and a set $\CL_i(\filter)\in\Phi^+$ of $i$-colored ceilings. Recall that $NN^{(k)}(W)$ has $\Cat^{(k)}(W)$ elements, $\Cat^{(k-1)}(W)$ elements with $n$ $k$-colored floors and $\Cat_+^{(k)}(W)$ elements corresponding to bounded chambers. Given a set $A\in\Phi^+$, let $S(A)=A\cap\Pi$ denote the set of simple roots in $A$.

Recall that $\Delta^{(k)}(W)$ is the $k$-divisible cluster complex on the set of colored almost-positive roots $\Pge^{(k)}$. Given $A\in\Delta^{(k)}(W)$, let $A=A^+\sqcup A^-$ denote its decomposition into positive colored roots $A^+$ and negative colored roots $A^-$.

Now, for each of the objects $NC_{(k)}(W)$, $NN^{(k)}(W)$ and $\Delta^{(k)}(W)$, we define a two variable generating function.

\begin{definition}
\label{def:triangles}
For each Weyl group $W$ and positive integer $k$, define the $M$-triangle, $H$-triangle and $F$-triangle as follows:
\begin{align*}
& M(x,y):= \sum_{\pi\leq \mu\,\in\,NC_{(k)}(W)} \mu(\pi,\mu)\cdot x^{n-\rk(\pi)}\cdot y^{n-\rk(\mu)} \\
& H(s,t):= \sum_{\filter\in NN^{(k)}(W)} s^{\abs{S(FL_k(\filter))}}\,\cdot\,t^{\abs{FL_k(\filter)}} \\
& F(p,q):= \sum_{A\in\Delta^{(k)}(W)} p^{\abs{A^+}}\,\cdot\,q^{\abs{A^-}}
\end{align*}
\end{definition}

In the case $k=1$, these generating functions were defined by Chapoton. 

A note on notation: $M$ here stands for ``M\"obius'', $H$ stands for ``height'' and $F$ stands for ``face''. The variables $s$ and $t$ correspond to ``simple'' roots and ``total'' roots, respectively. The variables $p$ and $q$ stand for ``positive'' and ``negative''. We hope that this notation will facilitate intuition.

Clearly each of the polynomials $M(x,y)$, $H(s,t)$ and $F(p,q)$ has degree $n$. The notation ``triangle'' indicates that all monomials with $x$-degree greater than $y$-degree; $s$-degree greater than $t$-degree; or sum of $p$- and $q$-degrees greater than $n$ have coefficient 0. That is, if we arrange the coefficients in an $n\times n$ matrix, the possibly-nonzero entries form a triangle. The $M$, $H$ and $F$ triangles encode all known occurrences of the Fuss-Catalan, Fuss-Narayana and positive Fuss-Catalan numbers in these three families of objects.

To follow our favorite example, consider $k=2$ and $W=A_2$. Observing Figures \ref{fig:nc_2(a2)}, \ref{fig:nn^2(a2)} and \ref{fig:2-associahedron}, the $M$, $H$ and $F$ triangles are given by
\begin{align*}
& M(x,y)= 5-12y+7y^2+6xy-6xy^2+x^2y^2, \\
& H(s,t)=  5+2t+4st+s^2t^2, \\
& F(p,q)= 1+6p+2q+7p^2+4qp+q^2,
\end{align*}
with corresponding matrices
\begin{equation*}
M=\left(\!\!\!\!\begin{array}{rrr} 5&& \\ -12&6& \\ 7&-6&1\end{array}\!\!\right), \,\,\,H=\begin{pmatrix} 5&& \\ 2&4& \\ 0&0&1 \end{pmatrix}, \,\,\,F=\begin{pmatrix} 1&6&7\\ 2&4& \\ 1&& \end{pmatrix}.
\end{equation*}
Notice, for example, that the bottom left corner of the $M$-triangle is $\Cat_+^{(2)}(A_2)=7$, and the diagonal entries of the $M$-triangle are the Fuss-Narayana numbers. The sum of the first column of the $H$-triangle is $\Cat_+^{(2)}(A_2)$, and the sums of its rows are Fuss-Narayana. The top right corner of the $F$-triangle is $\Cat_+^{(2)}(A_2)$ and the sum of its diagonal entries is $\Cat^{(2)}(A_2)=12$. There is really a lot going on here; the raison d'\^etre of the triangles is the following conjecture.

\begin{conjecture}
\label{conj:triangles}
For all rank $n$ Weyl groups $W$ and positive integers $k$, the following three equivalent statements hold:
\begin{align*}
&\textstyle{M(x,y) = (xy-1)^n F\left(\frac{1-y}{xy-1},\frac{1}{xy-1}\right) = (1-y)^n H\left(\frac{x}{x-1},\frac{y(x-1)}{1-y}\right)},\\
\vspace{.1in}
&\textstyle{H(s,t) = (1+(s-1)t)^n M\left(\frac{s}{s-1},\frac{(s-1)t}{1+(s-1)t}\right) = (t-1)^n F\left(\frac{1}{t-1},\frac{1+(s-1)t}{t-1}\right)},\\
\vspace{.5in}
&\textstyle{F(p,q) = q^n M\left(\frac{1+q}{q-p}, \frac{q-p}{q}\right) = p^n H\left(\frac{1+q}{1+p},\frac{1+p}{p}\right)}.
\end{align*}
\end{conjecture}

This conjecture underscores the numerological coincidences between the three Fuss-Catalan families.

There is also a natural ``duality'' for the triangles. In Definition \ref{def:triangles}, we have defined the $M$-triangle $M(x,y)$ in terms of the poset $NC_{(k)}(W)$ of $k$-delta sequences, but we could just as easily have used the dual poset $NC^{(k)}(W)$. Since $NC_{(k)}(W)$ is locally self-dual (Section \ref{sec:shifting}), it is easy to see that
\begin{equation*}
\sum_{\pi\leq \mu\,\in\,NC^{(k)}(W)} \mu(\pi,\mu)\cdot x^{n-\rk(\pi)}\cdot y^{n-\rk(\mu)} = (xy)^n\cdot M\left(\frac{1}{y},\frac{1}{x}\right).
\end{equation*}
Hence we might regard $(xy)^n\cdot M(1/y,1/x)$ as a ``dualization'' of the $M$-triangle. Transferring this duality to the $H$- and $F$-triangles via Conjecture \ref{conj:triangles}, we define the {\sf dual triangles}.

\begin{definition}
\label{def:dualtriangles}
Define the {\sf dual} $M$-, $H$- and $F$-triangle as follows:
\begin{align*}
&M^*(x,y) := (xy)^n M\left(\frac{1}{y}\,,\frac{1}{x}\right), \\
&H^*(s,t) := t^n H\left(1+(s-1)t\,, \frac{1}{t}\right),\\
&F^*(p,q) := (-1)^n  F\left(-1-p,-1-q\right).
\end{align*}
\end{definition}

In the case $k=1$, Conjecture \ref{conj:triangles} is the combination of Conjecture 1 in \cite{chapoton:one} and Conjecture 6.1 in \cite{chapoton:two}. In particular, Chapoton provides several heuristic arguments for the conjecture in \cite{chapoton:one}. Also in the case $k=1$, Athanasiadis recently proved Chapoton's conjecture for $M$ and $F$ \cite{athanasiadis:triangles}, and this proof is uniform when combined with recent results in \cite{ABMW}. Krattenthaler has recently proved our conjecture for the Fuss-generalized $M$ and $F$ triangles (see \cite{krattenthaler1,krattenthaler2,kratt-muller}). Tzanaki \cite{tzanaki:faces} gave an alternative proof, generalizing the $k=1$ result of Athanasiadis \cite{athanasiadis:triangles}. The case of the $H$-triangle is still not understood.

When $k=1$, it turns out that each of the triangles is self-dual, and hence Definition \ref{def:dualtriangles} is redundant. The fact that $M(x,y)=M^*(x,y)$ follows easily from the fact that the lattice $NC^{(1)}(W)$ is self-dual. Since the $F$-triangle $F(p,q)$ is really a generalization of the face polynomial of the complex $\Delta(W)$, the relation $F(p,q)=F^*(p,q)$ follows from the Dehn-Sommerville relations and the fact that $\Delta^{(1)}(W)$ is a simplicial sphere. It also follows that $H(s,t)=H^*(s,t)$ in the $k=1$ case, but here there is {\em no known duality} on $NN^{(1)}(W)$ to explain this. 

When $k>1$, the dual triangles really differ from the usual triangles. For example, in the case $k=2$ and $W=A_2$, we have
\begin{align*}
& M^*(x,y)= 1-16y+7y^2+6xy-12xy^2+5x^2y^2, \\
& H^*(s,t)=  1+4t+2st+2t^2+2st^2+s^2t^2, \\
& F^*(p,q)= 5+12p+4q+7p^2+4qp+p^2.
\end{align*}
Thus, it is reasonable to ask for combinatorial interpretations of $H^*(p,q)$ and $F^*(s,t)$. We have no suggestion for the dual $H$ triangle, but we have a conjectural interpretation\footnote{now proved by Krattenthaler \cite{krattenthaler1}.} of the dual $F$-triangle.

\begin{conjecture}
\label{conj:dualF}
For all finite Coxeter groups $W$ and positive integers $k$, we have
\begin{equation*}
F^*(p,q)= \sum_{A\in\Delta^{(k)}(W)} \frac{\Nar^{(k)}(W,\abs{A})}{\Nar^{(1)}(W,\abs{A})} \cdot p^{\abs{A^+}}\,\cdot\,q^{\abs{A^-}},
\end{equation*}
where $\Nar^{(k)}(W,i)$ are the Fuss-Narayana numbers as in Figure \ref{fig:fussnar}.
\end{conjecture}

Note that we have seen the numbers $\Nar^{(k)}(W,i)/\Nar^{(1)}(W,i)$ before, in Conjecture \ref{conj:hofh}. In general, they are not integers. However, the coefficients of the polynomial $F^*(x,x)$ {\em are} integers, and they are described as follows. Let $f_i$ denote the number of faces in $\Delta^{(k)}(W)$ with cardinality $i$. Fomin and Reading \cite{fomin-reading} call these the {\sf generalized Kirkman-Cayley numbers} corresponding to $W$ and $k$. Then Conjecture \ref{conj:dualF} implies that the coefficient of $x^i$ in $F^*(x,x)$ is equal to
\begin{equation*}
f^*_i:=\frac{\Nar^{(k)}(W,i)}{\Nar^{(1)}(W,i)}\cdot f_i.
\end{equation*}

There is a possible interpretation of these numbers.

The $f$-vector $(f_0,f_1,\ldots,f_n)$ of $\Delta^{(k)}(W)$ consists of the Kirkman-Cayley numbers, and the $h$-vector $(h_0,h_1,\ldots,h_n)$ consists of the Fuss-Narayana numbers, $h_i=\Nar^{(k)}(W,n-i)$. If we naively suppose that the reverse vector $(h_n,\ldots,h_1,h_0)$ is an ``$h$-vector'' and compute its corresponding ``$f$-vector'', we get $(f_0^*,f_1^*,\ldots,f_n^*)$. We will call this the {\sf $f^*$-vector} of the complex $\Delta^{(k)}(W)$. In general, {\em the $f^*$-vector of a complex $\Delta$ is not the $f$-vector of any complex}. However, there is a construction, called the {\sf canonical module} of $\Delta$, that has the $f^*$-vector as the coefficients in its Hilbert series (see \cite[Chapter 2.7]{stanley:green}). We do not know if this interpretation is significant here.

The main significance of the Chapoton triangles is that they encode the enumerative information about ``rank'' and ``positivity'' in the Fuss-Catalan objects. It is natural to wonder if the triangles can be generalized to encode finer information. For example, the notion of ``type'' refines that of ``rank''. There is also a natural way to refine ``positivity'': in the $H$-triangle, we may distinguish between roots based on height in the root poset, not just whether they are ``simple'' or ``not simple''.

\begin{problem}
Do there exist generalizations of the Chapoton triangles that encode a notion of ``type'', with a corresponding Conjecture \ref{conj:triangles}?. Do there exist multi-variable versions of the Chapoton triangles that encode finer degrees of ``positivity''?
\end{problem}

In general, the problem of relating the $F$- and $M$-triangles is currently much better understood than the $H$-triangle. This is explained by the fact that the set $NN^{(k)}(W)$ is defined only for Weyl groups, whereas $NC^{(k)}(W)$ and $\Delta^{(k)}(W)$ make sense for all finite Coxeter groups. One of the biggest open problems in the field is to extend the construction of $NN^{(k)}(W)$ to the noncrystallographic types. Conjecture \ref{conj:triangles} seems to suggest a way forward on this problem.

\begin{definition}
When $W$ is a noncrystallographic finite Coxeter group, we {\em define} the $H$-triangle by means of the relations in Conjecture \ref{conj:triangles}.
\end{definition}

This {\sf noncrystallographic $H$-triangle} encodes very refined enumerative information, but {\em corresponding to what object}? We are seeing a shadow, but we have no idea what is casting it.

\section{Future Directions}
\label{sec:future}

As a coda, we sketch three suggestions for future research.

\subsection{Noncrystallographic Root Poset}

Among the finite Coxeter groups $W$, ``most'' are Weyl groups. The noncrystallographic exceptions are the dihedral groups $I_2(m)$ (symmetries of the regular $m$-gon in $\reals^2$) for $m\not\in\{2,3,4,6\}$; the group $H_3$ (symmetries of the dodecahedron and icosahedron in $\reals^3$); and the group $H_4$ (symmetries of the $120$-cell and $600$-cell in $\reals^4$) (see \cite{coxeter:polytopes}).

As mentioned in Section \ref{sec:triangles}, Conjecture \ref{conj:triangles} allows us to define an ``$H$-triangle'' for these noncrystallographic types. For example, Figure \ref{fig:Htriangles} displays the $H$-triangles for the noncrystallographic types, in the case $k=1$. (Thanks to Fr\'ed\'eric Chapoton for this data.)

\begin{figure}
\vspace{.1in}
\begin{center}
\scalebox{0.8}{
\begin{tabular}{|c|l|}
\hline
$W$ & $H(s,t)$ \\
\hline
\hline
$I_2(m)$ & $1+2st+(m-2)t+s^2t^2$\\
\hline
$H_3$ & $1+3st+12t+3s^2t^2+4st^2+8t^2+s^3t^3$\\
\hline
$H_4$ & $1+4st+56t+6s^2t^2+19st^2+133t^2+4s^3t^3+5s^2t^3+9st^3+42t^3+s^4t^4$\\
\hline
\end{tabular}}
\end{center}
\caption{$H$-triangles for the noncrystallographic types, with $k=1$}
\label{fig:Htriangles}
\end{figure}

Notice that the coefficients here are nonnegative, so in principle they may be counting something. Morally, they should encode information about antichains in the ``root poset''. However, when we defined the root poset earlier (Definition \ref{def:rootposet}), we deliberately excluded the noncrystallographic types. Of course, one might extend Definition \ref{def:rootposet} by setting $\alpha\leq\beta$ for $\alpha,\beta\in\Phi^+$ whenever $\beta-\alpha$ is in the positive {\em real} span of the simple roots $\Pi$, but it turns out the result of this definition for the noncrystallographic types has completely the wrong properties. It seems that this is not the correct way to proceed.

What properties {\em should} a noncrystallographic root poset have? The following are known combinatorial properties of crystallographic root posets. Let $W$ denote a Weyl group of rank $n$.

\medskip

\noindent{\bf Property 1:}\, Recall the definition and properties of the exponents $m_1\leq m_2\leq\cdots\leq m_n$ of $W$ from Section \ref{sec:invariant}. Writing the exponents in reverse (weakly decreasing) order, we note that
\begin{equation*}
\lambda=(m_n,\ldots,m_2,m_1)\vdash N
\end{equation*}
is an integer partition of the number $N$ of reflections in $W$, where the largest part is $m_n=h-1$. The {\sf dual partition} $\lambda^*=(k_1,k_2,\ldots,k_{h-1})\vdash N$ is defined by setting
\begin{equation*}
k_i:=\#\left\{ m_j: m_j=i\right\}.
\end{equation*}

The following appears as Theorem 3.20 in Humphreys \cite{humphreys}. He attributes the result to A.~Shapiro (unpublished), with the first uniform proof by Kostant \cite{kostant}.

\begin{theorem}
Let $\Phi^+$ denote the crystallographic positive roots of type $W$. The number of positive roots of height $i$ is equal to $k_i$, as defined above. In particular, the highest root has height $h-1$.
\end{theorem}

That is, the integer partition of $N$ that is dual to the exponents gives the rank numbers of the root poset $(\Phi^+,\leq)$.

\medskip

\noindent{\bf Property 2:}\, Recall from Section \ref{sec:antichains} that the number of antichains in the root poset $(\Phi^+,\leq)$ is equal to the Coxeter-Catalan number $\Cat^{(1)}(W)$ \eqref{eq:fusscat}, and the number of antichains with cardinality $i$ is equal to the Narayana number $\Nar^{(1)}(W,i)$ \eqref{eq:fussnar}.

\medskip

\noindent{\bf Property 3:}\, Recall from Section \ref{sec:antichains} that the number of antichains in the root poset $(\Phi^+,\leq)$ consisting of nonsimple roots is equal to the positive Coxeter-Catalan number $\Cat_+^{(1)}(W)$ \eqref{eq:posfusscat}. Furthermore, it is believed that the numbers of antichains of nonsimple roots, refined by cardinality, coincide with the entries of the $h$-vector of the positive cluster complex $\Delta_+(W)$. (Athanasiadis and Tzanaki proved this result in the classical types \cite{athanasiadis-tzanaki:positive} and conjectured the result in general \cite[Conjecture 1.2]{athanasiadis-tzanaki:positive}.)

\medskip

\noindent{\bf Property 4:}\, The following statistic was defined by Chapoton in \cite{chapoton:two}.
\begin{definition}
Given a finite Coxeter group $W$ with exponents $m_1\leq\cdots\leq m_n$, define
\begin{equation*}
M(W):= \frac{nh}{\abs{W}} \prod_{i=2}^n (m_i-1).
\end{equation*}
\end{definition}
He proved \cite[Proposition 1.1]{chapoton:two} that this number counts a certain class of reflections in $W$.

\begin{theorem}[\cite{chapoton:two}]
The number of reflections in $W$ which do not occur in any standard parabolic subgroup is equal to $M(W)$.
\end{theorem}

In the case that $W$ is a Weyl group, it follows that $M(W)$ also counts the number of roots in $(\Phi^+,\leq)$ with {\sf full support} (we say a root has full support if it occurs above all simple roots in the root order; equivalently it contains all simple roots in its simple root expansion). Fomin and Reading \cite{fomin-reading} also studied this statistic.

\medskip

\begin{figure}
\vspace{.1in}
\begin{center}
\scalebox{.6}{
\input{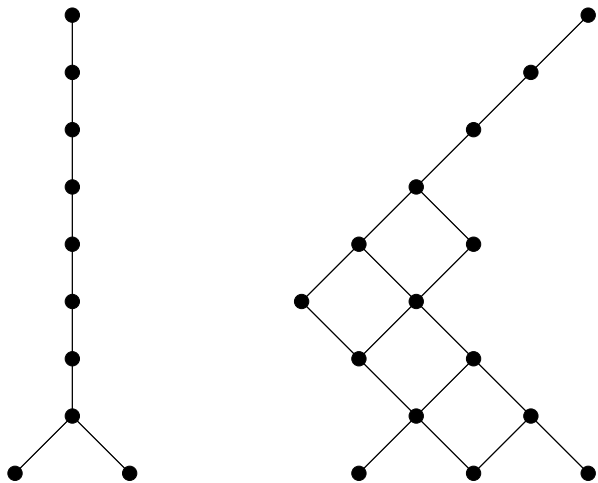}}
\end{center}
\caption{The root posets of types $I_2(m)$ and $H_3$?}
\label{fig:noncrystalroots}
\end{figure}

Based on these properties, we are led to guess the structure of the ``root posets'' of types $I_2(m)$ and $H_3$. One can verify that the posets in Figure \ref{fig:noncrystalroots} satisfy all four of the above properties, and their $H$-triangles agree with Figure \ref{fig:Htriangles}. We are unable at this moment to suggest a root poset of type $H_4$, since the calculations are much more complicated.

What are these posets? Where do they come from? It is remarkable that posets satisfying the desired properties exist at all. Something is going on here.

\begin{problem}
Define a ``root poset'' on the positive roots of a possibly-noncrystallographic root system. This poset should satisfy combinatorial properties 1--4 above and agree with Figures \ref{fig:Htriangles} and \ref{fig:noncrystalroots}. This should allow the definition of a poset $NN(W)$ that embeds in the partition lattice $\Pi(W)$ for noncrystallographic types, extending the results of Section \ref{sec:nonnesting}.
\end{problem}

Since the root poset of a crystallographic root system is an important structure in Lie theory and representation theory, the existence of a noncrystallographic root poset may have a wide impact.

\subsection{Cyclic Sieving}
In \cite{reiner-stanton-white}, Reiner, Stanton and White introduced and studied the {\em cyclic sieving phenomenon}; a mysterious connection between number theory and combinatorics.

Let $X$ be a finite set acted on by a cyclic group $C$ of order $n$, and let $X(q)$ be a polynomial with nonnegative integer coefficients and with the property that $X(1)=\abs{X}$. We also fix an isomorphism $\omega:C\hookrightarrow\complex^\times$ of $C$ with the complex $n$-th roots of unity.

\begin{definition}[\cite{reiner-stanton-white}]
We say that the triple $\left(X,X(q),C\right)$ exhibits the {\sf cyclic sieving property} if we have
\begin{equation*}
\left[ X(q) \right]_{q=\omega(c)} =\abs{\left\{x\in X:c(x)=x\right\}}.
\end{equation*}
for all $c\in C$.
\end{definition}

That is, the number of elements of $X$ fixed pointwise by the $d$-th power of a generator for $C$ is equal to $X(q)$ evaluated at the $d$-th power of a primitive $\abs{C}$-th root of unity.

They showed that this phenomenon occurs surprisingly often in connection with several well-known ``$q$-statistics'' such as the $q$-binomial coefficients. The theory of $q$-enumeration (usually over a finite field) has a rich history, so there is a large class of examples in which to search.

Athanasiadis and Garoufallidis have suggested the following $q$-analogue of the Fuss-Catalan number and conjectured that this is a polynomial in $q$ with nonnegative integer coefficients (see \cite[Problem 2.1]{armstrong:braids}).

\begin{definition}
\label{def:qfusscat}
Given a finite Coxeter group $W$ and positive integer $k$, define the {\sf $q$-Fuss-Catalan number}
\begin{equation}
\label{eq:qfusscat}
q\text{-}\Cat^{(k)}(W):=\prod_{i=1}^n \frac{[kh+d_i]_q}{[d_i]_q},
\end{equation}
where $[m]_q:=1+q+q^2+\cdots+q^{m-1}$ is the usual $q$-analogue of the integer $m$.
\end{definition}

We suggest\footnote{Conjectures about cyclic sieving are not hard to make. Priority for the Fuss-version ($k>1$) of Conjecture \ref{conj:sievingNC} technically rests with us since the object $NC_{(k)}(W)$ appears in this memoir for the first time. Conjecture \ref{conj:sievingclusters} was likely observed by several people independently.} that this statistic exhibits the cyclic sieving phenomenon in two distinct ways.

First, recall from Section \ref{sec:automorphisms} the definition of the automorphism $C^*$ on the poset $NC_{(k)}(W)$ of delta sequences (Definition \ref{def:C*}). This map generates a cyclic group of order $kh$ acting on the set $NC_{(k)}(W)$.

\begin{conjecture}
\label{conj:sievingNC}
The triple $\left(NC_{(k)}(W),q\text{-}\Cat^{(k)}(W),\langle C^*\rangle\right)$ exhibits the cyclic sieving property.
\end{conjecture}

Reiner, Stanton and White \cite{reiner-stanton-white} originally posed and solved this problem in the case $k=1$ and $W=A_{n-1}$. Bessis and Reiner \cite{bessis-reiner} recently posed and solved the problem for $k=1$ and all {\em well-generated complex reflection groups} $W$.

Second, recall the generalization $\tau^*$ of the tropical Coxeter element (Definition \ref{def:tropical*}), acting as a ``rotation'' on the $k$-divisible cluster complex $\Delta^{(k)}(W)$. Fomin and Reading showed that $\tau^*$ has order $(kh+2)/2$ if the longest element $w_\circ\in W$ is the antipodal map, and order $kh+2$ otherwise \cite[Lemma 5.2]{fomin-reading}. Instead of the vertices of $\Delta^{(k)}(W)$, we are interested in the action of $\tau^*$ on the maximal faces (colored clusters). Let $\Delta^{(k)}(W)^*$ denote the set of colored clusters.

\begin{conjecture}
\label{conj:sievingclusters}
The triple $\left( \Delta^{(k)}(W)^*,q\text{-}\Cat^{(k)}(W),\langle \tau^*\rangle\right)$ exhibits the cyclic sieving property.
\end{conjecture}

Again, the case $k=1$ and $W=A_{n-1}$ was posed and solved by Reiner, Stanton and White in their original paper \cite{reiner-stanton-white}. Eu and Fu \cite{eu-fu} have proved Conjecture \ref{conj:sievingclusters} case-by-case for all types, following a suggestion of Reiner.

It is interesting that the polynomial $q\text{-}\Cat^{(k)}(W)$ seems to have integer values when evaluated at $kh$-th roots of unity {\em and} when evaluated at $(kh+2)$-th roots (or $(kh+2)/2$-th roots) of unity.

\begin{problem}
What role do the polynomials $q\text{-}\Cat^{(k)}(W)$ play in the Fuss-Catalan combinatorics? Do these count some objects defined over a finite field of order $q$? Is there a natural definition of ``$q$-Fuss-Narayana polynomials'' $q\text{-}\Nar^{(k)}(W,i)$ generalizing \eqref{eq:fussnar}?
\end{problem}

\subsection{Diagonal Harmonics}
\label{sec:harmonics}

Our final suggestion connects with the theory of diagonal harmonics, as discussed by Haiman in \cite{haiman}. Let
\begin{equation*}
\rational[X,Y]:=\rational[x_1,\ldots,x_n,y_1,\ldots,y_n]
\end{equation*}
denote the polynomial ring over $\rational$ in the $2n$ variables $x_1,\ldots,x_n,y_1,\ldots,y_n$. The {\sf diagonal action} of the symmetric group $\frak{S}_n$ on the ring $\rational[X,Y]$ is defined by setting
\begin{equation*}
\sigma\cdot f(x_1,\ldots,x_n,y_1,\ldots,y_n)= f(x_{\sigma(1)},\ldots,x_{\sigma(n)},y_{\sigma(1)},\ldots,y_{\sigma(n)})
\end{equation*}
for all $\sigma\in \frak{S}_n$ and $f\in\rational[X,Y]$. That is, the symmetric group acts simultaneously on the sets $\{x_1,\ldots,x_n\}$ and $\{y_1,\ldots,y_n\}$ by permutations. 

\begin{definition}
If $I$ is the ideal in $\rational[X,Y]$ generated by $\frak{S}_n$-invariant polynomials of positive degree, we define the quotient ring
\begin{equation*}
R(n):=\rational[X,Y]/I.
\end{equation*}
\end{definition}

At first, one might suppose that $R(n)$ is just the product of the invariant rings of $\rational[X]$ and $\rational[Y]$. However, since there exist ``mixed invariants'', the structure of $R(n)$ is much more interesting, and much more difficult to study. There is also an important interpretation of $R(n)$ as a space of ``harmonics''. The study of the ring $R(n)$ was initiated by Garsia and Haiman and it is a thriving object of current research. Haiman surveyed the main features of this subject in \cite{haiman}. Stanley gave a more recent survey \cite[Section 3]{stanley:progress} after some of the main conjectures in \cite{haiman} had been proven.

We will describe only the features of the theory that are relevant to our current purpose. Since the ring $\rational[X,Y]$ is ``bigraded'' by $x$-degree and $y$-degree and the ideal  $I$ is homogeneous, the quotient $R(n)$ inherits this bigrading. Some of the following conjectures have now been proven and some have not. We refer to \cite[Section 3]{stanley:progress}.

\begin{conjecture}[\cite{haiman}]
The ring $R(n)$ has dimension $(n+1)^{n-1}$ as a real vector space.
\end{conjecture}

The diagonal action of $\frak{S}_n$ naturally induces an action on $R(n)$ that respects the bigrading; thus, each isotypic component of $R(n)$ is bigraded. Let $\Cat(n,t,q)$ denote the dimension of the $(t,q)$-bigraded component in the isotypic component of the sign representation.

Let $[n]_q!=[n]_q[n-1]_q\cdots[2]_q[1]_q$ denote the $q$-factorial and let $\qbin{n}{k}_q$ denote the standard $q$-binomal coefficient $\frac{[n]_q!}{[k]_q![n-k]_q!}$.

\begin{conjecture}[\cite{haiman}]
\hspace{.1in}
\begin{enumerate}
\item $\Cat(n,1,1)$ is equal to the Catalan number $\frac{1}{n+1}\binom{2n}{n-1}$.
\item $q^{\binom{n}{2}}\Cat(n,q^{-1},q)$ is equal to the $q$-Catalan number 
\begin{equation*}
q\text{-}\Cat^{(1)}(A_{n-1})=\frac{1}{[n+1]_q}\qbin{2n}{n}_q
\end{equation*}
as in Definition \ref{def:qfusscat}.
\end{enumerate}
\end{conjecture}

Thinking of $\frak{S}_n$ as the Weyl group of type $A_{n-1}$, one may interpret the ring $R(n)=R(A_{n-1})$ in the context of reflection groups. Since there is also a natural invariant ring for other Weyl groups (see Section \ref{sec:invariant}), Haiman wondered whether the diagonal action of $W$ would yield a ring $R(W)$ of coinvariants with similar combinatorial properties. In \cite[Section 7]{haiman}, he observed that the obvious definition of $R(W)$ does not yield the expected combinatorics, but that it is very close. He gave a conjecture \cite[Conjecture 7.2.3]{haiman} for what the ``correct'' definition of $R(W)$ should be. We reproduce his definition:

\begin{definition}[\cite{haiman}]
\label{def:quotient}
Let $\frak{g}$ be a semisimple complex Lie algebra, $\frak{h}$ a Cartan subalgebra, and $W$ the Weyl group. Let $\mathcal{C}_{\mathcal{N}}$ be the commuting nullcone, i.e. the set of pairs $(X,Y)\in\frak{g}\times\frak{g}$ such that $X$ and $Y$ commute and are nilpotent. Let $I(W)$ be the restriction to $\rational[U\oplus U]$ (which is the coordinate ring of $\frak{h}\times\frak{h}$) of $I(\mathcal{C}_{\mathcal{N}})$, and set $R(W):=\rational[U\oplus U]/I(W)$.
\end{definition}

Recently, Gordon used rational Cherednik algebras to construct a quotient with the desired combinatorial properties. He proved the following.

\begin{theorem}[\cite{gordon}]
There exists a quotient $R^*(W)$ of the coordinate ring of $\frak{h}\times\frak{h}$ satisfying the following properties:
\begin{enumerate}
\item $\dim R^*(W)= (h+1)^n$.
\item $R^*(W)$ is $\integers$-graded with Hilbert series $t^{-N}(1+t+\cdots+t^h)^n$.
\item The image of the polynomial ring $\complex[\frak{h}]$ in $R^*(W)$ is the classical coinvariant algebra, $\complex[\frak{h}]^{\rm{co}W}$.
\item If $\varepsilon$ is the sign representation of $W$, then $R^*(W)\otimes\varepsilon$ is isomorphic as a $W$-module to the permutation representation of $W$ on the reduction of the root lattice modulo $h+1$, written $Q/(h+1)Q$. 
\end{enumerate}
\end{theorem}

It is not clear at this point whether Gordon's ring $R^*(W)$ coincides exactly with Haiman's $R(W)$. It is also not immediately clear whether Gordon's ring carries a bigrading. However, the following is immediate.

\begin{theorem}
\label{th:aftergordon}
The dimension of the sign-isotypic component in the $W$-module $R^*(W)$ is equal to the Coxeter-Catalan number
\begin{equation*}
\Cat(W)=\frac{1}{\abs{W}}\prod_{i=1}^n (h+d_i).
\end{equation*}
\end{theorem}

\begin{proof}
Let $\chi$ denote the character of $W$ acting on $R^*(W)$, and let $\chi_\varepsilon$ denote the sign character (determinant character) of $W$; thus, the module $R^*(W)\otimes\,\varepsilon$ carries the character $\chi\cdot\chi_{\varepsilon}$. Since $\chi\cdot\chi_{\varepsilon}$ is a permutation character, we have $\chi(w)\cdot\chi_{\varepsilon}(w)=\abs{\Fix(w)}$, the number of elements of $Q/(h+1)Q$ fixed by $w\in W$. Considering inner products of characters, the dimension of the sign-isotypic component of $R^*(W)$ as a $W$-module is equal to
\begin{align*}
\langle \chi,\chi_{\varepsilon}\rangle &= \frac{1}{\abs{W}}\sum_{w\in W}\chi(w)\cdot\overline{\chi_{\varepsilon}(w)} \\
&= \frac{1}{\abs{W}}\sum_{w\in W}\chi(w)\cdot\chi_{\varepsilon}(w) \\
 &= \frac{1}{\abs{W}}\sum_{w\in W}\abs{\Fix(w)},
\end{align*}
which by the Burnside-Cauchy-Frobenius-Polya-Redfield counting lemma is equal to the number of orbits in $Q/(h+1)Q$ under the action of $W$. Haiman proved \cite[Theorem 7.4.4]{haiman} that this number is equal to $\Cat(W)$.
\end{proof}

Thus, the occurrence of the Coxeter-Catalan number suggests a link between the theory of diagonal harmonics and the Catalan combinatorics as described in this memoir. We expect that these two subjects will eventually become one.

\begin{problem}
\label{prob:theend}
Find a bigrading on Gordon's ring $R^*(W)$ and define $\Cat(W,t,q)$ to be the $(t,q)$-bigraded component of the sign-isotypic component of $R^*(W)$. In this case, we {\em should} have
\begin{equation*}
q^N \Cat(W,q^{-1},q)=q\text{-}\Cat(W)=\prod_{i=1}^n\frac{[h+d_i]_q}{[d_i]_q},
\end{equation*}
where $N$ is the number of reflections in $W$. Explain the significance of this bigrading for the noncrossing partitions $NC(W)$, the nonnesting partitions $NN(W)$ and the cluster complex $\Delta(W)$. Is there a $(t,q)$-bigraded Narayana number $\Nar(W,t,q,i)$ refining $\Cat(W,t,q)$, and such that $\Nar(W,1,1,i)=\Nar(W,i)$? Does there exist a bigraded $W$-module $R^{(k)}(W)$ generalizing Definition \ref{def:quotient}, with dimension $(kh+1)^n$ and having sign-isotypic component with dimension equal to the Fuss-Catalan number $\Cat^{(k)}(W)$ \eqref{eq:fusscat}?\footnote{Griffeth \cite{griffeth} has generalized Gordon's construction to the Fuss-Catalan case.} To what extent can the theory of diagonal harmonics be extended to noncrystallographic reflection groups?
\end{problem}


\backmatter

\bibliographystyle{amsalpha}

\end{document}